\documentclass[reqno]{amsart}

\usepackage[letterpaper]{geometry}
\usepackage[mathscr]{euscript}
\usepackage{amsmath,amsfonts,hyperref,amsthm,amsbsy,bm,thmtools,amssymb,tikz,cleveref,comment}
\usetikzlibrary{arrows} 
\DeclareMathOperator{\Adm}{Adm}
\DeclareMathOperator{\diam}{diam}

\DeclareMathOperator{\Imaginary}{Im}
\DeclareMathOperator{\Real}{Re}

\DeclareMathOperator{\vol}{vol}
\DeclareMathOperator{\mix}{mix}
\DeclareMathOperator{\Id}{Id}
\DeclareMathOperator{\dist}{dist}

\numberwithin{equation}{section}
\newtheorem{thm}{Theorem}[section]
\newtheorem{prop}[thm]{Proposition}
\newtheorem{lem}[thm]{Lemma}
\newtheorem{cor}[thm]{Corollary}

\newtheorem{rem}[thm]{Remark}
\theoremstyle{definition}
\newtheorem{definition}[thm]{Definition}

\newcommand{\td}{\mathtt d}
\newcommand{\te}{\mathtt e}
\newcommand{\cH}{\mathcal H}
\newcommand{\cU}{\mathcal U}

\DeclareMathOperator{\eq}{eq}

\newcommand{\mq}{q}
\newcommand{\tmix}{t_{\rm mix}}

\newcommand{\me}{\mathfrak e}
\newcommand{\md}{\mathfrak d}
\newcommand{\mX}{\mathfrak X}
\newcommand{\mF}{\mathfrak F}
\newcommand{\mC}{\mathfrak C}
\newcommand{\mA}{\mathfrak A}

\newcommand{\mL}{\mathfrak L}
\newcommand{\bfK}{{\bf K}}

\newcommand{\ONa}{\Omega^N_{a,b,c}}
\newcommand{\PNa}{\mathbb P^N_{a,b,c}}

\newcommand{\bu}{{\bm u}}
\newcommand{\bv}{{\bm v}}
\newcommand{\be}{{\bm e}}
\newcommand{\bd}{{\bm d}}

\DeclareFontFamily{U}{mathx}{}
\DeclareFontShape{U}{mathx}{m}{n}{<-> mathx10}{}
\DeclareSymbolFont{mathx}{U}{mathx}{m}{n}
\DeclareMathAccent{\widehat}{0}{mathx}{"70}
\DeclareMathAccent{\widecheck}{0}{mathx}{"71}

\makeatletter
\@addtoreset{equation}{section}

\makeatother



\newcommand{\bw}{{\bm w}}

\newcommand{\PNqa}{\mathbb P^{N,q}_{a,b,c}}

\title{Mixing times for Glauber dynamics of lozenge tilings of the hexagon}
\author{Amol Aggarwal}
\address{Stanford University, Department of Mathematics, 450 Jane Stanford Way, Stanford, CA 94305} \email{agga@stanford.edu}
\author{Fabio Toninelli}
\address{TU Wien, Institut f\"ur Stochastik und Wirtschaftsmathematik, Wiedner Hauptstra{\ss}e 8-10, A-1040 Vienna, Austria} \email{fabio.toninelli@tuwien.ac.at}

\begin{document}

\begin{abstract}
We prove that the continuous-time, single-flip Glauber dynamics for lozenge tilings of the size-$N$ hexagon mix
in time $N^{2+o(1)}$. This was predicted to hold on fairly general domains of diameter $N$ (on the basis of the ``Lifshitz law'' heuristic) but had previously only been established in domains such that the associated limit shape has no frozen facets. To access the hexagon, we introduce a multi-scale comparison argument between the height function of the random tiling evolving under the Glauber dynamics and the limit shape of a volume-tilted tiling (whose tilting parameter varies suitably in time).
 \end{abstract}

\maketitle
\tableofcontents

\section{Introduction}

\subsection{The lozenge tiling Glauber dynamic on the hexagon}
Random lozenge tilings, or equivalently the dimer model on the
honeycomb lattice, are a classical subject of study in statistical
mechanics, probability and combinatorics; see e.g., the monograph
\cite{RT}.  To introduce the problem, let $\mathbb T_N,N\in\mathbb N$,
be the infinite triangular grid with equilateral triangular faces of
side-length $1/N$, and let $D$ be a bounded domain of the plane that can be
tiled by lozenge-shaped tiles, obtained as the union of two adjacent
triangular faces. Call $\Omega^N$ the set of all
lozenge tilings of $D$ and $\mathbb P^N$ the uniform probability
measure on $\Omega^N$.

In addition to being of intrinsic interest, random tilings are archetypes of strongly correlated systems that are very sensitive to boundary conditions. They may be viewed as models of random surfaces, which are known \cite{CLP,CKP} to converge to a \emph{limit shape} as $N$ tends to $\infty$. A salient feature of these limit shapes is that they are highly inhomogeneous. They admit \emph{frozen regions}, where the surface is almost deterministically flat, and \emph{liquid regions}, where it appears more rough and random; the curve separating these phases is called an \emph{arctic boundary} \cite{RTT}. Across these different regions, the tiling admits considerably different asymptotic behaviors, which can be computed exactly due to the underlying determinantal structure of the model \cite{SDL,LSLD} that becomes explicit on specific polygonal domains, such as the hexagon \cite{NPRT}. For example, when $D$ is a regular hexagon, in the liquid region the tiling has Gaussian free field fluctuations \cite{AURT}; around the arctic boundary it has Kardar--Parisi--Zhang (KPZ) fluctuations \cite{CFPLG,SFFC,DP,ARTS}; and near points where the arctic boundary meets the hexagon, called tangency locations, the tiling has Gaussian unitary ensemble (GUE) fluctuations \cite{EM,TRM}. 

While the measure $\mathbb P^N$ is by now well studied,  local Markovian dynamics on  $\Omega^N$ (admitting $\mathbb P^N$ as stationary measure) and their mixing properties are still poorly understood. In this work, we consider the most natural such dynamic, that is  
the \emph{single-flip Glauber dynamic} on  $\Omega^N$, which is the
continuous-time Markov chain on $\Omega^N$ whose updates are the
elementary rotations depicted in Figure  \ref{fig:esagono} (right), both
having transition rate $1$. This dynamic  is routinely used\footnote{a freely available code implementing the Glauber dynamic can be found at  https://lpetrov.cc/lozenge-draw/ }, to numerically sample uniformly chosen random lozenge tilings of the domain $D$ (and can be implemented within the coupling-from-the-past framework; see  \cite{Propp-Wilson} and \cite[Chapter 25]{LevinPeres}).
We mention also that the single-flip Glauber dynamic is the zero-temperature limit of the Glauber dynamic of the three-dimensional Ising model with Dobrushin boundary conditions \cite{TSMF,CMT}. 

The uniform measure $\mathbb P^N$ on $\Omega^N$ is reversible for the Glauber dynamic and is
its  unique stationary measure. Since the state space is
finite, the total variation mixing time
$t_{\mix}^N$ is also finite. The question of quantifying $t_{\mix}^N$ as a function of
$N$ has an interesting history, and an incomplete status. It was first proven in \cite{LRS} that a \emph{non-local version} of the Glauber
dynamic, called the ``tower-move dynamic,'' mixes in polynomial (in $N$)  time. Via standard spectral
gap comparison tools, this implies polynomial time mixing for the single-flip
dynamic, as well. Subsequently, when $D$ is  the regular hexagon $\mX$ of side-length $1$ with sides parallel to the
axis of the triangular grid (that can be tiled by exactly $3N^2$ tiles), 
the mixing time of the tower-move
dynamic  was determined in \cite{Wilson} to be of
order $N^2\log N$ (with upper and lower bounds differing by a constant
prefactor). This result implies (via the above spectral gap comparison)
that $t^N_{\mix}\lesssim N^{6+o(1)}$ \cite{Tetali}. A more refined comparison argument (that combines the analysis of the tower-dynamic mixing time with the censoring inequality \cite{PW}) allows to prove $t^N_{\mix} \lesssim N^{4+o(1)}$ \cite{TSMF}; this was the best existing upper bound,
 for the mixing time of the Glauber dynamic in the hexagon,
before the present work.

However, on the basis of a general ``mean curvature motion heuristic'' (see also \Cref{sec:novelty} for more discussion), it was predicted that $t^N_{\mix}\approx N^2$ also for
the single-flip dynamic, possibly with logarithmic sub-leading
corrections.  This conjecture turned out to be rather elusive, for
several reasons. First, in contrast with the tower-move dynamic,
coupling arguments fail to provide any contraction for the Glauber
dynamic. Second, the dynamic has no perturbative parameter to tune
(such as the temperature for Ising models, or the number of colors and
graph degree for random colorings). Third, equilibrium correlations
for the system are not summable, so one cannot exploit strong spatial
mixing properties (in the spirit of \cite{Martinelli}) to upper bound the mixing time. Moreover, due
to the hard constraints intrinsic in the tiling model, configurations
can be blocked (unable to perform any local move) in arbitrarily large
regions of the graph.

Still, a first step towards the conjecture was done in
\cite{CMT}, where the bounds
\begin{equation}
  \label{eq:ltdmls}
N^2\lesssim t^N_{\mix}\lesssim N^{2+o(1)}  
\end{equation}
were proven  for the single-flip dynamic in a tilable domain $D$ of the plane (of diameter of order $1$), provided that the limit shape on  $D$  is affine. Later, \eqref{eq:ltdmls} was generalized  \cite{LTDMLS,MTTD} to the case where the limit shape is $C^\infty$. These partial results point out that the geometry of the domain $D$, and
in particular the regularity properties of the associated limit shape, play a subtle role in the convergence to equilibrium.

Arguably, the most natural and classical geometric domain where 
lozenge tilings can be studied is the hexagon $\mX$, and neither of the cited results  \cite{CMT,LTDMLS,MTTD}
applies there. Indeed, the limit shape on $\mX$ \cite{CLP} (as well as on other polygonal domains \cite{KO,Astala}) admits an arctic boundary, along which the limit shape is not smooth.

The main result of the present work  is a proof of the asymptotic 
$t_{\mix}^N\approx N^2$ for the Glauber dynamic on $\mX$, up to
sub-leading multiplicative corrections:
        \begin{thm}\label{thm:main}
        	For any $\delta>0$, there exist constants $\kappa > 0$ and $C = C(\delta)>0$ such that the mixing time of the single-flip Glauber dynamic on tilings of the hexagon $\mX$ satisfies
          \begin{eqnarray}
            \label{eq:boundstmix}
            \kappa N^2\le             t^N_{\rm mix}   \le  C N^{2+\delta}.
          \end{eqnarray}
        \end{thm}

        With minor modifications, one can also  prove the same 
        statement for a hexagon of fixed positive side-lengths $a,b,c$ (in which case the constants $\kappa,C$ would depend on $a,b,c$). However, we will not pursue this here, in order to avoid making the notation heavier.

        The   effort in proving \Cref{thm:main} is the upper bound: the argument for $t^N_{\rm mix} \ge \kappa N^2$, already given in
        \cite[Section 5.2]{CMT} for the dynamic in a domain $\Omega\subset \mathbb{T}_N$ where the
        limit shape is affine, works identically for the dynamic in
        $\mX$ and we will not repeat it. 
\subsection{Proof ideas}
\label{sec:novelty}

{
	The conjecture $t_{\mix} = N^{2+o(1)}$ is based on the ``Lifshitz law'' heuristic \cite{LODT} stating that, upon rescaling time by $N^2$ and in the limit $N\to\infty$, the limiting evolution of the tiling height function $\varphi$ under the Glauber dynamics should follow an anisotropic variant of mean curvature motion, of the form
	\begin{flalign}
		\label{tmotion} 
		\partial_t \varphi = \mu (\nabla \varphi) \cdot \mathcal{L} \varphi.
	\end{flalign} 
	
	\noindent In \eqref{tmotion}, $\mathcal{L}$ is an explicit semi-elliptic, non-linear operator (related to the surface tension functional of random tilings) and $\mu(\nabla \varphi)$ is a certain slope-dependent ``mobility coefficient'' governing the speed of the evolution; see \cite{Spohn} for further commentary. The limit \eqref{tmotion} has only been proven in the case of the tower-move dynamics \cite{LasTonhydro} (assuming the solution to \eqref{tmotion} is $C^2$), where the associated mobility coefficient $\mu$ is known explicitly. For the Glauber dynamics, there is no predicted expression for $\mu$. Moreover, under hexagonal boundary conditions, the solution $\varphi$ to \eqref{tmotion} should not be $C^2$, so even the well-posedness of \eqref{tmotion} is unclear. 
	
	At least the first issue was circumvented in \cite{CMT,LTDMLS,MTTD}, by instead proving that the evolution $H_t$ of the tiling height function under the Glauber dynamics likely remains between two deterministic ``barrier'' functions $\varphi^-_t$ and $\varphi^+_t$. The latter might be viewed as ``proxies'' to solutions of \eqref{tmotion} and were obtained by perturbing $\mathcal{H}$ by  functions $\psi^\pm_t$ whose Hessian dominates that of $\mathcal{H}$. The latter in particular necessitates that $\mathcal{H}$ is $C^2$, which does not hold if $\mathcal{H}$ is the limit shape for a random tiling of the hexagon. 
	
	To address this, we instead introduce different barriers, given by the limit shapes $\mathcal{H}_q$ of a random tiling of the hexagon, weighted by $e^{q \vol}$, where $\vol$ denotes the volume of the stepped surface associated with the tiling. Such ``$q$-deformed'' tiling models have been studied quite extensively over the past two decades \cite{LTEM,CFPLG,KO,DBPP,DLE}. Here, we use their limit shapes to study the Glauber dynamics of the original ($q=0$) tiling model. In particular, we prove (Proposition \ref{hqihqi1}) by induction on $t$ that
	\begin{flalign}
		\label{hqtqth}
		\mathcal{H}_{-q(t)} \le H_t \le \mathcal{H}_{q(t)},
	\end{flalign}
	
	\noindent holds with high probability, where $q(t)>0$ depends suitably on $t$ (and $N$), decreasing to $0$ as $t$ becomes large. We will eventually take $q(t)$ as small as $N^{-1}$ for 
        $t=N^{2+o(1)}$. This implies at such times that the height $H_t$ is so
close to the limit shape $\cH$ associated to the uniform measure
$\mathbb P^N$, that after another time-lag of $N^{2+o(1)}$ equilibrium is
reached. 

These $\mathcal{H}_q$ limit shapes are both fully explicit \cite{DBPP,DLE} and more closely model the behavior of the Glauber dynamics on the hexagon than the barrier functions $\varphi^\pm_t$ used in earlier works.\footnote{In fact, if one stops these Glauber dynamics before they mix, the associated height function looks numerically (see, e.g., https://lpetrov.cc/lozenge-draw/)  qualitatively quite similar to (though is not predicted to exactly coincide with) some $\mathcal{H}_q$.} For instance, they also admit arctic boundaries separating liquid from frozen phases. We will see that the mechanism under which the Glauber dynamics locally mixes will vary in mesoscopic regions between the bulk of the liquid phase and the arctic boundary, and that this must be taken into account when establishing \eqref{hqtqth}.
	
	To implement this idea, we begin by showing a ``base case'' for \eqref{hqtqth}, verifying it (Theorem \ref{thm:step1}) when $t = N^{2+o(1)}$ and $q(t) = q_0 > 0$ is a constant that is arbitrarily small but independent of $N$. First, we perturb the boundary conditions of $\mathcal{H}$ to make the associated limit shape smooth (Proposition \ref{thm:approx}), enabling us to apply the prior results \cite{LTDMLS,MTTD} to show that the Glauber dynamics under these perturbed boundary conditions mix in time $N^{2+o(1)}$. By comparing these to the original Glauber dynamics, we obtain (Proposition \ref{thm:perturb}) that $|H_t - \mathcal{H}| \lesssim q_0$ for $t = N^{2+o(1)}$ and fixed $q_0>0$. This is in the direction of \eqref{hqtqth} but does not quite verify it where $\mathcal{H}_{q_0}$ is frozen, as there we have $\mathcal{H}_{-q_0} = \mathcal{H}_{q_0}$. These regions are addressed by induction on $N$, involving ``embedding'' smaller hexagons (whose frozen regions contain those of $\mathcal{H}_{q_0}$) into the original one; see Section \ref{ProofHqt} for further details. 
	
        Next, we must show \eqref{hqtqth} for $q(t) < q_0$, which is
        the main task in this work.  A primary heuristic here is that,
        throughout the Glauber dynamics, one should expect all level
        lines of $H_t$ to move at a speed of order $N^{-2}$. However,
        unlike on domains admitting $C^2$ limit shapes, different
        level lines of $H_t$ on the hexagon should exhibit different
        fluctuations. Thus, it should take different amounts of time
        to ``detect'' whether the level line has moved; in this way,
        local mixing occurs at different rates across the hexagon. For
        example, in the bulk of the liquid region, the level lines
        should have $N^{o(1)-1}$ fluctuations, so their motion is
        detected after time $N^{1+o(1)}$; generically along the arctic boundary, they
        should have $N^{o(1)-2/3}$ fluctuations, so their motion is
        detected after time $N^{4/3+o(1)}$; and at the tangency locations, they should have $N^{o(1)-1/2}$ fluctuations, so their motion is detected only after time $N^{3/2+o(1)}$.
	
	To resolve this, we decompose the hexagon $\mathfrak{X}$ into ``annuli'' (Definition \ref{aq}), in which these level line fluctuations (and thus the rate of local mixing) are of the same order. We then establish the following two types of statements; in what follows, we focus on decreasing the upper bound on $H_t$ (increasing the lower bound is entirely analogous). The first, which we may refer to as ``liquid phase mixing,'' states that if the height function $H_t$ is upper bounded in an annulus not too close to the arctic boundary, then the upper bound on $H_t$ improves (that is, it decreases further) on inner annuli, after running the Glauber dynamics for a certain time (Lemma \ref{ak1akak1}). The second, which we may refer to as ``edge phase mixing,'' states that if $H_t$ is upper bounded in an annulus reasonably close to the arctic boundary, then this bound extends to the outer annuli (up to a negligible error) after some time (Lemmas \ref{akl0} and \ref{boundaryaj2}). By combining these two statements, we then prove \eqref{hqtqth} by induction on $t$. 
	
	To show these annulus mixing results, we first reduce them to their ``local versions'' (Propositions \ref{prop:hdecreasebulk}, \ref{htqhtq02}, and \ref{htq3}) that restrict to mesoscopic cells contained within the annuli. The benefit to doing this is that the limit shapes $\mathcal{H}_q$ are more regular in these mesoscopic cells. For example, although $\mathcal{H}_q$ is not $C^2$ near the arctic boundary,  it becomes smooth upon an explicit rescaling\footnote{This is sometimes called a ``blow-up'' procedure in the literature on obstacle problems \cite{TR}.} dependent on the cell's location (Proposition \ref{prop:formerassumption}). Given this, the proof of local bulk mixing (Proposition \ref{prop:hdecreasebulk}) uses an idea introduced in \cite{LTDMLS}. Specifically, we compare (the ``regular rescaling'' of) the limit shape $\mathcal{H}_{q(t)}$ to that $\mathcal{H}_{a,b,c}$ of an unweighted (that is, uniform) tiling on a certain $a \times b \times c$ hexagon, so that their gradients coincide at the cell's center but $\mathcal{H}_q$ is ``more concave'' than $\mathcal{H}_{a,b,c}$ (Proposition \ref{prop:last}). Using the added concavity of $\mathcal{H}_q$, with concentration bounds for the unweighted tiling around $\mathcal{H}_{a,b,c}$, we show that $H_{t'}$ must ``deflate'' to a smaller $\mathcal{H}_{q(t')}$ after some time $t'-t$; see Section \ref{sec61} for further details. 

	To show the local edge mixing results, we separate into cases, depending on whether the cell is away from or at the arctic boundary (Propositions \ref{htqhtq02} or \ref{htq3}), and whether it is away from or at the tangency location (Sections \ref{Proofhtqd0d} or \ref{Proofd}, and Sections \ref{sec:77caseI} or \ref{Proofh000}). The proofs in all of these cases are quite similar to each other, but they differ from the bulk case, since the limit shape $\mathcal{H}_q$ will not behave regularly in our mesoscopic cells close to the arctic boundary. However, the level lines of $\mathcal{H}_q$ will be more regular (Proposition \ref{prop:improvedconv}). So, to show the local edge mixing results, we develop an analog of the local bulk mixing reasoning that compares the convexity properties between the level lines of the tiling height functions (which are equivalent to non-intersecting random walks) approximating $\mathcal{H}_q$ and $\mathcal{H}_{a,b,c}$. See the beginning of Section \ref{Uq00h} for a more detailed exposition. }
	
 In order to implement the above strategy, we need to prove, among others,  several analytic statements concerning the volume-tilted limit shape $\cH_q$, that make precise the intuition that $\cH_q$ becomes more concave as $q$ grows, and that the curvature of its level lines changes in a specific way as a function of $q$. These statements are given in \Cref{sec:Vtilted} and proved (by direct, though tedious, Taylor expansions of the exact expressions for these limit shapes) in Appendices \ref{Shape} and \ref{appC}.

\subsection{Relation with existing literature}
\label{sec:literature}
The Lifshitz law heuristic mentioned above applies
to any reasonable reversible interface dynamics in $(d+1)$ dimensions, where $d$ refers to the space dimension and $1$ to time, corresponding in the case of lozenge tilings to $d=2$. This heuristic
leads to the conjecture that the inverse spectral gap is of order
$N^2$ and that the mixing time is of order $N^2$ times a logarithmic
(in $N$) factor, the latter originating from the fact that once the
interface is macroscopically close to the equilibrium one (which
happens after a time of the order of the inverse spectral gap), it
still takes some time before its law is close to the equilibrium one
in total variation distance.

The prior rigorous proofs of the
$\tmix\asymp N^2\log N$ asymptotic for interface dynamics had been been mostly limited to $(1+1)$-dimensional
random interface models \cite{CapLabLac,MartSinc,Lacoin} (see also
\cite{Funaki} and references therein).
A simplifying feature of one-dimensional
interface models is that their stationary measure, when viewed in
terms of the interface gradients, has an i.i.d. structure. In
contrast, for random tilings, equilibrium correlations have a slow
decay in space, of the order of the squared inverse distance.
Another case where such results were proven 
is  the Glauber dynamic of the $(d+1)$-dimensional Gaussian lattice free field \cite{GangGhe} (where the spectral gap can be computed exactly). In a different but related direction, \cite{BauBod} proved sharp results for the inverse spectral gap of Glauber dynamics of interface models on hierarchical lattices; the inverse gap has the same order as that of the corresponding Gaussian free field dynamic.
We point out that \cite{CapLabLac,GangGhe,Lacoin}
further prove that these dynamics exhibit the cutoff (in total variation) phenomenon.

Finally, let us mention the earlier  works \cite{CMTrsa,Green} on mixing of Glauber dynamics of random lozenge tilings, that take a point of view quite different from that of the present paper. These works study a biased version of the Glauber dynamic, where the two updates in \Cref{fig:esagono} have different transition rates. In this case, the phenomenology (both dynamic and equilibrium) is very different. First, the stationary measure of the biased dynamic is not uniform on the set of all tilings, but is instead the volume-tilted measure $e^{q\vol}$, where $q =  -\gamma N$ is of order $N$, and the parameter $\gamma\ne 0$ is related to the ratio of the two transition rates. Such a measure concentrates on configurations whose volume is of order $1/N$, while typical configurations sampled from $\mathbb P^N$ have volume of order $N^2$. Also, the Lifshitz law heuristic does not apply to the biased dynamic, because the drift  induced by the bias dominates the curvature effect; in fact, the mixing time and inverse spectral gap turn out to be of order $N$ and $O(1)$ respectively, in contrast with $N^{2+o(1)}$ and $N^2$ for the symmetric dynamic we study here.

\subsubsection{Notation}

\label{Notation} 

We say that an event $\mathscr{E}$ (that may depend on $N$) holds with \emph{overwhelming probability (w.o.p.)} (with respect to some sequence of probability measures $\mathbb P^N,N\in\mathbb N$) if the following holds. For any integer $L \ge 1$, there exists a constant $C > 1$ independent of $N$ (but possibly dependent on $L$ and other parameters involved) such that $\mathbb{P}^N[\mathscr{E}^{\complement}] \le CN^{-L}$ for all $N \ge 1$.

We use also the following standard notation: $f_N \lesssim g_N$ and $g_N \gtrsim f_N$ if there exists a positive constant $C>0$ (independent of $N$) such that $f_N \le C g_N$ (we may write $f_N \lesssim_Q g_N$ if $C$ depends on some parameter $Q$); $f_N \ll g_N$ and $g_N \gg f_N$ if $\lim_{N\rightarrow \infty} f_N / g_N = 0$; and $f_N \asymp g_N$ if $f_N \lesssim g_N$ and $g_N \lesssim f_N$. Given $D\subset \mathbb R^2$,    we let  ${\rm diam}(D)$ denote its diameter.

Finally, given $\gamma\in \mathbb{Z}_{\ge 0}^2$, we will let $\partial_\gamma$ denote the differential operator $\partial^{\gamma_1}_x\partial^{\gamma_2}_y$, that will be applied to smooth functions of $z=(x,y)$, while $D^2f$ denotes the Hessian of a function $f:\mathbb R^n\mapsto \mathbb R$.

\subsection*{Acknowledgments}
Amol Aggarwal heartily thanks Jiaoyang Huang for valuable conversations. The work of Amol Aggarwal was partially supported by a Clay Research Fellowship and a Packard Fellowship for Science and Engineering. The research of  Fabio Toninelli was funded by the Austrian Science Fund (FWF) 10.55776/F1002. For open access purposes, the authors have applied a CC BY public copyright license to any author accepted manuscript version arising from this submission.

\section{Some preliminaries}
	
\label{sec:notations}
\subsection{Height function and domains}

Given $a, b, c > 0$, let $\mathfrak{X}_{a, b, c} \subset \mathbb{R}^2$ denote the  $a \times b \times c$ hexagon with vertices $(0,0)$, $(a, 0)$, $(a+c, c)$, $(a+c, b+c)$, $(c, b+c)$, and $(0, b)$. We call these six corners the SW (southwest), S (south), SE (southeast), NE (northeast), N (north), and NW (northwest) corners of this hexagon, respectively. The sides of the hexagon are labeled W, SE, SW, E, NE, NW, with the obvious convention. Throughout, we assume implicitly that $aN,bN,cN$ are integers, so that the vertices of the hexagon are also vertices of $\mathbb T_N$. 
Whenever $(a,b,c)=(1,1,1)$, we drop them from the notation.

\begin{figure}[h]
 \begin{center} \includegraphics[width=9cm]{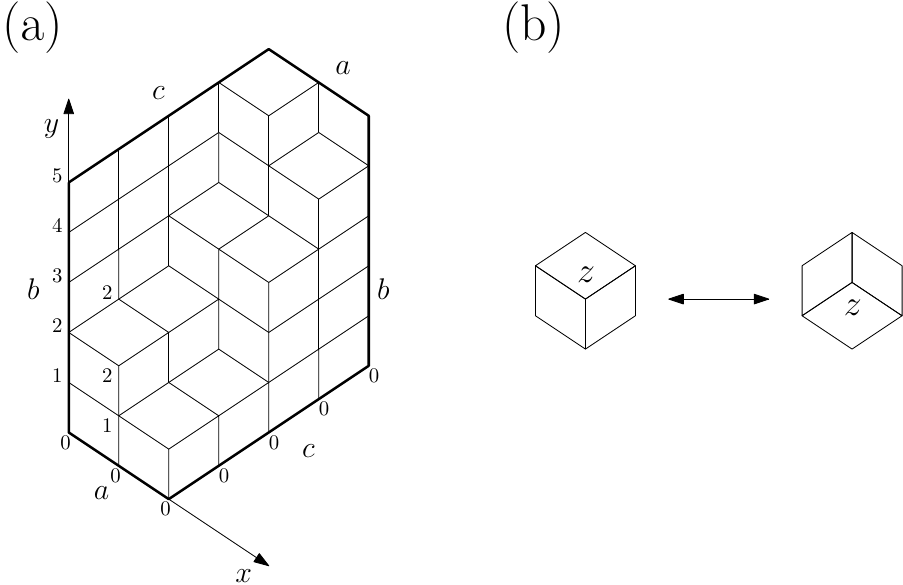}
  \caption{(a) The hexagon $\mathfrak X_{a,b,c}$. The height values are associated to lozenge vertices and height is measured with respect to the horizontal plane, if one interprets the tiling as a monotone stack of cubes of side-length $1/N$. The side-lengths $a,b,c$ are assumed to be multiples of $1/N$. (b) The elementary move $\eta\leftrightarrow \eta^{(z)}$ of the Glauber dynamic, both with rate $1$. $z$ is the common point of the three lozenges with different orientations.}
\label{fig:esagono}
\end{center}
\end{figure}

We adopt the conventions in Figure \ref{fig:esagono} for the coordinates $(x,y)$ in the hexagon. Since the $x$-axis is tilted in our figures, we sometimes refer to  $(1,0)$ as the SE direction, and to  $(1,1)$ as the NE direction.
\begin{definition}[$\ONa$ and $\mathbb P^N_{a,b,c}$]\label{def:Ona}
The set of lozenge tilings of $\mX_{a,b,c}$ is denoted by $\ONa$ and the uniform measure on $\ONa$ is denoted by $\PNa$.
More generally, given a finite tilable domain $D\subset \mathbb R^2$, 
we let $\Omega_D$ denote the set of lozenge tilings of $D$.  
\end{definition}
  To each $\eta\in\Omega_D$ we can associate a height function $H=H_\eta$ (indexed by the vertices of $\mathbb T_N$ contained in $D$), that takes values in $N^{-1}\mathbb Z$.
  More precisely:
\begin{definition}[Discrete height function]
\label{def:H}  Given a simply-connected, tilable domain $D$, a tiling $\eta\in\Omega_D$, and a point $z_0\in D\cap \mathbb T_N$, we define $H_{\eta} : D \cap \mathbb{T}_N \rightarrow \mathbb{R}$ as follows. First set $H_\eta(z_0)=0$; then, $H_{\eta}$ is determined by the relations 
  \begin{eqnarray}
    \label{eq:hfunctx}
    H_\eta(z')-H_\eta(z)=N^{-1}\cdot
    \begin{cases}
      -1+{\bf 1}_{(z,z') \text{ is the edge of a lozenge in $\eta$}} & \text{if } z'-z=(N^{-1},0) \\
      {\bf 1}_{(z,z') \text{  is the edge of a lozenge in $\eta$}} & \text{if } z'-z=(0,N^{-1})
    \end{cases}.
  \end{eqnarray}

\noindent We further extend $H_{\eta}$ to $D$ by linearity. In the specific case $D=\mX_{a,b,c}$ we always choose $z_0=(0,0)$.
\end{definition}
\begin{rem}
  The correspondence between $\eta$ and $H_\eta$ is a bijection. Therefore, we will implicitly identify a lozenge tiling with its height function. We will also often drop the index $\eta$ in the height and write just $H$.

  Note also that $H_\eta$ changes by either $-1/N$, $0$, or $1/N$ between adjacent vertices of $\mathbb T_N$. Lastly, note that when $D=\mX_{a,b,c}$,  $H_\eta$  is  zero on the SW, S, SE corners of the hexagon, it equals $b$ at the NE, N, NW corners and it is affine on each side of $\mX_{a,b,c}$.  See  Figure \ref{fig:esagono}.
\end{rem}

	The configuration space $\Omega_D$ has a natural partial order:
\begin{definition}[Partial order on $\Omega_D$]\label{def:partorder}
Given a simply-connected, tilable domain $D$ and $\eta,\eta'\in\Omega_D$, we write
$\eta\le \eta'$ if and only if $H_{\eta} (z) \le H_{\eta'} (z)$ for all $z \in D$. We let $\eta^\wedge\in \Omega_D$ (resp. $\eta^\vee\in \Omega_D$) denote the maximal (resp. minimal) configuration, that is the unique tiling such that $\eta\le \eta^\wedge$ (resp. $\eta\ge \eta^\vee$) for every $\eta\in\Omega_D$.
\end{definition}
The existence and uniqueness of $\eta^\wedge,\eta^\vee$ follows from the fact that, if $\eta,\eta'\in\Omega_D$, then $z\in D\mapsto \max(H_{\eta}(z),H_{\eta'}(z))$ and $z\in D\mapsto \min(H_{\eta}(z),H_{\eta'}(z))$ are height functions of tilings of $D$.

\begin{definition}\label{def:stochdom}
  Given a simply-connected, tilable domain $D$ and two probability measures $\mu,\nu$ on $\Omega_D$, we say that $\mu$ is stochastically dominated by $\nu$  (in formulas, $\mu\preceq \nu$) if there exists a coupling $(\eta,\eta')$ with $\eta\sim \mu,\eta'\sim\nu$ such that $H_\eta\le H_{\eta'}$ almost surely.
\end{definition}

\begin{figure}[h]
 \begin{center} \includegraphics[width=5cm]{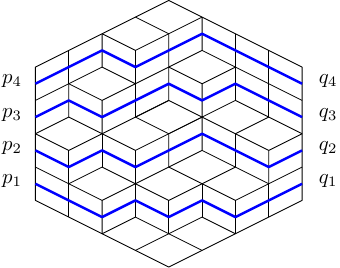}
  \caption{A lozenge tiling of $\mX_{a,b,c}$ with $a=b=c=1,N=4$ and the corresponding height level lines, that form a collection of  non-intersecting Bernoulli paths. Note that level lines are shifted vertically  by $1/(2N)$ with respect to the vertices of the triangular graph.}
\label{fig:paths}
\end{center}
\end{figure}

Lozenge tilings of $\mX_{a,b,c}$ are also in bijection with a collection of $bN$ non-intersecting simple-random walk paths. To set up this  correspondence, we first define level lines of the height function. First of all, given $\eta\in\ONa$, we extend its height function to a function $H_\eta:\mX_{a,b,c}\mapsto  [0,b]$ by extending $H_\eta$ of \Cref{def:H} linearly in each lozenge tile. Then, level lines are defined as follows:

        \begin{definition}[Level lines of the discrete height]\label{ttU}
          Given $a,b,c>0$ and the height function $H$ of a tiling $\eta$ of $\mX_{a,b,c}$ we define the level lines $\mathtt U_{a,b,c}^h:[0,a+c]\cap(1/N)\mathbb Z\mapsto[0,b+c]\cap(1/N)(\mathbb Z+1/2)$ with $h\in[0,b]\cap(1/N)(\mathbb Z+1/2)$ as
          \begin{eqnarray}
            \mathtt U^h_{a,b,c}(x)=y\in [0,b+c]\cap(1/N)(\mathbb Z+1/2):H(x,y)=h.
          \end{eqnarray}
        \end{definition}
Note from Figure \ref{fig:paths} that level lines are disjoint zig-zag paths whose endpoints are spaced by  $1/N$.
Later, we will need to consider collections of non-intersecting paths with different geometries. Therefore, we introduce the following:
 \begin{definition}[Bernoulli paths and Bernoulli path ensembles]
   \label{def:le}
Given $p=(x_p,y_p),q=(x_q,y_q)\in \mathbb T_N+(0,1/(2N))$ such that $x_p<x_q$ and $0 \le y_q-y_p \le x_q-x_p$, we call \emph{Bernoulli path} with endpoints $p,q$ any path $S$ that starts at $p$, ends at $q$ and 
 such that 
 each step of $S$  is either $(1/N,0)$  or $(1/N,1/N)$ (that is, steps are either in the SE or in the  NE direction).
 
    Given $n\in \mathbb N$ and distinct points
    $p_1,\dots,p_n,q_1,\dots,q_n$  of $\mathbb T_N+(0,1/(2N))$ such that all the $p_i$ have
    the same $x-$coordinate  $x_p$ (and the $q_i$ also have all the same $x$-coordinate $x_q>x_p$) and ordered so that the $y-$coordinate of $p_i$ is larger than that of  $p_j$ (resp. the $y$-coordinate of $q_i$
    is larger than that of $q_j$) for $i>j$, we let the \emph{Bernoulli path ensemble} with endpoints
    $p_1,\dots,p_n,q_1,\dots,q_n$ denote the uniform measure on
    $n-$ples     $(S_1,\dots,S_n)$ of non-intersecting (i.e., disjoint) Bernoulli paths, where $S_i$ has endpoints $(p_i,q_i)$. (We assume
    that the points $p_i,q_i,i\le n$ are such that the configuration
    space is nonempty).
  \end{definition}
  Note that the lines $\mathtt U_{a,b,c}^h$ are Bernoulli paths in the sense of \Cref{def:le}.

\subsection{The Glauber dynamics on lozenge tilings}
\label{sec:generalities}

\subsubsection{Definition of the dynamics and general facts}
In this section we  define the Glauber dynamic on lozenge tilings of a finite, simply-connected, tilable domain $D\subset \mathbb R^2$. First we introduce some notation.
Given $\eta\in\Omega_D$, we let $M_\eta$ be the set of points  $z\in D\cap \mathbb T_N$ such  that $z$ belongs to three tiles in $\eta$  of the three possible orientations (see Figure \ref{fig:esagono}(b)). Given $z\in M_\eta$,   we let $\eta^{(z)}$ be the tiling obtained by rotating the three tiles that contain $z$ by an angle $\pi$  around $z$ (see Figure \ref{fig:esagono} (b)).

   \begin{definition}
  [Glauber dynamics]
The lozenge tiling Glauber dynamic  on  $D$
  is the  continuous-time Markov process with state space  $\Omega_D$
 defined by the properties that the only allowed transitions  from  a configuration $\eta\in\Omega_D$ are the configurations $\{\eta^{(z)}\}_{z\in M_\eta}$ and that each such transition occurs according to an exponential clock of rate $1$.
   \end{definition}
   Since the transitions $\eta\to\eta^{(z)}$ and $\eta^{(z)}\to \eta$  have the same rate, the uniform measure on $\Omega_D$ is stationary and reversible.
  It is quickly verified that, when $D=\mX_{a,b,c}$, the Glauber dynamic is irreducible.

  We let $\eta_t$ denote the configuration at time $t$, and we write $H_t:=H_{\eta_t}$ for the associated height function. 

  \begin{definition}
    [$\mu_t^{\eta_0},\mathbf P_{\eta_0},\mathbf P_\nu$]\label{def:Pp}
    Given  $\eta_0\in \Omega_D$ and $t\ge0$, we let $\mathbf P_{\eta_0}$ denote the law of Glauber dynamic with initial condition $\eta_0$, and  $\mu_t^{\eta_0}$  the law of  $\eta_t$.
 If $\eta_0$ is random and sampled according to a distribution $\nu$ on  $\Omega_D$, then we write $\mathbf P_{\nu}$ for the law of the Glauber dynamic.
\end{definition}

We recall the definition of mixing time: for an irreducible Markov chain on a finite state space $\Omega$ and stationary measure $\pi$, one defines for $\epsilon\in(0,1)$
\begin{eqnarray}
  \label{eq:tmix}
  t_{\mix}(\epsilon)=\inf\{t\ge0:\max_{\eta\in\Omega}\|\mu_t^\eta-\pi\|\le \epsilon\}
\end{eqnarray}
where  $\|\cdot\|$ denotes the total variation distance. The following holds \cite[Section 4.5]{LevinPeres} for $0<\epsilon'\le\epsilon<1/2$:
\begin{eqnarray}
  \label{e:submult}
 t_{\mix}(\epsilon')\le t_{\mix}(\epsilon)\left\lceil\frac{|\log \epsilon'|}{|\log 2\epsilon|}
    \right\rceil.
\end{eqnarray}
When we do not specify the value of $\epsilon$, it is meant that $\tmix:=\tmix(1/4)$.

We will use repeatedly a
constrained version of the dynamic.

\begin{definition}[Constrained dynamic]
  \label{def:constrained}
  Fix $H_-,H_+:D\mapsto \mathbb R$ with $H_-\le H_+$ pointwise,
let $\Omega_{H_\pm}$ be the set of tiling configurations in $\Omega_D$ such that  ${H_-}\le H_\eta\le {H_+}$ (we assume $\Omega_{H_\pm}\ne \emptyset$) and fix an initial
condition $\eta_0\in\Omega_{H_\pm}$. The Glauber dynamic
with \emph{floor constraint} $H_-$ and \emph{ceiling constraint}
$H_+$ is the Markov chain on $\Omega_{H_\pm}$ that is defined just like the original Glauber dynamic, except that any update that would violate the constraint $H_-\le H_{\eta_t}\le H_+$ is discarded.
\end{definition}

{It is immediate to check that the uniform measure  $\pi_{H_\pm}$  on $\Omega_{H_\pm}$ is stationary and reversible for the dynamic constrained between $H_\pm$.}
The following stochastic ordering result holds:
\begin{prop}\label{prop:cancouple}
Given two
floors $H_-,H'_-:D\mapsto \mathbb R$ with $ H_-\le H'_-$, two ceilings $H_+,H'_+:D\mapsto \mathbb R$ with $H_+\le H'_+$ and two
initial conditions $\eta_0,\eta'_0\in\Omega_D$ satisfying $ \eta_0\le \eta'_0$, ${\eta_0}\in\Omega_{H_\pm}, {\eta'_0}\in\Omega_{H'_\pm}$, the Glauber  dynamic $\eta_t$ in $D$ constrained between $H_\pm$ and the dynamic $\eta'_t$ constrained between $H'_\pm$
can be coupled in such a way that $\eta_t\le \eta_t'$ for all times $t$, with probability $1$.
{In particular, it holds that $\pi_{H_\pm}\preceq \pi_{H'_\pm}$}.
\end{prop}
\begin{proof}
  Assign an independent Poisson clock of rate $2$ to each $z\in D \cap \mathbb{T}_N$. If the clock at $z$ rings at some time $t$, then flip an independent fair coin and:
  \begin{itemize}
  \item   if $z\not\in M_{\eta_t}$ then do nothing.
  \item if $z\in M_{\eta_t}$, then call $\eta_t^{(z,+)}$ (resp. $\eta_t^{(z,-)}$) the configuration among $\{\eta_t,\eta_t^{(z)}\}$ whose height function is maximal (minimal) at $z$.  If the coin gives head (resp. tail), then replace $\eta_t$ with $\eta^{(z,+)}$ (resp. $\eta^{(z,-)}$) provided that the new configuration satisfies the floor and ceiling constraint, otherwise do nothing.
  \end{itemize}
  If we use the same Poisson clocks and coin flips for all processes with initial conditions $\eta_0\in\Omega_D$ (compatible with their respective floor/ceiling),  it is immediately checked that ordering is preserved.

 { The statement $\pi_{H_\pm}\preceq \pi_{H'_\pm}$ follows in the limit $t\to\infty$, when the laws of  $\eta_t$ and $\eta'_t$  approach their stationary distributions.}
\end{proof}
In fact, the coupling described above provides  a global monotone coupling of  the evolutions started from any $\eta\in \Omega_D$
in such a way that if
$\eta_0\le \eta_0'$, then the corresponding evolutions satisfy
$\eta_t\le \eta'_t$ for every $t$, almost surely.   

The following stochastic domination statement similar to  \Cref{prop:cancouple} holds also for Bernoulli path ensembles (\Cref{def:le}); we omit its proof (as it follows from a very similar argument as used in that of  \Cref{prop:cancouple}; for the case without constraints $V^\pm,W^\pm$; see also \cite[Lemma 18]{CEP}).

  \begin{lem}\label{le:stochdomB}
  	
  	Let $S=(S_1,\dots,S_n)$ and $S'=(S'_1,\dots,S'_n)$ be two 
    Bernoulli path ensembles with  endpoints $(p_1=(x_p,y_{p_1}),\dots,p_n=(x_p,y_{p_n})),(q_1=(x_q,y_{q_1}),\dots,q_n=(x_q,y_{q_n}))$ for $S$,  and $(p'_1=(x_p,y'_{p'_1}),\dots,p'_n=(x_p,y'_{p'_n})),(q'_1=(x_q,y'_{q'_1}),\dots,q'_n=(x_q,y'_{q'_n}))$ for $S'$. Assume that $y_{p_i}\ge y'_{p'_i}$ and $y_{q_i}\ge y'_{q'_i}$
    for every $i=1,\dots,n$. Then, we can couple $S,S'$ in such a way that, almost surely, the $S_i$ is weakly above $S'_i$ for every $i=1,\dots,n$.

	More generally, fix functions $V^\pm_i,W^\pm_i:(x_p,x_q)\mapsto \mathbb R$  such that $V_i^\pm \ge W_i^\pm$. Condition  $S$ to the event that $S_i$ is above $V^-_i$ and below $V^+_i$ for every $i=1,\dots, n$, and condition $S'$ to the event that $S'_i$ is above $W^{-}_i$ and below $W^+_i$ for every $i=1,\dots,n$ (we assume that $(V_i^\pm , W_i^\pm)_{i=1,\dots,n}$ are such that these events are nonempty). Then, we can couple $S,S'$ in such a way that, almost surely, the $S_i$ is weakly above $S'_i$ for every $i=1,\dots,n$.
\end{lem}

        \begin{prop}[{\cite[Theorem 4.3]{CMT}}]
\label{prop:tmixconstrained}
Let $D \subset \mathbb{T}_N$ be a finite, simply-connected tileable domain. For the Glauber dynamic in $D$ constrained between $H_-$ and $H_+$, one has
          \begin{eqnarray}
            \label{eq:tmixconstrained}
            \tmix\le C (N {\rm diam}(D))^2\times (N H_{\max})^2 (\log N)^2, \quad H_{\max}:=\max_{z\in D}(H_{+}(z)-H_{-}(z))
          \end{eqnarray}
          for some universal constant $C$.
        \end{prop}
        Note that  $N {\rm diam}(D)$ and $N H_{\rm max}$ are the diameter of the domain and the maximal height distance
        between floor and ceiling (before the rescaling by $1/N$), respectively.
        As a corollary, the mixing time of the unconstrained dynamic satisfies
\begin{eqnarray}
  \label{eq:apriori}
  \tmix^N(\epsilon)\lesssim  C N^4(\log N)^2\log \frac1\epsilon.
\end{eqnarray}

We need a last general result, that uses the notation in \Cref{def:Pp} and \Cref{def:constrained}:
\begin{prop}[{\cite[Lemma 3.1]{MTTD}}]
  \label{rem:forall}
  Let $D$ be a tilable domain with $|D|$ denoting the cardinality of $D\cap \mathbb T_N$,  let $H_\pm:D  \mapsto\mathbb R $ be as in \Cref{def:constrained} and
  $A\subset\Omega_{H_\pm}$. For any 
fixed  $T \ge 1$, the following holds for the dynamic in $D$ constrained between $H_-$ and $H_+$:
  \begin{eqnarray}
    \label{eq:forall}
    \mathbf P_{\pi_{H_\pm}}(\exists t\le T: \eta_t\in A)\le 8 |D| T \,\pi_{H_\pm}(A).
  \end{eqnarray}
\end{prop}
We will use this to deduce 
 that for the stationary process, if  $\pi_{H_\pm}(A)$ is smaller than any inverse power of $N$, then the probability that there exists a time before $N^C$ (with $C$ any fixed constant) when $A$ occurs is still smaller than any inverse power of $N$.

\subsubsection{The Peres-Winkler censoring inequality}
\label{sec:PW}
The Peres-Winkler censoring inequality \cite{PW} states that, for certain monotone Markov dynamics under certain initial data, censoring (that is, ignoring) updates slows down convergence to equilibrium. 
        We give  in \Cref{prop:PW} the precise statement we need in  the specific case of the lozenge tiling dynamic.

        \begin{definition}[Restricted dynamic] \label{def:restricted} Given the Glauber dynamic in a tilable domain $D$ and a (not necessarily tilable) subset $D'\subset D$, the Glauber dynamic \emph{restricted to $D'$} is the one obtained by ignoring the Poisson clocks in $D\setminus D'$.
        \end{definition}

\begin{definition}
  \label{def:P+} 
  Given $H_+,H_-:\mX_{a,b,c}\mapsto \mathbb R$, with ${H_-}\le {H_+}$ pointwise,
  let $\mathcal P^+[H_-,H_+]$ denote the collection of  probability distributions $\mu$ on $\Omega_{H_\pm}$ that are increasing with respect to  $\pi_{H_\pm}$: 
  \begin{eqnarray}
    \label{eq:P+}
   \eta\ge \eta'\Rightarrow {\mu(\eta)}\ge \mu(\eta').
  \end{eqnarray}
  When $H_+$ and $H_-$ are the height functions of the maximal and minimal configurations in $\Omega^N_{a,b,c}$, respectively, we write $\mathcal P^+$ instead of $\mathcal P^+[H_-,H_+]$.
\end{definition}

 Examples of $\mu\in\mathcal P^+[H_-,H_+]$  are the uniform measure on $\Omega_{H_{\pm}}$ and the  Dirac measure on the maximal configuration in  $\Omega_{H_\pm}$.

  In the following proposition, $H_\pm$ and $\mathcal P^+[H_-,H_+]$ are as in \Cref{def:P+}.
  \begin{prop}[{\cite[Theorem 1.1]{PW}}]
  \label{prop:PW}
  Let $\mu\in\mathcal P^+[H_-,H_+],n\in \mathbb N$ and $t_0=0<t_1<\dots<t_n=T$. For any integer $i \in [1,n]$, let $\Lambda_i \subseteq \mathfrak{X}$ be a subset and $H^\pm_i: \mX\mapsto \mathbb R$ satisfy  $H^-_i\le H^+_i$. Let $(\tilde\eta_t)_{t\ge0}$ be the censored Glauber dynamics restricted to $\Lambda_i$ in the time interval $[t_{i-1},t_i)$ (as in Definition \ref{def:restricted}); with ceiling constraint  $\min(H_i^+,H_+)$ and floor constraint $\max(H_i^-,H_-)$ (as in \Cref{def:constrained}); and with initial condition distributed according to $\mu$.  Let $\tilde \mu_t$ be the law of $\tilde \eta_t$, and $\mu_t$ be the law of the unrestricted and uncensored Glauber dynamics at time $t$, with the same initial condition $\mu$. Then, for any $t\le T$:
  \begin{eqnarray}
     \tilde \mu_t\in\mathcal P^+[H_-,H_+]; \qquad \mu_t\preceq \tilde \mu_t; \qquad  \|\mu_t-\pi_{H_\pm}\|\le \|\tilde\mu_t-\pi_{H_\pm}\|.
  \end{eqnarray}
\end{prop}

\section{Volume-tilted measure and its limit shape}

\label{sec:Vtilted}

\subsection{Volume-tilted measure and limit shape} 
Even though the Glauber dynamic is reversible with respect to the uniform measure $\mathbb P^N_{a,b,c}$, we will frequently make use of a tilted version of this measure.
        Fixing real numbers
        $a, b, c > 0$ and $q\in\mathbb R$, we
         define the volume-tilted measure $\PNqa$ with parameter $q$ in $\mX_{a,b,c}$ as
         \begin{eqnarray}
           \label{eq:voltilt}
\PNqa(\eta)\propto e^{\mq \sum_{z\in\mX_{a,b,c}\cap\mathbb T_N} H_\eta(z)}, \quad \eta\in \ONa.           
         \end{eqnarray}
  Lozenge tilings of $\mX_{a,b,c}$ sampled according to $\PNqa$ exhibit the following limit shape phenomenon.\footnote{The parameters
$(r; \mathsf{q}; \mathsf{N}, \mathsf{T}, a_1)$ in \cite[Theorem 1.6]{DLE} are equal to
$(1; q; b, a+c, -c)$ here. Also, the limit shape in \cite[Theorem 1.6]{DLE} is shifted by $(0,-c)$ with respect to the one here.}
  \begin{thm}[{\cite[Theorem 1.6]{DLE}}]
    Fix $q\in\mathbb R, a,b,c>0$.   There exists a deterministic function $\cH_{q;a,b,c}:\mX_{a,b,c}\mapsto \mathbb R$ such that the height function $H_\eta$ of a random tiling sampled from $\mathbb P^{N,q}_{a,b,c}$ tends in probability to $\cH_{q;a,b,c}$ as $N\to\infty$:
    for every $\varepsilon>0$,
    \begin{eqnarray}
      \label{eq:limtheorem}
      \lim_{N\to \infty}\mathbb P^{N,q}_{a,b,c}(\exists z\in \mX_{a,b,c}:|H_\eta(z)-\cH_{q;a,b,c}(z)|\ge \varepsilon)=0.
    \end{eqnarray}
    \label{th:limitshape}
  \end{thm}
 { \begin{rem}\label{rem:stochdom}
    It is quickly verified that, if $q \ge q'$, then the measure $\mathbb P^{N,q}_{a,b,c}$ stochastically dominates $\mathbb P^{N,q'}_{a,b,c}$. Together with \eqref{eq:limtheorem}, this implies that $\cH_{q;a,b,c}\ge \cH_{q';a,b,c}$ everywhere in $\mX_{a,b,c}$.
  \end{rem}}
  We drop $q$ from the notation if $q = 0$ and $a, b, c$ if
  $a=b=c=1$. Thus, we may write for example
  $\mathcal{H}_{0; a, b, c} = \mathcal{H}_{a, b, c}$,
  $\cH_{q; 1, 1, 1} = \cH_q$,
  $\mathcal{H}_{0; 1, 1, 1} = \mathcal{H}$, and
  $\mathfrak{X}_{1, 1, 1} = \mathfrak{X}$.  For $q=0$,
  \Cref{th:limitshape} is due to \cite{CLP,CKP}. In that
    specific case, we provide a more quantitative concentration bound
    in \Cref{volumeheight} below.

  The limit shape $\cH_{q;a,b,c}$ is rather explicit, though involved (see \Cref{CurvatureDensity}) and  exhibits a \emph{phase separation}
        between liquid and frozen regions. In order to precisely define these regions, as well as the arctic boundary separating them, we start by noting that 
because of \Cref{def:H}, the gradient $\nabla \cH_{q;a,b,c}$ (wherever defined) belongs to the closure $\overline{\mathcal{T}}$ of the triangle
	\begin{flalign}
	\label{e:cT}	\mathcal{T}  = \big\{ (r, s) \in (-1, 0) \times (0, 1):  0 < s + r < 1 \}.
	\end{flalign}
        \begin{definition}\label{def:Adm}[Liquid and frozen region, arctic boundary]
  Given any open set $\Omega \subset \mathbb{R}^2$, we call a function $F : \overline{\Omega} \rightarrow \mathbb{R}$ \emph{admissible} if $F$ is $1$-Lipschitz and satisfies $\nabla F (z) \in \overline{\mathcal{T}}$ for Lebesgue almost every $z \in \Omega$. Let $\Adm (\Omega)$ denote the set of admissible functions $F : \overline{\Omega} \rightarrow \mathbb{R}$ and, for any function $f : \partial \Omega \rightarrow \mathbb{R}$, let $\Adm (\Omega; f)$ denote the set of $F \in \Adm (\Omega)$ satisfying $F |_{\partial \Omega} = f$.  Define the \emph{liquid region} $\mathfrak{L}_F \subseteq \Omega$ associated with an admissible height function $F \in \Adm (\Omega)$ by 
         \begin{flalign*}
         	\mathfrak{L}_F = \big\{ z \in \Omega : \nabla F(z) \in \mathcal{T} \big\},
         \end{flalign*}
         where the condition $\nabla F(z) \in \mathcal{T}$ in particular implies that $\nabla F(z)$ exists.

         In particular, we define the liquid region of $\cH_{q;a,b,c}$ as
         \[\mL_{q;a,b,c}:=\mL_{\cH_{q;a,b,c}}.\] The  boundary
        $\mathfrak A_{q;a,b,c}:=\partial\mL_{q;a,b,c}$ is called the
        \emph{arctic boundary} and the interior of
        $\mX_{a,b,c}\setminus\mL_{q;a,b,c}$ is called the \emph{frozen
        region}.
        \end{definition}

{In fact, for $q=0$, \Cref{th:limitshape} is a particular case of a general limit shape theorem (\cite{CKP}; see also \cite[Theorem 5.15]{RT}) that we recall  here for later convenience.
  \begin{definition}\label{def:harpoon}
    Let $D$ be a bounded open subset of $\mathbb R^2$ with piecewise smooth boundary and $f:\partial D\mapsto \mathbb R$ be such that $\Adm(D;f)\ne\emptyset.$ Let $(D_N)_{N\ge 1}$ be a sequence of tilable domains of $\mathbb T_N$. We write $D_N\stackrel{N\to\infty}\rightharpoonup (D,f)$ if,  as $N\to\infty$, $D_N\to D$ in Hausdorff metric, and at the same time $\lim_{N\to\infty}\max\{|h_N(z)-f(z)|:z\in\partial D_N\}=0$, where $h_N:\partial D_N\mapsto N^{-1}\mathbb Z$ is the boundary height function associated to $D_N$. 
  \end{definition}
  \begin{thm}\label{th:limshapeD}
    Let $D,f$ be as in \Cref{def:harpoon}. There exists $\phi\in\Adm(D;f)$ such that, for any sequence $(D_N)_{N\ge1}$ of tilable domains such that $D_N\stackrel{N\to\infty}\rightharpoonup (D,f)$ and  for every $\varepsilon>0$, it holds
    \begin{eqnarray}
      \lim_{N\to\infty}\pi_{D_N}(\exists z\in D_N:|H_\eta(z)-\phi(z)|>\varepsilon)=0
    \end{eqnarray}
    where $\pi_{D_N}$ is the uniform measure on $\Omega_{D_N}$. 
  \end{thm}
  The function $\phi$ is called \emph{limit shape} (in the domain $D$, with boundary condition $f$) and it satisfies the following monotonicity property and maximum principle:
  \begin{prop}\label{prop:maxprinc} Let $D$ be an open bounded subset of $\mathbb R^2$ with piecewise smooth boundary, let  $f_1,f_2:\partial D\mapsto \mathbb R$ be such that $\Adm(D;f_1),\Adm(D;f_2)\ne\emptyset$ and let $\phi_1,\phi_2$ be the associated limit shapes. Then, it holds:
    \begin{itemize}
    \item If $f_1(z)\le f_2(z)$ for all $z\in\partial D$, then $\phi_1(z)\le\phi_2(z)$ for all $z\in D$;
      
    \item The maximum of $\phi_1-\phi_2$ in $\bar D$ is realized on $\partial D$. 
    \end{itemize}
    
  \end{prop}
  \begin{proof}
    The first statement is proven for instance in \cite[Lemma 2.5]{Amol-Universality}. The second immediately follows from the first, with the fact that $\partial D$ is compact and that, if two boundary heights differ by an additive constant, the limit shapes do as well. 
  \end{proof}
}

           Several facts about the limit shape $\cH_{q;a,b,c}$ and arctic boundary of the volume-tilted measure are proven in  \Cref{Shape}. We begin with the following lemma about the arctic boundary $\mathfrak{A}_{q;a,b,c}$. This statement is known, but we are not aware of a proof written anywhere, so we provide one  in \Cref{Slope2} below.
	\begin{lem} 
		
		\label{axy0} 
		
		The arctic boundary $\mathfrak{A}_{q;a,b,c}$ is a smooth curve that is tangent to each edge of $\mathfrak{X}_{a,b,c}$ at exactly one point. 
	\end{lem}

\begin{figure}
 \begin{center} \includegraphics[width=5cm]{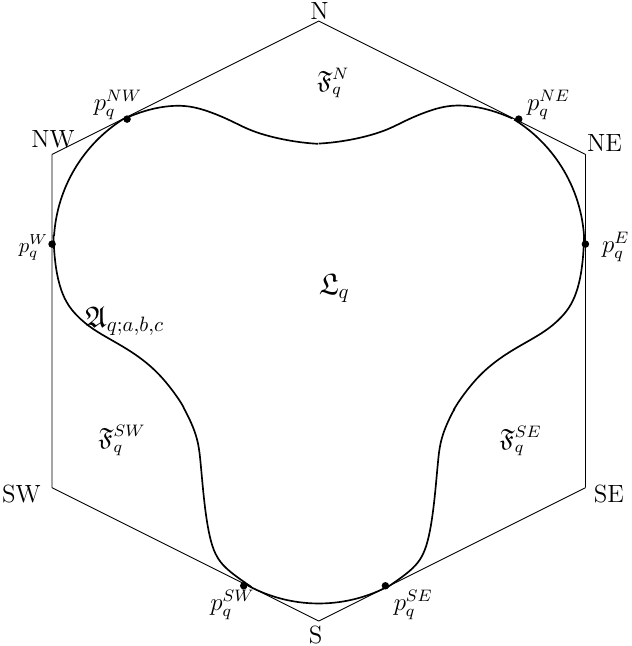}
   \caption{A pictorial view of the liquid region for $\mq>0$, in the case $a=b=c=1$ (so that we drop the indices $a,b,c$). For small enough $\mq$, the liquid region is  convex, for large $q$ it is not. In the three larger frozen regions $\mathfrak F^N_q,\mF^{SW}_q,\mF^{SE}_q$ the limit shape $\mathcal H_\mq$ coincides with the maximal height function. In the three other (smaller) frozen regions, it coincides with the minimal height function.}  \label{fig:liquid}
   
\end{center}
\end{figure}

	In view of \Cref{axy0}, we pose the following definition.
                      \begin{definition}\label{def:pa}Removing the six locations at which $\mathfrak{A}_{q;a,b,c}$ is tangent to $\partial \mathfrak{X}_{a,b,c}$ disconnects $\mathfrak{X}_{a,b,c} \setminus \overline{\mathfrak{L}}_{q;a,b,c}$ into six frozen regions $\mathfrak F^i_{q;a,b,c},i\in\{ N, NE, SE, S, SW, NW \}$ (see \Cref{fig:liquid}).
  We let $p^i_{q;a,b,c},i\in \{ NW, NE,E,  SE, SW, W \}$ denote the six tangency points; see Fig \ref{fig:liquid}, and 
   we let  $p^{SW}_{q;a,b,c}=(x^{SW}_{q;a,b,c},0)$ and $p^W_{q;a,b,c}= (0,y^W_{q;a,b,c})$.

  Moreover, we let $\mathfrak{A}_{q;a,b,c}^{i},i\in \{ N, NE, SE, S, SW, NW \}$  denote the portion of $\mathfrak A_{q;a,b,c}$ adjacent to the frozen region $\mathfrak F^i_{q;a,b,c}$. 
        \end{definition}
Note that                 $0<x^{SW}_{q;a,b,c}<a$ and  $0<y^W_{q;a,b,c}<b$.
        The following lemma determines $\nabla \mathcal{H}_{q;a,b,c}$ in each of the six frozen regions; we establish it in \Cref{Slope2} below.
	
	\begin{lem} 
		
		\label{hlx} 
		
		The following statements hold for any $z \in \mathfrak{X}_{a,b,c} \setminus \overline{\mathfrak{L}}_{q;a,b,c}$.
		\begin{enumerate} 
			\item If $z\in \mathfrak F^S_{q;a,b,c}\cup\mathfrak F^N_{q;a,b,c}$, then $\nabla \mathcal{H}_{q;a,b,c} = (0,0)$. 
			\item If $z\in\mathfrak F^{SW}_{q;a,b,c}\cup\mathfrak F^{NE}_{q;a,b,c}$, then $\nabla \mathcal{H}_{q;a,b,c} = (0, 1)$.
			\item If $z\in\mathfrak F^{SE}_{q;a,b,c}\cup\mathfrak F^{NW}_{q;a,b,c}$, then $\nabla \mathcal{H}_{q;a,b,c} = (-1,1)$. 
		\end{enumerate} 
	\end{lem} 
        For $q=0$, the arctic curve $\mathfrak A_{a,b,c}$
  is  the ellipse inscribed in $\mX_{a,b,c}$ \cite{CKP}. It is therefore natural to expect that   $\mathfrak A_{q;a,b,c}$ is convex for $q$ small. This is confirmed by the next result:
\begin{lem} 
	
\label{convex} 

There exists a constant $\varepsilon > 0$ such that, for any $q \in (-\varepsilon, \varepsilon)$ and $a, b, c \in (1-\varepsilon, 1 + \varepsilon)$, the arctic boundary $\mathfrak{A}_{q;a,b,c}$ is a smooth, convex curve, whose curvature is bounded below by a positive constant $\kappa=\kappa(a,b,c)>0$.
\end{lem}

In order to state more refined properties of the limit shape, we need some more geometric notations.

\begin{definition}
  [Directional distance]
  For any points $z, w \in \mathbb{R}^2$, sets $\mathcal{S}, \mathcal{S}' \subset \mathbb{R}^2$, and vector $V\in
\mathbb R^2$, we set $\dist_{V} (z, w) = \big| (z-w) \cdot V \big|$; $\dist_V (z,\mathcal{S}) = \inf_{s \in \mathcal{S}} \dist_V (z,s)$; and $\dist_V(\mathcal{S}, \mathcal{S}') = \inf_{s\in\mathcal{S}} \inf_{s' \in \mathcal{S}'} \dist_V(s,s')$. 
\end{definition}

	\begin{definition} [$\mathfrak{e}_{z;q},\md_z,\td_z,\te_{z;q}$]
		\label{lz} 
	 Given  $z \in \mathfrak{X}_{a, b, c}$ let $\ell_z $ denote the side of $\partial \mathfrak{X}_{a, b, c}$ closest to $z$ (if two such sides exist, we select one arbitrarily) and $\mathfrak{d}_z = \dist (z, \ell_z)^{1/2}$.  We let $\bm{v} = \bm{v}_z \in \big\{ (1, 0), (-1, 0), (0, 1), (0, -1), (1, 1), (-1, -1) \big\}$ denote the vector
        parallel to $\ell_z$; as for its orientation, we orient the boundary of the hexagon anti-clockwise, and we let $\bm{v}$ be oriented like the side $\ell_z$ that defines it. { We also let $\bm{v}^{\perp} = \bm{v}_z^{\perp}$ denote the vector obtained by a $\pi/2$ clockwise rotation of $\bm{v}$ (so that it is in particular orthogonal to $\bm{v}$, in the $(x,y)$ coordinate system).} Given $q\in\mathbb R$, we let $\bm{w} = \bm{w}_{z;q}$ denote the
unit        vector pointing from $z$ to the closest point of 
        $\mathfrak{A}_{q; {a}, {b},
          {c}}$.  Finally,   we let $\bm{u} = \bm{u}_{z;q}$ denote the
    unit   pointing from  $z$ and tangent to
        $\mathfrak{A}_{q; {a}, {b},
          {c}}$ (note that it is orthogonal to $\bm{w}$). The orientation of $\bm{u}$ is chosen by orienting  $\mathfrak{A}_{q; {a}, {b},
          {c}}$ anti-clockwise. If the points involved in the definition of $\bw,\bu$ are not unique, then choose one arbitrarily.\footnote{The distinction between the directions $\bu$ and $\bw$ will be only important  for points $z$ close to the arctic boundary, in which case the choice is unique.} 

        Moreover, we let (see Fig. \ref{fig:ed})
        \begin{eqnarray}
          \label{eq:eqz}
          \mathfrak{e}_{z;q}:=\dist_{\bm{v}} (z, \mathfrak{A}_{q;a,b,c}).
        \end{eqnarray}
        If $z$ is outside the liquid region, we set $\mathfrak{e}_{z;q}=0$. We further set 
   \begin{flalign}
     \mathtt{d}_z = \max \{ \mathfrak{d}_z, N^{-1/2} \}; \qquad \mathtt{e}_{z;q} = \max \{ \mathfrak{e}_{z;q}, \mathtt{d}_z^{-1/3} N^{-2/3} \}.
     \label{eq:tdte}
   \end{flalign}

  \begin{figure}[h]
 \begin{center} \includegraphics[width=5cm]{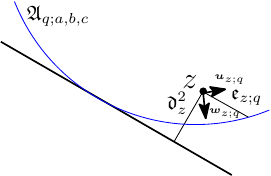}
  \caption{A point $z$ in $\mL_{q;a,b,c}$. The distance from $z$ to $\partial\mX_{a,b,c}$ is $\md_z^2$ (in this case, $\ell_z$ is the southwest side of the hexagon). The distance from $z$ to the arctic boundary $\mA_{q;a,b,c}$ in the direction parallel to $\ell_z$ is $\me_{z;q}$. The unit vector $\bw_{z;q}$ points to the closest point on the arctic boundary, and $\bu_{z;q}$ is perpendicular to it.
}
\label{fig:ed}
\end{center}
\end{figure}

\end{definition}

\begin{rem}
          \label{rem:prop}
          Let $q,a,b,c$ satisfy the conditions of \Cref{convex}. Elementary geometric considerations (using Lemmas \ref{axy0} and \ref{convex}) show that, if $z\in\mL_{q;a,b,c}$, then
          \begin{eqnarray}
            \label{eq:elegeo0}
            \me_{z;q}\lesssim \md_z
          \end{eqnarray}
          and
          \begin{eqnarray}
            \label{eq:elegeo}
             \dist_{\bm{w}_{z;q}} (z; \mathfrak{A}_{q;a,b,c}) \asymp \mathfrak{d}_z \mathfrak{e}_{z;q},\quad \dist_{\bm{u}_{z;q}} (z; \mathfrak{A}_{q;a,b,c}) \asymp (\mathfrak{d}_z \mathfrak{e}_{z;q})^{1/2}.
          \end{eqnarray}
        \end{rem}

\begin{rem}
	
	Let us briefly (and informally) comment on the reason for the definition \eqref{eq:tdte} of $\mathtt{d}_z$ and $\mathtt{e}_{z;q}$. Given $z \in \mathfrak{L}_q$, let $z' \in \mathfrak{A}_q$ denote the closest point to $z$. Then $|z-z'| \le \mathfrak{d}_z^2$ so, if $\mathfrak{d}_z < N^{-1/2}$, then $\mathcal{H}_q (z) - \mathcal{H}_q (z') \le N^{-1}$. Moreover, if $\mathfrak{e}_{z;q} < \mathtt{e}_{z;q}$, then we will show (see Lemma \ref{hde0}) that $|\mathcal{H}_q (z) - \mathcal{H}_q (z')| \asymp \mathfrak{d}_z^{1/2} \mathfrak{e}_{z;q}^{3/2} \lesssim N^{-1}$. Thus, the values of $\mathtt{d}_z$ and $\mathtt{e}_{z;q}$ are essentially the minimal ones to ``detect'' that a point $z$ is in the liquid region, up to an error of $N^{-1}$ in the height function $\mathcal{H}_q (z)$ evaluated at $z$.
	
\end{rem}

 \begin{definition}[Augmented liquid region]\label{def:augmented}
        Given  $\delta > 0$ and $N \ge 1$,
        define the \emph{augmented liquid region} $\mathfrak{L}^{+}_{q; a, b, c} (\delta) = \mathfrak{L}_{q; a, b, c} (\delta; N) \subset \mathfrak{X}_{a, b, c}$ by setting 
		\begin{flalign}
			\label{lqabc} 
			\mathfrak{L}^{+}_{q; a, b, c} (\delta) = \mathfrak{L}_{q; a, b, c} \cup \big\{ z \in \mathfrak{X}_{a, b, c} : \dist (z, \mathfrak{A}_{q; a, b, c}) \le \mathfrak{d}_z^{2/3} N^{\delta-2/3} \big\}.
		\end{flalign}
		
              \end{definition}
              
                \begin{rem}
                  \label{rem:lism}
                  
                  Observe that, if $z \in \mathfrak{L}_{q;a,b,c}^+ (\delta) \setminus  \mathfrak{L}_{q;a,b,c}$, then $\dist_{\bm{v}_z} (z, \mathfrak{A}_{q;a,b,c}) \lesssim \md_z^{-1/3}N^{\delta-2/3}$.
        \end{rem}

The following result gives precise bounds on the probability that, under the measure $\mathbb P^N_{a,b,c}$ (with $q=0$), the height function deviates from the limit shape. Its first part, which is due to \cite{AURT} (see also \cite{LTDMLS}), shows that the height function of a uniformly random tiling of a hexagon satisfies a nearly optimal concentration estimate with error $N^{\delta-1}$; its second part shows that this height function is frozen, in that it coincides with its limit shape, outside of the augmented liquid region. The latter is similar to \cite[Theorem 2.6]{ESC}, though there the augmentation of the liquid region was slightly larger than the one here (namely, there was no factor of $\mathfrak{d}_z^{2/3}$ as in \eqref{lqabc}).
 We provide a proof of \Cref{volumeheight} in \Cref{ProofHeightH}.

	\begin{prop}
		
		\label{volumeheight}
		
		For any real numbers $R > 1$ and $\delta \in (0,1)$, the following holds w.o.p.  with respect to $\mathbb P^N_{a,b,c}$. Fix real numbers $a, b, c \in [R^{-1}, R]$; sample a uniformly random tiling on a hexagon $\mathfrak{X}_{a, b, c}$ and denote its height function by $H : \overline{\mathfrak{X}}_{a, b, c} \rightarrow \mathbb{R}$. 
		
		\begin{enumerate} 
			\item \label{heightabc0} We have $\sup_{z \in \mathfrak{X}_{a, b, c}} \big| H (z) - \mathcal{H}_{a, b, c} (z) \big| \le N^{\delta-1}$. 
			\item \label{heightabc2} We have $H(z) = \mathcal{H}_{a, b, c} (z)$, for each $z \notin \mathfrak{L}_{a, b, c}^+ (\delta)$. 
                        \end{enumerate}
	\end{prop}

        \subsection{Properties of the limit shape and of its level lines}
        
        \label{sec:Taylor1}
        This section collects several results about the behavior of $\cH_{q;a,b,c}$ and of its level lines, defined below, near the arctic boundary.
        As mentioned in Section \ref{sec:Vtilted}, these limit shapes admit an explicit, though slightly intricate, formula in terms of elementary functions (provided in Appendix \ref{CurvatureDensity}). As such, the proofs of the limit shape results in this section will amount to a fairly direct (but tedious) Taylor expansion of those formulas.
        
          The following result, proven in \Cref{Proof00}, indicates that the derivatives of $\mathcal{H}$ exhibit a square root singularity as one nears the arctic boundary (with a correction that depends on $\mathfrak{d})$. Below, we recall the notation $x^{SW}_{q;a,b,c}$ from \Cref{def:pa}. 
	\begin{lem} 
	 \label{hde0} 
	 There exists a constant $C > 1$ such that the following holds, whenever $q \in (-\varepsilon, \varepsilon)$ and $a, b, c \in (1-\varepsilon, 1 + \varepsilon)$, for sufficiently small $\varepsilon \in (0, C^{-1})$. Fix $z_0 = (x_0,y_0) \in \mathfrak{L}_{q;a,b,c}$, and assume that $x_0\ge x^{SW}_{q;a,b,c}$ and that $\ell_{z_0}$ is the southwest edge of $\mathfrak{X}_{a,b,c}$. 
	 
	 \begin{enumerate} 
	 	\item We have that $C^{-1} (\mathfrak{d}_{z_0} \mathfrak{e}_{{z_0};q})^{1/2} \le -\partial_x \mathcal{H}_{q;a,b,c} ({z_0}) \le C (\mathfrak{d}_{z_0} \mathfrak{e}_{{z_0};q})^{1/2}$. 
	 	\item We have that $C^{-1} (\mathfrak{d}_{z_0}^{-1} \mathfrak{e}_{{z_0};q})^{1/2} \le \partial_y \mathcal{H}_{q;a,b,c} ({z_0}) \le C (\mathfrak{d}_{z_0}^{-1} \mathfrak{e}_{{z_0};q})^{1/2}$. 
		\item We have that $C^{-1} \mathfrak{d}_{z_0}^{1/2} \mathfrak{e}_{{z_0};q}^{3/2} \le \mathcal{H}_{q;a,b,c} ({z_0}) \le C \mathfrak{d}_{z_0}^{1/2} \mathfrak{e}_{{z_0};q}^{3/2}$. 
		\item Assume that $z_0 = (x_0,y_0) \in \mathfrak{L}_{q;a,b,c}$ and $\ell_{z_0}$ is the southwest edge of $\mathfrak{X}_{a,b,c}$, but now $ x_0\le x^{SW}_{q;a,b,c}$. Then $C^{-1} (\mathfrak{d}_{z_0} \mathfrak{e}_{{z_0};q})^{1/2} \le -\partial_x \mathcal{H}_{q;a,b,c} ({z_0}) \le C (\mathfrak{d}_{z_0} \mathfrak{e}_{{z_0};q})^{1/2}$ remains valid, but now $C^{-1} \le \partial_y \mathcal{H}_{q;a,b,c} ({z_0}) \le C$ and $C^{-1} \mathfrak{d}_{z_0}^2 \le \mathcal{H}_{q;a,b,c} ({z_0}) \le C \mathfrak{d}_{z_0}^2$. 
	\end{enumerate} 
      \end{lem}

In particular, under the assumptions of \Cref{hde0} we have
    \begin{eqnarray}
      \label{eq:lpa2}
C^{-1}(\cH_{q;a,b,c}({z_0}))^{1/3}\md_{z_0}^{-2/3}\le      \partial_y \cH_{q;a,b,c}({z_0})\le C  (\cH_{q;a,b,c}({z_0}))^{1/3}\md_{z_0}^{-2/3}.
    \end{eqnarray}

The conditions that $x_0 \ge x^{SW}_{q;a,b,c}$ and that $\ell_{z_0}$ is the southwest edge of $\mathfrak{X}_{a,b,c}$ essentially mean that ${z_0}$ is close to $\mathfrak A^S_{q;a,b,c}$, and closer to $p^{SW}_{q;a,b,c}$  than to  $p^{SE}_{q;a,b,c}$.

The next two lemmas compare the limit shapes (with $a=b=c=1$) for  two different values of $q$. They are shown in \Cref{HqEstimate}. 

            \begin{lem} 
            	\label{boundarydistance}
            		There exist constants $\varepsilon > 0$ and $C>1$ such that the following holds for any $-\varepsilon\le q'\le q\le \varepsilon$. Let $z = (x, y) \in \mathfrak{A}_q^S$ and $z' = (x', y) \in \mathfrak{A}_{q'}^S$ be such that $\ell_z = \ell_{z'}$ is the southwest edge of $\mathfrak{X}$. Then, we have $C^{-1} (q-q') \le x-x' \le C(q-q')$.
              \end{lem} 
              \begin{rem}
                \label{rem:bysymmetry}
                              By symmetry, an analogous statement holds, if we assume that $z \in \mathfrak{A}_q^{\mathrm{l}}$ and $z' \in \mathfrak{A}_{q'}^{\mathrm{l}}$ for some $\mathrm{l} \in \{ NW, N, NE, SE, S, SW \}$; that $\ell_z = \ell_{z'}$; and that $\mathfrak{d}_z = \mathfrak{d}_{z'}$, and then replace $x-x'$ by $\dist_{\ell_z} (z, z')$ above.
              \end{rem}

            \begin{lem}
                               \label{qlineq}
                               
                               There exist constants $\varepsilon > 0$ and $C>1$ such that the following holds for any $q,q' \in [-\varepsilon, \varepsilon]$, with $q \ge q'$, and $z \in \mathfrak{X}$. For each $\mathfrak{e}_z \in \{ \mathfrak{e}_{z;q}, \mathfrak{e}_{z;q'} \}$, we have 
                       \begin{multline}
                                    \label{eq:cambioq}
            \mathcal{H}_{q'} (z) + C^{-1} (\mathfrak{d}_z \mathfrak{e}_z)^{1/2} (q-q') - C \mathfrak{d}_z^{1/2} (q-q')^{3/2} \\ 
            \le \mathcal{H}_q (z) \le \mathcal{H}_{q'} (z) + C(\mathfrak{d}_z \mathfrak{e}_{z})^{1/2} (q-q') + C \mathfrak{d}_z^{1/2} (q-q')^{3/2}.                        
\end{multline}

\noindent In addition, we have 
    \begin{flalign}
    	\label{qlineq2} 
             \mathcal{H}_q (z) \ge \mathcal{H}_{q'} (z) + C^{-1} (\mathfrak{d}_z \mathfrak{e}_z)^{1/2} (q-q').                       
\end{flalign}

          \end{lem}

\begin{figure}[h]
 \begin{center} \includegraphics[width=5cm]{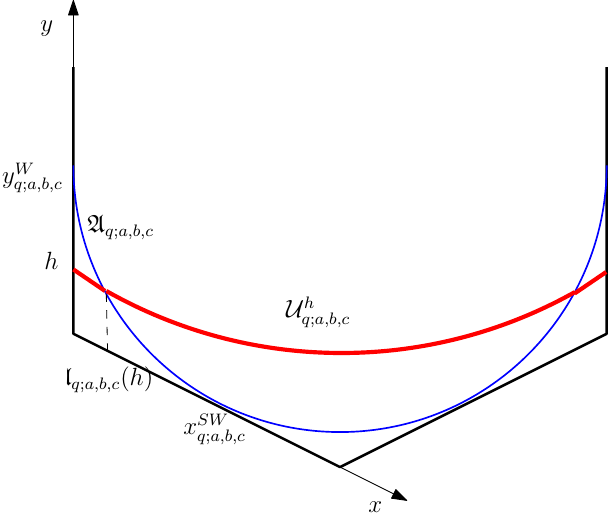}
  \caption{The arctic boundary (in blue) and the level line $\mathcal U^h_{q;a,b,c}$ (in red). For $x$ smaller than $\mathfrak l_{q;a,b,c}$, $\mathcal U^h_{q;a,b,c}(x)$ is constant equal to $h$. }
\label{fig:U}
\end{center}
\end{figure}

	Next, we need a precise control  on the  behavior of the level lines of $\cH_{q;a,b,c}$, as well as of the random level lines of the height function sampled from $\mathbb P^N_{a,b,c}$.
        In analogy with \Cref{ttU}, 
        the following definition (see also Figure \ref{fig:U}) provides level lines $\mathcal{U}$ of $\cH_{q;a,b,c}$,
        and the function $\mathfrak{l}$ describing $\mathfrak A^{SW}_{q;a,b,c}$ (recall \Cref{def:pa}). 	\begin{definition} [Level lines of the macroscopic shape]
		
		\label{ulevelnew} 
		For any $a, b, c >0$, $\mq \in \mathbb{R}$ and $0<h<b$, define $\mathcal{U}_{\mq; a, b, c}^h : [0,a+c]\mapsto [0,b+c] $ by setting
		\begin{equation} 
                  \mathcal{U}_{\mq; a, b, c}^h (x) = y\in [0,b+c]: \mathcal{H}_{\mq; a, b, c}(x,y)=h.
		\end{equation}
 If instead $h=0$ or $h=b$, then given $x$ there may exist multiple $y$ satisfying $ \mathcal{H}_{\mq; a, b, c}(x,y)=h$. In that case,  we define
 \begin{flalign}
   \label{e:geninversa}
   \begin{aligned}
                  \mathcal{U}_{\mq; a, b, c}^0 (x) & = \sup\{y\in [0,b+c]: \mathcal{H}_{\mq; a, b, c}(x,y)=0 \}; \\
                  \mathcal{U}_{\mq; a, b, c}^b (x) & = \inf\{y\in [0,b+c]: \mathcal{H}_{\mq; a, b, c}(x,y)=b\}.
                  \end{aligned}
                \end{flalign}
Also, define
                 $\mathfrak{l}_{\mq; a, b, c} : [0,y^W_{q;a,b,c}] \rightarrow [0,x^{SW}_{q;a,b,c}]$  by letting
                  $\mathfrak{l}_{\mq; a, b, c} (y)$ be the $x$ coordinate of the unique point of $\mathfrak A^{SW}_{q;a,b,c}$ with vertical coordinate $y$ (namely, so that $(\mathfrak{l}_{q;a,b,c}, y) \in \mathfrak{A}^{SW}_{q;a,b,c}$). 
                  Note that  $ \mathfrak{l}_{\mq; a, b, c} (0)=x^{SW}_{q;a,b,c}$ and $\mathfrak{l}_{\mq;a,b,c}(y^W_{q;a,b,c})=0$.
                See Figure \ref{fig:U}. As usual, if either $q=0$ or $a=b=c=1$, we drop these indices from the notation.
	\end{definition}

The next lemma is shown in \Cref{Slope2}. We refer to Figure \ref{fig:U} for a depiction.

	\begin{lem} 
		
		\label{ul} 
		
		There exists a constant $\varepsilon > 0$ such that, for any $q \in (-\varepsilon, \varepsilon)$,  $a, b, c \in (1 - \varepsilon, 1 + \varepsilon)$, and $h \in [0, \varepsilon]$, we have that $\mathcal{U}_{\mq; a, b, c}^{h} (x) = h$, for all $x \in [0, \mathfrak{l}_{\mq;a,b,c} (h)]$. 
		
                \end{lem}

    An immediate consequence of \Cref{hde0} is the following:
    \begin{cor}
      \label{newcor}
      
      { There exists a constant $C>1$ such that the following holds.} Fix $z = (x,y) \in \mathfrak{L}_{q;a,b,c}$, such that $x \ge x^{SW}_{q;a,b,c}$ and $\ell_z$ is the southwest edge of $\mathfrak{X}_{a,b,c}$. {Then, for any nonnegative $h \in [\mathcal{H}_{q;a,b,c} (z) - N^{-1}, \mathcal{H}_{q;a,b,c} (z) + N^{-1}]$ and $\textbf{V} \in \{ \bm{v}_z^{\perp}, \bm{w}_{z;q} \}$, we have $\dist_{\textbf{V}} (z, \mathcal{U}_{q;a,b,c}^h) \le C (\md_z \me_{z;q}^{-1})^{1/2} N^{-1}$.}  
    \end{cor}

    \begin{proof}[Proof of \Cref{newcor}]
    This follows quickly from \Cref{hde0}.
   	\end{proof}

The height concentration statements of \Cref{volumeheight} have the following immediate consequence on concentration of the level lines.
        \begin{cor}\label{cor:conclines}
          Let $H$ be the height function of a uniformly random tiling of   $\mX_{a,b,c}$, and let  $\delta>0$  be a small but $N$-independent constant. The following holds w.o.p.:
          \begin{enumerate}
          \item For every $h\in[0,b]\cap(1/N)(\mathbb Z+1/2)$, $\mathtt U^h_{a,b,c}$ is between $\mathcal U^{h+N^{\delta-1}}_{a,b,c}$ and $\mathcal U^{h
              -N^{\delta-1}}_{a,b,c}$.
            
          \item Let $z=(x,y)\in\mathfrak L_{a,b,c}$
and assume that $x \ge x^{SW}_{a,b,c}$ and that $\ell_z$ is the southwest edge of $\mathfrak{X}_{a,b,c}$. The distance along the direction $\bw_z$ between $z$ and
            $\mathtt U^h_{a,b,c}$ is at most $\td_z^{2/3}h^{-1/3}N^{-1+\delta}$, with $h=\cH_{a,b,c}(z)$.
          \end{enumerate}
        \end{cor}
        \begin{proof}
          The first statement of the corollary is a consequence of \Cref{volumeheight}. Indeed, any point $z$  above $\mathcal U^{h+N^{\delta-1}}_{a,b,c}$ and below  $\mathtt U^h_{a,b,c}$ by definition satisfies $\cH_{a,b,c}(z)\ge h+N^{\delta-1}$ and $H(z)\le h$. Thus, for such a point $z$ we have $|H(z) - \mathcal{H}_{a,b,c}| \ge N^{\delta-1}$. Since this w.o.p. happens nowhere in $\mX_{a,b,c}$, by the first item in \Cref{volumeheight}, no such point $z$ exists, meaning that $\mathtt{U}_{a,b,c}^h$ cannot go above $\mathcal{U}_{a,b,c}^{h+N^{\delta-1}}$. Similarly, $\mathtt U^h_{a,b,c}$ cannot go below $\mathcal U^{h-N^{\delta-1}}_{a,b,c}$.
          
          To show the second, observe from \Cref{hde0} that $h\asymp \md_z^{1/2}\me_{z}^{3/2}$, so  $\td_z^{2/3}h^{-1/3}N^{-1+\delta} \gg (\md_z/\me_{z})^{1/2}N^{\delta/2-1}.$ By definition, $z$ lies on the curve $\mathcal U^h_{a,b,c}$. Since \Cref{newcor} implies that the event $\{ \dist_{\bm{w}_z} (z, \mathtt{U}_{a,b,c}^h) \ge (\mathfrak{d}_z \mathfrak{e}_z^{-1})^{1/2} N^{\delta/2-1} \}$ 
          implies the event that $\mathtt U^h_{a,b,c}$  is not
          between $\mathcal U^{h+N^{-1+\delta/4}}_{a,b,c}$ and $\mathcal U^{h
              -N^{-1+\delta/4}}_{a,b,c}$, one concludes using the first statement of the corollary (with that $\delta$ replaced by $\delta/4$).        \end{proof}

\begin{figure}[h]
 \begin{center} \includegraphics[width=12cm]{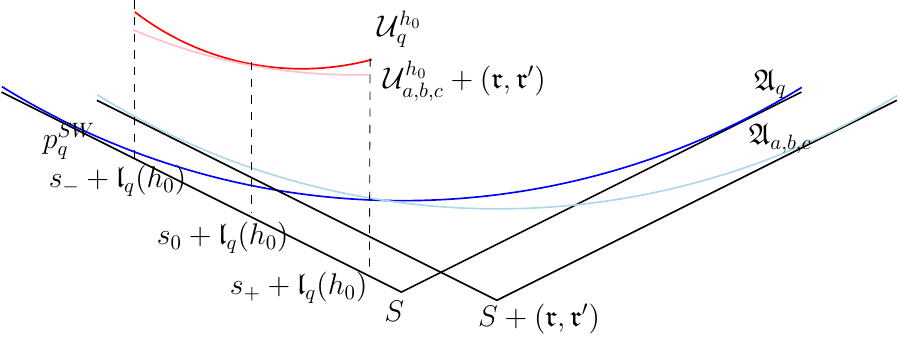}
  \caption{{A schematic illustration of \Cref{prop:improvedconv}. The figure focuses on the south corners of two hexagons (of side-lengths $(1,1,1)$ and $(a,b,c)$ respectively, and translated by $(\frak r,\frak r')$ one with respect to the other). The thick lines in blue/light blue are the corresponding arctic curves, while the thick lines in red/pink are the level lines at the same height $h_0$ of the limit shapes $\cH_q,\cH_{a,b,c}$. If $a,b,c,\frak r,\frak r'$ are chosen as in \Cref{prop:improvedconv}, the level line $\mathcal U^{h_0}_q$ is tangent to $\mathcal U^{h_0}_{a,b,c}+(\frak r,\frak r')$, the former is above the latter and the   curvature of the former is  $q s_0$ larger than that of the latter.}}
\label{fig:Prop327}
\end{center}
\end{figure}

          The next proposition (shown in \Cref{sec:lultima}) will be
          used to compare the level lines of $\cH_q$ and of
          $\cH_{a,b,c}$.  Equation \eqref{eq:meccia2} indicates the
          level lines of the $q$-deformed limit shape are ``more
          convex'' than those of the non-deformed limit shape, upon
          approximately matching (after a shift in the hexagon, given
          by $\mathfrak{r}$ and $\mathfrak{r}'$ in the $x$ and $y$
          coordinates, respectively) their gradients through
          \eqref{eq:meccia} and \eqref{eq:r'}; \eqref{e:primamail}
          allows a specific inequality between these derivatives. See \Cref{fig:Prop327}.  A heuristic explanation for this increased convexity is as follows. Suppose one fixes the values of a level line at two $x$-coordinates, $x^-$ and $x^+$. Then increasing $q$ would, in order to maximize the volume of the associated height function, cause this level line to extend further to the right in the interval $(x^-, x^+)$, thereby increasing its convexity.
          
          Simultaneously enabling
          an approximate match between these gradients but a gain in
          convexity when comparing these limit shapes can require a
          slight change in the dimensions of the hexagon, from $1
          \times 1 \times 1$ for the $q$-deformed limit shape to some $a
          \times b \times c$ for the non-deformed one. This gain \eqref{eq:meccia2}
          in convexity is not imposed too close to the left edge of
          the arctic boundary, where we instead only require that the
          level lines of $\mathcal{H}_q$ are above those of
          $\mathcal{H}_{a,b,c}$, by \eqref{eq:weaker}. The estimate \eqref{heighthabcq} lower bounds the height function $\mathcal{H}_{a,b,c}$ in terms of $\mathcal{H}_q$ (observe that this corresponds to upper bounding the level lines $\mathcal{U}_{a,b,c}$ of $\mathcal{H}_{a,b,c}$ in terms of the $\mathcal{U}_q$, which the estimates \eqref{eq:meccia2} and \eqref{eq:weaker} do for the $h_0$-level line).

\begin{prop}
                \label{prop:improvedconv} 
                
                There exist constants $\varepsilon >0$ and $C>1$ such that the following holds for any $q \in [0, \varepsilon]$. Let $h_0\in (0,\varepsilon]$ and $0 \le s_- \le s_0 \le s_+$ be real numbers such that
                \begin{eqnarray}
                  \label{eq:condizs}
                 \mathfrak{l}_q (0)-\mathfrak{l}_q (h_0) \le  s_0
                  \le s_+ \le 1 + \varepsilon-\mathfrak{l}_q (h_0) \\
                  0 \le s_+ - s_- \le C^{-1} s_0^{1/2}
                  \label{eq:condiz1}.
                \end{eqnarray}
                Then there exist real numbers
                $a,b,c\in [1-C \varepsilon,1+C\varepsilon]$ and $\mathfrak{r}$ such that
                                  \begin{eqnarray}
                                    \label{eq:meccia}
                                \partial_s \cU^{h_0}_{a,b,c}  (\mathfrak{l}_q (h_0) +s_0+\mathfrak r) = \partial_s \cU^{h_0}_{q} (\mathfrak l_q(h_0)+s_0),
                                  \end{eqnarray}
                                  and, defining $\mathfrak{r}'$ via
                                  \begin{eqnarray}
                                    \label{eq:r'}
                                  \mathcal{U}_{a,b,c}^{h_0} (\mathfrak{l}_q(h_0) + s_0 + \mathfrak{r} 
                                  )+   \mathfrak{r}' = \mathcal{U}_q^{h_0} (s_0+\mathfrak l_q(h_0)),  
                                  \end{eqnarray}
                                 \noindent the following statements hold.
                                \begin{enumerate}

                                \item \label{314item1} For all $s \in [ \max(s_-,C^{-1} h_0^{1/2}), s_+]$, we have 
                                  \begin{eqnarray}
                                    \label{eq:meccia2}
                                    \begin{aligned} 
					\mathcal{U}_q^{h_0} & (\mathfrak{l}_q (h_0) + s) - \mathcal{U}_q^{h_0} (\mathfrak{l}_q (h_0) + s_0) \\
					& \ge \mathcal{U}_{a,b,c}^{h_0} (\mathfrak{l}_q (h_0) + s + \mathfrak{r}) - \mathcal{U}_{a,b,c}^{h_0} (\mathfrak{l}_q (h_0) + s_0 + \mathfrak{r}) + C^{-1} qs_0 |s-s_0|^2.
					\end{aligned} 
                                  \end{eqnarray}
                                  
                                \item \label{314item1bis} For all $s \in [s_-,C^{-1} h_0^{1/2}]$, we have 
                                  \begin{eqnarray}
                                    \label{eq:weaker}
                                    \mathcal U_q^{h_0}(\mathfrak l_q(h_0)+s)\ge
                                    \mathcal U^{h_0}_{a,b,c}(\mathfrak{l}_q (h_0) +s+\mathfrak r)+\mathfrak r'.
                                  \end{eqnarray}
                                                                      \item \label{314item2}
                                                                        For all $h\in[0,h_0]$, we have
                                                                         \begin{equation}
                                                                          \label{e:primamail}
                                                                         \cU_q^{h_0}(\mathfrak l_q(h_0)+s_0)-\cU_q^{h}(\mathfrak l_q(h_0)+s_0) \ge \cU_{a,b,c}^{h_0}(\mathfrak{l}_q (h_0) +s_0+\mathfrak r)- \cU_{a,b,c}^{h}(\mathfrak{l}_q (h_0) +s_0+\mathfrak r).
                                  \end{equation}
                                  
                                \item   \label{314item3} For all $s\in[ C^{-1} h_0^{1/2},s_+]$ and $h\in[0,h_0]$, we have
                                  \begin{equation}
                                    \label{eq:primamail'}
                                    \big| \big( \cU_q^{h_0}(\mathfrak l_q(h_0)+s) - \cU_q^{h}(\mathfrak l_q(h_0)+s) \big)
                                   - \big( \mathcal U_{a,b,c}^{h_0}(\mathfrak{l}_q (h_0) + s + \mathfrak r) - \mathcal U_{a,b,c}^{h}(\mathfrak{l}_q (h_0) + s + \mathfrak r) \big) \big|  \le Cqs_+^{5/3} h_0^{2/3}.
                                  \end{equation}
                                  
                                  \item \label{heighthabcq2} For any point $z = (x,y)$ with $x - \mathfrak{l}_q (h_0) \in [s_-, s_+]$ and $\mathcal{H}_q (z) \le h_0$, we have  that 
                                  \begin{flalign}
                                  	\label{heighthabcq} 
                                  	\mathcal{H}_{a,b,c} (x + \mathfrak{r}, y + \mathfrak{r}') - \mathcal{H}_q (x, y) \ge -Cq \mathfrak{d}_z h_0.
                                  \end{flalign}

                      \end{enumerate}
              \end{prop}
              
\subsection{Rescaled limit shape}

\begin{figure}[h]
 \begin{center} \includegraphics[width=5cm]{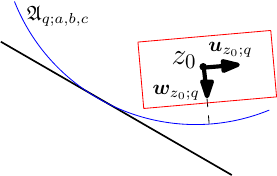}
   \caption{In the definition \eqref{e:funzrescaled} of the rescaled limit shape,
   the coordinates around a point $z_0$ are rescaled differently in direction $\bu_{z_0;q}$ and $\bw_{z_0;q}$. The rectangular region centered at $z_0$ and with sides parallel to $\bu_{z_0;q},\bw_{z_0;q}$ is the cell $\mathfrak Y$ (\Cref{def:bulkcells}, with $w_0$ there equal to $z_0$ here).}
\label{fig:rescaledH}
\end{center}
\end{figure}

In this section, we state a result that compares the macroscopic shapes $\cH_q$ and $\cH_{a,b,c}$,  up to a suitable rescaling. As in \Cref{sec:Taylor1}, the proofs of the statements in this section are based on an elementary Taylor expansion of the formulas for these limit shapes. First, while \Cref{hde0} suggests that second derivatives of $\mathcal{H}$ diverge in the (inverse of the) distance to the arctic boundary, the following proposition (see also \Cref{fig:rescaledH}) provides a rescaling of $\mathcal{H}$ that remains smooth.

\begin{prop}
          \label{prop:formerassumption}          
          The following holds for $\varepsilon > 0$ sufficiently small. Fix ${a}, {b}, {c} \in (1 - \varepsilon, 1 + \varepsilon)$, $q \in (-\varepsilon, \varepsilon)$ and a point $z_0=(x_0,y_0) \in {\mathfrak{L}}_{q; a,b,c}$,
          and assume that $x_0 \ge x^{SW}_{q;a,b,c}$ and that $\ell_{z_0}$ is the southwest edge of $\mathfrak{X}_{a,b,c}$.
          Define the rescaled limit shape as
            \begin{flalign} \label{e:funzrescaled}
                                                                  \widehat{\mathcal{H}}_{q;a,b,c} (\hat z) = \mathfrak{d}_{z_0}^{-1/2} \mathfrak{e}_{z_0;q}^{-3/2}  \cH_{q;a,b,c}(z_0+(\mathfrak{d}_{z_0} \mathfrak{e}_{z_0;q})^{1/2} \bu_{z_0;q} \hat x + (\mathfrak{d}_{z_0} \mathfrak{e}_{z_0;q})\bw_{z_0;q} \hat y ), \quad \hat z=(\hat x,\hat y).
              \end{flalign}

For every $B>1$ there exists a constant $C > 1$ (depending only on $B$) such that  for any differential operator $\partial_{\gamma}$ {(recall the notation in \Cref{Notation})} with $|\gamma| \in [1, B]$ and any $\hat z\in \mathbb{R}^2$ with $|\hat{z}| \le C^{-1}$, we have $|\partial_{\gamma} \widehat{\mathcal{H}}_{q; {a}, {b}, {c}} (\hat{z})| \le C$. 
\end{prop}

Note that the condition $|\hat{z}| \le C^{-1}$, together with \Cref{rem:prop}, guarantees that the argument $z$ of $\cH_{q;a,b,c}$ in \eqref{e:funzrescaled} is within $\mL_{q;a,b,c}$ (and that $\dist_{\bm{u}_{z_0;q}} (z, \mathfrak{A}_{q;a,b,c}) \asymp (\mathfrak{d}_{z_0} \mathfrak{e}_{z_0;q})^{1/2}$ and $\dist_{\bm{w}_{z_0;q}} (z, \mathfrak{A}_{q;a,b,c}) \asymp \mathfrak{d}_{z_0} \mathfrak{e}_{z_0;q}$ for such $\hat{z}$). 

The next proposition will enable us to compare the rescaled limit shapes associated with $\mathcal{H}_q$ and $\mathcal{H}_{a,b,c}$. Through \eqref{e:thereexist}, it essentially states (the plausible notion) that, fixing their  gradients, the latter is ``more concave'' than the former, when $q \ge 0$.
As in \Cref{prop:improvedconv}, simultaneously matching the gradients but  gaining in concavity when comparing these limit shapes can require slightly changing the dimensions of the hexagon, from $1 \times 1 \times 1$ for the $q$-deformed limit shape to some $a \times b \times c$ for the non-deformed one. A similar notion (in the case when $q=0$ and $\mathfrak{d}_{z_0}, \mathfrak{e}_{z_0;q} \asymp 1)$
was witnessed and used in \cite[Section 7]{MTTD}.
              \begin{prop}
                \label{prop:last} 
                There exist constants $\varepsilon > 0$ and $C>1$ 
                such that the following holds whenever $q \in [0,\varepsilon]$. Fix a point $z_0 \in {\mathfrak{L}}_{q}$, and define $\widehat{\mathcal{H}}_q = \widehat{\mathcal{H}}_{q;1,1,1}$ as in \eqref{e:funzrescaled}.
         There exist $a,b,c\in [1-Cq, 1+Cq]$  and $z_1 \in \mathfrak{L}_{ {a}, {b}, {c}}$ such that, defining 
                  \begin{flalign}
                    \label{eq:checkH}
                    \widecheck \cH_{a,b,c}(\hat z)=\mathfrak{d}_{z_0}^{-1/2} \mathfrak{e}_{z_0;q}^{-3/2}  \cH_{a,b,c}(z_1+(\mathfrak{d}_{z_0} \mathfrak{e}_{z_0;q})^{1/2} \bu_{z_0;q} \hat x + (\mathfrak{d}_{z_0} \mathfrak{e}_{z_0;q})\bw_{z_0;q} \hat y ), \quad \hat z=(\hat x,\hat y),
                  \end{flalign}
                  we have 
                  \begin{equation}
                                 \nabla \widehat{\mathcal{H}}_{q} (0,0) = \nabla \widecheck{\mathcal{H}}_{a,b,c} (0,0), 
                  \qquad D^2\widecheck{\mathcal{H}}_{a,b,c}(0,0)-D^2\widehat{\mathcal{H}}_{q}(0,0) \ge C^{-1} \mathfrak{d}_{z_0} q \cdot \Id.
                                   \label{e:thereexist}
                                 \end{equation}
                                 Moreover, for every $B\ge 0$ there exists $C'>1$ such that 
		 for any $\gamma$ with $1 \le |\gamma| \le B$ we have 
		\begin{flalign*}
			\big| \partial_{\gamma} \widehat{\mathcal{H}}_{q} (0,0) - \partial_{\gamma} \widecheck{\mathcal{H}}_{a,b,c} (0,0) \big| \le C' \mathfrak{d}_{z_0} q.
		\end{flalign*}
                
              \end{prop}

                Note that, in the definition of  $\widecheck \cH_{a,b,c}$, the quantities $\md_{z_0},\me_{z_0;q}$ are intentionally computed at $z_0$ and not at $z_1$. The reason is that, otherwise, the equality $\nabla \widehat{\mathcal{H}}_{q} (0,0) = \nabla \widecheck{\mathcal{H}}_{a,b,c} (0,0)$ would not imply equality $\nabla\cH_{q}(z_0)=\nabla\cH_{a,b,c}(z_1)$  of the gradients of the unrescaled functions.

\section{Proof of Theorem \ref{thm:main}: main steps}
\label{sec:mainsteps}
\subsection{Proof of Theorem \ref{thm:main}}

We aim at proving $\tmix^N(1/4)\le C(\delta)N^{2+\delta}$ for every $\delta>0$.
Let
$t^{\mq,N}_{\mix}(\epsilon)$ be the mixing time {(cf. \eqref{eq:tmix})} of the Glauber dynamic in $\mX$,
constrained between floor $\mathcal{H}_{{-q}}$ and ceiling
$\mathcal{H}_{q}$, in the sense of \Cref{def:constrained}.
The first step, shown in Section \ref{sec:thm53}, is to prove:
\begin{thm}
  \label{thm:therecursion}
  
  For any real numbers $q_0>0, \delta > 0$, there exist constants $r = r(q_0) \in (0, 1)$ and $\mC = \mC(q_0, \delta) > 1$ such that the following holds. 
 Define
  \[
     i_{max}=\left\lfloor\frac 9{10}\frac{\log N}{\log(1/(1-r))}\right\rfloor
    \] and
  the sequence
\[
  N_i=\lfloor N (1-r)^i\rfloor,\quad 0\le i\le i_{max}\] that
decreases from $N_0=N$ to $N_{i_{max}}\asymp N^{1/10}$. Define also
the sequence \[
\epsilon_i=10^{-i}N^{-100},\quad 0\le i\le i_{max}
.\]
 For every $i\le i_{max}-1$, we have
   \begin{flalign}\label{eq:recursion}
     \tmix^{N_i}(\epsilon_i) \le \mC{N_i}^{2+\delta} + t_{\rm mix}^{N_{i+1}}(\epsilon_{i+1})+t^{\mq_0,N_i}_{\rm mix}\left(\frac   {\epsilon_i}4\right).
   \end{flalign}
 \end{thm}
 
{ \begin{rem}\label{rem:whoisL}
   Note that $\epsilon_i$ 
decreases from $\epsilon_0= N^{-100}$ to
$\epsilon_{i_{max}}=N^{-L}$ for some $L=L(q_0)>100$. 
\end{rem}}

{For the following result, for $i\le i_{max}$ and $t>0$, let $H^{N_i}_t$ denote
the height at time $t$ of the tiling evolved according to the Glauber dynamic in $\mX$, with $N$ equal $N_i$.}
 The main estimate leading to Theorem \ref{thm:therecursion}, shown in Section
 \ref{ProofHqt}, is:
\begin{thm}
  \label{thm:step1}
  
  Adopt the notation of \Cref{thm:therecursion}. For every $i\le i_{max}-1$ and for every initial condition, the following holds off of an event of probability at most  $N^{-2L}+2\epsilon_{i+1}$ (with $L$ the same constant as in  \Cref{rem:whoisL}):
  for all times $t$ in the interval
  \begin{eqnarray}
    \label{eq:intervallo}
I^N_i=    [\mC N_i^{2+\delta}+t_{\rm mix}^{N_{i+1}}(\epsilon_{i+1}),N^5],
  \end{eqnarray}
 we have that { $\mathcal{H}_{-q_0} (z) \le H^{N_i}_t (z) \le \mathcal{H}_{q_0} (z)$} for all $z \in \mathfrak{X}$.
 
\end{thm}
{In \Cref{thm:step1} and in the following, $N^5$ could have been  replaced by any power of $N$ with exponent larger than $2$, since we want this time to be much larger than our target for the mixing time of the Glauber dynamic.}
The next step, shown in \Cref{Proofq0n}, is:
\begin{thm}
  \label{thm:fN}
  
  There exists $\varepsilon > 0$ so that the following holds. For any real numbers $q_0 \in (0, \varepsilon)$, $\delta > 0$, and $\mathsf{L} > 1$, there is a {constant $C > 1$} such that
  \begin{eqnarray}
  t^{q_0,N}_{\rm mix}(N^{-\mathsf L})\le C N^{2+\delta}   .
  \end{eqnarray}
\end{thm}

Putting together the two results, we immediately have Theorem \ref{thm:main}.
\begin{proof}[Proof of  Theorem \ref{thm:main}, assuming Theorems \ref{thm:therecursion} and \ref{thm:fN}]
  
For $N_{i_{max}}\asymp N^{1/10}$ we know from \eqref{eq:apriori} that
{\begin{eqnarray}
  \label{eq:iniz}
  \tmix^{N_{i_{max}}}(\epsilon_{i_{max}})\lesssim (N_{i_{max}})^4(\log N_{i_{max}})^2|\log(\epsilon_{i_{max}})|\lesssim \sqrt N\le N.
\end{eqnarray}}
 On the other hand, thanks to \Cref{rem:whoisL} there exists a constant $L=L(q_0) > 0$ such that $\epsilon_i / 4 \ge  N_i^{-L}$. Hence, 
\[
  t^{ q_0,N_i}_{\rm mix}\left(\frac {\epsilon_i}4\right)\le t^{q_0,N_i}_{\rm mix}(N_i^{-L})\le C N_i^{2+\delta/2}
\]
where the second inequality holds by Theorem \ref{thm:fN}. Then, \eqref{eq:recursion} yields 
\[
 \tmix^{N_i}(\epsilon_i)\le (C+\mC) N_i^{2+\delta/2}+ t_{\rm mix}^{N_{i+1}}(\epsilon_{i+1}).
  \]
  Solving the iteration and using the final condition \eqref{eq:iniz}, we get the desired bound
  \begin{eqnarray*}
    \tmix^{N}(1/4)=\tmix^{N_0}(1/4)\le  \tmix^{N_0}(\epsilon_0)\lesssim N^{2+\delta}+N\lesssim N^{2+\delta},
  \end{eqnarray*}
  
  \noindent thereby establishing the theorem.
\end{proof}

The rest of the proof is organized as follows. In \Cref{sec:thm54} we will prove Theorem \ref{thm:step1}. In Section
\ref{sec:thm53} we show how this implies Theorem
\ref{thm:therecursion}. Finally, the proof
of Theorem \ref{thm:fN} is contained in Sections
\ref{TimeBarrier}.

\subsection{Proof of Theorem \ref{thm:step1}}
\label{sec:thm54}

\subsubsection{Perturbing the boundary condition}
The first step is to prove that, for any $\delta,\omega>0$, after time
of order $N^{2+\delta}$, the height function is w.o.p.  within $\omega$ from the macroscopic height. If the limit
shape had no frozen regions, this would follow from the main result in
\cite{LTDMLS}. In our case this assumption fails, but the problem can
be resolved by a suitable perturbation of the boundary height.  Let
$f:\partial \mathfrak X \mapsto \mathbb R$ denote the restriction to $\partial \mX$ of the height
function of any lozenge tiling of $\mX$. Note that $f$ is affine
along each side of the hexagon.

\begin{prop}
  \label{thm:perturb}
   Assume that there exists a sequence of  boundary  height functions $f^{(n)} : \partial \mathfrak{X} \rightarrow \mathbb{R}$, for each $n \ge 1$, such that $\Adm(\mX;f^{(n)})$ is not empty (recall \Cref{def:Adm}) and
  \begin{enumerate}
  \item [(A)] $\lim_{n\to\infty}\sup_{z\in \partial\mathfrak X
    }|f^{(n)}(z)-f(z)|=0$
    
  \item [(B)] the limit shape  $\mathcal H^{(n)}$ in $\mathfrak X$ with boundary condition $f^{(n)}$ {(cf. \Cref{th:limshapeD})}
    has the property that $\nabla \mathcal H^{(n)}\in \mathcal T$ and $\dist(\nabla \mathcal H^{(n)},\partial\mathcal T)\ge 1/(2n)$ hold almost everywhere in $\mathfrak{X}$, where $\mathcal T$ was defined in \eqref{e:cT}.
  \end{enumerate}
  Then, for the Glauber dynamic in $\mathfrak X$, for every $\delta,L,\omega>0$ there {exists $C<\infty$} such that, for every $N$,
  \begin{eqnarray}
    \label{eq:fromCMP}
    \max_{\eta\in \Omega^N}\mathbf P_\eta\Big(\exists t: C N^{2+\delta}\le t\le N^{L}, z\in \mathfrak X, |H_t(z)-\mathcal H(z)|\ge \omega\Big)\le CN^{-L}.
  \end{eqnarray}
\end{prop}

\begin{proof}
  Let
  $\mathfrak X'$ denote the hexagon with the same center as $\mathfrak X$ and  with side-lengths  $1-\omega/4 $ instead of $1$, and let $\psi^{(n)}$ be the restriction of $\cH^{(n)}$ to $\partial \mathfrak X'$. Let $(D_N)_{N\ge1}$ be a sequence of tilable domains of $\mathbb T_N$ such that $D\stackrel{N\to\infty}\rightharpoonup (\mX',\psi^{(n)})$ (cf. \Cref{def:harpoon}). Observe that the limit shape in $\mX'$ with boundary height $\psi^{(n)}$  is just the restriction of $\cH^{(n)}$ to $\mX'$.

 Consider now the Glauber dynamic $(\eta'_t)_{t\ge 0}$ in $D_N$. Thanks to assumption (B) above about the gradient of $\mathcal H^{(n)}$, and since  $\mathfrak X'$
is at positive distance  $\omega/4$ from $\partial \mathfrak X$,
it follows from \cite[Theorem 2.5]{MTTD} that its mixing time is
$O(N^{2+\delta'})$ for every $\delta'>0$. Using the equilibrium height
fluctuation statement of \cite[Theorem 3.7 and Corollary 3.8]{MTTD}, and the
submultiplicativity property \eqref{e:submult}, together with
\Cref{rem:forall}, we conclude that, if $\delta>\delta'$ then, w.o.p.,
\begin{eqnarray}\label{eq:lp}
  H_{\eta'_t}(z)\le \mathcal H^{(n)}(z)+\omega/4\quad \forall t\in[N^{2+\delta},N^L], \quad \forall z\in D_N,
\end{eqnarray}
 for any choice of the initial condition $\eta'_0\in\Omega_{D_N}$.  

Choose $n$ large
enough so that
$\sup\{|f^{(n)}(z)-f(z)|:z\in\partial
\mX\}\le \omega/4$, { which is possible by assumption (A) above. By the first item of \Cref{prop:maxprinc}}, this implies that
$|\mathcal H^{(n)}-\mathcal H|\le \omega/4$ everywhere in
$\mathfrak X$.  Note that, since the height is $1$-Lipschitz, the
maximal height function in the hexagon $\mathfrak X$, restricted to
$\partial \mathfrak X'$, is smaller than
$\mathcal H^{(n)}|_{\partial \mathfrak X'}+\omega/2$. Then, it follows from \eqref{eq:lp} together with \Cref{prop:cancouple} that the following holds w.o.p., for the
Glauber {dynamic $(\eta_t)_{t\ge0}$} in $\mathfrak X$: for any initial condition $\eta_0\in\Omega_{\mathfrak{X}}$,
\begin{eqnarray}
  \label{eq:lp2}
  H_t(z):=H_{\eta_t}(z)\le \mathcal H(z)+\omega, \quad \forall t\in [N^{2+\delta},N^L],\quad \forall z\in D_N.
\end{eqnarray}
 On the other
hand, in $\mathfrak X\setminus D_N$, {$H_{t}$}  is
deterministically smaller than $\mathcal H+\omega/2$, again because
it is $1$-Lipschitz and equals $\cH$ on $\partial\mX$ (and $\dist (\partial \mathfrak{X}, \partial \mathfrak{X}') \le \omega/2$). Altogether, {$H_{t}$} is w.o.p. smaller than $\mathcal H+\omega$ everywhere in
$\mathfrak X$, in the whole time interval and for any initial
condition. The statement that the height is larger than
$\mathcal H-\omega$ is proven analogously.
\end{proof}

 The next important observation is that such sequence
 $\{ f^{(n)} \}_{n \geqslant 1}$ indeed exists:
 \begin{prop}
   \label{thm:approx}
   There exists a sequence of boundary height functions $\{f^{(n)}\}_{n\ge1}$ on  $\partial\mathfrak{X}$ that satisfies conditions (A) and (B) of Proposition \ref{thm:perturb}.
 \end{prop}
 \begin{proof}
    We show how a small perturbation of the
   hexagon boundary height function leads to its limit shape having no
   frozen regions. We define
   $\phi^{(n)} : \mathfrak{X} \mapsto \mathbb{R}$ as
\[ \phi^{(n)} (x, y) = \left\{\begin{array}{lll}
                             y\big(1-\frac1n\big)-\frac{x}  n & \text{if} & 0 \leqslant x \leqslant 1, 0\le y\le 1\\
                             & & \\
     -x\big(1-\frac2n\big)  + y\big(1-\frac1n\big)+ \big(1-\frac3n\big) & \text{if} &
                                                                      1 \leqslant x \leqslant 2, y \leqslant x\\
                             & & \\
    \frac {2y}n- \frac{x }{n}+\big(1-\frac3n)&\text{if} &
     1 \leqslant y \leqslant 2, y \geqslant x
                           \end{array}\right. \]
                       and $f^{(n)}:\partial \mathfrak{X}\mapsto \mathbb R$ as the restriction of $\phi^{(n)}$ to $\partial\mX$. Note that $f^{(n)}$ is an admissible boundary height, that is, $\Adm(\mX;f^{(n)})$ in \Cref{def:Adm} is nonempty, since the gradient of $\phi^{(n)}$ belongs to $\overline{\mathcal T}$ (for sufficiently large $n$) wherever defined.
Moreover, $\phi^{(n)}$ clearly satisfies condition (A) of Proposition \ref{thm:perturb}.

Call $\mathcal{H}^{(n)}$ the limit shape  associated to  $f^{(n)}$. To prove that $f^{(n)}$ also satisfies condition (B), we need to show that for every $r> 0$ one has
\begin{multline}\label{eq:lp3}
  {r/n} \leqslant \mathcal{H}^{(n)} (x, y + r) - \mathcal{H}^{(n)} (x, y)
  \leqslant (1 - 1/n) r, \\{r/n} \leqslant \mathcal{H}^{(n)} (x - r, y) -
  \mathcal{H}^{(n)} (x, y) \leqslant (1 - 1/n) r  \\
  {r/n} \leqslant \mathcal{H}^{(n)} (x + r, y + r) - \mathcal{H}^{(n)} (x,
  y) \leqslant (1 - 1/n) r. 
\end{multline}
whenever the argument of $\cH^{(n)}$ belongs to $\mX$.
We only discuss the proof of the upper bound in the first line, since those of the others are
entirely analogous.

Fix $r > 0$, and let $(x, y)\in \mX$ be such that
$\mathcal{H}^{(n)} (x, y + r) - \mathcal{H}^{(n)} (x, y)$ is
maximal. We first claim that either $(x, y + r)$
or $(x, y)$ can be assumed to be on the boundary of the hexagon
$\mathfrak X$.  Indeed, note that both $u(\cdot):=\cH^{(n)}(\cdot)$
and $v(\cdot):=\cH^{(n)}(\cdot+(0,r))$ are { limit shapes of the tiling model} in $\mX$ and $\mX-(0,r)$, respectively. In
particular, they are both minimizers in $\mX\cap(\mX-(0,r))$ (with
different boundary data). Let $(x,y),(x,y+r)\in\mX$, so that
$(x,y)\in\mX\cap(\mX-(0,r))$. Then,
        \begin{multline}
          \cH^{(n)}(x,y+r)-\cH^{(n)}(x,y)=v(x,y)-u(x,y)\le v(z_1,z_2)-u(z_1,z_2)\\=\cH^{(n)}(z_1,z_2+r)-\cH^{(n)}(z_1,z_2)
        \end{multline}
        for some $(z_1,z_2)\in\partial\left(\mX\cap(\mX-(0,r))\right)$, by {the second item of \Cref{prop:maxprinc}}. That is, either $(z_1,z_2)$ or $(z_1,z_2+r)$ belongs to $\partial\mX$ and the claim follows.

        By symmetry, to prove the upper bound in the first of \eqref{eq:lp3} we may restrict to the case when $(x, y)$ is on
        either the SW or SE side of $\mathfrak{X}$ (it could be also
        on a vertical side, in which case the bound
        $\mathcal{H}^{(n)} (x, y + r) - \mathcal{H}^{(n)} (x, y)
        =\phi^{(n)}(x,y+r)- \phi^{(n)}(x,y)\leqslant (1 - 1/n) r$
        holds with equality); let us assume it is on the SE side, as
        the proof is very similar if it is on the SW one.  Let $G$ be
        the affine function on $\mathfrak{X}$ of slope
        $(-1+2/n, 1-1/n)$, such that
        $G(x, y) = \mathcal H^{(n)}(x, y)$.  It is quickly verified
 that $G\ge f^{(n)}$ on
        $\partial\mX$, with equality on the whole SE and E sides.
        From {\Cref{prop:maxprinc}}, it follows that
        $G \ge \mathcal H^{(n)}$ holds also in the interior of
        $\mathfrak{X}$. In particular,
\[\mathcal H^{(n)} (x, y+r) -\mathcal H^{(n)} (x, y) =\mathcal H^{(n)}
 (x, y+r) - G(x, y) \le G(x, y+r) - G(x, y) =
(1-1/n) r,\] from which the above upper bound follows.   
 \end{proof}

 \subsubsection{Reducing to $\mathcal H_{-\mq_0}\le H_t\le \mathcal H_{\mq_0}$ with $\mq_0$ small}
 
 \label{ProofHqt} 

 In this section we prove Theorem \ref{thm:step1}. {Throughout we will implicitly assume that $q_0$ is sufficiently small (so that all of the results in \Cref{sec:Vtilted} hold for $\mathcal{H}_{q_0}$). This is allowed, since  the statement we want to prove is weaker for larger $q_0$.}
 Recall that we are considering the Glauber dynamic in $\mathfrak X$, with the original $N$ equal to $N_i$ (so that the tile size is $1/N_i$) for $i\le i_{max}$, so that $N_i \in [N^{1/10}, N]$. The statements of Proposition \ref{prop:regioneliquida}, Proposition \ref{prop:astuzia} and Lemma \ref{lem:astuzia} below refer to that situation, even if not explicitly stated. As before, $H_t^{N_i}$ denotes the height function of the configuration at time $t$; we abbreviate $H_t^{N_i} = H_t$ throughout this section. 

         By symmetry and monotonicity of the Glauber dynamics, it is enough to start from the maximal configuration $\eta^\wedge$ (see \Cref{def:partorder})  and to prove that $H_t\le\mathcal H_{\mq_0}$ off of an event of probability at most
         \begin{equation}
           \label{eq:xiNi}
\xi^N_i:=\frac12 N^{-2L}+3\epsilon_{i+1},
         \end{equation}
for all times in the interval $I^N_i$ defined in \eqref{eq:intervallo}.  Fix
$\omega>0$. Due to Proposition \ref{thm:perturb} and Proposition
\ref{thm:approx}, for 
$  t\in I^N_i$,
the height function $H_t$  is within $\omega$ from
$\mathcal{H}$, except with probability
$(1/4)N^{-2L}$. Note that $H_t\le \mathcal H+\omega$ does not  imply that $H_t \le \mathcal{H}_{\mq_0}$
everywhere in $\mathfrak X$, even for small $\omega$ and large $\mq_0$ (since $\mathcal{H} = \mathcal{H}_{q_0}$ in the regions where both of these height functions are frozen).

Recall that the frozen region of $\mathcal H_{\mq_0}$ is composed of {
six connected regions labeled $\mathfrak F^i_{q_0},i\in\{ N, NE, SE, S, SW, NW \}$.}
 In $\mathfrak F^N_{q_0}, \mathfrak F^{SW}_{q_0}, \mathfrak F^{SE}_{q_0}$,  the inequality
$H_t \le \mathcal{H}_{\mq_0}$ holds since $\mathcal{H}_{\mq_0}$ is
maximal there by \Cref{hlx}. Call  {$\mathfrak T_{\mq_0}=\mathfrak X\setminus (\mathfrak F^N_{q_0}\cup \mF^{SE}_{q_0}\cup \mF^{SW}_{q_0})$} (see \Cref{fig:liquid2}): it includes the liquid region $\mL_{q_0}$ of $\mathcal{H}_{\mq_0}$ and
the three frozen  $\mF^S_{q_0},\mF^{NE}_{q_0},\mF^{NW}_{q_0}$, where $\mathcal{H}_{\mq_0}$ is minimal
instead (also by \Cref{hlx}). For $j \in \{ S, NE, NW \}$, let $\mathfrak B^j_{q_0} = \{ z \in \mathfrak{X} : \dist (z, \mF^j_{q_0}) < c_{\omega} \}$ denote the $c_\omega$-neighborhood of 
$\mF^j_{q_0}$, for some $c_\omega>0$ to be chosen later. We first claim:
\begin{prop}
  \label{prop:regioneliquida}
  
  There exists a real number $\delta_0 = \delta_0 (q_0)>0$ such that the following holds. For $\omega < \delta_0$, there exists a constant $c_{\omega} > 0$ such that $\lim_{\omega \rightarrow 0} c_{\omega} = 0$ and the inequality
           \[
             \mathcal H+\omega\le \mathcal H_{\mq_0}, \qquad \text{holds on $\mathfrak T_{\mq_0}\setminus(\cup_{j\in\{S,NE,NW\}}\mathfrak B^j_{q_0})$}. 
           \]
           In particular, off of an event of probability at most $(1/4)N^{-2L}$, we have that $H_t \le \mathcal{H}_{\mq_0}$ in $\mathfrak T_{\mq_0}\setminus(\cup_{j\in\{S,NE,NW\}}\mathfrak B^j_{q_0})$, for every $t\in I^N_i$ (recall \eqref{eq:intervallo}).
         \end{prop}

   The next result deals  with the three regions   $\mathfrak B^j_{q_0},j\in\{S,NE,NW\}$.
           \begin{prop}\label{prop:astuzia}
             Off of an event of probability at most $\epsilon_{i+1}+N^{-2L}/12$, we have that $H_t \le \mathcal{H}_{q_0}$ on $\cup_{j\in\{S,NE,NW\}}\mathfrak B^j_{q_0}$, for each $t \in I_i^N$ (recall \eqref{eq:intervallo}).             
             \label{prop:regionesolida}
           \end{prop}

	\begin{proof}[Proof of Theorem \ref{thm:step1}]
		By Propositions \ref{prop:regioneliquida} and \ref{prop:regionesolida}, with a union bound, we deduce that $H_t^{N_i} (z)=H_t(z) \le \mathcal{H}_{q_0} (z)$ for all $z \in \mathfrak{X}$ and $t \in I_i^N$, with probability at least $1-(N^{-2L}/3 + \epsilon_{i+1})$. By symmetry (or entirely analogous reasoning), we also have  $H_t^{N_i} (z) \ge \mathcal{H}_{-q_0} (z)$ for all such $(z,t)$, with the same probability. Thus, the theorem follows from a union bound. 
	\end{proof}

  \begin{proof}[Proof of Proposition \ref{prop:regioneliquida}]
  	
    Throughout this proof, we refer to Figure \ref{fig:liquid2} and use
    the notations depicted there.  Since $q_0 > 0$, we have by
    {\Cref{rem:stochdom}} that
    $\mathcal H_{\mq_0} \ge \mathcal{H}$ everywhere on $\mathfrak{X}$;
by \Cref{qlineq} that
    $\mathcal{H}_{q_0} > \mathcal{H}$ on
    $\mathfrak{L}_{q_0} \setminus (\partial\mF^{NE}_{q_0} \cup
    \partial \mF^{NW}_{q_0} \cup \partial \mF^{S}_{q_0})$;
     {and by \Cref{hlx}} that
    $\mathcal{H}_{q_0} = \mathcal{H}$ on
    $\partial \mF^{NE}_{q_0} \cup \partial \mF^{NW}_{q_0} \cup
    \partial \mF^{S}_{q_0}$.
    By continuity of {$z\mapsto \cH(z)$ and
      $z\mapsto \mathcal{H}_{q_0}(z)$}, 
    it follows that for sufficiently small
    $\omega > 0$ there exists a constant $c_\omega>0$ (which we can
    take to satisfy $\lim_{\omega \rightarrow 0} c_{\omega} = 0$) for
    which $\mathcal H_{\mq_0} \ge \mathcal H+\omega$ holds on
    $\mL_{\mq_0}\setminus(\cup_{j=S,NE,NW}\mathfrak B^j_{q_0}) =
    \mathfrak{T}_{\mq_0}\setminus(\cup_{j=S,NE,NW}\mathfrak B^j_{q_0})$. This shows the first statement of the proposition; the second follows from the previously mentioned fact that $H_t \le \mathcal{H} + \omega$ with probability at least $1 - N^{-2L}/4$. 
          \begin{figure}[h]
 \begin{center} \includegraphics[width=5cm]{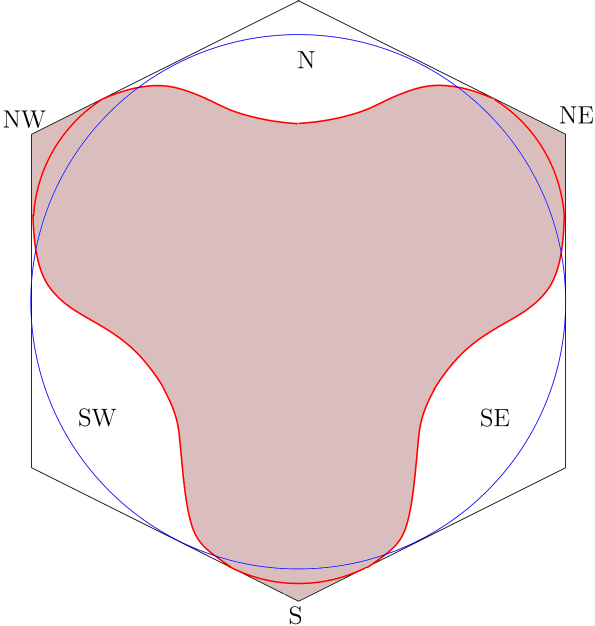}
   \caption{The arctic curve of $\mathcal H_{\mq_0}$ (in red) enclosing the liquid region $\mL_{\mq_0}$ and the arctic curve of $\mathcal H$ (in blue), enclosing the liquid region $\mL$. The shaded region is $\mathfrak
   T_{q_0}$} \label{fig:liquid2}
\end{center}
\end{figure}
        \end{proof}

\begin{figure}[h]
 \begin{center} \includegraphics[width=9cm]{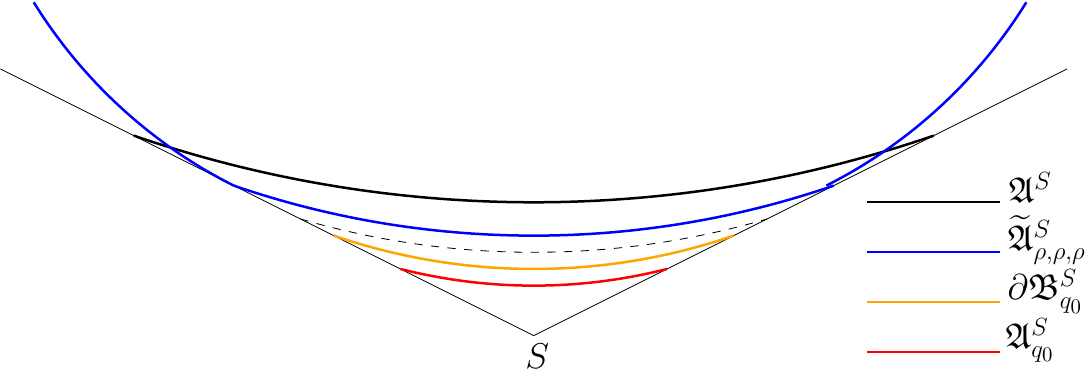}
   \caption{ A zoom-in around the south corner of 
     $\mathfrak X$. The line $\mA^S$  (respectively $\widetilde{\mA}^S_{\rho,\rho,\rho}$) is the south
     arc of the arctic boundary of the limit shape in
     $\mathfrak X$ (respectively $\widetilde{\mathfrak X}_{\rho,\rho,\rho}$). The dashed line is the boundary of the region at distance $c_\omega$ below $\widetilde{\mA}^S_{\rho,\rho,\rho}$.
     The line $\mA^S_{\mq_0}$ (respectively
     $\partial \mathfrak B^S_{q_0}$) is the south arc of the arctic
     boundary of $\mathcal H_{\mq_0}$ (respectively, the boundary of
     $\mathfrak B^S_{q_0}$). For $c_\omega,r,\mq_0$ small, these four curves are
     close to each other.
   } \label{fig:lines}
\end{center}
\end{figure}

        \begin{proof}[Proof of Proposition \ref{prop:astuzia}]
          By symmetry, it suffices to prove the statement in $\mathfrak B^S_{q_0}$, namely, that $H_t (z) = 0$ for $z \in \mathfrak{B}_{q_0}^S$ and $t \in I_i^N$. Let us briefly outline how we will proceed. We will first embed a hexagon $\widetilde{\mathfrak{X}}_{\rho,\rho,\rho}$ of side length $\rho \approx 1$ within $\mathfrak{X}$ and use monotonicity to bound the Glauber dynamics $(H_t)$ in $\mathfrak{X}$ from above by those $(G_t)$ in $\widetilde{\mathfrak{X}}_{\rho,\rho,\rho}$, in a region containing $\mathfrak{B}_{q_0}^S$. After time $t \ge t_{\mix}^{N_{i+1}}$, the latter dynamics will have mixed. By the second part of Lemma \ref{volumeheight}, we will have w.o.p. that $G_t = 0$ in the associated south frozen region, which we will see contains $\mathfrak{B}_S^{q_0}$. Hence, monotonicity will yield w.o.p. that $H_t = 0$ on $\mathfrak{B}_S^{q_0}$, proving the proposition. 
          
          To implement this, let us fix some notation. While $\mathfrak X$ as usual
          has side length $1$, we let {$\widetilde{\mathfrak X}_{\rho,\rho,\rho}$ denote
          the smaller hexagon $\mX_{\rho,\rho,\rho}$ of side length $\rho=1-r$, translated so that it shares 
          the south corner of $\mathfrak X$}. The constant $r$ is the same as in
          \Cref{thm:therecursion} and it will be chosen
          small enough later. We assume throughout that $(1-\rho)^j N \in \mathbb{Z}$ if $j \le i+1$, for notational convenience.
          \begin{rem}            \label{eq:ordine}
            We will choose $r$ sufficiently small (once $q_0$ is given), then $\omega$ sufficiently small (once $r$ is given). None of these parameters depends on $N$.
          \end{rem}
          \noindent Further let $\mathfrak X^{N_i}$ and
          $\widetilde{\mathfrak {X}}_{\rho,\rho,\rho}^{N_i}$ denote the intersection of
          $\mathfrak X$ and $\widetilde{\mathfrak X}_{\rho,\rho,\rho}$ with $\mathbb{T}_{N_i}$ (the triangular grid of mesh
          $1/N_i$). In accordance with \Cref{def:pa}, {let $\mathfrak A^S$ (resp. $\widetilde{\mathfrak A}^S_{\rho,\rho,\rho}$)} denote the  south arc of the arctic boundary of $\cH$ (resp. of $\widetilde{\mathcal H}_{\rho,\rho,\rho}$).  See Figure \ref{fig:lines}.
          Note that $\widetilde{\mathcal H}_{\rho,\rho,\rho}$ is a simple linear rescaling of
          $\mathcal H$.

          Let $(G_t)_{t\ge0}$  denote the height function associated to the  Glauber dynamic on $\widetilde{\mathfrak {X}}_{\rho,\rho,\rho}$ starting from the maximal
          height function in $\Omega_{\widetilde{\mathfrak {X}}_{\rho,\rho,\rho} }$  and recall that $(H_t)_{t\ge0}$ is the height function of the Glauber  dynamic on $\mathfrak {X}$ starting from the maximal
          height function in $\Omega_{\mX}$.  In both cases, the tile size is $1/N_i,i\le i_{\rm max}$.
{We couple the two processes as described in the proof of \Cref{prop:cancouple} (except that for the process in $\widetilde{\mathfrak {X}}_{\rho,\rho,\rho}$, the Poisson clocks and coin flips in $\mX\setminus\widetilde{\mathfrak {X}}_{\rho,\rho,\rho}$ induce no update).} We will show at end of the proof  that the  following holds except with probability $N^{-3L}$:
\begin{eqnarray}
  \label{lem:astuzia}
H_{\mC N_i^{2+\delta}+t} (z) \le G_t (z),\; \text{ for every } t\le N^5 \text{ and }
          z\in \mathfrak X\; \text{ with } \dist (z, \mathfrak A^S) \le N_i^{-1}. 
\end{eqnarray}
 {It follows easily from the definition of the monotone  coupling in \Cref{prop:cancouple} that,}
          except with probability $N^{-3L}$,  
          $H_{\mC N_i^{2+\delta}+t} (z) \le G_t (z)$ for all $z$ below $\mathfrak A^S$ 
          and all $t \le N^5$.  The mixing time in
          $\widetilde{\mathfrak {X}}_{\rho,\rho,\rho}$ is $\tmix^{N_{i+1}}$, since  $\widetilde{\mathfrak {X}}^{N_i}_{\rho,\rho,\rho}$ is nothing but
          $\mathfrak X^{N_{i+1}}$, up to a rescaling by
          $\rho=(1-r)$. After time
          $t_{mix}^{N_{i+1}} (\epsilon_{i+1})$ we have (by the
          concentration statement in Proposition \ref{volumeheight}) that
          $G_t (z) = \widetilde{\mathcal H}_{\rho,\rho,\rho} (z) = 0$ for $z$ at distance at least  $c_\omega$ below $\widetilde{\mathfrak A}^S_{\rho,\rho,\rho}$, except with probability
          $\epsilon_{i+1}+N^{-3L}$.  It follows that, except with probability
          $\epsilon_{i+1}+2 N^{-3L}\le (1/3)\xi^N_i$ (recall \eqref{eq:xiNi}), for all $t\le N^5$ we have
          $H_{\mC N_i^{2+\delta}+t} (z) = 0$ for all $z$ at distance {at least $c_\omega$}
          below $\widetilde{\mathfrak A}^S_{\rho,\rho,\rho}$.
          On the other hand, $c_\omega$ sufficiently small (depending only on $\mq_0$ and $r$) the region at distance at least $c_\omega$ below $\widetilde{\mathfrak A}^S_{\rho,\rho,\rho}$  includes the region $\mathfrak B^s_{q_0}$ (see also Figure \ref{fig:lines}). {In fact, $\partial \mathfrak B^S_{\mq_0}$ is at distance  $c_\omega$ from  $\mA^S_{\mq_0}$  (by definition of $\mathfrak B^S_{\mq_0}$); $\mA^S_{\mq_0}$ is strictly below 
            $\mA^S$ (by \Cref{boundarydistance}), at a distance depending only on $q_0$; the distance between $\mA^S$ and $\widetilde{\mA}^S_{\rho,\rho,\rho}$ tends to zero as $r\to0$ (since the two are related via a scaling by $1-r$). Choosing the parameters as in \Cref{eq:ordine}, it follows that the region of distance at least $c_{\omega}$ below $\widetilde{\mathfrak{A}}_{\rho,\rho,\rho}^S$ contains $\mathfrak{B}_{q_0}^S$; see \Cref{fig:lines}.  The desired result then follows.}

It remains to prove \eqref{lem:astuzia}. Let $(F_t)_{t\ge0}$ denote the  height function of the stationary Glauber dynamic  on 
            $\widetilde{\mathfrak {X}}_{\rho,\rho,\rho}$, {that is, the one started from the uniform measure on $\Omega_{\widetilde{\mX}_{\rho,\rho,\rho}}$.}  {Under the coupling of \Cref{prop:cancouple}},  $G_t \ge F_t$  almost surely. Therefore,  it suffices to show that
            \begin{equation}
              \label{eq:daprovare}
\text{except with probability $N^{-3L}$}, \quad              H_{\mC N_i^{2+\delta}+t} (z)
            \le F_t(z)\quad \text{for every}\quad  t\le N^5
            \end{equation} and
            for every $z \in \mathfrak{X}$ such that $\dist (z, \mA^S) \le N_i^{-1}$. For
            lightness, we say that $z$ is \emph{along $\mA^S$} if $\dist (z, \mA^S) \le N_i^{-1}$.
            By the concentration bounds for
            the height under the stationary measure in $\widetilde{\mX}_{\rho,\rho,\rho}$ (Theorem \ref{volumeheight}) and \Cref{rem:forall}, the following two statements hold w.o.p.:
\begin{enumerate}
            \item [(i)]
              $F_t(z) > \widetilde{\mathcal H}_{\rho,\rho,\rho}(z) - \omega$ for all $z \in
              {\mathfrak X_{\rho,\rho,\rho}}$ and $t\le N^5$.
            \item [(ii)]
               $F_t(z) = \widetilde{\mathcal H}_{\rho,\rho,\rho}(z)$ for all $t\le N^5$ and for {all $z$ in  $\widetilde{\mF}^S_{\rho,\rho,\rho}$, with $\dist (z, \widetilde{\mA}^S_{\rho,\rho,\rho}) \ge c_\omega$}.
            \end{enumerate}
            Let us first prove \eqref{eq:daprovare}, for all $z$
             along $\mA^S$ that are sufficiently close to the boundary of $\mathfrak X$, namely, for $z$ satisfying $\dist (z, \partial \mathfrak{X}) \le b(r)$, for some sufficiently small constant $b(r) > 0$ (to be chosen later). In fact, observe in this case that $z$ is {in $\widetilde{\mF}^{SW}_{\rho,\rho,\rho}\cup\widetilde{\mF}^{SE}_{\rho,\rho,\rho} $
              and bounded away from $\widetilde{\mA}_{\rho,\rho,\rho}$} (see Figure \ref{fig:lines}).  Thus, point (ii) above and the fact that {we are taking $\omega$ sufficiently small (recall \Cref{eq:ordine})  so that $c_{\omega} < b(r)$}   implies that $F_t (z) =
             \widetilde{\mathcal H}_{\rho,\rho,\rho} (z)$ w.o.p..  {In $\widetilde{\mF}^{SW}_{\rho,\rho,\rho}\cup\widetilde{\mF}^{SE}_{\rho,\rho,\rho} $, the limit shape  $\widetilde{\cH}_{\rho,\rho,\rho}$ coincides with the maximal height function, so that 
             $\widetilde{\mathcal H}_{\rho,\rho,\rho}(z) \ge H_{\mC N_i^{2+\delta}+t}
             (z)$ deterministically}. This yields \eqref{eq:daprovare} for all such $z$.

 Next, we must  verify \eqref{eq:daprovare} for all $z$ along
 $\mA^S$ with $\dist (z, \partial \mathfrak{X}) \ge b(r)$.   Observe by Proposition \ref{thm:perturb} that w.o.p. we have $H_{\mC N_i^{2+\delta}+t} \le \mathcal H +
 \omega$ for every $t\le N^5$. Since $\mathcal H(z) \le N_i^{-1} \le \omega$ for $z$ along $\mA^S$, it follows that $H_{\mC N_i^{2+\delta}+t}(z) \le 2\omega$.   Since $\mA^S$ is above  $\widetilde{\mA}^S_{\rho,\rho,\rho}$ by a distance $c(r)>0$, it follows that $\widetilde{\mathcal H}_{\rho,\rho,\rho} (z) > \psi(r)$ for some constant $\psi(r) > 0$ (dependent on $r$ but not on $\omega$) for $z$ along $\mA^S$ with $\dist (z, \partial \mathfrak{X}) \ge b(r)$.  Together with item (i) above, we conclude that w.o.p.
 \[F_t(z)\ge \widetilde{\mathcal H}_{\rho,\rho,\rho}(z)-\omega\ge \psi(r)-\omega\ge 2 \omega\ge H_{\mC N_i^{2+\delta}+t} \quad \text{for all}\quad t\le N^5,\]
 and for all $z$ along $\mA^S$ with $\dist (z, \partial \mathfrak{X}) \ge b(r)$,
 where the third inequality holds provided by taking $\omega$ {sufficiently small}, as in  \Cref{eq:ordine}.
  \end{proof}
        
\subsection{Proof of  \Cref{thm:therecursion}}
\label{sec:thm53}

  Theorem \ref{thm:step1} implies the following. Let as usual $(\eta_t)_{t\ge0}$ denote the Glauber dynamic in $\mX$, started from some $\eta_0\in\Omega_\mX$, and  
 let $$(\tilde \eta_t)_{t\ge \tau}, \quad \tau:=\mC N_i^{2+\delta}+t_{\rm mix}^{N_{i+1}}(\epsilon_{i+1})$$ be the Glauber dynamic in $\mX$, constrained (in the sense of \Cref{def:constrained}) between floor  
$\mathcal H_{-\mq_0}$ and ceiling $\cH_{q_0}$  and with initial condition $\tilde \eta_{ \tau}=\eta_{ \tau}$. If the two processes are coupled using the same Poisson clocks and coin flips (as in \Cref{prop:cancouple}), then off of an event of probability $N^{-2L}+2\epsilon_{i+1}$, we have $\eta_t=\tilde \eta_t$ for all $t\in I^N_i$.
 Note also that, because of Proposition \ref{prop:tmixconstrained}, \eqref{e:submult}, and the facts that $\min_{i\le i_{max}}\epsilon_i=N^{-L}$ and $\max_{i\le i_{max}}N_i =N$, 
\[
T^N_i:=\mC N_i^{2+\delta}+\tmix^{N_{i+1}}(\epsilon_{i+1})+t_{\mix}^{\mq_0,N_i}\left(\frac{\epsilon_i}4\right)\le N^{5}.
\]
  Therefore, uniformly in the initial condition $\eta_0$, at time $T^N_i$ {the law of  $\eta_t$}  has total variation distance at most
  $N^{-2L}+2\epsilon_{i+1}+\epsilon_i/4$ from the invariant measure of
  the constrained process $\tilde \eta$, that is denoted $\pi_{\cH_{\pm q_0}}$, using the notation of \Cref{def:P+}.
Next, note by \Cref{volumeheight} (and the bound $N_i\ge N^{1/10}$ for every $i\le i_{max}$) that the total variation distance between the uniform measure $\mathbb P^{N_i}$ on $\Omega^{N_i}$ and $\pi_{\cH_{\pm q_0}}$ is at most $N^{-2L}$. Hence, {the total variation distance between  $\mathbb
  P^{N_i}$ and the distribution of $\eta_{T^N_i}$ } is at most \[\xi^N_i:=
  2N^{-2L}+2\epsilon_{i+1}+\frac{\epsilon_i}4=\left(\frac1{5}+\frac14\right)\epsilon_i+2N^{-2L}\le
  \epsilon_i\]
  (where the last inequality holds because $\epsilon_i\ge N^{-L}\gg N^{-2L}$).
  This implies the desired bound
    \eqref{eq:recursion}.

	\section{Estimates for the height reduction} 
	\label{TimeBarrier} 

        The goal of this section is to prove \Cref{thm:fN}. Throughout this section, $\varepsilon_0>0$ is chosen small
        enough (but independently of $N$). Also,  $\delta \in (0, 10^{-1000})$ is fixed throughout to be a sufficiently small real number. Throughout, we will fix a constant $A > 10^{10}$. Moreover, $M,K,D>0$ will be large constants that are much larger one with respect to the other (and $A$), in a certain order. The values of the constants $(M,K,D)$ need not be the same at each occurrence, but we will always impose that they satisfy 
        \begin{eqnarray}
        	\label{eq:costanti}
        	M>10^5 K>10^{10}A>10^{15}D>10^{20}.
        \end{eqnarray}

	\subsection{Annulus decomposition of the hexagon} 
	
	\label{AProperty}
	
\begin{figure}[h]
 \begin{center} \includegraphics[width=7cm]{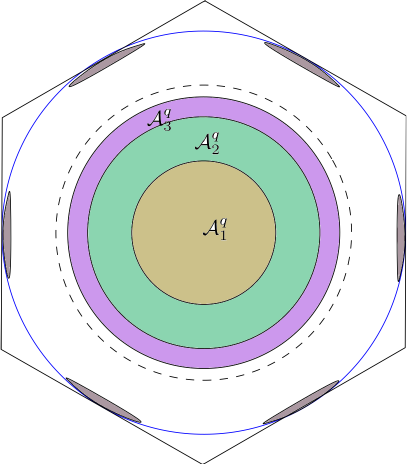}
   \caption{A schematic view of the  decomposition of
     $\mL_q$ into $\mathcal A^q_\ell$. In the figure, $q$ is assumed to be very small, so that
     the liquid region is essentially a disk. For $\ell=1$,
     $\mathcal A_\ell^q$ is approximately a disk and for $\ell>1$ the
     sets $\mathcal A_\ell^q$ are approximately concentric annuli,
     closer and closer to the arctic boundary. However, for $\ell$
     sufficiently large (such that $4^{-\ell}\lesssim N^{-2/3}$) $\mathcal A_\ell^q$ consists of six disjoint
     regions, one per tangency point (this can be checked via  \eqref{eq:tdte}). For  $\ell > \log N$, the annuli are empty,
     because by definition $\td_z\ge N^{-1/2}$ and $\te_{z;q}\ge N^{-2/3}$.  The sets $\mathcal B^q_\ell$ are
     similar.} \label{fig:annuli}
\end{center}
\end{figure}

	We begin with the following definition of annuli; see Figure \ref{fig:annuli} (that explains why we call these regions ``annuli'').
	
	\begin{definition}[Annuli]
		
		\label{aq} 
	For any real numbers $q, u \ge 0$ and integer $\ell \ge 1$, define the following subsets of the liquid region:
	\begin{flalign*}
		\mathcal{A}_{\ell}^{q} = \big\{ z \in \mathfrak{L}_q: 4^{-\ell} \le \mathtt{d}_{z} \mathtt{e}_{z;q} \le 4^{1-\ell} \}; \qquad  \mathcal{B}_{\ell}^q = \big\{ z \in \mathfrak{L}_q : 4^{1/2-\ell} \le \mathtt{d}_z  \mathtt{e}_{z;q} \le 4^{1 - \ell} \big\}.
	\end{flalign*}
as well as
	\begin{flalign*}
		\mathcal{A}_{\ell}^{q} (u) = \big\{ z \in \mathcal{A}_{\ell}^q: \mathtt{d}_z^{1/2} \mathtt{e}_{z;q}^{3/2} \ge uN^{-1} \big\}; \qquad \mathcal{B}_{\ell}^q (u) = \big\{ z \in \mathcal{B}_{\ell}^q : \mathtt{d}_z^{1/2} \mathtt{e}_{z;q}^{3/2} \ge uN^{-1} \big\}.
	\end{flalign*}
	
	\end{definition} 
          
\begin{rem}
	The relevance of the parameter $\mathtt{d}_z \mathtt{e}_{z;q}$ (in the definitions of $\mathcal{A}_{\ell}^q$ and $\mathcal{B}_{\ell}^q$) is due to \Cref{qlineq}, which indicates that it dictates how the limit shape $\mathcal{H}_q (z)$ varies upon changing $q$. The relevance of $\mathtt{d}_z^{1/2} \mathtt{e}_{z;q}^{3/2}$ (in the definitions of $\mathcal{A}_{\ell}^q (u)$ and $\mathcal{B}_{\ell}^q (u)$) is due to the third statement of \Cref{hde0}, which indicates that it dictates the value of $\mathcal{H}_q (z)$.
	
\end{rem}

For any   $q>0$, we will treat differently the  \emph{bulk  liquid phase}, consisting of 
	 the points  $z\in\mL_{q}$ such that	\begin{eqnarray}
		\label{eq:distanza}
		\md_z\me_{z;q}\ge\mq^{-2/3}N^{-2/3+A\delta}
	\end{eqnarray}
	and the \emph{edge phase}, consisting of the $z\in\bar\mX$ such that
	\begin{eqnarray}
		\label{eq:stavolta}
		\md_z\me_{z;q} <4^{20}\mq^{-2/3}N^{-2/3+A\delta},
	\end{eqnarray}
	\noindent as the mechanisms for mixing in these two regions will be different (we will explain the reason for the definition \eqref{eq:distanza} of the bulk liquid phase in Remark \ref{yde} below). Note that the bulk liquid and edge phases  overlap on $20$ ``scales'' of the form $4^{-\ell}$. Let $\ell_0 = \ell_0 (q) $ denote the tenth such scale, defined to satisfy
	\begin{flalign}
		\label{adelta4}
		2^{-21} q^{2/3} N^{2/3 - A\delta} \le 4^{\ell_0(q)} < 2^{-19} q^{2/3} N^{2/3 - A\delta}.
	\end{flalign}
        \begin{rem}
          \label{rem:utile}
          
          Observe from the first bound in \eqref{adelta4} that $\ell_0(q) \gtrsim \log N \gg 1$ for $q > N^{2A\delta-1}$ (in the discussion below, we will always take $q \ge N^{M\delta-1}$ for a large constant $M \ge 2A$ since, for all practical purposes, $q \in [0, N^{M\delta-1}]$ will be indistinguishable from $q=0$). 
        \end{rem}
    
    The next lemma lower bounds $\mathcal{H}_q$ in $\mathcal{A}_k^q$ if $k$ is not too large, and the following lemma shows containment properties for the $\mathcal{A}$ and $\mathcal{B}$ annuli. 
	 \begin{lem}
		
		\label{hqa} 
		Let $0<q\le \varepsilon_0$ be a real number; $k \le \ell_0 (q) + \delta \log N$ be a positive integer; and $z \in \mathcal{A}_k^q$ be a point satisfying $\dist (z, \mathfrak{A}_q) = \dist (z, \mathfrak{A}_q^S)$. Then, $\mathcal{H}_q (z) \ge N^{(5/4)A\delta-1}$ and $(\mathfrak{d}_z, \mathfrak{e}_{z;q}) = (\mathtt{d}_z, \mathtt{e}_{z;q})$.
		
	\end{lem}

	\begin{proof} 
		We have 
		\begin{flalign}
			\label{dze0}
			\mathtt{d}_z^{1/2} \mathtt{e}_{z;q}^{3/2} \ge (\mathtt{d}_z \mathtt{e}_{z;q})^{3/2} \ge 4^{-(3/2)\ell_0 (q)- (3/2)\delta \log N} & \gg 2^{27} N^{(3A/2-3)\delta} (qN)^{-1} \ge N^{5A\delta/4-1},
		\end{flalign} 
		\noindent where in the first inequality we used the fact that $\mathtt{d}_z \le 1$; in the second we used the fact that $z \in \mathcal{A}^q_k$ for some $k \le \ell_0 (q) + \delta \log N$; and in the last two we used \eqref{adelta4}, with the fact that $A \ge 30$ and $q\le 1$.
In particular, $\te_{z;q}\gg N^{-2/3}\td_z^{-1/3}$, which implies $\te_{z;q}=\me_{z;q}$. Moreover, since \eqref{dze0} and $\me_{z;q}\lesssim \md_z\le \td_z$ ({where we used \eqref{eq:elegeo0} in the first bound}) imply that $\mathtt{d}_z \gg N^{-1/2}$, we also have $\td_z=\md_z$. Together with \eqref{dze0}, these yield $\mathfrak{d}_z^{1/2} \mathfrak{e}_{z;q}^{3/2}\gg N^{5A\delta/4-1}$, so the lemma follows from \Cref{hde0}.             
	\end{proof} 
	
	\begin{lem}
		
		\label{aqaq} 
		Fix $\delta>0$.
		Let $k \ge 1$ be a positive integer and $q', q \in [0, \varepsilon_0]$, and $u \ge 0$ be real numbers such that $|q-q'| \le 2^k N^{-1-4\delta} u$. We have
		\begin{flalign}
			\label{bk1u} 
		\mathcal C^q_k(u):=	\left[\mathcal{B}_{k+1}^q (u) \cup \bigcup_{j=1}^k \mathcal{A}_j^q (u)\right] \subseteq \bigcup_{j=1}^{k+1} \mathcal{A}_j^{q'} \big( (1 - N^{-\delta}) u \big),
		\end{flalign} 
and
		\begin{flalign}
			\label{dze2} 
			1 - 2N^{-2\delta} \le \mathtt{e}_{z;q'} \mathtt{e}_{z;q}^{-1} \le 1 + 2 N^{-2\delta}, \qquad \text{for every $z\in \mathcal C^q_k(u)$}.
		\end{flalign}
              \end{lem}
	\begin{proof} 		
		Fix $z \in \mathcal C^q_k(u)$ {and assume, by symmetry, that $\ell_z$ is the southwest edge of $\mX$}. To verify \eqref{bk1u}, it suffices to show that 
		\begin{flalign} 
			\label{dzqe}
			\mathtt{d}_z \mathtt{e}_{z;q'} \ge 4^{-k-1}, \quad \text{and} \quad \td_z^{1/2}\te_{z;q'}^{3/2}  \ge (1 - N^{-\delta}) uN^{-1}.
		\end{flalign} 
		\noindent To show the former, observe that
		\begin{flalign}
			\label{dzqe2}
			\mathtt{d}_z \mathtt{e}_{z;q'} = \mathtt{d}_z \mathtt{e}_{z;q} + \mathtt{d}_z (\mathtt{e}_{z;q'} - \mathtt{e}_{z;q}) \ge 4^{-k-1/2} - 2^{k} N^{-1-3\delta} u \mathtt{d}_z,
		\end{flalign}
		where the second statement follows from the fact that $z \in \mathcal C^q_k(u)$ and $\dist_{\bm{v}} (\mathfrak{A}_q, \mathfrak{A}_{q'}) \lesssim  |q-q'| \le 2^{k} N^{-1-4\delta} u$ (by \Cref{boundarydistance}\footnote{Recall that \Cref{rem:bysymmetry} states how to modify the statement of \Cref{boundarydistance}  if $\ell_z$ is not the southwest edge of $\mX$.}).
 Further observe that it is enough to prove \eqref{dzqe} when $\mathtt e_{z;q'}$ is minimal for a given value of $\mathtt d_z$, that is, when $z\in \mathcal B_{k+1}^q(u)\cup \mathcal A_k^q(u)$. In that case,  $8^{1-k} \ge (\mathtt{d}_z \mathtt{e}_{z;q})^{3/2} \ge \mathtt{d}_z u N^{-1}$. Therefore, $4^{-k-1} \ge 2^{k-2} N^{-1} u \mathtt{d}_z$, which with \eqref{dzqe2} confirms the first statement in \eqref{dzqe}. It also (using \eqref{dzqe2}) implies that 
		\begin{flalign*} 
			1 - 2N^{-2\delta} \le \mathtt{d}_z \mathtt{e}_{z;q'} ( \mathtt{d}_z \mathtt{e}_{z;q})^{-1} \le 1 + 2 N^{-2\delta},
		\end{flalign*} 
		
		\noindent which verifies \eqref{dze2}. Since $\mathtt{d}_z^{1/2} \mathtt{e}_{z;q}^{3/2} \ge uN^{-1}$, this gives $\mathtt{d}_z^{1/2} \mathtt{e}_{z;q'}^{3/2} \gtrsim (1 - 2N^{-2\delta})^{3/2} u N^{-1}$, which verifies the second statement in \eqref{dzqe}.	
	\end{proof}

	\subsection{Height reduction across annuli}
	
	\label{ReduceAnnuli0} 
	
	 In this section, we state three variants of the following statement. If a height function $H_0$ initially satisfies w.o.p. an upper bound of the type $H_0\le \cH_q$ (possibly with a small additive error) in a certain region $R\subseteq \mX$, then after a suitable amount of time under the Glauber dynamics, $H_t$ satisfies w.o.p. the improved upper bound $H_t\le \cH_{q'}$, with $q'<q$, in a modified region $R' \subseteq \mathfrak{X}$. Lemma \ref{ak1akak1} addresses regions in the bulk liquid phase, while \Cref{akl0} and \Cref{boundaryaj2} address regions in the edge phase (close to or at the arctic boundary, respectively).
	 
	  In all cases, the time $t$ at which this improved bound holds will be at least $N^2 (q-q')$ (and will attain this bound, up to powers of $N^{\delta}$, if $q-q'$ is made as large as permitted). This suggests that the parameter $q$ defining the upper barrier $\mathcal{H}_q$ decreases ``at speed'' $N^{-2-o(1)}$. As such, $q$ should reduce from $\varepsilon_0$ to $0$ in time $N^{2+o(1)}$, which is consistent with the notion that the Glauber dynamics mix in time $N^{2+o(1)}$. Furthermore, in all of these statements, we must run for time\footnote{In Lemma \ref{boundaryaj2}, one may identify $2^k$ with $\mathtt{d}^{-1/3} N^{1/3+o(1)}$, as one can verify that this is approximately the value of $(\mathtt{d}_z \mathtt{e}_{z;q'})^{-1/2}$ there.} $2^k N$ to ``witness'' a height reduction on the annulus $\mathcal{A}_k$. To heuristically see why this should be the case, observe that one expects $q$ to change by $2^k N^{-1}$ after time $2^k N$. Hence, Lemma \ref{qlineq} indicates that the height function $(H_t)$ on $\mathcal{A}_k^{q'}$ reduces by $(\mathfrak{d}_z \mathfrak{e}_{z;q'})^{1/2} 2^k N^{-1}$ in this time. Since $(\mathfrak{d}_z \mathfrak{e}_{z;q'})^{1/2} \sim 2^{-k}$ on $\mathcal{A}_k^{q'}$, this reduction is around $N^{o(1)-1}$, which is the minimal that can be ``detected'' (as it is the minimal height difference).

	 The following lemma, to be shown in \Cref{ProofAnnulus00}, states the bulk liquid phase mixing bound. It indicates that, if a height function $H_0$ is initially bounded by some $\mathcal{H}_q$, then running the Glauber dynamics $(H_t)_{t \ge 0}$ decreases this upper bound to $\mathcal{H}_{q'}$ at points $z$ in the bulk liquid phase (see the bound $k \le \ell_0 (q) + 5$ in \eqref{eq:assu}), for some $q' < q$. 

\begin{figure}[h]
 \begin{center} \includegraphics[width=11cm]{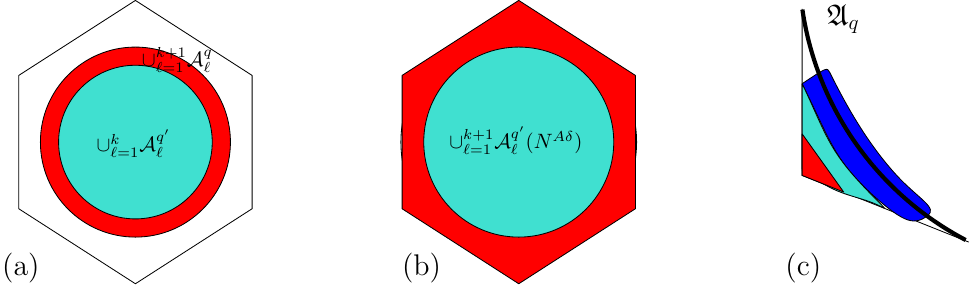}
  \caption{         {A pictorial illustration of Lemmas 
      \ref{ak1akak1}, \ref{akl0} and \ref{boundaryaj2}. (a): in \Cref{ak1akak1}, we assume an upper bound on the height function at time zero in $\cup_{\ell=1}^{k+1}\mathcal A^q_\ell$ (union of red and blue regions) and we deduce an improved upper bound at time $\mathtt t$, in the smaller blue region $\cup_{\ell=1}^{k}\mathcal A^{q'}_\ell$. (b): in \Cref{akl0}, we assume at time zero an upper bound in the whole hexagon, and we deduce an improved upper bound at time $\mathtt t$, in the blue region $\cup_{\ell=1}^{k+1}\mathcal A^{q'}_\ell(N^{A\delta})$. (c): in \Cref{boundaryaj2} we assume at time zero an upper bound in the whole hexagon and a stronger upper bound in the complement of $\mL_q^+(K\delta),K>1$ (red region), and we obtain at time $\mathtt t$  improved upper bounds, both in the complement of $\mL_q^+(\delta)$ (union of light blue and red regions), and in the region  (dark blue) close to the arctic boundary, identified by condition \eqref{dzhq}. In (c),   we zoom on a corner of the hexagon, to make the figure understandable, and the thick curve denotes the arctic boundary $\mA_q$.} }
\label{fig:55}
\end{center}
\end{figure}

See Fig. \ref{fig:55}.
           
	\begin{lem} 
		
		\label{ak1akak1}
		
		Let $A,D$ satisfy \eqref{eq:costanti}. Also let $q$ be a real number and $k \ge 0$ be an integer satisfying 
		\begin{equation}
			\label{eq:assu}
			q \in [N^{A\delta-1}, \varepsilon_0], \quad k \le \ell_0 (q) + 5.
		\end{equation}
		
		\noindent Further let $q' \ge 0$ satisfy
		\begin{flalign}
			\label{qkq} 
			q' \in [q - 2^k N^{(D+2)\delta-1}, q - 2^{k}  N^{(D-2) \delta-1}]; \qquad \mathsf{t} = q'^{-1} k 2^k N^{1+6A\delta}.
		\end{flalign} 
		
		\noindent Consider the Glauber dynamics $(H_t)_{t \ge 0}$ starting from an initial condition $\mu \in \mathcal P^+$ (recall Definition \ref{def:P+}) and assume that, w.o.p. under $\mu$, 
		\begin{flalign}
			\label{htzq0} 
			\begin{aligned} 
				& H_0 (w) \le \mathcal{H}_{q} (w) + N^{\delta-1}, \qquad  \text{for each $w \in \bigcup_{\ell=1}^{k+1} \mathcal{A}_{\ell}^{q}$ }. 
			\end{aligned}
		\end{flalign}
		
		\noindent Then we have w.o.p. that
		\begin{flalign}
			\label{htzq1}
			H_t (z) \le \mathcal{H}_{q'} (z), \qquad \text{for each $z \in \bigcup_{\ell = 1}^k \mathcal{A}_{\ell}^{q'}$ and $t \in [\mathsf{t}, N^5]$}.
		\end{flalign}  
		
	\end{lem}

	The next lemma, to be shown in Section \ref{ProofAnnuli20}, states the following. Suppose a height function initially satisfies an upper bound on the full hexagon (observe that this assumption differs from that in \Cref{ak1akak1}, which is indicative of the eventually different mechanisms of mixing in the bulk liquid and edge phases). Then, it satisfies an improved upper bound on a union of annuli, potentially within the edge phase. See Fig. \ref{fig:55}.

	\begin{lem} 
		
		\label{akl0} 
		Let $A,K,D$ satisfy \eqref{eq:costanti}. Also let $k \ge 1$ be an integer and $q \in [2^k N^{K \delta-1}, \varepsilon_0]$ be a real number. Further let $q' \ge 0$ satisfy
		\begin{flalign}
			\label{qdq}
			\begin{aligned} 
				& q' \in [q - 2^k N^{(D+1)\delta-1}, q - 2^k N^{(D-1)\delta - 1}]; \qquad 	\mathsf{t}:=q^{-1} k 2^{k+3} N^{1+11K\delta}.
			\end{aligned} 
		\end{flalign}
		
		\noindent  Consider the Glauber dynamics $(H_t)_{t \ge 0}$ starting from an initial condition $\mu \in \mathcal P^+$ (recall Definition \ref{def:P+}) and assume that, w.o.p. under $\mu$, 
		\begin{flalign}
			\label{hszq2} 
			\begin{aligned}
				H_0 (w) \le \mathcal{H}_{q} (w) + N^{\delta-1}, \qquad  \text{for each $w \in \mathfrak{X}$}.
			\end{aligned} 
		\end{flalign}
		
		\noindent Then,  w.o.p., we have
		\begin{flalign}
			\label{htz2} 
			H_t (z) \le \mathcal{H}_{q'} (z), \qquad \text{for each $z \in \bigcup_{j=1}^{k+1} \mathcal{A}_j^{q'} (N^{A\delta})$ and $t \in [ \mathsf{t}, N^5]$.}
		\end{flalign}
		
	\end{lem}

The next lemma, to be shown in \Cref{Proofqn2}, essentially states the following. Assume that a height function is initially bounded above by $\mathcal{H}_q + N^{K\delta-1}$ on all of $\mathfrak{X}$ and by $\mathcal{H}_q$ outside of the augmented liquid region $\mathfrak{L}_q^+ (K\delta)$ (recall \eqref{lqabc}). Then this height function will reduce under the Glauber dynamics to $\mathcal{H}_{q'} + N^{\delta/4-1}$ on all of $\mathfrak{X}$, and to $\mathcal{H}_{q'}$ outside of $\mathfrak{L}_{q'}^+ (\delta)$, for $q' \le q$ satisfying the third bound in \eqref{dzhq}. In addition to reducing $q$, this effectively reduces the error $K$ from a potentially large constant to $1$. 

Let us provide a bit more explanation for the second condition in \eqref{dzhq}. It implies that
$\mathtt{e}_{z;q'} \le \mathtt{d}_z^{-1/3}
N^{2K\delta/3-2/3}$. Since
$\mathtt{e}_{z;q'} \ge \mathtt{d}_z^{-1/3} N^{-2/3}$, this
indicates that $\mathtt{e}_{z;q'}$ is almost as small as
possible (recall \eqref{eq:tdte}), meaning that $z$ is essentially at the arctic
boundary. Thus, the below lemma indicates how height functions reduce very close to the arctic boundary. See Fig. \ref{fig:55}.

\begin{lem} 
	
	\label{boundaryaj2}
	
	Fix a real number $K$ satisfying \eqref{eq:costanti}; $\mathtt{d} \in [N^{-1/2}, 1)$; and real numbers $q \ge q' \ge 0$ with $\varepsilon_0\ge q \ge N^{3500K\delta-2/3} \mathtt{d}^{-1/3}$. Consider the Glauber dynamics $(H_t)_{t \ge 0}$ starting from an initial condition $\mu\in\mathcal P^+$ and assume that, w.o.p., 
	\begin{flalign}
		\label{hsz2}
		\begin{aligned} 
			& H_0 (w) \le \mathcal{H}_q (w) + N^{K\delta-1}, \qquad \text{for all $w \in \mathfrak{X}$}; \\
			& H_0 (w) \le \mathcal{H}_q (w), \qquad \qquad \qquad \text{for all $w \notin \mathfrak{L}_q^+ ({K\delta})$}. 
		\end{aligned} 
	\end{flalign} 
	
	\noindent Denoting $\mathsf{t} = q^{-1} \mathtt{d}^{-1/3} N^{4/3+15000K\delta}$, the following holds w.o.p. for any $t \in [\mathsf{t}, N^5]$. For any $z \in \mathfrak{X}$ satisfying
	\begin{flalign}
		\label{dzhq}
		\mathtt{d}_z \ge N^{-K\delta} \mathtt{d}; \qquad \mathtt{d}_z^{1/2} \mathtt{e}_{z;q'}^{3/2} \le N^{K\delta-1}, \qquad \text{and} \qquad N^{K\delta} (q-q') \le N^{-2/3} \mathtt{d}_z^{-1/3},
	\end{flalign}
	\noindent we have
	\begin{flalign}
		\label{htzq2} 
		\begin{aligned}
			& H_t (z) \le \mathcal{H}_{q'} (z) + N^{\delta/4-1}, \qquad \text{and} \qquad H_t (z) \le \mathcal{H}_{q'} (z),  \quad \text{if $z \notin \mathfrak{L}_{q'}^+ (\delta)$}.
		\end{aligned} 
	\end{flalign} 
\end{lem}

\subsection{Proof of \Cref{thm:fN}} 

\label{Proofq0n}

To show \Cref{thm:fN}, we will use the results from Section \ref{ReduceAnnuli0} in an inductive way, to successively improve the upper bound on the height function under the Glauber dynamics (in suitably chosen time intervals), to reduce it from the initial $\cH_{q_0}$ to $\cH_{q_{i_{\max}}}$, with $q_{i_{\max}}$ so small that $\cH_{q_{i_{\max}}}$ is barely distinguishable from the equilibrium height. The three results from that section that we will use are \Cref{ak1akak1}, which applies in the bulk liquid phase; \Cref{akl0}, which applies in the edge phase not at the arctic boundary; and \Cref{boundaryaj2}, which applies at the arctic boundary. 

All of these decrease the value of $q$ in the function $\mathcal{H}_q$ upper bounding the Glauber dynamics. It will be beneficial to decrease $q$ by as much as possible under each application of the above results; this is limited by \Cref{ak1akak1}. Ignoring factors of $N^{\delta}$, after taking $k = \ell_0 (q) + 5$ to the maximum permitted by \eqref{eq:assu}, this decreases $q$ by about $2^{\ell_0(q)} N^{-1}$, which is by \eqref{adelta4} nearly $q^{1/3} N^{-2/3}$. This suggests the following definition for a sequence of $(q_i)$ that will serve as the parameters in the $\mathcal{H}_q$ upper bounding the Glauber dynamics.

\begin{definition} 

\label{qi} 
Fix constants $M,K$ satisfying \eqref{eq:costanti}. Define the sequence of real numbers $q_0, q_1, \ldots $ by, for each $i \ge 0$, setting $q_0 = \varepsilon_0$ and
\begin{flalign}
	\label{qdelta} 
	q_{i+1}= q_i - \Delta_i, \qquad \text{where} \qquad \Delta_i = q_i^{1/3} N^{-2/3-M\delta}.	
\end{flalign}

\noindent We will always work with values of $i$ for which $q_i>0$.

\end{definition}

Now, recall that we seek for $\mathcal{H}_q$ to be ``indistinguishable'' from $\mathcal{H}$ for some sufficiently small $q$; in our context, this means $|\mathcal{H}_q - \mathcal{H}|$ is at most $N^{-1}$, again ignoring factors of $N^{\delta}$. Since $\mathcal{H}_q$ behaves differently across $\mathfrak{X}$, the value of $q$ for which this ``indistinguishability,'' or ``stabilization,'' occurs varies across $\mathfrak{X}$. By Lemma \ref{qlineq}, this is governed by the annuli from Definition \ref{aq}; on $\mathcal{A}_k$, this ``stabilizing'' value of $q$ is about $N^{-1} 2^k$. This quantity ranges from $N^{-1}$ (in the center of the hexagon) to $N^{-1/2}$ (at the tangency locations). 

Thus, we will analyze the Glauber dynamics on a sequence of times $(\mathsf{T}_{i,j})$. Here, $j$ will denote the index of the annulus that we are attempting to stabilize during the time sequence $(\mathsf{T}_{i,j})_i$; as such, it will range from $1$ to the integer $\mathfrak{j}$, for which $N^{-1} 2^{\mathfrak{j}}$ is approximately its maximal value of $N^{-1/2}$. Since the outermost annuli $\mathcal{A}_k$ (with the largest values of $2^k N^{-1}$) will stabilize first, one should imagine $j$ as decreasing from $\mathfrak{j}$ to $1$. Then $i$ will index the $q_i$ for which the height function of our Glauber dynamics $(H_t)$ is bounded by $\mathcal{H}_{q_i}$ at times $t \ge \mathsf{T}_{i,j}$ (see \Cref{hqihqi1}). It will range from $0$ to the value $a_j$, for which $q_{a_j}$ is approximately $2^j N^{-1}$. Lemma \ref{ak1akak1} implies that the amount of time it takes for $q_i$ to decrease to $q_{i+1}$ is (at least in the bulk liquid region) around $q_i^{-1} (q_i - q_{i+1}) N^2$, which is nearly $q_i^{-2/3} N^{4/3}$. As such, we will have that $\mathsf{T}_{i+1,j} - \mathsf{T}_{i,j}$ is about $q_i^{-2/3} N^{4/3}$. 

These notions are made specific through the following definition.

\begin{definition} 

\label{tij} 

Let $\mathfrak{j}:=\lfloor (\frac12+\delta)\log_2 N\rfloor$ denote the largest integer such that $2^{\mathfrak{j}} \le N^{1/2+\delta}$. For each $j \in [0, \mathfrak{j}]$, let $a_j$ denote the minimal integer with $q_{a_j} < \max \{2^j N^{5000K\delta-1}, N^{M\delta-1} \}$; observe that the sequence $(a_j)$ is non-increasing in $j$. For each $j \in [0, \mathfrak{j}]$ and $i \in [0, a_j]$, define the sequence $\mathsf{T}_{i,j}$ inductively (starting at $(i,j) = (0, \mathfrak{j})$, and then increasing in $i \in [0, a_j]$ and decreasing in $j \in [0, \mathfrak{j}]$) by setting 
\begin{flalign*} 
	\mathsf{T}_{0,\mathfrak{j}} = N^{2 + M^2 \delta}; \quad \mathsf{T}_{i+1,j} = \mathsf{T}_{i,j} + q_i^{-2/3} N^{4/3 + M^2\delta}; \quad \mathsf{T}_{0,j} = \mathsf{T}_{a_{j+1},j+1},
\end{flalign*} 
observing that
\begin{eqnarray}
	\label{eq:obsthat}
	\mathsf{T}_{i+1,j} = \mathsf{T}_{i,j} + \frac{q_i-q_{i+1}}{ q_i} \cdot N^{2+ (M^2+M) \delta}.  
\end{eqnarray}
\end{definition}

\begin{rem}

\label{qideltai} 

Observe that $q_i$ and $\Delta_i$ are decreasing in
$i$, and that $N^{M \delta} \Delta_i \le q_i$ for each
$i \in [1, a_0]$ (as $q_i \ge N^{M\delta-1}$). Hence,
$q_{i+1} \ge (1 - N^{-M\delta}) q_i$, meaning that
$\Delta_{i+1} \ge (1 - N^{-M \delta}) \Delta_i$.  Note
also that $\mathsf{T}_{i,j}$ is increasing with $i$ but
nonincreasing with $j$ (specifically, $\mathsf{T}_{i,j}\ge \mathsf{T}_{i',j+1}$
for every $i\le a_j$ and $i'\le a_{j+1}$). \end{rem}

\begin{lem} 

\label{qiestimate} 
We have for each $j \in [0, \mathfrak{j}]$ and $i \in [0, a_j-1]$ that 
\begin{flalign} 
	\label{m0qm0} 
	\begin{aligned} 
		& (1 - N^{-M\delta}) q_{i} \le q_{i+1} \le  q_{i}; \qquad \quad q_{a_j} \ge 2^{j-1} N^{5000K\delta-1}; \\
		& a_j \le a_0 \le N^{2/3+2M\delta}; \qquad \qquad \qquad \quad \mathsf{T}_{a_j, j} \le N^{2+2M^2 \delta}. 
	\end{aligned}
\end{flalign} 

\end{lem}

\begin{proof}
The first statement follows from \Cref{qideltai}. This in particular yields $q_{a_j} \ge q_{a_j-1}/2$; together with the fact that $q_{a_j-1} \ge 2^j N^{5000K\delta-1}$ (by the minimality of $a_j$ in its definition), this implies the second bound in \eqref{m0qm0}. The fact that the $(a_j)$ are non-increasing implies the first estimate in the third bound of \eqref{m0qm0}.  

A heuristic for the second estimate $a_0 \le N^{2/3 + 2M\delta}$ there is that one may view the $q_i$ as an approximate solution, at time $i$, of the differential equation $q'(t)=-q(t)^{1/3} \cdot N^{-2/3-M\delta}$, with $q(0)=\varepsilon_0$. Since this equation satisfies $q(t) \ll N^{-1}$ for $t \gg N^{2/3+M\delta}$, it is plausible that $a_0$, the first integer $i$ for which $q_i <  N^{5000K\delta-1}$, should be smaller than $ N^{2/3+2M\delta}$. 

To prove this more carefully, first observe for each integer $k \ge 0$ that there exist at most $2^{1-2k/3} N^{2/3+M\delta}$ indices $i$ (which are consecutive, as the $q_i$ are decreasing) for which $q_i \in [2^{-k-1}, 2^{-k}]$. Indeed, assume to the contrary that this is false. Then there exist indices $i_1 \le i_2$ with $i_2 - i_1 \ge 2^{1-2k/3} N^{2/3+M\delta}$ and $q_{i_1}, q_{i_1+1}, \ldots , q_{i_2} \in [2^{-k-1}, 2^{-k}]$. This yields
\begin{flalign*}
	0 \le q_{i_2} \le q_{i_1} - (i_2-i_1) \Delta_{i_2} \le q_{i_1} - 2^{1-2k/3} q_{i_2}^{1/3} \le 2^{-k} - 2^{2/3-k} < 0,
\end{flalign*}

\noindent which is a contradiction (where in the second inequality we used \eqref{qdelta}; in the third we used \eqref{qdelta} and the fact that $i_2 - i_1 \ge 2^{1-2k/3} N^{2/3+M\delta}$; and in the fourth we used that $q_{i_1}, q_{i_2} \in [2^{-k-1}, 2^{-k}]$). Hence, $a_0 \le N^{2/3+M\delta} \sum_{k=0}^{\infty} 2^{1-2k/3} \ll N^{2/3+2M\delta}$.

To confirm the fourth bound in \eqref{m0qm0}, observe that
\begin{flalign*} 
	\mathsf{T}_{a_j, j} & \le \mathsf{T}_{0,\mathfrak{j}} + N^{4/3+M^2\delta} \displaystyle\sum_{j=0}^{\mathfrak{j}} \displaystyle\sum_{i=0}^{a_j} q_i^{-2/3} \\
	& \le  N^{2+M^2\delta} + (\mathfrak{j}+1) N^{4/3 + M^2\delta} \displaystyle\sum_{k=0}^{N^{\delta}} 2^{1-2k/3} N^{2/3+M\delta} \cdot 2^{2k/3} \le N^{2+2M^2 \delta},
\end{flalign*} 	
where the first bound follows from the definition of $\mathsf{T}_{i,j}$; the second from applying the above fact, for any integer $k \ge 1$, there are at most $2^{1-2k/3} N^{2/3 + M\delta}$ many indices $i$ for which $q_i \in [2^{-k}, 2^{1-k}]$; and the third follows from performing the sum (together with the fact that $\mathfrak{j} \lesssim \log N$).
\end{proof} 

After each annulus stabilizes, it will be convenient to shift the Glauber dynamics height function by $N^{M\delta-1}$. Since there are at most $\log N \ll N^{\delta}$ annuli, this produces an overall shift of at most $N^{(M+1)\delta-1}$, which will be negligible for our purposes.

\begin{definition}
\label{def:Hs}

Let $(H_s)_{s \ge 0}$ denote the height function at time $s$ for the Glauber dynamic on $\mX$, with floor constraint $\cH_{-q_0}$ and ceiling constraint $\cH_{q_0}$ (as in \Cref{def:constrained}), started from the maximal configuration compatible with these constraints (recall \Cref{def:partorder}). Then, for each $j \in [0, \mathfrak{j}]$ and $w \in \mathfrak{X}$, let 
\begin{flalign}
	\label{eq:Hsj}
	H_{s,j} (w) = H_s (w) + (j-\mathfrak{j}) N^{M\delta-1}.
\end{flalign} 
\end{definition}

Under the above notation we have the following inductive upper bounds on $H_{s,j}$, to be shown in \Cref{ProofHm}, that quickly imply \Cref{thm:fN}.
\begin{prop} 

\label{hqihqi1} 

Fix integers $j_0 \in [0, \mathfrak{j}]$ and $i_0 \in [0, a_{j_0}]$. Assume that, for any integers $j \in [j_0, \mathfrak{j}]$ and $i \in [0, a_j]$ such that either $j > j_0$ or $i < i_0$, the following holds w.o.p.: 
\begin{flalign}
	\label{htm}
	\begin{aligned}
		& H_{\mathsf T_{i,j},j} (w) \le \mathcal{H}_{q_{i}} (w) + N^{\delta-1}, \qquad \text{for each $w \in \mathfrak{X}$}; \\
		& H_{\mathsf T_{i,j},j} (w) \le \mathcal{H}_{q_{i}} (w), \qquad \qquad  \quad  \text{for each $w \notin \mathfrak{L}_{q_i}^+ (\delta)$}.
	\end{aligned}
\end{flalign}

\noindent Then w.o.p., we have for any $t \in [\mathsf{T}_{i_0, j_0}, N^5]$ that 
\begin{flalign}
	\label{hm1z} 
	\begin{aligned}
		& H_{t,j_0} (z_0) \le  \mathcal{H}_{q_{i_0}} (z_0) + N^{\delta-1}, \qquad \text{for each $z_0 \in \mathfrak{X}$}; \\ 
		& H_{t,j_0} (z_0) \le \mathcal{H}_{q_{i_0}} (z_0), \qquad \qquad \quad    \text{for each $z_0 \notin \mathfrak{L}_{q_{i_0}}^+ (\delta)$}.
	\end{aligned} 
\end{flalign}

\end{prop} 
\begin{rem}
\label{rem:fza}

By \Cref{prop:PW}, the law of $H_s $ is in $\mathcal P^+$, for any $s \ge 0$. Thus \Cref{prop:PW} implies that, if the bounds \eqref{htm} hold w.o.p. at time $\mathsf T_{i,j}$, then the same inequalities hold w.o.p. for every $t\in [\mathsf T_{i,j}, N^5]$.
\end{rem}

\begin{proof}[Proof of \Cref{thm:fN}] Recall from \Cref{def:Hs} that $(H_t)_{t\ge0}$ denotes the dynamic constrained between $\cH_{\pm q_0}$.
By \Cref{hqihqi1} and \Cref{rem:fza}, induction on increasing $i\ge0$ and decreasing $j\le \mathfrak j$, \eqref{m0qm0}, and the definition \eqref{eq:Hsj} of $H_{t,j}$, we have w.o.p. that $(H_t)_{t \ge 0}$ satisfies
\begin{eqnarray}
	H_{t} (z) \le \mathcal{H}_{q_{a_0}} (z) + \mathfrak{j} N^{M\delta-1} + N^{\delta-1} \le \mathcal{H}_{q_{a_0}} (z) + N^{(M+1)\delta-1}
\end{eqnarray}
for every $z\in\mX$ and every
$t\in [\mathsf T_{a_0,0},N^5]$.
Note that
$q_{a_0} \le N^{M\delta-1}$ (by definition of $a_j$) and
$\mathcal{H}_{q_{a_0}} (z) \le \mathcal{H} (z) + N^{\delta}
q_{a_0}\le \mathcal{H} (z) + N^{(M+1)\delta-1}
$ (by \Cref{qlineq}).
In conclusion, w.o.p. the height function $H_t$ satisfies $H_t \le \cH^+:=\min\{ \cH+ 2N^{(M+1)\delta-1},\cH_{q_0} \}$ for all $t\in [\mathsf T_{a_0,0},N^5]$. The analogous lower bound $H_t\ge \cH^-:=\max \{ \mathcal{H} - 2N^{(M+1)\delta-1},\cH_{-q_0} \}$ follows by symmetry. Recall also that $\mathsf T_{a_0,0}\le N^{2+2M^2\delta}$ from \eqref{m0qm0}.
For some sufficiently large constant $C>1$, we have then the following facts:
\begin{enumerate}
	\item $\cH^-\le H_t\le \cH^+$, w.o.p., for all $t\in [N^{2+2M^2\delta},N^5]$. This allows to couple (w.o.p.) the dynamic $(H_t)_{t \ge 0}$ with the Glauber dynamics with ceiling constraint $\mathcal{H}^+$ and floor constraint $\mathcal{H}^-$, in the whole interval $t\in [N^{2+2M^2\delta},N^5]$.
	
	\item the mixing time of the Glauber dynamics, with ceiling constraint $\mathcal{H}^+$ and floor constraint $\mathcal{H}^-$, is at most $CN^{2+2(M+2)\delta}$, by \Cref{prop:tmixconstrained};
	
	\item the mixing time $t^{q_0,N}_{\rm mix}$ of the Glauber dynamics, with ceiling constraint $\mathcal{H}_{q_0}$ and floor constraint $\mathcal{H}_{-q_0}$, is at most $CN^{4+\delta}\ll N^5$, also by \Cref{prop:tmixconstrained}.
\end{enumerate}

From items (1)--(3) and the submultiplicativity \eqref{e:submult} we deduce that, at time $t=N^5$, the law of the configuration $H_t$ is, on the one hand, at total variation distance at most $N^{-2}$ from the uniform measure $\pi_{\cH^{\pm}}$ and, on the other hand, at total variation distance at most $N^{-2}$ from the uniform measure $\pi_{\cH_{\pm q_0}}$ (recall \Cref{def:P+}). It follows that $\|\pi_{\cH^{\pm}}-\pi_{\cH_{\pm q_0}}\|\le N^{-1}$. Then, item (2) (with item (1)) implies that 	$t^{q_0,N}_{\rm mix}(1/4+N^{-1})\le N^{2+4M^2\delta}+CN^{2+2(M+2)\delta}\le  2 N^{2+4M^2\delta}$. 	Finally, applying once more the submultiplicativity \eqref{e:submult} of the mixing time we conclude 
\begin{equation}
	t^{q_0,N}_{\rm mix}(N^{-\mathsf L})\lesssim N^{2+4M^2\delta}|\log (N^{-\mathsf L})|  \lesssim N^{2+5M^2\delta}
\end{equation}
which proves the theorem. 
\end{proof}

\subsection{Proof of \Cref{hqihqi1}}

\label{ProofHm}

Recall from \Cref{qi} that $q_0=\varepsilon_0$ and that the sequence $(q_i)_i$ is decreasing.  We may assume that $i_0$ is sufficiently large (and in particular $i_0 \ge 1$) so that 
$q_{i_0} \le \varepsilon_0 / 2$, since for  $q_{i_0} \ge \varepsilon_0/2$ the statement we want to prove follows from \Cref{thm:step1} (with the fact that $H_{s,j} (w) \le H_s (w)$ for all $j \in [0,\mathfrak{j}]$ and $w \in \mathfrak{X}$).

Now fix $z_0 \in \mathfrak{X}$, and set $\mathfrak{q} = q_{i_0}$ and $\ell_0 = \ell_0 (\mathfrak q)$ (recall \eqref{adelta4}).

\begin{lem}
If $j_0 < \mathfrak{j}$ and $\mathfrak{q}\not \in [2^{j_0-1} N^{5000K\delta-1}, 2^{j_0+1} N^{5000K\delta-1}]$,  then the statement of \Cref{hqihqi1} holds.
\end{lem}
\begin{proof}
Note that $H_{s,j}$ is increasing in $j$, and that $\mathsf T_{i,j}$ is decreasing in $j$, so the statement \eqref{hm1z} is weaker for smaller values of $j$ (for the same $i$). Therefore (recalling also \Cref{rem:fza}), if $j_0<\mathfrak j$ \eqref{hm1z} follows from \eqref{htm} applied at $(i,j) = (i_0, j_0+1)$, except if the condition $i \le a_j$ is violated when we increase $j$ from $j_0$ to $j_0+1$, namely, if $i_0 > a_{j_0+1}$. 
Since $q_i$ is decreasing in $i$, this implies that $\mathfrak{q} = q_{i_0} < q_{a_{j_0+1}}$.  
Since $i \le a_{j_0}$, we always have $\mathfrak{q} \ge q_{a_{j_0}}$, so this means 
$q_{a_{j_0}} \le \mathfrak{q} < q_{a_{j_0+1}},$ 
which upon recalling the definition of $a_j$ (and the first statement in \eqref{m0qm0}) implies that $\mathfrak q \in [2^{j_0-1} N^{5000K\delta-1}, 2^{j_0+1} N^{5000K\delta-1}]$.  
\end{proof}
From now on, we assume therefore that
\begin{eqnarray}
\label{eq:fromn}
j_0=\mathfrak j, \quad \text{or} \quad \mathfrak q \in [2^{j_0-1} N^{5000K\delta-1}, 2^{j_0+1} N^{5000K\delta-1}].
\end{eqnarray}

The following lemma verifies \Cref{hqihqi1} if $z_0$ is close to the arctic boundary.

\begin{lem} 

\label{dz0z0h} 

If $\mathtt{d}_{z_0}^{1/2} \mathtt{e}_{z_0;\mathfrak{q}}^{3/2} \le N^{K\delta-1}$, then the statement of \Cref{hqihqi1} holds.

\end{lem}

\begin{proof} 

We distinguish between whether $\mathfrak{q} \le N^{4000K\delta-2/3} \mathtt{d}_{z_0}^{-1/3}$ or $\mathfrak{q} > N^{4000K\delta-2/3} \mathtt{d}_{z_0}^{-1/3}$. In the former case, we will show that \eqref{hm1z} follows from \eqref{htm} at a larger value of $j > j_0$. To that end, let us first verify that $j_0 < \mathfrak{j}$. Indeed, we have $N^{4000K\delta-2/3} \mathtt{d}_{z_0}^{-1/3} \ge \mathfrak{q} \ge q_{a_{j_0}} \ge 2^{j_0-1} N^{5000K\delta-1}$, so that $\mathtt{d}_{z_0}^{-1/3} \ge 2^{j_0-1} N^{1000K\delta-1/3}$. Since $2^{\mathfrak{j}} \ge N^{1/2}$ and $\mathtt{d}_{z_0} \ge N^{-1/2}$, it follows that $j_0 < \mathfrak{j}$ and thus $\mathfrak{q} \ge 2^{j_0-1} N^{5000K\delta-1}$ because of \eqref{eq:fromn}.

Now, set $j_0' = j_0+1$, and observe that $q_{a_{j_0'}} \le 2^{j_0+1} N^{5000K\delta-1}$. Hence, we have w.o.p. that
\begin{flalign*}
	H_{t,j_0} (z_0) & = H_{t,j_0+1} (z_0) - N^{M\delta-1} \\
	& \le \mathcal{H}_{q_{a_{j_0'}}} (z_0) -
	\frac12 \cdot N^{M\delta-1} \\
	& \le \mathcal{H}_{\mathfrak{q}} (z_0) + N^{\delta} (q_{a_{j_0'}} - \mathfrak{q}) (\mathtt{d}_{z_0} \mathtt{e}_{z_0;\mathfrak{q}})^{1/2} + N^{\delta} \mathtt{d}_{z_0}^{1/2} q_{a_{j_0'}}^{3/2} -\frac12 \cdot N^{M\delta-1}\\ 
	& \le \mathcal{H}_{\mathfrak{q}} (z_0) + N^{5100K\delta-1} 2^{j_0+1} (\mathtt{d}_{z_0}  \mathtt{e}_{z_0;\mathfrak{q}})^{1/2} +  N^{7600K\delta-1}  - \frac12 \cdot N^{M\delta-1},
\end{flalign*}		
where in the first statement we used the definition of $H_{t,j}$; in the second we used \eqref{htm} together with \Cref{rem:fza}; in the third we used \Cref{qlineq};
and in the fourth we used the above fact that $q_{a_{j_0'}} \le 2^{j_0+1} N^{5000K\delta-1}$ and that $q_{a_{j_0'}} \le 4\mathfrak{q} \le N^{5000K\delta-2/3} \mathtt{d}_{z_0}^{-1/3}$. Since $\mathtt{d}_{z_0}^{1/2} \mathtt{e}_{z_0;\mathfrak{q}}^{3/2} \le N^{K\delta-1}$ and $\mathtt{d}_{z_0}^{-1/3} \ge 2^{j_0-1} N^{1000K\delta-1/3}$, we have $\mathtt{d}_{z_0}^{1/2} \mathtt{e}_{z_0}^{1/2} \le 2^{1-j_0} N^{-900K\delta}$, and so
\begin{flalign*}
	H_{t,j_0} (z_0) \le \mathcal{H}_{\mathfrak{q}} (z_0) + N^{8000K\delta-1} -\frac12 \cdot N^{M\delta-1} \le \mathcal{H}_{\mathfrak{q}} (z_0),
\end{flalign*}
which verifies the lemma if $\mathfrak{q} \le N^{4000K\delta-2/3} \mathtt{d}_{z_0}^{-1/3}$ (recall that $M>10^5 K$).

Now instead suppose that $\mathfrak{q} \ge N^{4000K\delta-2/3} \mathtt{d}_{z_0}^{-1/3}$; we will deduce \eqref{hm1z} from \eqref{htm} at a smaller value of $i < i_0$, with \Cref{boundaryaj2}. Set 
\begin{flalign} 
	\label{rmqm1} 
	r = \mathfrak{q} + N^{-2/3-K\delta} \mathtt{d}_{z_0}^{-1/3}\le 2\mathfrak q;
\end{flalign}
the bound holds since $\mathfrak{q} \ge N^{4000K\delta-2/3} \mathtt{d}_{z_0}^{-1/3}$.
\noindent Further let $m \le i_0$ be the smallest integer such that $q_{m} < r$. Observe that such an integer $m$ exists since $r < \mathfrak{q}  +o(1) \le \varepsilon_0/2 + o(1) < \varepsilon_0$ and $m < i_0$ since $r - \mathfrak{q} = N^{-K\delta-2/3} \mathtt{d}_{z_0}^{-1/3} \gg \Delta_{m-1}$ (by the definition \eqref{qdelta} of $\Delta_{m-1}$).

We will apply \Cref{boundaryaj2} with the $(\mathtt{d}; z; q, q'; H_s)$ there equal to $(\mathtt{d}_{z_0}; z_0; r, \mathfrak{q}, H_{s+\mathsf{T}_{m,j_0}})$ here. These parameters satisfy \eqref{dzhq}, where \eqref{htm} (with the $(i_0, j_0)$ there equal to $(m,j_0)$ here) verifies \eqref{hsz2}. Therefore, \Cref{boundaryaj2} implies w.o.p. that $H_{t,j_0} (z_0) \le \mathcal{H}_{\mathfrak{q}} (z_0) + N^{\delta/4-1}$ and $H_{t,j_0} (z_0) \le \mathcal{H}_{\mathfrak{q}} (z_0)$ if $z_0 \notin \mathfrak{L}_{\mathfrak{q}}^+ (\delta)$ for $\mathsf{T}_{m,j_0} + \mathfrak{q}^{-1} \mathtt{d}_{z_0}^{-1/3} N^{4/3+15000K\delta} \le t \le N^5$. Since $\mathsf{T}_{i_0,j_0} - \mathsf{T}_{m,j_0} \ge (i_0-m) r^{-1} N^{2 + M^2 \delta} \Delta_{i_0} \ge (i_0 - m) (2 \mathfrak{q})^{-1} N^{2+M^2 \delta} $  and $i_0-m \ge (r-\mathfrak{q}) \Delta_m^{-1}  \ge (r - \mathfrak{q}) (2\Delta_{i_0})^{-1}$ (as $r \le 2 \mathfrak{q}$ by \eqref{rmqm1}), we have that $\mathsf{T}_{i_0,j_0} - \mathsf{T}_{m,j_0} \gg \mathfrak{q}^{-1} \mathtt{d}_{z_0}^{-1/3} N^{4/3+ 15000K \delta}$, and so this establishes \eqref{hm1z}.
\end{proof}

The next lemma verifies \Cref{hqihqi1} if $z_0$ is in the bulk liquid phase.

\begin{lem}

\label{z0k}

If $z_0 \in \bigcup_{j=1}^{\ell_0+4} \mathcal{A}_j^{\mathfrak{q}}$, then the statement of \Cref{hqihqi1} holds.

\end{lem}   

\begin{proof}

Recalling $D$ from \eqref{eq:costanti}, set $\mathfrak{q}' = q_{i_0-1} - 2^{\ell_0} N^{D\delta-1}$. We will apply \Cref{ak1akak1}, with the $(k;q,q';H_s)$ 
there equal to $(\ell_0+5; q_{i_0-1} ,\mathfrak{q}'; H_{s+\mathsf{T}_{i_0-1,j_0}})$ 
here. These parameters satisfy \eqref{qkq}, and \eqref{htm} implies \eqref{htzq0}. Hence, \Cref{ak1akak1} yields w.o.p. that $H_t (z) \le \mathcal{H}_{\mathfrak{q}'} (z)$ for all $z \in \bigcup_{j=1}^{\ell_0+5} \mathcal{A}_j^{\mathfrak{q}'}$ and $ \mathsf{T}_{i_0-1,j_0} + (\mathfrak{q}')^{-1} (\ell_0+5) 2^{\ell_0+5} N^{1+10A\delta}\le t\le N^5$. Note that $\mathfrak{q} \ge \mathfrak{q}'$, since 
\begin{flalign*} 
	\mathfrak{q} - \mathfrak{q}' = 2^{\ell_0} N^{D\delta-1} - q_{i_0-1}^{1/3} N^{-2/3-M\delta} \ge \mathfrak{q}^{1/3} N^{-2/3-A\delta/2} - q_{i_0-1}^{1/3} N^{-2/3-M\delta} > 0,
\end{flalign*}
where in the first statement we used the
definitions of $\mathfrak{q}$ and $\mathfrak{q}'$; in
the second we used \eqref{adelta4}; and in the third
we used the fact from \Cref{qideltai} that
$2 \mathfrak{q} > q_{i_0-1}$, together with the condition $M \ge 10^6 A$ (by \eqref{eq:costanti}).

This implies that
$H_t (z) \le \mathcal{H}_{\mathfrak{q}} (z)$ for all
$z \in \bigcup_{j=1}^{\ell_0+5}
\mathcal{A}_j^{\mathfrak{q}'}$ and thus for all
$z \in \bigcup_{j=1}^{\ell_0+4}
\mathcal{A}_j^{\mathfrak{q}}$ (by \eqref{bk1u} and
\Cref{hqa}, where we used the fact that
$\mathfrak{q} - \mathfrak{q}' \le 2^{\ell_0} N^{D\delta-1} \le
2^{\ell_0+5} N^{A\delta - 4\delta - 1}$). Thus, to show the lemma, it suffices to verify that $\mathsf{T}_{i_0-1,j_0} + (\mathfrak{q}')^{-1} (\ell_0+5) 2^{\ell_0+5} N^{1+10A\delta} \le \mathsf{T}_{i_0,j_0}$. This holds since 
\begin{flalign*}
	\mathfrak{q}'^{-1} (\ell_0 + 5) 2^{\ell_0 + 5}
	N^{1+10A\delta} & \le 2 q_{i_0-1}^{-1} (\ell_0 + 5) 
	\mathfrak{q}^{1/3} N^{4/3+10A\delta} \\
	& \le q_{i_0-1}^{-2/3}
	N^{4/3+7A\delta} \le \mathsf{T}_{i_0,j_0} -
	\mathsf{T}_{i_0-1,j_0},
\end{flalign*}

\noindent by \eqref{adelta4} and the bounds
$q_{i_0-1} \le 2\mathfrak{q}'$ (as $q_{i_0-1} \ge \mathfrak{q} \ge 2^{\ell_0+1} N^{D\delta-1}$, where the last inequality holds since $(\mathfrak{q}N)^{1/3} \ge 2^{\ell_0+1} N^{A\delta/2}$, by \eqref{adelta4} and \eqref{eq:assu}), $\mathfrak{q}\le q_{i_0-1}$, and $M \ge 10^{10} A$.
\end{proof}

\begin{proof}[Proof of \Cref{hqihqi1}]

By Lemmas \ref{dz0z0h} and \ref{z0k}, it remains to consider the case when $z_0 \in \mathcal{A}_k^{\mathfrak{q}} (N^{K\delta})$ for some $k \ge \ell_0 + 5=\ell_0(\mathfrak q)+5$. The proof will be similar to that of \Cref{dz0z0h}; as there, we will distinguish between whether $\mathfrak{q} \le 2^k N^{K\delta-1}$ or $\mathfrak{q} \ge 2^k N^{K\delta-1}$.

{\bf Case 1:  $\mathfrak{q} \le 2^k N^{K \delta-1}$.} We will verify that \eqref{hm1z} follows from \eqref{htm} at a larger value of $j > j_0$. To that end, we first confirm that $j_0 < \mathfrak{j}$. Observe that $2^{j_0-1} N^{5000K\delta-1} \le \mathfrak{q} \le 2^k N^{K\delta-1} \le 2 N^{K\delta-1} (\mathtt{d}_{z_0} \mathtt{e}_{z_0})^{-1/2} \le N^{K\delta-1/2}$, where in the first inequality we used \eqref{m0qm0} and the fact that $\mathfrak{q} = q_{i_0}$; in the second we used our assumption that $\mathfrak{q} \le 2^k N^{K\delta-1}$; in the third we used the fact that $z_0 \in \mathcal{A}_k^{\mathfrak{q}}$; and in the fourth we used the fact that $\mathtt{d}_{z_0} \mathtt{e}_{z_0;\mathfrak{q}} \gtrsim  \mathtt{d}_{z_0}^{1/2} \mathtt{e}_{z_0;\mathfrak{q}}^{3/2} \ge N^{K\delta-1}$ (by \eqref{eq:elegeo0} and the fact that $z_0 \in \mathcal{A}_k^{\mathfrak{q}} (N^{K\delta})$). Hence, $2^{j_0} \le N^{1/2-4900K\delta}$ meaning since $2^{\mathfrak{j}} \asymp N^{1/2+\delta}$ that $j_0 \le \mathfrak{j} - 1$.

Now, set $j' = j_0+1$. We claim (and prove below) that
\begin{equation}
	\label{eq:qle4q}
	q_{a_{j'}} \le 2^{j_0+1} N^{5000K\delta-1} \le 4\mathfrak{q} .    
\end{equation}
Thus, we have w.o.p. that 
\begin{flalign*}
	H_{t,j_0} (z_0) & = H_{t,j'} (z_0) - N^{M\delta-1} \\
	& \le \mathcal{H}_{q_{a_{j'}}} (z_0) + N^{\delta-1} - N^{M\delta-1} \\
	& \le \mathcal{H}_{\mathfrak{q}} (z_0) + N^{\delta} (q_{a_{j'}} - \mathfrak{q}) (\mathtt{d}_{z_0} \mathtt{e}_{z_0;\mathfrak{q}})^{1/2} + N^{\delta} \mathtt{d}_{z_0}^{1/2} q_{a_{j'}}^{3/2} - N^{M\delta-1} \\ 
	& \le \mathcal{H}_{\mathfrak{q}} (z_0) + N^{2\delta} 2^{3-k} \mathfrak{q} + 8 N^{\delta} \mathtt{d}_{z_0}^{1/2} \mathfrak{q}^{3/2} - N^{M\delta-1} \\
	& \le \mathcal{H}_{\mathfrak{q}} (z_0) +  N^{2K\delta-1} - N^{M\delta-1} \le \mathcal{H}_{\mathfrak{q}} (z_0),
\end{flalign*}
where in the first statement we used the definition of $H_{t,j}$; in the second we used the $(i,j) = (a_{j'}, j')$ case of \eqref{htm} (together with \Cref{rem:fza}); in the third we used \Cref{qlineq}; in the fourth we used \eqref{eq:qle4q}, with the fact that $z \in \mathcal{A}_k^{\mathfrak{q}}$; in the fifth we used the facts that $\mathfrak{q} \le 2^k N^{K\delta-1}$ and $\mathtt{d}_{z_0}^{1/2} \le 2^{-3k/2} N^{1/2-K\delta/2}$ (as $\mathtt{d}_{z_0}^{1/2} \mathtt{e}_{z_0;\mathfrak{q}}^{3/2} \ge N^{K\delta-1}$ and $\mathtt{d}_{z_0} \mathtt{e}_{z_0} \le 4^{-k}$ for $z_0 \in \mathcal{A}_k^{\mathfrak{q}}$); and in the sixth we used the fact that $M \ge 10^5 K$. 

To establish the proposition when $\mathfrak{q} \le 2^k N^{K\delta-1}$, it remains to verify \eqref{eq:qle4q}.
Recall that $a_{j'}$ is the minimum integer such that $q_{a_{j'}}<\max \{ 2^{j'} N^{5000K \delta-1},N^{M\delta-1} \}$. If the maximum is realized by the first alternative, then we directly have $q_{a_{j'}}<2^{j'} N^{5000K \delta-1}\le 4\mathfrak q$ using \eqref{m0qm0}. So, assume instead that $2^{j'} N^{5000K\delta-1} \le q_{a_{j'}} < N^{M\delta-1}$. Then, since $j' = j_0+1$, we have  $2^{j_0} N^{5000K\delta-1} \le q_{a_{j'}} \le N^{M\delta-1}$. Therefore, since $a_{j'}$ is the minimal integer for which $q_{a_{j'}} < N^{M\delta-1}$, it follows that $a_{j'} = a_{j_0}$, so $q_{a_{j_0}} = q_{a_{j'}}$, meaning that 
\begin{equation}
	\mathfrak{q} = q_{i_0} \ge q_{a_{j_0}} = q_{a_{j'}},  
\end{equation}
where the first statement holds by the definition of $\mathfrak{q}$; the second since $q_{i}$ is decreasing in $i$; and the third holds by the above. This verifies the claim.

{\bf Case 2: $\mathfrak{q} \ge 2^k N^{K \delta-1}$.} We will use
\eqref{htm} applied at a smaller value of $i < i_0$, together with 
\Cref{akl0}.  More specifically, recalling $D$ from \eqref{eq:costanti}, let $m < i_0$ be any integer such
that
$q_{m} - \mathfrak{q} \in [2^k N^{D\delta-1}, 2^{k+2}
N^{D\delta-1}]$. In particular, $\mathfrak q\ge 2^k N^{K \delta-1}\ge q_m/2$, since $K > 10^{15}$ by \eqref{eq:costanti}. To see that such an integer $m$ exists, note first of all 
that
\begin{eqnarray}
	\label{eq:trt}
	q_{i_0-1} - \mathfrak{q}=   q_{i_0-1} -    q_{i_0}\asymp
	\mathfrak q^{1/3}N^{-2/3-M\delta}\ll 2^k N^{-1+D\delta}
\end{eqnarray}

\noindent by the definition \eqref{qdelta} of the $(q_i)$ (with \Cref{qideltai}), and the assumption $k>\ell_0(\mathfrak q)+4$ (with \eqref{adelta4}). This, together with the fact that $q_j-q_{j+1} \asymp q_{i_0-1}-q_{i_0} = q_{i_0-1} - \mathfrak{q}$ whenever $\mathfrak{q} \le q_j \le 2 \mathfrak{q}$ (which holds by \eqref{qdelta}) yields the existence of such an integer $m$.

Let us apply \Cref{akl0}, with the $(q, q'; H_s)$ there equal to $(q_{m}, \mathfrak{q}; H_{s+\mathsf{T}_{m,j_0}})$ here, where we observe that \eqref{hszq2} is verified by the $(i,j) = (m, j_0)$ case of \eqref{htm}. This implies that, w.o.p.,  $H_{s,j_0} (z_0) \le \mathcal{H}_{\mathfrak{q}} (z_0)$ for all \[\mathsf{T}_{m,j_0} + q_m^{-1} k^2 2^{k+3} N^{1+11K\delta} \le s \le N^5.\] The proof is then concluded by noting that   $\mathsf{T}_{i_0,j_0} - \mathsf{T}_{m,j_0} \gg k^2 2^{k+3} q_{m}^{-1} N^{1+11K\delta}$, as
\begin{flalign*} 
	\mathsf{T}_{i_0, j_0} - \mathsf{T}_{m, j_0} \gtrsim  (q_{m} - \mathfrak{q}) \mathfrak q^{-1} N^{2+M\delta} \ge \mathfrak{q}^{-1} 2^{k} N^{D\delta-1} \cdot N^{2+M\delta} \gg q_{m}^{-1} k^2 2^{k+3} N^{1+11K\delta}.
\end{flalign*}
where we used \eqref{eq:obsthat} and the bounds $q_m \le 2\mathfrak{q}$ and $M>10^5 K$.
\end{proof}

\section{Height reduction in the bulk liquid phase} 

\label{HeightReduce00}

\subsection{Proof of \Cref{ak1akak1}} 

\label{ProofAnnulus00}

\begin{figure}[h]
 \begin{center} \includegraphics[width=9cm]{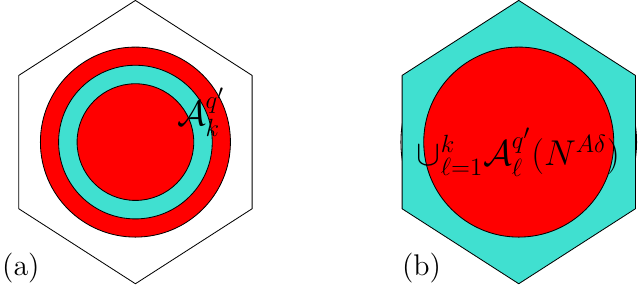}
  \caption{{(a): a schematic illustration of \Cref{aj0}. We assume an upper bound on the height function at time zero in  $\mathcal C^q_k$ (union of red and blue regions) and we deduce a stronger upper bound at time $\mathtt t$ in the annulus $\mathcal A^{q'}_k$ (blue region). (b): a schematic illustration of \Cref{ajnear}.  We assume an upper bound on the height function at time zero in the whole hexagon, and a better upper bound in the red region $\bigcup_{j=1}^k \mathcal A^{q'}_j(N^{A\delta})$ (red region) and we deduce an improved upper bound at time $\mathtt t$ in the blue region.}}
\label{fig:61}
\end{center}
\end{figure}

In this section we establish \Cref{ak1akak1} using the below lemma, to be shown in Section \ref{sec61prec}. It explains how the height function under the Glauber dynamics reduces in a single annulus (corresponding to the outermost one in the union of annuli considered in Lemma \ref{ak1akak1}) in the bulk liquid phase. Observe that the reduction $q-q'$ of $q$ and time $\mathsf{t}$ given by \eqref{q2q} below are similar to those in \eqref{qkq}. See \Cref{fig:61}.

\begin{lem} 
	\label{aj0}
	Let $A,D$ satisfy \eqref{eq:costanti}; let $q \in (0, \varepsilon_0)$ and $k\in \mathbb N$ be such that $k \le \ell_0 (q)+5$; and set
	\begin{flalign} 
		q' \in [q - 2^k  N^{D \delta-1}, q]; \qquad \mathsf{t} = 2^k q^{-1} N^{1+6A\delta}.
		\label{q2q}
	\end{flalign} 
	Consider the Glauber dynamics $(H_t)_{t \ge 0}$ starting from an initial condition $\mu\in\mathcal P^+$ (recall \Cref{def:P+}) and assume that, w.o.p. under $\mu$,
	\begin{flalign}
		\label{h0zb} 
		H_0 (w) \le \mathcal{H}_q (w) + N^{\delta-1}, \quad \text{for each $w \in \mathcal C^q_k:=\mathcal{B}_{k+1}^q \cup \bigcup_{j=1}^k \mathcal{A}_j^q$}. 
	\end{flalign}
	
	\noindent Then we have w.o.p. that, for any $t \in [\mathsf{t}, N^{5}]$,   
	\begin{flalign} \label{eq:follo}
		H_t (z_0) \le \mathcal{H}_{q'} (z_0), \qquad \text{for each $z_0 \in \mathcal{A}_k^{q'}$}.
	\end{flalign}  
	
\end{lem}

\begin{proof}[Proof of \Cref{ak1akak1}]
	To establish the lemma, we induct on $k \in [0, \ell_0+5]$.	The statement \eqref{htzq1} is empty for $k=0$; fix $k \in [1, \ell_0+5]$. Assume that  \eqref{qkq} and \eqref{htzq0} imply \eqref{htzq1} w.o.p., with the $k$ there replaced by $k-1$; we will then show that  \eqref{qkq} and \eqref{htzq0} w.o.p. imply \eqref{htzq1}. 
	To that end, since $\mu\in\mathcal P^+$, 	it suffices by \eqref{h0zb} and \Cref{prop:PW} to show \eqref{htzq1} when the Glauber dynamics $(H_t)_{t \ge 0}$ are replaced by the Glauber dynamics restricted  (in the sense of \Cref{def:restricted}) the region $\bigcup_{\ell=1}^{k+1}\mathcal A_\ell^q$, and with ceiling constraint given by  $\cH_q+N^{\delta-1}$ (as in Definition \ref{def:constrained}).  Thus, \Cref{aj0} (applied with $D$ there equal to $D+2$ here) yields, for any $r \in [q', q]$,
	\begin{flalign} 
		\label{htzsi2} 
		H_t (z) \le \mathcal{H}_{q'} (z), \qquad \text{for each $t \in [q'^{-1}  2^k N^{1+6A\delta}, N^5]$ and $z \in \mathcal{A}_k^r$},
	\end{flalign}
	w.o.p. (where we used the fact that $q\ge q'$). Set $\mathtt{s} = q'^{-1} 2^k N^{1+6A\delta}$. 
	
	By the $r = q'$ case of \eqref{htzsi2}, it remains to verify that $H_s (z) \le \mathcal{H}_{q'} (z)$ holds w.o.p., for each $s \in [\mathsf{t}, N^5]$ and $z \in \bigcup_{j=1}^{k-1} \mathcal{A}_{j}^{q'}$. This will follow from two applications of the inductive hypothesis. More specifically, set $\mathsf{q}_j = q - j(q-q') / 2$ for each $j \in \{ 0, 1, 2 \}$. Then, we have that $\mathsf{q}_0 = q$ and $\mathsf{q}_2 = q'$; that $\mathsf{q}_i \in \big[ \mathsf{q}_{i-1} - 2^{k-1} N^{(D+2)\delta-1}, \mathsf{q}_{i-1} - k 2^{k-1} N^{(D-2)\delta-1} \big]$ for $i \in \{ 1, 2 \}$ (where the inclusion follows from \eqref{qkq}). By \eqref{htzq0}, the inductive hypothesis \eqref{htzq1} (with the $(k, q, q')$ there given by $(k-1, \mathsf{q}_0, \mathsf{q}_1)$ here) yields, w.o.p.,
	\begin{flalign} 
		\label{htzsi3} 
		H_t (z) \le \mathcal{H}_{\mathsf{q}_1} (z), \quad \text{for $t \in \big[ \mathtt{s} +  \mathsf{q}_1^{-1} (k-1) 2^{k-1} N^{1+6A\delta}, N^5 \big]$ and $z \in \bigcup_{j=1}^{k-1} \mathcal{A}_{j}^{\mathsf{q}_1}$}.
	\end{flalign}
	Now by \eqref{htzsi2} (with the $r$ there equal to $\mathsf{q}_1$ here) and \eqref{htzsi3}
	the inductive hypothesis \eqref{htzq1} (with the $(k, q, q')$ there given by $(k-1, \mathsf{q}_1, \mathsf{q}_2)$ here) yields w.o.p. that $H_t (z) \le \mathcal{H}_{\mathsf{q}_2} (z) = \mathcal{H}_{q'} (z)$, for each $z \in \bigcup_{j=1}^{k-1} \mathcal{A}_j^{q'}$ and each $t \in [\mathtt{s}', N^5]$. Here, we have denoted 
	\begin{flalign*} 
		\mathtt{s}' = \mathtt{s} + \mathsf{q}_1^{-1} (k-1) 2^{k-1} N^{1+6A \delta} + \mathsf{q}_2^{-1} (k-1) 2^{k-1} N^{1+6A\delta} \le  q'^{-1} 2^k k N^{1+6A\delta} = \mathsf{t},
	\end{flalign*} 
	where in the last inequality we recalled the definitions of $\mathtt{s}$ and $\mathsf{t}$, and we also used the facts that $\mathsf{q}_1, \mathsf{q}_2 \ge q'$. This establishes the lemma.
	Note that, in applying the inductive hypothesis the second time, we have implicitly assumed that $k-1\le \ell_0(\mathsf{q}_1)+5$, which is needed in \eqref{eq:assu}.
	The claim follows if $\ell_0(\mathsf{q}_1)\ge \ell_0(q)-1$, since $k$ satisfies $k\le \ell_0(q)+5$. 
	First observe by \eqref{adelta4} and \eqref{eq:assu}
	\begin{flalign*}
		q \ge  (qN)^{1/3} N^{A\delta/2-1} \ge 2^{\ell_0 (q)+5} N^{A\delta-1} \ge 2^k N^{A\delta-1}.
	\end{flalign*}
	Hence, it follows from \eqref{qkq} that $\mathsf{q}_1 \le q \le 2\mathsf{q}_1$.
	Using again \eqref{adelta4} we deduce $\ell_0(\mathsf{q}_1)> \ell_0(q)-2$ which implies $\ell_0(\mathsf{q}_1)\ge \ell_0(q)-1$,
	since $\ell_0$ is an integer.
\end{proof}

        \subsection{Mixing in bulk liquid phase cells}
        
        \label{sec61prec}

In this section we show \Cref{aj0}, by reducing it to a ``more local'' height reduction statement on suitably defined mesoscopic rectangular cells. Throughout this section, we let $q \in (0, \varepsilon_0)$ and $z_0 \in \mathfrak{L}_{q}$ satisfy
\begin{eqnarray}
  \label{eq:distanzabis}
  \md_{z_0}\me_{z_0;q}\ge\mq^{-2/3}N^{-2/3+A\delta}.
\end{eqnarray}

\noindent Recall from \Cref{AProperty} that this assumes $z_0$ is in the bulk liquid phase.  Since $\me_{z_0;q}\md_{z_0}\le 1$,
\eqref{eq:distanzabis} implies
\begin{eqnarray}
  \label{eq:fnt}
  q\ge \frac{N^{-1+(3/2)A \delta}}{ (\md_{z_0}\me_{z_0;q})^{3/2}}\ge \frac{N^{-1+(3/2)A \delta}}{ (\md_{z_0}\me_{z_0;q})^{1/2}}.
\end{eqnarray}

The following definition provides the mesoscopic cells in which we will locally analyze the Glauber dynamics in the bulk liquid phase. The reasoning behind the specific choice of dimensions in \eqref{gothiccells} will be explained in \Cref{rem:whythecell} below. See also \Cref{fig:rescaledH}.
\begin{definition}\label{def:bulkcells}
Define  the cell $\mathfrak Y \subset \bar{\mathfrak{X}}$ centered at $z_0$ by
  \begin{flalign}
    \label{gothiccells}
    \begin{aligned}
  \mathfrak{Y} = \bigg\{ z \in \overline{\mathfrak{X}} : & \dist_{\bm{u}_{z_0;q}} (z, z_0) \le q^{-1/2} N^{(\frac34 A-\frac12)\delta -1/2} (\md_{z_0}\me_{z_0;q})^{-1/4},\\
  &  \dist_{\bm{w}_{z_0;q}} (z, z_0) \le q^{-1/2} N^{(\frac34A-\frac12)\delta -1/2} (\md_{z_0}\me_{z_0;q})^{1/4} \bigg\}.
  \end{aligned}
	\end{flalign} 

\noindent Further define $\mathfrak{Y}^- \subset \mathfrak{X}$ by  $\mathfrak{Y}^-=(1/L)\cdot(\mathfrak{Y}-z_0)+z_0$ (namely, to be its shrinking by a factor $L$), where $L$ is a sufficiently large constant that depends only on the constant $C$ in \Cref{prop:last}.  
\end{definition}

\begin{rem}

\label{yde} 

By \eqref{eq:distanzabis} and \eqref{gothiccells} (with the fact that $A \ge 30$), we have that 
\begin{flalign*}
	\dist_{\bm{u}_{z_0;q}} (z_0, \partial \mathfrak{Y}) \le N^{-\delta/2} (\md_{z_0}\me_{z_0;q})^{1/2}; \qquad \dist_{\bm{w}_{z_0;q}} (z_0, \partial \mathfrak{Y}) \le N^{-\delta/2} \md_{z_0}\me_{z_0;q}.
\end{flalign*}

\noindent By \eqref{eq:elegeo} and the smoothness of $\mathfrak{A}$ (\Cref{convex}), it quickly follows that, for any $z \in \mathfrak{Y}$, we have 
\begin{flalign}
	\label{ed0} 
	(1 - N^{-\delta/4}) \mathfrak{e}_{z_0;q} \le \mathfrak{e}_{z;q} \le (1+N^{-\delta/4}) \mathfrak{e}_{z_0;q}; \qquad 	(1 - N^{-\delta/4}) \mathfrak{d}_{z_0} \le \mathfrak{d}_{z} \le (1+N^{-\delta/4}) \mathfrak{d}_{z_0}. 
\end{flalign}

\noindent In particular, the cell $\mathfrak{Y} \subseteq \mathfrak{L}_q$ is contained in the liquid region. Observe that the derivation of \eqref{ed0} used (and in fact required) the condition \eqref{eq:distanzabis}, which is why we are using the latter to define the bulk liquid phase.

\end{rem} 

\begin{figure}[h]
 \begin{center} \includegraphics[width=5cm]{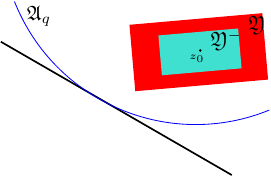}
   \caption{Proposition \ref{prop:hdecreasebulk} assumes the bound $H\le \cH_q+N^{\delta-1}$ in the cell $\mathfrak Y$ (the union of the blue and red regions), and deduces w.o.p. the stronger bound $H_t\le \cH_{\bar q}$ in the smaller cell $\mathfrak Y^-$ (the blue region), after time $\mathtt t$.}
\label{fig:Prop58}
\end{center}
\end{figure}

The following proposition, to be shown in Section \ref{sec61}, provides a rate of decrease for the height function of tilings under Glauber dynamics in $\mathfrak{Y}$. See \Cref{fig:Prop58}. Observe below that $\bar q\ge0$, by \eqref{eq:fnt} and \eqref{qdecreaseq}.

 \begin{prop}[Bulk mixing]\label{prop:hdecreasebulk}
   The following holds if the constant $L$ in \Cref{def:bulkcells} is large enough. Fix real numbers $q \ge \overline{q}$ and a point $z_0 \in \mathfrak{L}_{q}$ satisfying \eqref{eq:distanzabis} and  
\begin{flalign}
	\label{qdecreaseq} 
	0 \le q - \overline{q} < N^{5A\delta/4-1} (\mathfrak{d}_{z_0} \mathfrak{e}_{z_0;q})^{-1/2}.
\end{flalign}

\noindent Further assume that $H_0 (z) \le
        \mathcal{H}_{q} (z)+N^{\delta-1}$ for each $z \in 
        \mathfrak{Y}$. 
        Under the Glauber dynamics restricted to $\mathfrak Y$, with the ceiling constraint $\cH_q+N^{\delta-1}$ on $\mathfrak{Y}$, 
        we have w.o.p. 
        \begin{eqnarray}
          \label{eq:bulkstat}
H_t (z) \le \mathcal{H}_{\bar q} (z),  
        \end{eqnarray}
        for all
        \begin{equation}
          \label{eq:therange}
          z\in \mathfrak{Y}^-,\qquad t \in[\mathtt t,N^{5}
            ],  \quad \mathtt t:= q^{-1} N^{1+6 A\delta}(\md_{z_0}\me_{z_0;q})^{-1/2}.
        \end{equation}
      \end{prop}

      \begin{rem} 

\label{tq} 

Let us briefly, and heuristically, explain the choices of $\bar{q}$ and $\mathsf{t}$ in \Cref{prop:hdecreasebulk} (omitting factors of $N^{\delta}$ in this informal discussion). The choice \eqref{qdecreaseq} of $\bar{q}$ arises since it, with \Cref{qlineq}, implies that $0 \le \mathcal{H}_q (z) - \mathcal{H}_{\bar{q}} (z) \le N^{o(1)-1}$ for $z$ near $z_0$. Thus, by imposing a floor constraint of $\mathcal{H}_{\bar{q}}$ and a ceiling constraint of $\mathcal{H}_q + N^{o(1)-1}$ for the Glauber dynamics, the factor of $H_{\max}$ appearing in \eqref{eq:tmixconstrained} will be $N^{o(1)-1}$, which is essentially minimal. Then the mixing time in \eqref{eq:tmixconstrained} will be $N^{2+o(1)} \diam (\mathfrak{Y})^2 \le q^{-1} N^{1+o(1)} (\mathfrak{d}_{z_0} \mathfrak{e}_{z_0;q})^{-1/2}$, by \eqref{gothiccells}; this is consistent with the choice of $\mathsf{t}$ in \eqref{eq:therange}.

\end{rem} 

Lemma \ref{aj0} is a quick consequence of \Cref{prop:hdecreasebulk}.

\begin{proof}[Proof of \Cref{aj0}]
	To show \eqref{eq:follo}, fix a point $z_0 \in \mathcal{A}_k^{q'}$, and let $\mathfrak{Y}$ denote the cell associated to $z_0$ as in
	\eqref{gothiccells}, with the constant $A$ in the definition of the cell chosen to be $D+2/3$ here. Observe that $\mathfrak{Y} \subseteq \mathcal{C}_k^q$, as for any $z \in \mathfrak{Y}$ we have  
	\begin{flalign*} 
		\mathtt{d}_z \mathtt{e}_{z;q} \ge (1-N^{-\delta/4}) \mathtt{d}_{z_0} \mathtt{e}_{z_0;q} \ge (1-2N^{-\delta/4})^2 \ge 4^{-k-1/2},
	\end{flalign*} 
	
	\noindent where the first bound holds by \eqref{ed0}, the second holds by \eqref{dze2} (in \Cref{aqaq} at $u = N^{(D+5)\delta}$, which applies by \Cref{hqa}), and the third holds since $z_0 \in \mathcal{A}_k^{q'}$. 
	
	Since $\mu\in\mathcal P^+$, 
	it suffices by \eqref{h0zb} and \Cref{prop:PW} to show \eqref{eq:follo} when the Glauber dynamics $(H_t)_{t \ge 0}$ are replaced by the Glauber dynamics restricted  (in the sense of \Cref{def:restricted}) to $\mathfrak Y$, and with ceiling constraint given by $\cH_{q}+N^{\delta-1}$ (as in Definition \ref{def:constrained}). Hence, \eqref{h0zb} verifies w.o.p. the assumptions in \Cref{prop:hdecreasebulk}. The lemma then follows from that proposition with \eqref{q2q}, using the fact that the $\mathfrak{d}_{z_0} \mathfrak{e}_{z_0;q}$ in \Cref{prop:hdecreasebulk} is equal (by \Cref{hqa}) to $\mathtt{d}_{z_0} \mathtt{e}_{z_0;q} \in [4^{-k}, 4^{1-k}]$.
\end{proof}

\subsection{Proof of \Cref{prop:hdecreasebulk}} 
\label{sec61}

In this section we show \Cref{prop:hdecreasebulk}, whose proof will be similar to that of \cite[Claim 6.3]{LTDMLS}, upon replacing the use of \cite[Theorem 4.4]{LTDMLS} in the proof of \cite[Proposition 6.7]{LTDMLS} with \Cref{prop:last} here. Specifically (omitting factors of $N^{\delta}$ in this informal discussion), we will first use \Cref{prop:last} to match at $z_0$ the values of $\mathcal{H}_q$ and its gradient with those of some $q=0$ limit shape $\mathcal{H}_{a,b,c}$, so that the latter is ``less concave'' than $\mathcal{H}_q$ at $z_0$. This difference in concavity will ensure that $\mathcal{H}_{a,b,c} - \mathcal{H}_q \gg N^{-1}$ on the boundary $\partial \mathfrak{Y}$ of the cell containing $z_0$. Hence, monotonicity will imply that the stationary measure for $(H_t)$ on $\mathfrak{Y}$, with boundary data $\mathcal{H}_q$ on $\partial \mathfrak{Y}$, is lower than $\mathcal{H}_{a,b,c}$ by at least $N^{-1}$. Since $\mathcal{H}_{a,b,c}$ and $\mathcal{H}_q$ are shifted to match at $z_0$, this will imply that $H_t$ decreases by at least $N^{-1}$ at $z_0$. By \Cref{qlineq}, this corresponds (as in \eqref{qdecreaseq}) to decreasing $q$ by $N^{-1} (\mathfrak{d}_{z_0} \mathfrak{e}_{z_0;q})^{-1/2}$, in agreement with \eqref{eq:bulkstat}. We next implement this in detail.

\begin{proof}[Proof of \Cref{prop:hdecreasebulk}]

	Let $z_1,$ $a$, $b$, and $c$ be as in \Cref{prop:last} and, for lightness of notation, let $\cH':=\cH_{a,b,c}$, $\bu:=\bu_{z_0;q}$, and $\bw:=\bw_{z_0;q}$.  We let $\widehat \cH_q$ and $\widecheck \cH'$ be the macroscopic shapes $\cH_q$ and $\cH'$ rescaled according to  \eqref{e:funzrescaled} and \eqref{eq:checkH}, respectively, namely, 
	\begin{flalign}
		\label{Hwh}
		\begin{aligned}
			\widehat\cH_q(\hat x,\hat y) & =\md_{z_0}^{-1/2}\me_{z_0;q}^{-3/2} \cdot \cH_q(z_0+(\md_{z_0}\me_{z_0;q})^{1/2}\bu\hat x+ \md_{z_0}\me_{z_0;q} \bw\hat y)\\
			\widecheck\cH'(\hat x,\hat y) & =\md_{z_0}^{-1/2}\me_{z_0;q}^{-3/2} \cdot \cH_{a,b,c}(z_1+(\md_{z_0}\me_{z_0;q})^{1/2}\bu\hat x+ \md_{z_0}\me_{z_0;q} \bw\hat y).
		\end{aligned}
	\end{flalign}
	If the argument of $\cH_q$ in \eqref{Hwh} belongs to the 
	cell $\mathfrak{Y}$ from \eqref{gothiccells}, then the rescaled coordinates $\hat x,\hat y$ run
	in the range
	\begin{equation}
		\label{eq:range}
		|\hat x|,|\hat y|\le q^{-1/2}N^{(\frac34 A-\frac12)\delta-1/2}(\md_{z_0}\me_{z_0;q})^{-3/4} =:R\le N^{-\delta/2}\ll1, 
	\end{equation}
	(where we used assumptions \eqref{eq:distanzabis}). Thus, in these coordinates,  $\mathfrak{Y}^-$ and $\mathfrak Y$  are square $\mathfrak S^-$ and $\mathfrak S$ of side-lengths $2L^{-1} R$ and $2R$, respectively.
	
	By Proposition \ref{prop:last}, we have that $\nabla \widehat{\mathcal{H}}_q (0,0) = \nabla \widecheck{\mathcal{H}}' (0,0)$ and that $|\partial_\gamma \widehat{\mathcal{H}}_q - \partial_\gamma \widecheck{\mathcal{H}}'| \lesssim \mathfrak{d}_{z_0} q$ on $\mathfrak{S}$, whenever $|\gamma| \ge 2$. Therefore, for any $\rho \in [0, 1]$, we have 
	\begin{equation}
		\label{eq:diffHhat}
		\displaystyle\sup_{|w| \le \rho R} |\widehat\cH_q (w) -\widecheck\cH' (w) -\Delta \md_{z_0}^{-1/2}\me_{z_0;q}^{-3/2}|\lesssim \rho^2 R^2\md_{z_0} q\le \rho^2 N^{(\frac32A-1)\delta-1}\me_{z_0;q}^{-3/2}\md_{z_0}^{-1/2},
	\end{equation}
	where $\Delta=\cH_q(z_0)-\cH'(z_1)$. The value of $\Delta $ will play no role in the proof, so without loss of generality (and to lighten notations) we assume that it is zero. As such, $\widehat{\mathcal{H}}_q (0,0) = \widecheck{\mathcal{H}}' (0,0)$.
	The inequality \eqref{eq:diffHhat}, with \eqref{Hwh}, implies that the supremum norm difference between the two original functions $\cH_\mq$ and $\cH'$ in the cell $\mathfrak Y$
	is bounded as 
	\begin{eqnarray}
		\label{eq:HH'}
		\displaystyle\sup_{z \in \mathfrak{Y}} |\cH_q (z) -\cH' (z+z_1-z_0)| = \mathfrak d_{z_0}^{1/2}\mathfrak e_{z_0;q}^{3/2} \displaystyle\sup_{w \in \mathfrak{S}} |\widehat\cH (w) -\widecheck\cH' (w) | \lesssim N^{(\frac32A-1)\delta-1}.
	\end{eqnarray}
	
	We next estimate $\widecheck{\mathcal{H}}' (w) - \widehat{\mathcal{H}} (w)$ for $w\in\partial\mathfrak S$. To that end, Taylor expanding $ \widecheck{\cH}'$ yields for any integer $K \ge 2\delta^{-1}$ and point $w \in \partial \mathfrak{S}$ that
	\begin{equation*}
		\widecheck{\cH}'(w) - \widecheck{\cH}' (0,0)  =  w \cdot \nabla \widecheck{\cH}' (0,0) + \displaystyle\frac{w}{2} \cdot D^2          \widecheck{\cH}' (0,0) \cdot w + \sum_{3 \le |\gamma| \le K} w^{\gamma} \partial_{\gamma} \widecheck{\cH}' (0,0) + O ( |w|^{K+1}),
	\end{equation*}
	where $\gamma \in \mathbb{Z}_{\ge 0}^2$ is a multi-index and $w^{\gamma} \partial_{\gamma} \widecheck{\cH}' (0,0) $ symbolically denotes the corresponding term in the Taylor expansion. By  
	\Cref{prop:formerassumption} and  \Cref{prop:last} (and recalling that $q'=0$) we deduce for $w \in \mathfrak{S}$ and a sufficiently large constant $C>1$ that
	\begin{flalign}
		\label{e:ASb}
		\begin{aligned}
			\widecheck{\cH}'(w) - \widecheck{\cH}' (0,0) & \ge w \cdot \nabla \widehat{\cH} (0,0) + \displaystyle\frac{w}{2} \cdot D^2\widehat{\cH} (0,0) \cdot w + \frac 1 {2C} \cdot \mathfrak d_{z_0} q R^2\\  
			& \qquad + \displaystyle\sum_{3 \le |\gamma| \le K} w^{\gamma} \partial_{\gamma} \widehat{\cH} (0,0) + O \big( \mathfrak d_{z_0}q R^3 \big) + O \big( R^{K+1} \big) \\
			& \ge \widehat{\cH}(w) - \widehat{\cH}(0,0) + \frac1 {4C} \cdot \mathfrak d_{z_0}^{-1/2}\mathfrak e_{z_0;q}^{-3/2} N^{(\frac32A - 1)\delta-1} + O( R^{K+1})\\ 
			& \ge \widehat{\cH} (w) - \widehat{\cH}(0,0) + \frac 1{8C} \cdot \mathfrak d_{z_0}^{-1/2}\mathfrak e_{z_0;q}^{-3/2} N^{(\frac32A-1)\delta-1},
		\end{aligned}
	\end{flalign}
	where in the second and third inequalities we absorbed
	$O \big( \mathfrak d_{z_0} q R^3 \big)$
	into $(2C)^{-1} \mathfrak d_{z_0}q R^2$ by
	changing $(2C)^{-1}$ to $(4C)^{-1}$, and then we absorbed $O( R^{K+1}) = O(N^{-\delta K/2})$ by
	changing $(4C)^{-1}$ to $(8C)^{-1}$ (using the facts that $R \le N^{-\delta/2}$ and that $K \ge 2\delta^{-1}$). 
	Going back to the unrescaled limit shapes using \eqref{Hwh}, \eqref{e:ASb} implies for $z\in \partial\mathfrak Y$ that
	\begin{eqnarray}
		\label{eq:going}
		\cH'(z+z_1-z_0)\ge   \cH_{\mq}(z)+\frac1{8C} \cdot N^{(\frac32A-1)\delta-1}.
	\end{eqnarray}
	Therefore, by the assumptions of the proposition, the initial condition of the dynamic satisfies 
	\begin{flalign} 
		\label{h0h2} 
		H_0 (z) \le \mathcal{H}_{q}(z)+N^{\delta-1} \le \mathcal{H}'(z+z_1-z_0) - \frac 1{10C} \cdot N^{(\frac32A-1)\delta-1}, \quad \text{for $z\in{\partial \mathfrak{Y}}$},
	\end{flalign}
	where the error term $N^{\delta-1}$ has been absorbed by changing the constant  $8$ to $10$.
	Let us now upper bound $H_t(z)$ for $z\in \mathfrak Y^-$. Recall by assumption that the dynamics $(H_t)_{t \ge 0}$ are restricted to $\mathfrak Y$ (in the sense of \Cref{def:restricted}; in particular, the height function $H_t$ coincides with $H_0$ at all times, outside of $\mathfrak Y$) and  have the ceiling constraint $H_+:=\mathcal{H}_{q}+N^{\delta-1}$ in $\mathfrak Y$. To upper bound $H_t$, we may by monotonicity (\Cref{prop:PW}) further impose on $(H_t)_{t \ge 0}$ the floor constraint $H_-:=\mathcal{H}' - N^{\frac32A\delta-1}$  in $\mathfrak Y$. The distance
	between this floor and ceiling is by \eqref{eq:HH'}  at most 
	$CN^{\frac{3}{2} A \delta-1}$, for some constant $C>1$. So, the mixing time for these constrained dynamics
	is by \Cref{prop:tmixconstrained} at most
	\begin{eqnarray}
		\label{eq:isatmost}
		N^{2+\delta/2} \diam (\mathfrak{Y})^2 \cdot N^{3A\delta} \le q^{-1}  N^{5 A \delta+1}(\md_{z_0}\me_{z_0;q} )^{-1/2}\le N^{-A\delta} \mathtt t.          
	\end{eqnarray}
	
	To conclude the proof, let us introduce some notation. Let $\pi$ denote the uniform measure on tilings $\eta$ of the hexagon such that $H_{\eta} (z)=H_0 (z)$ for $z \notin \mathfrak{Y}$ and $H_- (z) \le H_{\eta} (z) \le H_+ (z)$ for $z \in \mathfrak Y$. Also let $(H^{\eq}_t)_{t\ge \mathtt t}$ denote the Glauber dynamics restricted to $\mathfrak Y$, with ceiling constraint $H_+$ and floor constraint $H_-$, with initial data $H_{\mathtt{t}}^{\eq}$ at $t=\mathtt t$ sampled under $\pi$; observe that this process is stationary. Further let  $\bar \eta$ be the tiling of $\mX$ whose height function is maximal, subject to the constraints that  $H_{\bar \eta}(z)\le \cH'(z+z_1-z_0)- (10C)^{-1} \cdot N^{(\frac32 A-1)\delta-1}$ for all $z\in\partial \mathfrak Y$. Additionally let $\mu$  denote the uniform measure on tilings of  $\mX$ that coincide with $\bar\eta$ outside of $\mathfrak Y$,  and let $\mu^+$ denote the measure obtained by conditioning $\mu$ on the event that $H(z)\ge H_-(z)$ for all $z \in \mathfrak Y$.
	
	Thanks to \eqref{eq:isatmost} and \eqref{e:submult}, at time $t=\mathtt t$ the law of $H_t$  is of total variation distance smaller than any inverse power of $N$ from  $\pi$. Hence, we can couple the processes $(H_t)_{t\ge \mathtt t}$ and $(H^{\eq}_t)_{t\ge \mathtt t}$ so that, w.o.p., they coincide for all $t\ge \mathtt t$. Proposition \ref{rem:forall} implies that the statement \eqref{eq:bulkstat} holds  w.o.p. simultaneously for all $t\in[\mathtt t,N^{5}]$, provided that
	\begin{eqnarray}
		\label{eq:dispe}
		H(z)\le \cH_{\bar q}(z)\quad \forall z\in\mathfrak Y^-,
	\end{eqnarray}
	holds w.o.p., where $H$ is sampled under $\pi$.

	To prove \eqref{eq:dispe} under $\pi$, note that the distribution of the height function $H$ in $\mathfrak Y$ sampled from $\pi$  is stochastically dominated (\Cref{def:stochdom}), because of the ceiling constraint, by that of the height in $\mathfrak Y$, sampled from $\mu^+$. Next, the height fluctuation bound\footnote{Recall that $\cH'=\cH_{a,b,c}$ and therefore the right-hand side of \eqref{h0h2}, as a function of $z$, is the limit shape in $\mX_{a,b,c}$, up to a translation $z_1-z_0$ of the hexagon, and to a global additive shift of the height. Therefore, we apply \Cref{volumeheight} to the hexagon $\mX_{a,b,c}$.} \Cref{volumeheight} implies on one hand that $\|\mu-\mu^+\|$ is smaller than any inverse power of $N$ and on the other that, w.o.p. under  $\mu$, 
	\[\big| H (z) - \mathcal{H}' (z+z_1-z_0) +
	\frac1{10C} \cdot N^{(\frac32A - 1)\delta-1} \big| \le N^{\delta/4-1}\; \text{ for every }
	z \in \mathfrak{Y}.\]
	
	Putting everything together we deduce that, w.o.p. under $\pi$ and
	for every 
	$z \in \mathfrak{Y}$,
	\begin{flalign*} 
		{H} (z) & \le \mathcal{H}' (z+z_1-z_0) + N^{\delta/4-1} - \displaystyle\frac{1}{10C} \cdot N^{(\frac32 A - 1)\delta-1} \\
		& \le \mathcal{H}' (z+z_1-z_0) - \frac{1}{16C} \cdot N^{(\frac32A-1) \delta-1}.
	\end{flalign*}
	
	Now, if $z\in\mathfrak Y^-$ then by \eqref{eq:diffHhat} (with the $\rho$ there to be $L^{-1}$ defining $\mathfrak{Y}^-$ here) and \eqref{Hwh}, we have $|\cH'(z+z_1-z_0)-\cH_q(z)|\le N^{(\frac32A-1)\delta-1}/(32C)$, if $L$ is sufficiently large. This implies  $ H (z)\le \cH_q(z)-N^{(\frac32A-1)\delta-1}/(32C)$. Finally, the fact that  $3A/2-1 \ge 5A/4$ (as $A \ge 30$), together with \Cref{qlineq}, \eqref{qdecreaseq}, and \eqref{ed0} (with \eqref{dze0} and \eqref{eq:fnt} to verify that $\mathfrak{e}_{z;q} \gg q-\bar{q}$, thereby enabling us to neglect the $\mathfrak{d}_z^{1/2} (q-\bar{q})^{3/2}$ error terms in \eqref{eq:cambioq}) imply that $\cH_q(z)-N^{(\frac32A-1)\delta-1}/(32C)\le \cH_{\bar q}(z)$ whenever \eqref{qdecreaseq} holds; 
	this establishes \eqref{eq:dispe} and hence the proposition. 
\end{proof} 

\begin{rem}
	\label{rem:whythecell} The dimensions of the cell $\mathfrak Y$ have been chosen in \eqref{gothiccells} exactly in such a way that the supremum norm difference between $\cH_q$ and $\cH'$ in $\mathfrak Y$ is of order $N^{o(1)-1}$ (\eqref{eq:HH'} provides the upper bound, and \eqref{eq:going} provides the lower bound).
\end{rem}

\section{Height reduction in the edge phase}

\subsection{Proof of Lemma \ref{akl0}} 

\label{ProofAnnuli20} 

	In this section we show Lemma \ref{akl0} using the following lemma, to be shown in \Cref{sec:across}. It essentially states the following. Assume that a height function is initially bounded above by $\mathcal{H}_q$ on all of $\mathfrak{X}$, and by $\mathcal{H}_{q'}$ up to an annulus $\mathcal{A}_k^{q'} (N^{A\delta})$ in the edge phase (see the condition $k \ge \ell_0$ below), for some $q' \le q$ satisfying \eqref{q2qnear} (which is similar to the condition \eqref{qdq} from Lemma \ref{akl0}). Then this height function will reduce under the Glauber dynamics to $\mathcal{H}_{q'} + N^{\delta-1}$, on all of $\mathfrak{X}$. In this way, the stronger height upper bound on inner annuli ``extends'' {(up to a small additive term $N^{\delta-1}$)} to outer ones, unlike in the bulk mixing phase result \Cref{aj0}, where the weaker height upper bound ``improves'' on inner annuli (indeed, to show \Cref{akl0}, we will first use \Cref{ak1akak1} to obtain a stronger height bound on the bulk liquid phase annuli, and then ``extend them out'' using \Cref{ajnear}). See \Cref{fig:61}.

\begin{lem} 
	
	\label{ajnear}
	
	Let $K,  D $ satisfy \eqref{eq:costanti}. 
	Let $0<q\le \varepsilon_0$, 
	denote $\ell_0 = \ell_0 (q)$ and let $k \ge \ell_0$ be an integer such that $q \ge 2^k N^{K \delta-1}$. Set 
	\begin{flalign}
		\label{q2qnear}
		\begin{aligned}
			& q' \in [q - 2^k N^{D \delta-1}, q]; \qquad \mathsf{t} = 2^k q^{-1} N^{1+11K\delta}
		\end{aligned}
	\end{flalign}
	and denote $\ell_0' = \ell_0 (q')$. Consider the Glauber dynamics $(H_t)_{t \ge 0}$ starting from an initial condition $\mu\in\mathcal P^+$ (recall \Cref{def:P+}) and assume that, w.o.p. under $\mu$,
	\begin{flalign}
		\label{h0z2new} 
		\begin{aligned} 
			& H_0 (w) \le \mathcal{H}_{q'} (w), \qquad \text{for each $w \in \bigcup_{j=1}^{k} \mathcal{A}_j^{q'} (N^{A\delta})$}; \\
			& H_0 (w) \le \mathcal{H}_q (w), \qquad \text{for each $w \in \mathfrak{X}$}.
		\end{aligned} 
	\end{flalign}
	Then we have w.o.p. that, for any $t \in [\mathsf{t}, N^{5}]$, 
	\begin{flalign}
		\label{hqtz}
		H_t (z) \le \mathcal{H}_{q'} (z) + N^{\delta-1}, \qquad \text{for each $z \in \mathfrak{X}$}.
	\end{flalign} 
	
\end{lem}

\begin{proof}[Proof of \Cref{akl0}]
	We induct on $k$, under the assumption  (less restrictive than \eqref{qdq}) that 
	\begin{flalign} \label{dqdbis}
		q' \in I_k:= \big[ q - (N^{\delta} - k) 2^k N^{(D+2)\delta-1}, q - (k+1) 2^k N^{(D-2)\delta-1} \big].
	\end{flalign}  
	Note that we can assume that $k<N^\delta$, as otherwise $\mathsf{t} \ge N^5$ and the claim of the lemma is empty.
	\noindent Setting $\ell_0 = \ell_0 (q)$, first observe by \Cref{ak1akak1} that the lemma holds for $k \le \ell_0 + 4$. So, fix $k \ge \ell_0 + 5$, and assume that \eqref{htz2} (with the $k$ there replaced by $k-1$) holds for all $q'\in I_{k-1}$ and $t\in [\mathsf{t}_{k-1},N^5]$, where $\mathsf{t}_{k-1} = q^{-1} (k-1) 2^{k+2} N^{1+11K\delta}$. 
	
	Now assume that $q'\in I_k$. Set $\mathsf{q}_0 = q' - 2^{k+1} N^{D\delta/2-1}$ and $q'' = (q+\mathsf{q}_0) / 2$, observing that $\mathsf{q}_0 \le q' \le q'' \le q \le 2\mathsf{q}_0$; also let $\mathsf{r} = 2^{k+1} q^{-1} N^{1+11K\delta}$.  Observe that $q-q'' = q'' - \mathsf{q}_0 \in \big[ k 2^{k-1} N^{(D-2)\delta-1}, (N^{\delta}-k+1) 2^{k-1} N^{(D+2)\delta-1} \big]$. Thus, the inductive hypothesis (with the $(q, q')$ there equal to the $(q, q'')$ here) yields w.o.p. that $H_t (z) \le \mathcal{H}_{q''} (z)$ for each $z \in \bigcup_{j=1}^k \mathcal{A}_j^{q''} (N^{A\delta})$ and $t \in [\mathsf{t}_{k-1}, N^5]$. Therefore, \Cref{ajnear} (with the $(q, q')$ there equal to $(q, q'')$ here) yields w.o.p. that $H_t (z) \le \mathcal{H}_{q''} (z) + N^{\delta-1}$ for all $z \in \mathfrak{X}$ and $t \in [\mathsf{t}_{k-1} + \mathsf{r}, N^5]$. Applying the inductive hypothesis again (now with the $(q, q')$ there equal to the $(q'', \mathsf{q}_0)$ here) yields w.o.p. that $H_t (z) \le \mathcal{H}_{\mathsf{q}_0} (z)$ for all $z \in \bigcup_{j=1}^{k} \mathcal{A}_j^{\mathsf{q}_0} (N^{A\delta})$  and $t \in [2\mathsf{t}_{k-1} + \mathsf{r}, N^5]$. 
	Again applying \Cref{ajnear} (with the $(q,q')$ there equal to $(q'',\mathsf{q}_0)$ here, and replacing $H_t$ by $H_t - N^{\delta-1}$) then gives $H_t (z) - N^{\delta-1} \le \mathcal{H}_{\mathsf{q}_0} (z) + N^{\delta-1}$ for all $z \in \mathfrak{X}$ and $t \in [2\mathsf{t}_{k-1} + 2 \mathsf{r}, N^5]$. 
	
	Since $2\mathsf{t}_{k-1} + 2 \mathsf{r} 
	\le \mathsf{t}$, it suffices to show that $\mathcal{H}_{\mathsf{q}_0} (z) + 2N^{\delta-1} \le \mathcal{H}_{q'} (z)$ for each $z \in \bigcup_{j=1}^{k+1} \mathcal{A}_j^{q'} (N^{A\delta}) \subseteq \mathfrak{L}_{q'}$. To do so, we consider two cases. First suppose that $z \in \mathfrak{L}_{\mathsf{q}_0}$. Then,   
	\begin{flalign*}
		\mathcal{H}_{q'} (z) - \mathcal{H}_{\mathsf{q}_0} (z) & \gtrsim  (q'-\mathsf{q}_0) (\mathtt{d}_z \mathtt{e}_{z;q'})^{1/2} \gtrsim N^{D\delta/2-1} \gg N^{\delta-1},
	\end{flalign*}
	
	\noindent where in the first bound we used \eqref{qlineq2}, in the second we used the facts that $q'-\mathsf{q}_0 = 2^{k+1} N^{D\delta/2-1}$ and $\mathtt{d}_z \mathtt{e}_{z;q'} \ge 4^{-k-1}$ (as $z \in \bigcup_{j=1}^{k+1} \mathcal{A}_j^{q'}$), and in the third we used that $D>5$. If $z \notin \mathfrak{L}_{\mathsf{q}_0}$, then since $z \in \bigcup_{j=1}^{k+1} \mathcal{A}_j^{q'} (N^{A\delta})$, we have $\mathcal{H}_{q'} (z) - \mathcal{H}_{\mathsf{q}_0} (z) \gtrsim \mathtt{d}_z^{1/2} \mathtt{e}_{z;q'}^{3/2} - N^{\delta -1} \gtrsim N^{A\delta-1}$, where for the second inequality we used \Cref{hde0} (and the fact that, if $(\mathfrak{d}_z, \mathfrak{e}_{z;q'}) \ne (\mathtt{d}_z, \mathtt{e}_{z;q'})$, then $\mathcal{H}_{q'} (z) \lesssim N^{-1}$, which is quickly verified also by \Cref{hde0}). Hence, again $\mathcal{H}_{q'} (z) \ge \mathcal{H}_{\mathsf{q}_0} (z) + 2N^{\delta-1}$, which confirms the lemma.
\end{proof}

\subsection{Mixing in edge phase cells}

\label{CellBoundary}

                As in \Cref{sec61prec}, we will show \Cref{ajnear} and \Cref{boundaryaj2} using ``more local'' height reduction results on certain mesoscopic cells; in this section we state those results (whose proofs will appear in Sections \ref{ProofEstimate0} and \ref{ProofEstimate1}). We begin with some notation. Throughout, we fix  real numbers  $B>A\ge 30$, $L \ge 10$, and
                \begin{equation}
                  \label{eq:condizq}
                q \in (N^{5B\delta-1}, \varepsilon_0)
              \end{equation}
              where we recall that $\varepsilon_0$ is chosen small enough, but independent of $N$ (so that the statements of \Cref{ul} and \Cref{thm:fN} with $q_0=\varepsilon_0$ hold, using \Cref{thm:step1}). 
               We further let $\mathfrak{d}, \omega > 0$ be real numbers satisfying  (the below conditions and choices will be explained in \Cref{omegaq0})
	\begin{eqnarray}
		\label{condizd}
	&	\mathfrak{d} \in [N^{-1/2}, 1) \\\label{domega} 
             &   N^{-2A\delta} \max\{(qN)^{1/3}, \mathfrak{d}^{-1}\} \le \omega \le  \displaystyle\min \{ \mathfrak{d}^{-1/3} N^{1/3 - A\delta}, qN^{1-3B\delta} \}.
	\end{eqnarray}
        Also set 
	\begin{flalign}
		\label{e0} 
		\bd = \max \big\{ \mathfrak{d}, L^{10} N^{B\delta} \omega^{1/2} (qN)^{-1/2} \big\}, \quad \be = \min \big\{ N^{2B\delta} \omega (q \mathfrak{d} N)^{-1}, L^{-10} N^{B\delta} \omega^{1/2} (qN)^{-1/2} \big\},
	\end{flalign} 
	\noindent which satisfy
	\begin{flalign}
		\label{ede} 
		\be \ge L^{-20}\mathfrak{d}^{-1/3}  N^{-2/3}, \qquad \bd \be = N^{2B\delta} \omega (qN)^{-1}.
	\end{flalign} 
	The second statement in \eqref{ede} holds since $\bd$ takes the first alternative in \eqref{e0} if and only if $\be$ does. To see that the first statement in \eqref{ede} also holds, note that  when  $\be = N^{2B\delta} \omega (q\mathfrak{d} N)^{-1}$, we have 
$\mathfrak{d}^{1/3} \be = N^{2B\delta} \omega (qN)^{-1} \mathfrak{d}^{-2/3} \ge N^{2B\delta} \omega (qN)^{-1}$,
since $\mathfrak{d} \le 1$. This, together with the lower bound on $\omega$ in \eqref{domega} (using the fact that $B > A$ and $q \le 1$), implies that $\mathfrak{d}^{1/3} \be \ge N^{-2/3}$. 
In the alternative case, we have from the minimum defining $\be$ that 
$
\be = L^{-10} N^{B\delta} \omega^{1/2} (qN)^{-1/2} 
\le N^{2B\delta} \omega (q\mathfrak{d} N)^{-1}
= (L^{10} \be)^2 \mathfrak{d}^{-1}.
$
Thus,
$
\be \ge L^{-20} \mathfrak{d},
$
and so 
$
\be \mathfrak{d}^{1/3} \ge L^{-20} \mathfrak{d}^{4/3}\ge L^{-20} N^{-2/3},
$
where we have used the fact that $\mathfrak{d} \ge N^{-1/2}.$ 

\begin{definition}
  \label{def:w0}
  Given $(\md,q,\omega,L)$ as above, we let $w_0$ be a point in $\mL_q$ satisfying the following conditions: 
  \begin{itemize}
  \item When $\bd=\md$, let $w_0 \in \mathfrak{L}_q$ be any point (which is not necessarily unique) satisfying 
	\begin{flalign}
		\label{w00} 
		\mathfrak{d}_{w_0} =  \bd; \qquad L^{-1} \be \le \mathfrak{e}_{w_0;q} \le L \be
	\end{flalign} 
	so that, by \eqref{ede}, we have
        \begin{flalign}
          \label{e:dew0}
        L^{-1} N^{2B\delta} \omega (qN)^{-1} \le \mathfrak{d}_{w_0} \mathfrak{e}_{w_0;q} \le L N^{2B\delta} \omega (qN)^{-1}.  
        \end{flalign}
      \item  When $\bd>\md$, let $w_0 \in \mathfrak{L}_q$ denote the point with
        $\md_{w_0}=\bd$, such that the orthogonal projection of $w_0$ to $\ell_{w_0}$ is the tangency point $\mathfrak{A}_q \cap \ell_{w_0}$. 
        
  \end{itemize}
\end{definition}

\begin{figure}[h]
 \begin{center} \includegraphics[width=6cm]{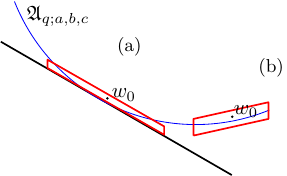}
   \caption{The cell $\mathcal Y$ of \Cref{yy}, with the point $w_0$ marked as a dot.  In the case $\md<\bd$ (a), the long sides of the cell are parallel to the boundary of the hexagon. In the case where $\md=\bd$ (b) they are parallel to $\bu_{w_0;q}$.  The short sides are in both cases vertical (that is, orthogonal  to the boundary of the hexagon, in the $(x,y)$ coordinate system). Note that the cell is not entirely contained in the liquid region; in particular, in case (a) the short sides are in the frozen region. }
\label{fig:cellY}
\end{center}
\end{figure}

	The distinction of cases in \eqref{e0}  will be relevant in defining the following cells $\mathcal{Y}$ and $\mathcal{Y}^-$ around $w_0$, on which we will analyze the Glauber dynamics. See also \Cref{fig:cellY}.

	\begin{definition} 
		
		\label{yy} 
	  Given $(\md,q,\omega,L)$ as above, and a point $w_0$ satisfying the properties in \Cref{def:w0},
          we define two subsets $\mathcal{Y} \subseteq \overline{\mathfrak{X}}$ and $\mathcal{Y}^- \subseteq \overline{\mathfrak{X}}$ as follows. Set
   	\begin{flalign}\label{otherwise}
          & \mathcal{Y} = \big\{ z \in \overline{\mathfrak{X}} : \dist_{\bm{v}_{w_0}} (z, w_0) \le L^2 (\mathfrak{d}_{w_0} \mathfrak{e}_{w_0;q})^{1/2}, \dist_{\bm{w}_{w_0}} (z, w_0) \le L^5 \mathfrak{d}_{w_0} \mathfrak{e}_{w_0;q} \big\}, \qquad \text{if $\bd = \mathfrak{d}$}; \\
          & \mathcal{Y} = \big\{ z \in \overline{\mathfrak{X}} :
          \dist_{\bm{v}_{w_0}} (z, w_0) \le L^2 N^{B\delta} \bd,\mathfrak{d}_z^2 \le L^5 N^{2B\delta} \bd^2
          \}, \qquad \qquad \qquad \qquad \qquad \qquad \quad
            \text{otherwise}.
            \label{otherwise2}
	\end{flalign}

	\noindent Then define $\mathcal{Y}^- \subseteq \mathcal{Y}$ by 
	\begin{flalign*}
		& \mathcal{Y}^- = (1/L) \cdot (\mathcal{Y} - w_0) + w_0, \qquad \qquad \qquad \qquad \qquad \quad \text{if $\bm{d} = \mathfrak{d}$}; \\
		& \mathcal{Y}^- = \big\{ z \in \overline{\mathfrak{X}} :
		\dist_{\bm{v}_{w_0}} (z, w_0) \le L \bd,\mathfrak{d}_z^2 \le L^4 \bd^2
		\}, \qquad 	\text{otherwise}.
              \end{flalign*}
	\end{definition}

       Note that, while the sides of the cell $\mathfrak Y$ are parallel to $\bw_{w_0;q},\bu_{w_0;q}$, those of $\mathcal Y$ are parallel to $\bv^\perp_{w_0},\bu_{w_0;q}$. This different choice is made for technical convenience. Note also that the angle between $\bv^\perp_{w_0}$ and $\bw_{w_0;q}$ is of order $\mathfrak{d}_{w_0} \lesssim 1$, so this modification will make little qualitative difference.

        \begin{rem}
		
		\label{omegaq0} 
		
		Let us briefly explain some of the choices and bounds above, by indicating how we will eventually proceed; in this informal discussion, we omit powers of $N^{\delta}$. We will
		view $\omega N^{-1}$ as the distance, in the $\bm{v}_{w_0}$
		direction, between the $N^{-1}$-quantiles of the height
		function $\mathcal{H}_q$ around some point
		$z \in \mathfrak{L}_q$ in the edge phase; we will also interpret
		$\mathfrak{d}$ as being of order $\mathfrak{d}_z$. Then, by the first statement of \Cref{hde0}, we have $\omega \sim (\mathfrak{d}_z \mathfrak{e}_{z;q})^{-1/2}$. 
		
		Since $z$ is in the edge phase, we have by \eqref{eq:stavolta} that $\omega \sim (\mathfrak{d}_z \mathfrak{e}_{z;q})^{-1/2} \gtrsim q^{1/3} N^{1/3}$, as in the first lower bound on $\omega$ in \eqref{domega}. Since $\mathfrak{e}_{z;q} \lesssim \mathfrak{d}_z$, we also have $\omega \sim (\mathfrak{d}_z \mathfrak{e}_{z;q})^{-1/2}  \gtrsim \mathfrak{d}^{-1}$, as in second lower bound in \eqref{domega}. Since $\mathtt{e}_{z;q} \gtrsim \mathfrak{d}^{-1/3} N^{-2/3}$, we obtain $\omega \sim (\mathfrak{d}_z \mathfrak{e}_{z;q})^{-1/2} \lesssim \mathfrak{d}^{-1/3} N^{1/3}$, as in the first upper bound  in \eqref{domega}. Since the dimensions of our cells in \Cref{yy} are of order $(\mathfrak{d}_{w_0} \mathfrak{e}_{w_0;q})^{1/2} \times (\mathfrak{d}_{w_0} \mathfrak{e}_{w_0;q})$, and the second upper bound  in \eqref{domega} implies $\bm{d} \bm{e} \sim \omega (qN)^{-1} \ll 1$, which ensures that these cells are mesoscopic. 
		
		These cell dimensions $(\mathfrak{d}_{w_0} \mathfrak{e}_{w_0;q})^{1/2} \times (\mathfrak{d}_{w_0} \mathfrak{e}_{w_0;q})$ are fixed so that, by \eqref{eq:elegeo}, the cell is both partially inside and partially outside of the liquid region $\mathfrak{L}_q$. To explain the choice \eqref{e:dew0} of $\mathfrak{d}_{w_0} \mathfrak{e}_{w_0;q}$, observe that the ``long dimension'' (in the $\bm{u}_{w_0}$ direction) of these cells is $(\mathfrak{d}_{w_0} \mathfrak{e}_{w_0;q})^{1/2}$. Upon identifying $\omega \sim (\mathfrak{d}_z \mathfrak{e}_{z;q})^{-1/2}$, \eqref{e:dew0} ensures that this is of order $(qN)^{-1/2} (\mathfrak{d}_z \mathfrak{e}_{z;q})^{-1/4}$, consistent with the long dimension of the cells $\mathfrak{Y}$ from \eqref{gothiccells}. 
			
			The parameter $\bd$ distinguishes between whether $w_0$ is reasonably distant ($\bd = \mathfrak{d}$) or too close ($\bd > \mathfrak{d}$) to a tangency location. In particular, the specific choice \eqref{e0} of $\bm{d}$ is fixed to ensure that, in the first case ($\bm{d} = \mathfrak{d}$) of \Cref{yy}, we have $\mathfrak{d}_z \asymp \mathfrak{d}$ for all $z \in \mathcal{Y}$; see \Cref{l:inthecell} below.
		
		\end{rem}
		
	The following lemma verifies that certain points in the liquid region, whose values of $\mathfrak{d}$ coincide with those of $w_0$, are contained in the above cells (observe that the condition \eqref{eq:no10pm} is in agreement with the beginning of \Cref{omegaq0}).
	
        \begin{lem}
          \label{lem:no10pm} Let $z \in \mathfrak{L}_q$ be any point in the liquid region with
$\mathfrak{d}_{z} = \mathfrak{d}_{w_0}=\md$ and $\ell_z = \ell_{w_0}$. In the case $\md=\bd$, assume furthermore that
\begin{equation}
  \label{eq:no10pm}
\omega \le  (\mathfrak{d}_{w_0} \mathfrak{e}_{w_0; q})^{-1/2}.
          \end{equation}
Then, if $L$ in \Cref{yy} is large enough, $z \in \mathcal{Y}^-$.
        \end{lem}
        \begin{proof}
        	
        	First suppose that $\bm{d} < \mathfrak{d}$. Since $\mathfrak{d}_{z} = \mathfrak{d}_{w_0}$ and $\dist (z, w_0) \lesssim \mathfrak{d}$ (as $z, w_0 \in \mathfrak{L}_q)$, we have from \Cref{yy} that $z \in \mathcal{Y}^-$. So, suppose instead that $\bm{d} = \mathfrak{d}$. 
        	
          Without loss of generality, assume that the closest tangency
          point to $w_0 = (x,y) = (x,\mathfrak{d}^2)$ is the SW one $(x^{SW}_q,0)$ and that $w_0$ is to the right
          of it (namely, $x \ge x^{SW}_q$).  Note that $w_0$ and $z$ have the same $y$-coordinate $\mathfrak{d}^2$; let $\ell$ be the (horizontal) line containing these points, and let $\mathbf{l} = (l,\mathfrak{d}^2)$ and $\mathbf{r} = (r, \mathfrak{d}^2)$ denote the left and right intersections of
          $\ell$ with $\partial \mathcal Y^-$, respectively. Denoting $z = (x', y) = (x', \mathfrak{d}^2)$, we then must show that $x' \in [l, r]$. From \Cref{yy} of $\mathcal Y^-$ and \eqref{eq:elegeo}, it follows (for sufficiently large $L$)  that $\mathbf{r} \notin \mathfrak{L}_q$, so $x' < r$ since $z \in \mathfrak{L}_q$. It therefore suffices to verify that $x \ge l$, or equivalently that $\mathfrak{e}_{z;q} \le \mathfrak{e}_{w_0;q}$. Since $\mathfrak{d}_{z} = \mathfrak{d}_{w_0}$, we thus must show $\mathfrak{d}_{z} \mathfrak{e}_{z;q} \le \mathfrak{d}_{w_0} \mathfrak{e}_{w_0;q}$, which holds since $\mathfrak{d}_{z} \me_{z;q}\le \omega^{-2}$ (from \eqref{eq:no10pm}), since $L^{-1} N^{2B\delta} \omega (qN)^{-1} \le \mathfrak{d}_{w_0} \mathfrak{e}_{w_0;q}$ (from \eqref{e:dew0}), and by \eqref{domega}. 
\end{proof}

\begin{figure}[h]
 \begin{center} \includegraphics[width=6cm]{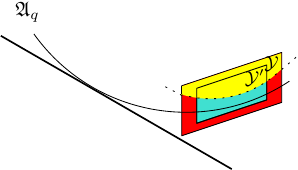}
   \caption{Proposition \ref{htqhtq02} assumes the weaker bound $H \le \cH_q$ in the (union of the blue and red) region in $\mathcal{Y}$ below the dashed line (corresponding to points $z$ where $(\mathfrak d_z \mathfrak e_{z;q'})^{1/2}=\omega^{-1}$), and it  assumes the stronger bound $H \le \mathcal{H}_{q'}$ in the (yellow) region in $\mathcal{Y}$ above the dashed line. It then deduces w.o.p. the improved bound $H_t\le \cH_{q'}+N^{\delta-1}$ in the (blue) region in $\mathcal Y^-$ below the dashed line, after time $t \ge q^{-1}\omega N^{1+10B\delta}$.}
\label{fig:Prop514}
\end{center}
\end{figure}

\begin{figure}[h]
 \begin{center} \includegraphics[width=6cm]{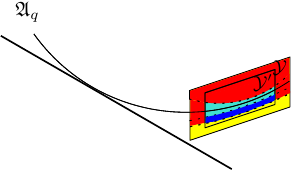}
   \caption{From top to bottom, the three dashed lines indicate the set of $z$ satisfying $(\td_z\te_{z;q})^{1/2}=N^{R\delta}\omega^{-1}$; the boundary of $\mL_q^+(\delta)$; and the boundary of $\mL_q^+(R\delta)$.	Proposition \ref{htq3} assumes the weaker bound $H \le \mathcal{H}_q  + N^{R\delta-1}$ in the (union of the red and blue) regions in $\mathcal{Y}$ above the bottommost dashed line and the stronger bound $H \le \mathcal{H}_q$ in the (yellow) region below this dashed line. It then deduces w.o.p. in $\mathcal{Y}^-$ the improved bounds $H_t\le \cH_{q}+N^{\delta-1}$ in the (light blue) region (between the top two dashed lines), and $H_t\le \cH_q$ in the (dark blue) region (between the bottom two dashed lines), after time $t \ge q^{-1}\omega N^{1+10B\delta}$.}
\label{fig:Prop515}
\end{center}
\end{figure}

Under this notation, the following results state how the tiling height function reduces in $\mathcal{Y}^-$, under the Glauber dynamics on $\mathcal{Y}$; see \Cref{fig:Prop514} and \Cref{fig:Prop515} for depictions.  Proposition \ref{htqhtq02} and Proposition \ref{htq3}  will be proven in Sections \ref{ProofEstimate0} and \ref{ProofEstimate1}, respectively. 

Proposition \ref{htqhtq02} essentially states that, assuming some weaker height upper bound \eqref{eq:cond1} in the full cell $\mathcal Y$ and stronger upper bound \eqref{eq:cond3} in a part of it, we are able to ``extend'' this stronger upper bound to a larger region \eqref{eq:statmesoe} within the smaller cell $\mathcal Y^-$ (which is analogous to \Cref{ajnear}). While  \Cref{htqhtq02} considers the edge phase in general, Proposition \ref{htq3} specifically considers the arctic boundary, through the condition \eqref{eq:q'2} (which is analogous to the second bound in \eqref{dzhq}). It in fact improves the bounds there, by effectively reducing the constant $R$ governing the upper bounds to $1$ (which is analogous to  \Cref{boundaryaj2}). The reason for the choices of $\mathsf{t}$ (and $q'$, for \Cref{htqhtq02}) in the below propositions are similar to those explained in Remark \ref{tq}. 

	The statements of \Cref{htqhtq02} and \Cref{htq3}, have been set up to be as analogous to each other as possible, so that their proofs will be quite similar. However, this will come at the expense of having to verify their hypotheses whenever we apply them in \Cref{sec:across}.

		\begin{prop}
			
                  \label{htqhtq02}
                  The following holds if the constant $L$ entering \Cref{yy} is chosen to be sufficiently large. Let the constants $A,B,R$ satisfy
			\begin{eqnarray}
				\label{eq:condizioni}
				B \ge 1000 R \ge 1000A \ge 30000. 
			\end{eqnarray}
			Let $q$ satisfy \eqref{eq:condizq}; let $\md,\omega > 0$ satisfy \eqref{condizd} and \eqref{domega}; and let the point $w_0$ around which the cell $\mathcal Y$ is centered satisfy the assumptions in \Cref{def:w0}. Also let $q'$ satisfy
			\begin{align}
				\label{eq:q'}
				\begin{aligned} 
					&  q' \in [q - N^{R\delta-1} \omega, q].
				\end{aligned}
			\end{align}
			Further assume:
			\begin{enumerate}
				\item [(i)]   For all $w \in \mathcal Y$, we have 
				\begin{eqnarray}
					\label{eq:cond1}
					H_0 (w) \le \cH_\mq (w).
				\end{eqnarray}
				
				\item [(ii)]  For all $w \in \mathcal Y$ such that $(\mathtt{d}_w \mathtt{e}_{w;q'})^{1/2} \ge \omega^{-1}$ and $\mathtt{d}_w^{1/2} \mathtt{e}_{w;q'}^{3/2} \ge N^{A\delta-1}$, we have
				\begin{eqnarray} 
					\label{eq:cond3}
					H_0 (w) \le \cH_{\mq'} (w) .
				\end{eqnarray}
			\end{enumerate}
			
			\noindent Apply the Glauber dynamics restricted to  $\mathcal Y$ with {ceiling constraint $\cH_{q'}$ in the region  $\mathcal R:=\{w:(\mathtt{d}_w \mathtt{e}_{w;q'})^{1/2} \ge \omega^{-1}, \mathtt{d}_w^{1/2} \mathtt{e}_{w;q'}^{3/2} \ge N^{A\delta-1}\}$ and ceiling constraint  $\cH_q$ in $\mathcal{Y} \setminus \mathcal R$.} Then, w.o.p., it holds
			\begin{eqnarray}
				\label{eq:statmesoe}
                          H_t (z) \le \cH_{\mq'} (z) + N^{\delta-1}
			\end{eqnarray}
			 for each $t \in [\mq^{-1} N^{1+10B\delta} \omega,N^{5}]$  and $ z \in \mathcal Y^-$ such that $ \mathtt{d}_z \in (L^{-2} \mathfrak{d}, L^2 \mathfrak{d})$.
		\end{prop}

		\begin{prop} 
			
			\label{htq3}

		The following holds if the constant $L$ entering \Cref{yy} is chosen to be sufficiently large. Let the constants $A,B,R$ satisfy
			\begin{eqnarray}
				\label{eq:condizioni2}
				B \ge 1000R \ge 1000A  \ge 30000. 
			\end{eqnarray}
			Let $q$ satisfy \eqref{eq:condizq}; let $\md,\omega > 0$ satisfy \eqref{condizd} and \eqref{domega}; and let the point $w_0$ around which the cell $\mathcal Y$ is centered satisfy the assumptions in \Cref{def:w0}. Also suppose that
			\begin{align}
				\label{eq:q'2}
				\begin{aligned} 
                                        \omega \ge \mathfrak{d}^{-1/3} N^{1/3-R\delta}.
				\end{aligned}
			\end{align}
			Further assume:
			\begin{enumerate}
				\item [(i)]   For all $w \in \mathcal{Y}$, we have 
				\begin{eqnarray}
					\label{eq:cond12}
					H_0 (w) \le \cH_\mq (w) + N^{R\delta-1}.
				\end{eqnarray}
				
				\item [(ii)] For all $w \in \mathcal{Y} \setminus  \mathfrak{L}_q^+ (R\delta)$, we have 
				\begin{eqnarray}
					\label{eq:cond32}
					H_0 (w) \le \mathcal{H}_q (w).
				\end{eqnarray}
			\end{enumerate}
			
			\noindent Apply the Glauber dynamics restricted to $\mathcal Y$, with ceiling constraint $\cH_q$ on $\mathcal{Y} \setminus  \mathfrak{L}_q^+ (R\delta)$ and $\cH_q+N^{R\delta-1}$ on $\mathcal{Y}\cap \mathfrak{L}_q^+ (R\delta)$. Then, w.o.p., the following two statements hold for each $t \in [\mq^{-1} N^{1+10 B\delta} \omega,N^{5}]$.
			\begin{enumerate} 
				\item  \label{66.1} For each $z \in \mathcal{Y}^-$ with $(\mathtt{d}_z \mathtt{e}_{z;q})^{1/2} \le N^{R\delta} \omega^{-1}$ and $\mathtt{d}_z \in (L^{-2} \mathfrak{d}, L^2 \mathfrak{d})$, we have 
				\begin{eqnarray}
					\label{eq:statmesoe3}
					H_t (z) \le \cH_q (z) + N^{\delta-1}.
				\end{eqnarray}
				
				\item  \label{66.2} For each $z \in \mathcal{Y}^- \setminus \mathfrak{L}_q^+ (\delta)$ with $\mathtt{d}_{z} \in (L^{-2} \mathfrak{d}, L^2 \mathfrak{d})$, we have
				\begin{flalign*}
					H_t (z) \le \mathcal{H}_{q} (z).
				\end{flalign*}
			\end{enumerate} 
		\end{prop} 
		
\subsection{Proof of \Cref{boundaryaj2}}

\label{Proofqn2} 

          As in the proof of \Cref{aj0}, it suffices by Proposition \ref{prop:PW} to show \eqref{htzq2}, with the Glauber dynamic $(H_t)_{t \ge 0}$ replaced by the Glauber dynamics with ceiling constraint $\cH_q+N^{K\delta-1}$ on the whole $\mX$ and $\cH_q$ on the complement of $\mL_q^+({K\delta})$. So, we let $(H_t)_{t \ge 0}$ denote the latter dynamic for the remainder of this section.
		
		Fix a point $z_0 \in \mathfrak{X}$ satisfying \eqref{dzhq}. We must show w.o.p. that, for all $t \in [\mathsf{t}, N^5]$, we have $H_t (z_0) \le \mathcal{H}_{q'} (z_0) + N^{\delta/4-1}$ and $H_t (z_0) \le \mathcal{H}_{q'} (z_0)$ if $z_0 \notin \mathfrak{L}_{q'}^+ (\delta)$. It suffices to show w.o.p. that, for all $t \in [\mathsf{t}, N^5]$, we have 
		\begin{flalign}
			\label{htz0q} 
			H_t (z_0) \le \mathcal{H}_q (z_0) + \frac{1}{2} \cdot N^{\delta/4-1}, \qquad \text{and} \qquad H_t (z_0) \le \mathcal{H}_q (z_0), \quad \text{if $z_0 \notin \mathfrak{L}_{q}^+ ( \delta / 2)$}.
		\end{flalign}
		
Indeed, the sufficiency of the first bound in \eqref{htz0q} follows from the bounds 
		\begin{flalign}
			\label{hqz03} 
			\mathcal{H}_q (z_0) - \mathcal{H}_{q'} (z_0) & \lesssim (q-q') (\mathtt{d}_{z_0} \mathtt{e}_{z_0;q'})^{1/2} +(q-q')^{3/2}\td_{z_0}^{1/2} \lesssim N^{K\delta/3-1/3} \mathtt{d}_{z_0}^{1/3} (q-q') \le N^{-1},
		\end{flalign}
                where the first statement follows from \Cref{qlineq}; the second and the third from \eqref{dzhq}. The sufficiency of the second bound in \eqref{htz0q} follows from the fact that the distance between the arctic boundaries $\mathfrak{A}_{q'}$ and $\mathfrak{A}_{q}$ on the level $\{ z : \mathtt{d}_z = \mathtt{d}_{z_0} \}$ is at most $N^{\delta} (q-q') \le \mathtt{d}_{z_0}^{-1/3} N^{-2/3-\delta}$, again by \Cref{boundarydistance}. Therefore, by Definition \ref{def:augmented} of augmented liquid region (together with \Cref{rem:lism}), for $z_0\not\in\mathfrak{L}_{q}^+ ( \delta / 2)$, we have $\cH_q(z_0)=\cH_{q'}(z_0)$, so the second statement of \eqref{htzq2} follows from that of \eqref{htz0q}.
                   
               To show \eqref{htz0q} (and thus \Cref{boundaryaj2}), we will use \Cref{htq3}, whose hypotheses are verified through the below lemma.    
                  
                   \begin{lem}\label{lem:checking1}
                   	The assumptions of \Cref{htq3} are satisfied with the parameters $(\mathfrak{d}, \omega; A, B, R; \delta)$ there given by $(\mathtt{d}_{z_0}, N^{-K\delta/2} (\mathtt{d}_{z_0} \mathtt{e}_{z_0;q})^{-1/2}; 2K, 4000K, 4K; \delta / 4)$ here.  
                   \end{lem}

\begin{proof}[Proof of \Cref{boundaryaj2}]
	
	By \Cref{lem:checking1}, we may apply \Cref{htq3}; taking $w_0$ as in \Cref{def:w0} and
the associated cells $\mathcal{Y}$ and $\mathcal{Y}^-$ as in \Cref{yy}, it
is quickly verified via \Cref{lem:no10pm} that
$z_0 \in \mathcal{Y}^-$. Since the initial condition $\mu$ belongs to
$\mathcal P^+$, by Proposition \ref{prop:PW} it suffices to show \eqref{htzq2}, with the Glauber dynamic $(H_t)_{t \ge 0}$ replaced by the Glauber dynamics restricted (in the sense of \Cref{def:restricted}) to $\mathcal Y$. Then, \Cref{htq3} yields
\eqref{htz0q} w.o.p. for
$q^{-1} N^{1+10B\delta/4} \omega \le t \le N^5$. Since
$B = 4000K$ and
$\omega \le (\mathtt{d}_{z_0} \mathtt{e}_{z_0; q})^{-1/2} \le
\mathtt{d}_{z_0}^{-1/3} N^{1/3}$, and we also have
$\mathtt{d}_{z_0} \ge N^{-K \delta} \mathtt{d}$, we have from the
definition of $\mathsf{t}$ that
$\mathsf{t} \ge q^{-1} N^{1+10000K\delta} \omega$, which implies the
lemma.
	\end{proof}

              \begin{proof}[Proof of \Cref{lem:checking1}]
                We need to check that
		\begin{flalign}
			\label{kndeltaomega}
			\begin{aligned}  
				\max \big\{ N^{-K\delta} (qN)^{1/3}, N^{-K\delta} \mathtt{d}_{z_0}^{-1}, \mathtt{d}_{z_0}^{-1/3} N^{1/3-K\delta} \big\} & \le N^{-K\delta/2} (\mathtt{d}_{z_0} \mathtt{e}_{z_0;q})^{-1/2} \\
				& \le \min \{ \mathtt{d}_{z_0}^{-1/3} N^{1/3-K\delta/2}, qN^{1-3000K\delta} \}.
			\end{aligned} 
		\end{flalign}
		\noindent First observe that $\mathtt{d}_{z_0}^{-1/3} N^{1/3-K\delta} \le N^{-K\delta/2} \mathtt{d}_{z_0}^{-1/2} \mathtt{e}_{z_0;q}^{-1/2} \le \mathtt{d}_{z_0}^{-1/3} N^{1/3-K\delta/2}$ follows from the facts that $\mathtt{e}_{z_0;q} \ge \mathtt{d}_{z_0}^{-1/3} N^{-2/3}$ (by its definition) and $\mathtt{d}_{z_0}^{1/2} \mathtt{e}_{z_0;q}^{3/2} \le N^{K\delta-1}$ by \eqref{dzhq}.             
 We also have $N^{-K\delta/2} (\mathtt{d}_{z_0} \mathtt{e}_{z_0;q})^{-1/2} \ge N^{-K\delta} \mathtt{d}_{z_0}^{-1}$, as $\mathtt{e}_{z_0;q} \lesssim \mathtt{d}_{z_0}$ (by \eqref{eq:elegeo0}).   We have that $N^{-K\delta} (qN)^{1/3} \le N^{-K\delta/2} (\mathtt{d}_{z_0} \mathtt{e}_{z_0;q})^{-1/2}$ since $q \le 1$ and $(\mathtt{d}_{z_0} \mathtt{e}_{z_0;q})^{-1/2} \ge (\mathtt{d}_{z_0}^{1/2} \mathtt{e}_{z_0;q}^{3/2})^{-1/3} \ge N^{1/3-K\delta/3}$ (the first since $\mathtt{d}_{z_0} \le 1$ and the second by \eqref{dzhq}). Finally, $N^{-K\delta/2} (\mathtt{d}_{z_0} \mathtt{e}_{z_0;q})^{-1/2} \le qN^{1-3000K\delta}$ since $q \ge \mathtt{d}_{z_0}^{-1/3} N^{3500K\delta-2/3}$ and $\mathtt{d}_{z_0}^{1/2} \mathtt{e}_{z_0;q}^{3/2} \ge N^{-1}$. This verifies \eqref{kndeltaomega}.
              \end{proof}

	\subsection{Proof of \Cref{ajnear}}
	\label{sec:across}

	Since $\mu\in\mathcal P^+$, to prove \Cref{ajnear} it suffices (as in the proof of \Cref{aj0}) by \Cref{prop:PW} to show \eqref{hqtz} when the Glauber dynamics $(H_t)_{t \ge 0}$ are replaced by the Glauber dynamics with ceiling constraint $\cH_q$ on all of
	$\mX$, and the stronger ceiling constraint $\cH_{q'}$ in $\bigcup_{j=1}^{k} \mathcal{A}_j^{q'}(N^{A \delta})$. Therefore, for the remainder of this section, we will let $(H_t)_{t \ge 0}$ denote the latter dynamics. In what follows, we fix a point  $z_0 \in \mathfrak{X}$. 
	
	The following lemma directly addresses several ``boundary cases'' for $z_0$.

	\begin{lem} 
		
		\label{lq00} 
		
		The bound \eqref{hqtz} holds for any $t \in [\mathsf{t}, N^5]$, w.o.p., either if $z_0 \in \overline{\mathfrak{L}}_{q'}$ and $\mathtt{d}_{z_0} \mathtt{e}_{z_0;q'} \le N^{-6D\delta} 4^{-k}$, or if $z_0 \in \bigcup_{j=1}^k \mathcal{A}_j^{q'} (N^{A\delta})$. Additionally, if w.o.p. \eqref{hqtz} holds for any $t \in [\mathsf{t}, N^5]$ and $z_0 \in \overline{\mathfrak{L}}_{q'}$, then w.o.p. \eqref{hqtz} holds for any $t \in [\mathsf{t}, N^5]$ and $z_0 \in \mathfrak{X}$.

	\end{lem} 
	
	\begin{proof} 
		
		The fact that \eqref{hqtz} holds if $z_0 \in \bigcup_{j=1}^k \mathcal{A}_j^{q'} (N^{A\delta})$ follows on the ceiling constraint for $(H_t)_{t \ge 0}$ mentioned above. Next, if $z_0 \in \overline{\mathfrak{L}}_{q'}$ and $\mathtt{d}_{z_0} \mathtt{e}_{z_0;q'} \le N^{-6D\delta} 4^{-k}$, then
		\begin{flalign*}
			H_t (z_0) \le \mathcal{H}_q (z_0) & \le \mathcal{H}_{q'} (z_0) + N^{\delta} (q-q') (\mathtt{d}_{z_0} \mathtt{e}_{z_0;q'})^{1/2} + N^{\delta} \mathtt{d}_{z_0}^{1/2} (q-q')^{3/2} \\
			& \le \mathcal{H}_{q'} (z_0) + N^{\delta - D\delta-1} + 2^{3k/2} N^{4D\delta - 3/2} \mathtt{d}_{z_0}^{1/2},
		\end{flalign*}
		
		\noindent where in the first bound we used the ceiling constraint mentioned above on $(H_t)_{t \ge 0}$; in the second we used \Cref{qlineq}; in the third we used the facts that $0 \le q-q' \le 2^k N^{2D\delta-1}$ and $\mathtt{d}_{z_0} \mathtt{e}_{z_0} \le N^{-6D\delta} 4^{-k}$. Since $\mathtt{d}_{z_0}^{1/2} \mathtt{e}_{z_0;q'}^{3/2} \ge N^{-1}$ and $\mathtt{d}_{z_0}^{1/2} \mathtt{e}_{z_0;q'}^{1/2} \le N^{-3D\delta} 2^{-k}$, we have that $\mathtt{d}_{z_0} \le 8^{-k} N^{1-9D\delta}$, and so insertion into the above inequality yields $H_t (z_0) \le \mathcal{H}_{q'} (z_0) + N^{-1-\delta} + N^{-1-D\delta/2} \le \mathcal{H}_{q'} (z_0) + N^{\delta-1}$, confirming \eqref{hqtz}.
		
		Next, suppose that $z_0 \notin \overline{\mathfrak{L}_{q'}}$. If $z_0$ is
		in a connected component of
		$\mathfrak{X} \setminus \overline{\mathfrak{L}}_{q'}$ for which
		$\mathcal{H}_{q'} (z_0)$ is maximal, then we have
		$H_t (z_0) \le \mathcal{H}_{q'} (z_0) < \mathcal{H}_{q'}
		(z_0) + N^{\delta-1}$. Otherwise, if
		$z_0 \notin \mathfrak{L}_{q'}$, then there is a point
		$z_0' \in \overline{\mathfrak{L}}_{q'}$ (obtained by intersecting $\mathfrak{A}_{q'}$ with the line through $z_0$ orthogonal to $\ell_{z_0}$) for which
		$\mathcal{H}_{q'} (z_0) = \mathcal{H}_{q'} (z_0')$ and
		$H (z_0) \le H(z_0')$ for any admissible height function
		$H : \mathfrak{X} \rightarrow \mathbb{R}$. Hence, if
		$H_t (z_0') \le \mathcal{H}_{q'} (z_0')+N^{\delta-1}$ holds w.o.p., then we would have
		$H_t (z_0) \le H_t (z_0') \le \mathcal{H}_{q'} (z_0') +
		N^{\delta-1} = \mathcal{H}_{q'} (z_0) + N^{\delta-1}$ w.o.p., confirming \eqref{hqtz}.
	\end{proof} 
	
	In view of Lemma \ref{lq00}, for the remainder of this section, we will assume that 
	\begin{flalign}
		\label{z0estimateaj} 
		z_0 \in \overline{\mathfrak{L}}_{q'}; \qquad z_0 \notin \bigcup_{j=1}^k \mathcal{A}_j^{q'} (N^{A\delta}); \qquad N^{-6D\delta} 4^{-k} \le \mathtt{d}_{z_0} \mathtt{e}_{z_0;q'}.
	\end{flalign} 
	
	\noindent Now set 
	\begin{flalign}
		\label{omega0d}
		\mathfrak{d} =  \mathtt{d}_{z_0}; \qquad \omega = \min \{ 2^k, \mathfrak{d}^{-1/3} N^{1/3-A\delta} \}.
	\end{flalign}

	\begin{lem}
		\label{lem:verification2}
		The parameters $(q,\md,\omega)$ satisfy the inequalities \eqref{domega}, with the $(A,B)$ there equal to $( A, K / 3)$ here. 
	\end{lem}
	\begin{proof}
		The bound $\omega \le 2^k \le qN^{1-K\delta}$ follows from our assumed lower bound on $q$; the bound $\omega \le \mathfrak{d}^{-1/3} N^{1/3-A\delta}$ follows from our definition of $\omega$;              
		the bound $N^{-2A\delta} \mathfrak{d}^{-1} \le \omega$ follows from the facts that $\mathfrak{d} \ge N^{-1/2}$ and that $N^{-2A\delta} \mathtt{d}_{z_0}^{-1} \lesssim N^{-2A\delta}  (\mathtt{d}_{z_0} \mathtt{e}_{z_0;q'})^{-1/2} \le N^{(3D-2A)\delta} 2^k \ll 2^k$ (the first by \eqref{eq:elegeo0}, the second by the last bound in \eqref{z0estimateaj}, and the third by \eqref{eq:costanti}); and the bound $q^{1/3} N^{1/3-2A\delta} \le \omega$ follows from the bounds $q^{1/3} N^{1/3-2A\delta} \le N^{1/3-2A\delta} \le \mathfrak{d}^{-1/3} N^{1/3-A\delta}$ and $2^k \ge 2^{\ell_0} \gtrsim q^{1/3} N^{1/3 - A\delta/2}$ by the fact that $k \ge \ell_0$ and \eqref{adelta4}, with the fact that $\omega = \min \{ 2^k, \mathfrak{d}^{-1/3} N^{1/3-A\delta} \}$. 
	\end{proof}

	Now let $L \ge 10$ be a sufficiently large real number so that Lemma \ref{lem:no10pm} holds. Further let $w_0$ be a point in  $\mathfrak{L}_{q}$ satisfying the properties in \Cref{def:w0}. Then define the rectangular cells $\mathcal{Y}$ and $\mathcal{Y}^-$ with respect to $(\mathfrak{d},q,\omega,L)$ and $w_0$, as in \Cref{yy}. 
	\begin{lem}
		\label{lem:cellinclude1}
		
		We have that  $z_0 \in \mathcal{Y}^-$.
	\end{lem}

	\begin{proof}
		 The claim holds by  \Cref{lem:no10pm} if either we are in the second case of \Cref{yy} ($\mathfrak{d}_{w_0} \ne \bm{d}$), or if
		$\mathtt{d}_{z_0} \mathtt{e}_{z_0;q} \le  4^{-k}$  (as, in the latter case, \eqref{eq:no10pm} holds by our definition \eqref{omega0d} of $\omega$). Therefore, we can assume that we are in the first case of \Cref{lem:no10pm}, i.e., $\bm{d} = \mathfrak{d}_{w_0} = \mathfrak{d}$. Observe that
		\begin{flalign}
			\label{hqz00}
			N^{A\delta-1} \ge \mathtt{d}_{z_0}^{1/2} \mathtt{e}_{z_0;q'}^{3/2} \gtrsim \mathfrak{d}^{1/2} \mathfrak{e}_{z_0;q}^{3/2},
		\end{flalign}
		where the first bound holds by the second statement in \eqref{z0estimateaj} and the second holds by \eqref{dze2} (applied with the $u$ there equal to $N^{A\delta}$ here, using \eqref{q2qnear}).
		
		As in the proof of \Cref{lem:no10pm}, it is enough to show  that $\me_{z_0;q} \le \me_{w_0;q}$. To do so, observe by \eqref{e:dew0} (at $B=K/3$) that 
		\begin{eqnarray*}
			\me_{w_0;q} \ge L^{-1} N^{2K/3\delta} \omega (q \mathfrak{d} N)^{-1}.
		\end{eqnarray*}
		\noindent Since \eqref{hqz00} implies that $\mathfrak{e}_{z_0;q} \lesssim \mathfrak{d}^{-1/3} N^{2A\delta/3 - 2/3}$, it suffices to show for a sufficiently large constant $C>1$ that
		\begin{eqnarray}
			\omega \ge CL q \mathfrak{d}^{2/3} N^{1/3 + \delta(2A/3-2K/3)}.   
		\end{eqnarray}
		Since $\omega \ge N^{-2A\delta} (qN)^{1/3}$ by \eqref{domega}, this holds by \eqref{eq:costanti} and bound $d, q \lesssim 1$. 
	\end{proof}

		By \Cref{lem:cellinclude1} (and the fact that $\mu \in \mathcal{P}^+$), to show \Cref{ajnear} it suffices by \Cref{prop:PW} to prove \eqref{hqtz} (with the $z$ there equal to $z_0$ here) when the dynamics $(H_t)_{t \ge 0}$ are further restricted to $\mathcal Y$ (as in Definition \ref{def:restricted}). Therefore, in the remainder of this section, we will let $(H_t)_{t \ge 0}$ denote the latter dynamics. Then the next lemma verifies \eqref{hqtz} for certain $z$ (satisfying \eqref{htzqkdelta} below).
		
		\begin{lem} 
			
			\label{htzt}
			 
			Setting $\mathsf{t}' = \mathsf{t}/2$, we have w.o.p. that, for any $z \in \mathfrak{X}$ and $t \in [\mathsf{t}', N^5]$,
			\begin{flalign}
				\label{htzqkdelta}
				H_t (z) \le \mathcal{H}_{q'} (z) + N^{\delta-1}, \quad \text{if either $q-q' \le \mathtt{d}_z^{-1/3} N^{A\delta-2/3}$ or $\mathtt{d}_z^{1/2} \mathtt{e}_{z;q'}^{3/2} \ge N^{4A\delta-1}$.}
			\end{flalign}
			
		\end{lem} 
		
		\begin{proof}

			It suffices to show the lemma at $z = z_0$. This would follow from \Cref{htqhtq02} (with the $R$ there equal to $2A$ here), upon verifying the assumptions there; by \eqref{eq:costanti} and \Cref{lem:verification2}, we must confirm \eqref{eq:q'}, \eqref{eq:cond1}, and \eqref{eq:cond3}. Since $\omega^{-1} \ge 2^{-k}$ by \eqref{omega0d}, the bounds \eqref{eq:cond1} and \eqref{eq:cond3} follow from \eqref{h0z2new}. Moreover, \eqref{eq:q'} holds whenever $q-q' \le N^{2A\delta-1} \omega$. If $\mathtt{d}_{z}^{-1/3} N^{1/3-A\delta} \ge 2^k$, then $\omega=2^k$, so by \eqref{q2qnear} we have $q-q' \le 2^k N^{D\delta-1} \le N^{2A\delta-1} \omega$, confirming \eqref{eq:q'} and thus the lemma. 
			
			So, assume instead that $\omega = \mathtt{d}_{z_0}^{-1/3} N^{1/3-A\delta} < 2^k$. If $q-q' \le \mathtt{d}_z^{-1/3} N^{A\delta-2/3}$, then $q-q' \le N^{2A\delta-1} \omega$, again confirming \eqref{eq:q'} and thus the lemma; this addresses the first case in \eqref{htzqkdelta}. To address the second, it suffices to show that $\mathcal{H}_{q'} (z_0) < N^{4A\delta-1}$ (when $\omega < 2^k$).  To that end, since $z_0 \notin \mathcal{A}_j^{q'} (N^{A\delta}) $ for any $j \le k$ (by \eqref{z0estimateaj}), we may assume that $z_0 \in \mathcal{A}_j^{q'}$ for some $j > k$. In that case, $\mathtt{d}_{z_0} \mathtt{e}_{z_0; q'} \le 4^{-k} \le \mathtt{d}_{z_0}^{2/3} N^{2A\delta -2/3}$ (as $2^k > \mathtt{d}_{z_0}^{-1/3} N^{1/3-A\delta}$), so $\mathtt{d}_{z_0}^{1/2} \mathtt{e}_{z_0;q'}^{3/2} \ll N^{4A\delta-1}$, confirming that the second case in \eqref{htzqkdelta} does not hold. 
		\end{proof} 
		
		If $q-q' > \mathtt{d}_z^{-1/3} N^{A\delta-2/3}$ (and $\mathtt{d}_z^{1/2} \mathtt{e}_{z;q'}^{3/2} < N^{A\delta-1}$), then $q-q'$ is too large to directly satisfy \eqref{eq:q'}, so we cannot immediately apply \Cref{htqhtq02}. In this case, we will inductively increase the ``range'' for $q-q'$, by a certain increment $u(z)$, at which we can obtain improved estimates. 
		
		More precisely, given $z \in \mathfrak{X}$, let $u(z) = N^{-25A\delta-2/3} \mathtt{d}_z^{-1/3}$ if $q-q' \le N^{A\delta-2/3}\mathtt{d}_z^{-1/3} $. If instead $q-q' > N^{A\delta-2/3} \mathtt{d}_z^{-1/3}$ then  let $u(z) \in [N^{-25A\delta-2/3} \mathtt{d}_z^{-1/3}/2, N^{-25A\delta-2/3} \mathtt{d}_z^{-1/3}]$ be any real number satisfying $(q-q')/u(z) \in \mathbb{Z}$. Under this notation, we have the following lemma that quickly implies \Cref{ajnear}; it will be shown by induction.

		\begin{lem} 
			
		\label{htzku} 
		
		For any point $z \in \mathfrak{X}$ and integer $m \in [0, (q-q') / u(z)]$, the following two statements hold w.o.p., for each $t \in [\mathsf{t}'+mN^{10K\delta+2} u(z), N^5]$. First,
		\begin{flalign}
			\label{1htzk} 
			H_{t} (z) \le \mathcal{H}_{q-m u(z)} (z) + N^{\delta-1}.
		\end{flalign}
		
		\noindent Second, we have
		\begin{flalign}
			\label{2htzk} 
			H_{t} (z) \le \mathcal{H}_{q-mu(z)} (z), \qquad \text{if $z \notin \mathfrak{L}_{q-mu(z)}^+ (\delta)$}.
		\end{flalign}
		\end{lem}

		\begin{proof}[Proof of \Cref{ajnear}]
			
			Recall we have fixed $z_0 \in \mathfrak{X}$. If $q-q' \le u(z_0)$, then \eqref{hqtz} follows from \eqref{htzqkdelta}. If instead $z_0 \in \mathfrak{X}$ satisfies $q-q' \ge u(z_0)$, then taking $m = (q-q') / u(z_0) \in \mathbb{Z}$ in \eqref{1htzk} yields \eqref{hqtz}, since $\mathsf{t}' + (q-q') N^{10K\delta+2} < \mathsf{t}$.
		\end{proof} 
		
		It remains to show \eqref{1htzk} and \eqref{2htzk}. They hold when $m=0$ because we imposed the ceiling $\cH_q$ everywhere in $\mX$. So, it suffices to show them for $m>0$; we also know them if $q-q' \le N^{A\delta-2/3} \mathsf{d}_z^{-1/3}$ (and $m>0$), by \eqref{htzqkdelta}.
		Let $\mathcal{S}$ denote the set of real numbers of the form $\mathsf{t}' + mN^{10K\delta+2} u(z)$ that are less than $N^5$, where $m>0$ ranges over all positive integers and $z \in \mathfrak{X} \cap \mathbb{T}_N$ ranges across all points such that $q-q' > N^{A\delta-2/3} \mathsf{d}_z^{-1/3}$. 
		Observe that $|\mathcal{S}| \le N^{11}$, since $(q-q') / u(z) \in \mathbb{Z}$ and $u(z) \ge N^{-25A\delta-2/3} \mathsf{d}_z^{-1/3}/2 \ge N^{-6/7}$. 
		
		It suffices to show that, for any $s \in \mathcal{S}$, \eqref{1htzk} and \eqref{2htzk} hold w.o.p. for all $t \in [s, N^5]$, where the $m$ in those equations is given by $\lfloor (s-\mathsf{t}') / (N^{10K\delta+2} u(z)) \rfloor$.
		To that end, induct on $s$. More specifically, we fix $s \in \mathcal{S}$ with $s = \mathsf{t}' + m N^{10K\delta+2} u(z)$, for some point $z \in \mathfrak{X}$ and integer $m \in [1, (q-q')/u(z)]$. We then prove \eqref{1htzk} and \eqref{2htzk} for all $t \in [s, N^5]$, assuming they hold for all $(m',z', t)$ such that $\mathsf{t}' + m' N^{10K\delta+2} u(z') < s$ and $t \in [\mathsf{t}' + m' N^{10K\delta+2} u(z'), N^5]$.
		
		Set 
		\begin{flalign*} 
			q'' = q- (m-1) u(z); \qquad s' = s - N^{10K\delta+2} u(z).
		\end{flalign*} 
		
		\noindent By symmetry, we may assume that the closest tangency location, between $\mathfrak{A}_{q''}$ and $\partial \mathfrak{X}$, to $z$ is the SW one $p_{q''}^{SW}$; we may also assume that $z \in \mathfrak{L}_{q''}^+ (2\delta)$ (as otherwise \eqref{2htzk} would follow from applying it to the point where $\mathfrak{L}_{q''}^+ (2\delta)$ meets the line through $z$ parallel to $\ell_z$). 
		
			\begin{lem}\label{lem:verify2}
			The assumptions of \Cref{htq3} are verified if the $(A,R,B,\delta)$ there equal to $(2A,10A,K/3,\delta/2)$ here; the $(q,\mathfrak{d},\omega)$ there equal to $(q'',\mathtt{d}_z, \mathtt{d}_z^{-1/3} N^{1/3-2A\delta})$ here; and the $H_0$ there equals to $H_{s'}$ here. Moreover, $z$ is in the cell $\mathcal Y^-$ of the proposition.
		\end{lem}	
		
	\begin{proof}[Proof of \Cref{htzku}]
		
		Applying \Cref{htq3} using \Cref{lem:verify2}, we obtain by its second statement that $H_t (z) \le \mathcal{H}_{q''} (z)$, for any $t \in [\mathsf{t}' + mN^{10K\delta+2} u(z), N^5]$ and $z \notin \mathfrak{L}_{q''}^+ (\delta/2)$. Also, since $\mathfrak{L}_{q''}^+ (\delta/2) \subseteq \mathfrak{L}_{q''-u(z)}^+ (\delta)$ (by \Cref{boundarydistance} and \Cref{rem:lism}, with \Cref{rem:bysymmetry}, as $u(z) \le \mathsf{d}_z^{-1/3} N^{-25A\delta-2/3}\ll \mathsf d_z^{-1/3} N^{\delta/2-2/3}$), we have for $z \notin \mathfrak{L}_{q''-u(z)}^+ (\delta)$ that $\mathcal{H}_{q''-u(z)}(z) = \mathcal{H}_{q''}(z)$. Hence, $H_t (z) \le \mathcal{H}_{q''} (z) = \mathcal{H}_{q''-u(z)} (z) = \mathcal{H}_{q-mu(z)}(z)$, so \eqref{2htzk} holds.
		
		To show \eqref{1htzk}, we may assume that $\mathtt{d}_z^{1/2} \mathtt{e}_{z;q''-u(z)}^{3/2} \le N^{4A\delta-1}$, as otherwise \eqref{2htzk} would follow from \eqref{htzqkdelta}. Then, we have $\mathcal{H}_{q''-u(z)} (z) \ge \mathcal{H}_{q''} (z) - N^{-1}$ by \Cref{qlineq}, as $u(z) \le \mathsf{d}_z^{-1/3} N^{-25A\delta-2/3}$ and $\mathtt{d}_z^{1/2} \mathtt{e}_{z;q''-u(z)}^{3/2} \le N^{4A\delta-1}$ implies that $u(z) (\mathfrak{d}_z \mathfrak{e}_{z;q''-u(z)})^{1/2} + u(z)^{3/2} \mathfrak{d}_z^{1/2} \ll N^{-1}$. Therefore, we would deduce \eqref{1htzk} from \eqref{eq:statmesoe3} if we have $(\mathtt{d}_z \mathtt{e}_{z;q''})^{1/2} \le N^{4A\delta} \omega^{-1}$. Recalling that $\omega = \mathtt{d}_z^{-1/3} N^{1/3-2A\delta}$, this is equivalent to $\mathtt{d}_z^{1/2} \mathtt{e}_{z;q''}^{3/2} \le N^{18A\delta-1}$, which holds since $\mathtt{d}_z^{1/2} \mathtt{e}_{z;q''-u(z)}^{3/2} \le N^{4A\delta-1}$. This establishes \eqref{1htzk} and \eqref{2htzk}, and hence the lemma, in both cases.
	\end{proof}

	\begin{proof}[Proof of \Cref{lem:verify2}]
		
		Observe under this choice of parameters that \eqref{eq:q'2} holds, and as does \eqref{domega}, where the last bound $qN^{1-3B\delta} \ge \omega$ in \eqref{domega} holds since $\omega \le N^{1+30A\delta} u(z) \le N^{1+30A\delta}(q-q') \le 2^k N^{31A\delta} \le N^{1-K\delta/2} q$, by \eqref{eq:costanti} and our assumptions on $(q,q')$ (and the estimates in \eqref{domega} are quickly verified). It remains to verify \eqref{eq:cond12} and \eqref{eq:cond32}. For the former, observe by \eqref{htzqkdelta} that $H_{s'} (w) \le \mathcal{H}_{q'} (w) + N^{5A\delta-1} \le \mathcal{H}_{q''} (w) + N^{5A\delta-1}$ holds w.o.p. for all $w \in \mathfrak{X}$; this confirms \eqref{eq:cond12} under our parameter choices.
		
		Next we verify \eqref{eq:cond32}. By the inductive hypothesis \eqref{2htzk} (applied with that $z$ equal to $w$ here) we have that $H_{s'} (w) \le \mathcal{H}_{q'''} (w)$ w.o.p. whenever $w \notin \mathfrak{L}_{q'''}^+ (\delta)$, where $q''' = q - \lfloor (s' - \mathsf{t}') N^{-10K\delta-2} / u(w) \rfloor u(w)$. Since $s' - \mathsf{t}' = s - \mathsf{t}' - N^{10K\delta+2} u(z) = (m-1) N^{10K\delta+2} u(z)$ and $(m-1) u(z) = q-q''$, we have that 
		\begin{flalign*} 
			(s'-\mathsf{t}') N^{-10K\delta-2} u(w)^{-1}  = (m-1) u(z) u(w)^{-1} = (q-q'') u(w)^{-1}.
		\end{flalign*} 
		
		\noindent As such, $q''' \le q'' + 2u(w)$, so $H_{s'} (w) \le \mathcal{H}_{q''+2u(w)} (w)$. By \Cref{boundarydistance} (and \Cref{rem:bysymmetry}), for some constant $C>1$, the distance in the $\bm{v}_w$-direction between $\mathfrak{A}_{q''+2u(w)}$ and $\mathfrak{A}_{q''}$ is at most $Cu(w) \le C \mathtt{d}_w^{-1/3} N^{-25A\delta-2/3}$. It quickly follows that, if $w \notin \mathfrak{L}_{q''}^+ (10A\delta)$, then $w \notin \mathfrak{L}_{q''+2u(w)}^+ (\delta)$. So, $\mathcal{H}_{q''+2u(w)} (w) = \mathcal{H}_{q''}(w)$ for all $w \notin  \mathfrak{L}_{q''}^+ (10A\delta)$. Therefore, $H_{s'} (w) \le \mathcal{H}_{q''+2u(w)} (w) = \mathcal{H}_{q''} (w)$ for $w \notin \mathfrak{L}_{q''}^+ (10A\delta)$, confirming \eqref{eq:cond32} under our parameter choices.	This shows the first statement of the lemma. The second follows from \Cref{lem:cellinclude1}.
	\end{proof}

	\section{Proof of Proposition \ref{htqhtq02}}
	
	\label{ProofEstimate0} 
	
	In this section, we prove \Cref{htqhtq02}. Recall from \Cref{sec61} that the proof of  \Cref{prop:hdecreasebulk} was based on the notion (justified by \Cref{prop:last}) that the normalized limit shape $\widehat{\mathcal{H}}_q$ from \eqref{e:funzrescaled} is ``more concave'' than $\widehat{\mathcal{H}}_{a,b,c}$. This reasoning does not seem to directly apply to show \Cref{htqhtq02}, since $\mathfrak{e}_{z;q'}$ is not approximately constant within the cell $\mathcal{Y}^-$. Thus the scaling required to define this renormalized height function $\widehat{\mathcal{H}}_{q'}$ will vary considerably, thereby skewing the notion of concavity, across $\mathcal{Y}^-$. 
	
	To remedy this, we will instead study level lines of $\mathcal{H}_{q'}$, using \Cref{prop:improvedconv} (in place of \Cref{prop:last}) to realize that those of $\mathcal{H}_{q'}$ are ``more convex'' than are those of $\mathcal{H}_{a,b,c}$. We will discuss this in more detail in the beginning of \Cref{Uq00h} below.

\begin{rem}
  \label{rem:wlog}
  Because of the definition \eqref{e0} of $\be$ and of the upper bound
  \eqref{domega} on $\omega$, we have $\be\le
  N^{-B\delta/2}$. The same considerations show that the product
  $\md_{w_0}\me_{w_0;q}$ in \eqref{e:dew0} is at most $N^{-B\delta}$. Therefore, for $N$ large the whole cell $\mathcal Y$ (recall \Cref{yy}) is within a distance $o(1)$ (with respect to $N$) 
  from the arctic boundary.  Without loss of generality, we will
  prove Propositions \ref{htqhtq02} and \ref{htq3} only in the case where
  $w_0$ is closer to $p^{SW}_q$ than to any of
  the other tangency points of $\mA_q$.  In fact, all other cases can be reduced to this one by
  symmetry.  In addition, we will only consider the case where the $x$-coordinate
  of $w_0$ is at least  $x^{SW}_q$ (recall \Cref{def:pa}). The proof of the propositions in
  the alternative case is entirely analogous.

  In view of this discussion, the cell $\mathcal Y$ is close to the frozen
  region $\mF^S_q$, where the limit
  shape $\mathcal H_q$ is zero.  By \Cref{hde0}, we
  have then, for $z\in\mathcal Y\cap \mL_q$ whose $x$-coordinate is at least $x^{SW}_q$, that
  \begin{equation}
    \label{e:wlog}
 \cH_q(z)\asymp \md_z^{1/2}\me_{z;q}^{3/2}.    
  \end{equation}
\end{rem}

Note for future reference that
\begin{eqnarray}
  \label{eq:qq'q''}
\left|\frac{q-q'}q\right|\le N^{-B \delta}
\end{eqnarray}
as follows from conditions \eqref{eq:q'} and \eqref{eq:condizioni}, together with  the upper bound in \eqref{domega}.
{Recall that, as in the statement of \Cref{htqhtq02}, the constant $L$ is assumed to be large enough (but independent of $N$).
}

\subsection{Case I:  \texorpdfstring{$\boldsymbol{d}=\mathfrak{d} $}{}}  

\label{Proofhtqd0d} 

In this section we show \Cref{htqhtq02} when $\mathfrak{d} = \bm{d}$. 

\subsubsection{Properties of cells and level lines} 

We start with a few preparatory lemmas. 

\label{ProofEstimated0} 

\begin{lem}
  \label{l:inthecell} For every $z=(x,y)\in \mathcal Y$, we have
  $\md_{w_0} / 2 \le \mathfrak{d}_{z}\le 2 \mathfrak{d}_{w_0}$ and
  $x\ge x^{SW}_q$.  Letting $h=\cH_q(z)$, we have
 { $s:=x-\mathfrak l_{q}(h)\asymp \md$}. If $q'$ satisfies \eqref{eq:q'} and $h':=\cH_{q'}(z)>0$, then   $s':=x-\mathfrak l_{q'}(h')\asymp \md$.
\end{lem}

\begin{proof}Write $w_0=(x_0,y_0)$. 
    Recall that
    \begin{eqnarray}
      \label{eq:y0large}
    y_0=\bd^2=\md_{w_0}^2\ge L^{20}N^{2B\delta}\frac\omega{q N},  
    \end{eqnarray}
    where the bound follows from \eqref{e0} and the assumption $\md=\bd$. As for $x_0$, note that $\me_{w_0;q}\le L\be \le L^{-19}\md_{w_0}$ (using \eqref{w00} and \eqref{e0}), which
    implies that
    \begin{eqnarray}
      \label{eq:x0large}
    x_0-x^{SW}_q\ge \md L^{-1},
    \end{eqnarray}
     if $L$ is sufficiently large, since the point of $\mA_q$ closest to, and with the same $y$-coordinate as, $w_0$ has $x$-coordinate larger than $x^{SW}_q$ by $c\mathfrak{d}$, for some  constant $c>0$.

    Since the slope of the vector $\bu_{w_0;q}$ is of order $\md$, the diameter of $\mathcal Y$ in the $y$ direction is of order
    $$L^5 \md_{w_0}\me_{w_0;q}+ \mathfrak{d} L^2 (\md_{w_0}\me_{w_0;q})^{1/2} \le L^6 N^{2B\delta}\frac\omega{qN}+L^{5/2}\bd N^{B\delta}\sqrt{\frac\omega{qN}}\le \displaystyle\frac{\mathfrak{\bm{d}^2}}{2L} =\displaystyle\frac{y_0}{2L},$$
    where the first bound uses \eqref{e:dew0} and the second uses \eqref{e0} (with the fact that $\bm{d} = \mathfrak{d}$). Together with \eqref{eq:y0large},
    this proves the first claim of the lemma for $L$ sufficiently large. 
    
    Similarly, the diameter
    of $\mathcal Y$ in the $x-$direction is, for some constant $C>1$, at most
    $C L^2\sqrt{\md_{w_0}\me_{w_0;q}}\le
    CL^{5/2}N^{B\delta}\sqrt{\omega/(qN)}\le CL^{-5}\md$. Together with \eqref{eq:x0large}, the claim
    $x\ge x^{SW}_q$ follows.  The statement $x-\mathfrak l_q(h)\asymp \md$ follows from the fact that $x-x^{SW}_q\asymp \md$ (as proven above) and $x^{SW}_q-\mathfrak l_q(h)\asymp \sqrt h\lesssim \md$ (the first bound follows from the definition of $\mathfrak l_q(h)$, together with the fact that the curvature of the arctic curve is bounded below, by \Cref{convex}, and the second from \eqref{e:wlog}).
    { As for the statement about $s'$, note that $s'\ge x-x^{SW}_{q'}\ge x-x^{SW}_{q}\asymp \md$, where the second bound follows from \Cref{boundarydistance} and the third one was proven just above. Moreover, if $\cH_{q'}(z)>0$ then $z\in\mL_{q'}$ and $x-x^{SW}_{q'}\lesssim \md$ follows from the fact that the curvature of  $\mA_{q'}$ is uniformly positive. Finally, $x^{SW}_{q'}-\mathfrak l_{q'}(h')\asymp \sqrt{h'}\le \sqrt{h}\lesssim \md$, where the first inequality follows from $q'\le q$ and the second was proven above. The bound $s'\asymp\md$ then follows.
}\end{proof}

Next, we require the following geometric definitions. Below, $\mathbf{S}$ prescribes a strip containing the cell $\mathcal{Y}$; we will focus on decreasing the height function $H_t$ in the part in $\mathbf{S}$ below the $k$-th level line $\mathcal{U}_{q'}^{k/N}$ of $\mathcal{H}_{q'}$, where $k \sim (\omega^3 \mathfrak{d})^{-1}$ (see \eqref{e:h0}) as, above this level line, the hypotheses of \Cref{htqhtq02} quickly imply that it holds (with the underlying reason being that, if $(\mathtt{d}_z \mathtt{e}_{z;q'})^{1/2} \le \omega^{-1}$, then by \Cref{hde0} we have $\mathcal{H}_q (z) \asymp \mathtt{d}_z^{1/2} \mathtt{e}_{z;q'}^{3/2} \lesssim (\omega^3 \mathfrak{d})^{-1}$, meaning that $z$ lies below $\mathcal{U}_{q'}^{k/N}$; see \Cref{rem:splits} for further details). Here, as always, $w_0=(x_0,y_0)$ is the point in $\mL_q$ that defines the center of the cell $\mathcal Y$. 

\begin{definition}\label{def:Lpm}
  Let ${\bf L}^{\pm}$ denote the line parallel to the $y$-axis and with $x$-coordinate $x^\pm:=x_0\pm L^2(\md_{w_0}\me_{w_0;q})^{1/2}$. Let $\bf S$ denote the strip between $\bf L^-$ and $\bf L^+$. Given $z\in
    \mathcal Y$, we let ${\bf L}^z$ denote the line parallel to the $y$-axis and containing $z$.    Note that the vertical sides of $\mathcal Y$ are subsets of $\bf L^-,\bf L^+$.

Moreover, as in \Cref{ulevelnew}, let\footnote{With a minor abuse of notation,  ``the level line $\mathcal U_{q'}^{k/N}$'' refers to the graph of the function $x\mapsto \mathcal U_{q'}^{k/N}(x)$.}
$\mathcal U_{q'}^{k/N}$ be the level line of $\cH_{q'}$ at height $k/N$, where
\begin{eqnarray}
	\label{e:h0}
	\displaystyle\frac{k}{N}=\frac{L^3}{\omega^3\md}.  
\end{eqnarray}

\end{definition}
  \begin{figure}[h]
 \begin{center} \includegraphics[width=7cm]{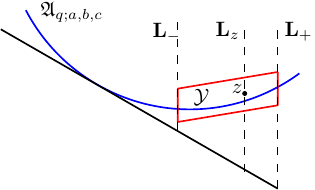} 
  \caption{{The cell $\mathcal Y$ with the vertical lines ${\bf L}_\pm,{\bf L}_z$ where $z\in \mathcal Y$. The strip $\bf S$ is contained between ${\bf L}_-$ and ${\bf L}_+$.}}
\label{fig:Def83}
\end{center}
\end{figure}

\noindent See \Cref{fig:Def83}. Note that, because of the bounds  $\omega\ge (qN)^{1/3}N^{-A \delta}$, $\md\ge N^{B\delta}\sqrt{\omega/Nq}$ and $Nq\ge N^{5B\delta}$ together with the fact that $B$ is sufficiently larger than $A$, we have $(k/N)\ll1$.

\begin{lem}\label{prop:Ucontained}
  The intersection $\mathcal U_{q'}^{k/N}\cap \bf S$ is contained in $\mathcal Y$.
\end{lem}

\begin{proof}[Proof of \Cref{prop:Ucontained}]
	
Let $z \in \mathcal{Y}$ denote any point on the top boundary of $\mathcal{Y}$. It suffices to show that $\mathcal{H}_{q'} (z) \ge kN^{-1}$. To do so, let us first compare $\mathcal{H}_{q'}(z)$ to $\mathcal{H}_q (z)$, to which end we compare $\mathfrak{e}_{z;q'}$ to $\mathfrak{e}_{z;q}$. By \Cref{hde0} (with \Cref{l:inthecell}), we have that $\mathfrak{e}_{z;q'} \asymp (\mathfrak{d}^{-1/2} \mathcal{H}_{q'} (z))^{2/3} \asymp L^2 \mathfrak{d}^{-1} \omega^{-2}$. We also have that $|q-q'| \le N^{R\delta-1} \omega,$ by \eqref{eq:q'}. Thus, by \eqref{domega} it follows that $|q-q'| \le N^{R\delta} \mathfrak{e}_{z;q'}$,  meaning (via \Cref{boundarydistance}) that $\mathfrak{e}_{z;q'} \le\mathfrak{e}_{z;q} \le N^{R\delta} \mathfrak{e}_{z;q'}$. Applying \Cref{hde0} again, this implies that $\mathcal{H}_q (z) \ll N^{2R\delta} \mathcal{H}_{q'} (z)$, so it suffices to show that $\mathcal{H}_q (z) \gg N^{2R\delta-1} k$. 

By \Cref{hde0} (with \Cref{l:inthecell}), we have $\mathcal{H}_q (z) \asymp \mathfrak{d}^{-1} (\mathfrak{d} \mathfrak{e}_{z;q})^{3/2}. $
Next, by \eqref{eq:elegeo} and the fact that $\dist (z; \mathfrak{A}_q) \gtrsim \mathfrak{d}_{w_0} \mathfrak{e}_{w_0;q}$ (by the definition  of $\mathcal{Y}$, with the smoothness of $\mathfrak{A}_q$), we have 
$\mathfrak{d} \mathfrak{e}_{z;q} \gtrsim \mathfrak{d}_{w_0} \mathfrak{e}_{w_0;q}$.
To lower bound the latter, observe that the lower bound on $\omega$ from \eqref{domega} yields that $\omega^3 \ge N^{1-6A\delta} q$. Together with \eqref{e:dew0} and the fact that $B \ge 1000A$, it follows that 
$\mathfrak{d}_{w_0} \mathfrak{e}_{w_0;q} \ge L^{-1} N^{2B\delta} \omega (qN)^{-1} \gg L N^{B\delta} \omega^{-2}$. Combining these bounds (and using the fact that $B \ge 1000R$) yields $\mathcal{H}_q (z) \gg N^{B\delta} \mathfrak{d}^{-1} \omega^{-3} \gg N^{2R\delta-1} k$, verifying the lemma.
\end{proof}

\begin{rem}
  \label{rem:splits}
As a consequence of \Cref{prop:Ucontained}, $\mathcal U_{q'}^{k/N}\cap \bf S$ splits $\mathcal Y$ into a region $\mathcal R^\uparrow$ above it where
$\cH_{q'}\ge \frac kN$ and the complement $\mathcal R^\downarrow$ where $\cH_{q'}< \frac kN$.
The statement of \Cref{htqhtq02} is obvious in the region $\mathcal R^\uparrow$ (defined in \Cref{rem:splits}), since $\mathcal R^\uparrow\subset \mathcal R$ (the latter being defined in \Cref{htqhtq02}). Indeed,
 because of Remark \ref{rem:wlog} and \Cref{l:inthecell} and the upper bound in \eqref{domega}, for { $w\in \mathcal R^\uparrow$}  one has $\me_{w;q'}^{3/2}\md_w^{1/2}\gtrsim L^{3}(\omega^3\md)^{-1}\asymp  L^{3}(\omega^3\md_w)^{-1} \ge N^{3A\delta-1}$ so that, in particular, $(\te_{w;q'}\td_w)^{1/2}\ge (\me_{w;q'}\md_w)^{1/2}\ge \omega^{-1}$ and $\mathtt{d}_w^{1/2} \mathtt{e}_{w;q'}^{3/2} \ge N^{A\delta-1}$, that is, $w\in\mathcal R$. \label{rem:obvious}  

Conversely, note that  $\mathfrak{d}_z^{1/2} \mathfrak{e}_{z;q'}^{1/2} \lesssim L \omega$ for all $z \in \mathcal{R}^{\downarrow}$. Indeed, by \eqref{e:wlog} and \eqref{e:h0}, we have that $\mathfrak{d}_{z'}^{1/2} \mathfrak{e}_{z';q'}^{1/2} \lesssim (k\mathfrak{d}/N)^{1/3} \lesssim L \omega$ for all $z' \in \mathcal{U}_{q'}^{k/N}$, and for every $z \in \mathcal{R}^{\downarrow}$ there exists some $z' \in \mathcal{U}_{q'}^{k/N}$ above $z$, with the same $x$-coordinate, so that $\mathfrak{d}_z^{1/2} \mathfrak{e}_{z;q'}^{1/2} \le \mathfrak{d}_{z'}^{1/2} \mathfrak{e}_{z';q'}^{1/2} \lesssim L \omega$.

Note also that, in the case $\bd=\md$, the condition $\mathtt d_z\in(L^{-2}\md,L^2\md)$ is automatically satisfied in the whole cell $\mathcal Y$ {for $L$ large enough}, in view of \eqref{condizd} and \Cref{l:inthecell}.
\end{rem}

\subsubsection{Comparison of level lines}

\label{Uq00h}

For the remainder of the proof of Proposition \ref{htqhtq02} when $\mathfrak{d} = \bm{d}$, we fix a point $\bar z=(\bar x,\bar y)\in \mathcal Y^-$; we then verify the bound \eqref{eq:statmesoe} w.o.p. at $\bar{z}$. 

Let us briefly outline how we will proceed. We will eventually (see the beginning of the proof of \Cref{htqhtq02} in \Cref{Proof000}) censor the Glauber dynamics $(H_t)$ outside of $\mathcal{R}^{\downarrow}$ and run them inside of $\mathcal{R}^{\downarrow}$ (with a floor constraint) until they mix, which will by \Cref{prop:tmixconstrained} take time at most $q^{-1} N^{1+5B\delta} \omega$. Since the hypotheses of \Cref{htqhtq02} imply that $H_t$ is bounded above by $\mathcal{H}_{q'}$ along the upper boundary of $\mathcal{R}^{\downarrow}$ and by $\mathcal{H}_q$ along its left and right boundaries, monotonicity implies that $H_t$ is bounded above by the uniform measure on tilings of $\mathcal{R}^{\downarrow}$ with these boundary conditions. To analyze this stationary measure, we will compare to random tilings of a suitably chosen hexagon. 

Specifically, we will first apply Proposition \ref{prop:improvedconv} to find $(a,b,c)$, such that the following holds. The $k$-th level line $\mathcal{U}_{a,b,c}^{k/N}$ of $\mathcal{H}_{a,b,c}$ is (after a shift by $(\mathfrak{r}, \mathfrak{r}')$) tangent to, and less convex than, the $k$-th level line $\mathcal{U}_{q'}^{k/N}$ of $\mathcal{H}_{q'}$ at $\bar{z}$. This difference in convexity will imply that the level lines of $\mathcal{H}_{q'}$ are above those of $\mathcal{H}_{a,b,c}$ along the boundary of $\mathcal{R}^{\downarrow}$; see \Cref{lem:rRlL}\footnote{The idea of comparing two random tilings by tuning the curvatures of their limit shape level lines was also used in \cite[Section 3]{aggarwal2025edge} (though in the distinct context of accessing local statistics).}. Together with monotonicity and the concentration estimate \Cref{volumeheight} for random tilings of $\mathfrak{X}_{a,b,c}$, this will imply that the height function for the uniform measure described above is w.o.p. bounded above by $\mathcal{H}_{a,b,c}$ (up to an error of $N^{\delta-1}$). Using \eqref{e:primamail}, it will be verified (see \Cref{lem:UbarGbar}) that the latter is at most $\mathcal{H}_{q'} (\bar{z})$, showing the proposition. 

To implement this, we must first verify that \Cref{prop:improvedconv} applies here. Below, given $x\in[x^-,x^+]$, set $s_x:=x-\mathfrak l_{q'}(k/N)$. 

\begin{lem}
	Let $s_\pm:=s_{x^\pm},s_0:=s_{\bar x}$. Then, the assumptions of \Cref{prop:improvedconv} are verified, with $q,\varepsilon,h_0$ there equal to $q',\varepsilon_0,(k/N)$ here.
  \label{lem:hereandthere}
\end{lem}

\begin{proof}[Proof of \Cref{lem:hereandthere}]
	
	The condition $k/N\le \varepsilon_0$ holds since $k/N \ll 1$, as noted just after \eqref{e:h0}. Also, condition $q'\in[0, \varepsilon_0]$ holds thanks to \eqref{eq:q'}, \eqref{eq:condizq} and \eqref{eq:qq'q''}.

  The first bound in \eqref{eq:condizs} holds because $\bar x\ge x^{SW}_q$ (by \Cref{l:inthecell}) and therefore $x\ge x^{SW}_{q'}$ as well (by \Cref{boundarydistance}). The last bound in \eqref{eq:condizs} holds  because $\bar z$ is closer to the $p^{SW}_q$ than to $p^{SE}_q$ (recall \Cref{rem:wlog}), and parameters in the definition of the cells ensure that their size tends to zero as $N\to\infty$ (recall \Cref{omegaq0}).

  As for condition \eqref{eq:condiz1},
  the upper bound on $s_+-s_-$ is satisfied in our case (if $L$ in the definition of $\mathcal Y$ is large enough) since the width of the cell $\mathcal Y$ in the $x$-direction is at most $L^{-7}\md_{w_0}\lesssim L^{-7}s_0\lesssim L^{-7}s_0^{1/2}$, as follows from \eqref{e0}, \eqref{w00},  \eqref{otherwise}. For the bound $\md_{w_0}\lesssim s_0$, we used \Cref{l:inthecell}.
\end{proof}

Thanks to \Cref{prop:improvedconv}, we can find $a,b,c\in [1-C\varepsilon_0,1+C\varepsilon_0]$
and $\mathfrak r,\mathfrak r'$ such that the statements of the proposition hold.
\begin{rem}\label{rem:nor}
    The specific value of the constants $\frak r,\frak r'$ plays no role in the subsequent arguments. Therefore, to lighten formulas, we will drop them (assume they equal zero).
  \end{rem}
In particular, we have by \eqref{eq:meccia} and \eqref{eq:r'} that 
 \begin{eqnarray}
    \label{eq:tangent}
  \mathcal U^{k/N}_{q'}(\bar x)=\mathcal U^{k/N}_{a,b,c}(\bar x), \quad  \partial_x\mathcal U^{k/N}_{q'}(\bar x)=\partial_x\mathcal U^{k/N}_{a,b,c}(\bar x) 
  \end{eqnarray}
  
  \noindent and by \eqref{eq:meccia2}, \Cref{l:inthecell} and  \eqref{eq:qq'q''} that, for $x \in [x^-, x^+]$, 
\begin{equation}
  \label{eq:tangent2}
  \displaystyle
  \big(  \mathcal{U}_{q'}^{k/N} (x) - \mathcal{U}_{a,b,c}^{k/N} (x)  \big) \ge C^{-1} q' \mathfrak{d} |x-\bar{x}|^2 \ge (2C)^{-1} q \mathfrak{d} |x - \bar{x}|^2.
\end{equation}

\begin{definition}
  \label{def:lLrR}
 For $j\le k$, with $k$ as in \eqref{e:h0}, let $y^{\pm,j}_q:=\mathcal U^{j/N}_q(x^\pm)$, $y^{\pm,j}_{q'} :=\mathcal U^{j/N}_{q'} (x^\pm)$, and $y^{\pm,j}_{a,b,c}:= \mathcal U^{j/N}_{a,b,c}(x^\pm)$. 
\end{definition}

Note that these are the $y$-coordinates of the intersections of $\bf L^\pm$ with the $(j/N)$-level line of $\cH_q$ and $\cH_{a,b,c}$ (the latter up to a translation of the hexagon). The next lemma compares the boundary data of the level lines of $\mathcal{H}_{q'}$ and $\mathcal{H}_{a,b,c}$ on $\mathcal{R}^{\downarrow}$.

\begin{lem}
  \label{lem:rRlL}
   For $1\le j\le k$, we have that $y^{\pm,j}_q- y^{\pm,j}_{a,b,c}\ge{L^2} N^{2B\delta-1}\omega\md$. Moreover, $\mathcal{U}^{k/N}_{q'} (x) \ge \mathcal{U}^{k/N}_{a,b,c} (x) $ for all $x \in [x^-, x^+]$. 
\end{lem}

\begin{proof}
  First of all, we claim that
  \begin{eqnarray}
    \label{eq:claimk}
      y^{\pm,k}_{q'}-y^{\pm,k}_{a,b,c}\ge(100 C)^{-1} L^3 N^{2B\delta-1}\omega\md.
  \end{eqnarray}
  In fact, this follows from \eqref{eq:tangent}, \eqref{eq:tangent2}, and the fact that  $\min(|\bar x-x^+|,|\bar x-x^-|)\ge \frac14(x^+-x^-)$ (because  $\bar{z} \in\mathcal Y^-$ and $L$ is large) and lastly 
$(x^+-x^-)^2=
4L^4 (\md_{w_0}\me_{w_0;q})\ge 4 L^3 N^{2B\delta}\omega(qN)^{-1}$ (by \eqref{otherwise} and \eqref{e:dew0}). All these facts, together, imply \eqref{eq:claimk}.
Next, we claim that
  \begin{eqnarray}
    \label{eq:claimj}
      y^{\pm,j}_{q'}-y^{\pm,j}_{a,b,c}\ge(200 C)^{-1} L^3 N^{2B\delta-1}\omega\md
  \end{eqnarray}
for every $j\le k-1$. This follows if we prove that
  \begin{eqnarray}
  |(y^{\pm,k}_{a,b,c}-y^{\pm,j}_{a,b,c})-(y^{\pm,k}_{q'}-y^{\pm,j}_{q'})|\ll  L^3 N^{2B\delta-1}{\omega\md}
  \end{eqnarray}
for
  $1\le j<k-1$.
  Using \eqref{e:primamail}, \eqref{eq:primamail'} and recalling from \Cref{l:inthecell} that
  $s_0\asymp \md$, we see
  that
  \begin{multline}
    \label{eq:imaa}
    |(y^{\pm,k}_{a,b,c}-y^{\pm,j}_{a,b,c})-(y^{\pm,k}_{q'}-y^{\pm,j}_{q'})|\lesssim q \md^{5/3}(k/N)^{2/3}=q\md \Big(\md \frac kN\Big)^{2/3}\\= L^2 \md\omega\times \frac q{\omega^3}\le L^2 N^{6A\delta-1}{\omega\md}
  \end{multline}
where in the last inequality  we used the lower bound \eqref{domega} for $\omega$. The claim \eqref{eq:claimj} then follows since $6A<2B$.

Finally, the claim of the lemma follows from \eqref{eq:claimj} since, 
by \Cref{qlineq} and the second part of \Cref{hde0}, we have that $| y^{\pm,j}_{q}- y^{\pm,j}_{q'}|\asymp\md(q-q')\lesssim \md\omega N^{R\delta-1}\ll  L^3
  N^{2B\delta-1}{\omega\md}$
 (the last bound holds because of \eqref{eq:condizioni}). 
\end{proof}

\subsubsection{Proof of \Cref{htqhtq02} if $\mathfrak{d} = \bm{d}$}

\label{Proof000} 

\begin{figure}[h]
 \begin{center} \includegraphics[width=7cm]{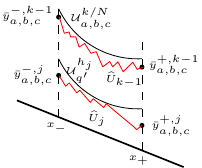}
  \caption{{A schematic illustration of \Cref{def:lineensembles}. The top Bernoulli path $\widehat U_{k-1}$ is weakly below $\mathcal U^{k/N}_{a,b,c}$ and the Bernoulli path $\widehat U_j,j<k-1$, is weakly below $\mathcal U^{h_j}_{q'}$, where $h_j:=j/N+N^{-1+(6A+R)\delta}$.}}
\label{fig:Def810}
\end{center}
\end{figure}

As mentioned in the beginning of \Cref{Uq00h}, we can (after running the Glauber dynamics until they mix on $\mathcal{R}^{\downarrow}$) bound $(H_t)$ above by the uniform measure on random tilings of $\mathcal{R}^{\downarrow}$, with upper level line given by $\mathcal{U}_{q'}^{a,b,c}$ and left and right level line endpoints given by $(y_q^{\pm,j})$. By \Cref{lem:rRlL} and monotonicity, we may replace the upper boundary condition by $\mathcal{U}_{a,b,c}^{k/N}$ and the left and right ones by $\mathcal{U}_{a,b,c}^{j/N}$. The following definition provides notation for the associated Bernoulli path ensemble (here, $\widehat{U}$ involves an additional floor condition to be briefly explained below \Cref{def:lineensembles} ). 
 See also \Cref{fig:Def810}.
 
\begin{definition}[$\widehat U,\widecheck U$]
  \label{def:lineensembles} For $j\le k-1$, let $(x^\pm,\bar y^{\pm,j}_{a,b,c})$ be a point in $\mathbb T_N+(0,1/(2N))$ at minimal distance from $(x^\pm, y^{\pm,j}_{a,b,c})$. We
  let $\widehat U=\{\widehat U_j\}_{1\le j\le k-1}$ denote the Bernoulli path ensemble
  with endpoints given by
  $\{(x^\pm,\bar y^{\pm,j}_{a,b,c})\}_{1\le j\le k-1}$, conditioned to
  the event that the top path $\widehat U_{k-1}$ is weakly {below}
  $\mathcal U^{k/N}_{a,b,c}$ and that each $\widehat U_j,j\le k-1$
  is weakly below $ \mathcal U^{(j/N+N^{-1+(6A+R)\delta})}_{q'}$.
  Similarly, we let $\widecheck U=\{\widecheck U_j\}_{1\le j\le k-1}$ denote the
  Bernoulli path ensemble defined like $\widehat{U}$, but without the constraint
  that $\widecheck U_j$ is below
  $ \mathcal U^{(j/N+N^{-1+(6A+R)\delta})}_{q'}$.
\end{definition}

Note that the condition that each $\widehat U_j$ is below $\mathcal U^{N^{-1}(j/N+N^{-1+(6A+R)\delta})}_{q'}$ is equivalent to the condition that the height function associated to the Bernoulli path ensemble $\widehat U$ is {above} the floor 
  $\cH_{q'}-N^{-1+(6A+R)\delta}$. This floor constraint will be a convenience when applying \Cref{prop:tmixconstrained}, as it will indicate that the quantity $H_{\max}$ from \eqref{eq:tmixconstrained} is nearly its minimal value $N^{o(1)-1}$ (see \eqref{q0q0} below). Since the concentration bound \Cref{cor:conclines} applies to $\widecheck{U}$ and not $\widehat{U}$, the following lemma, stating that this constraint is irrelevant, will be useful.

\begin{lem}
  \label{lem:whpcoupling}
  The Bernoulli path ensembles  $\widehat U,\widecheck U$ can be coupled in such a way that, w.o.p., $\widehat U=\widecheck U$.
\end{lem}
\begin{proof}[Proof of \Cref{lem:whpcoupling}]
  The claim follows once we prove
  \begin{eqnarray}
    \label{eq:lpa1}
\text{w.o.p., } \widecheck U_j \text{  is weakly below }
 \mathcal U^{(j/N+N^{-1+(6A+R)\delta})}_{q'}, \text{ for every } j\le k-1.
  \end{eqnarray}
In agreement with \Cref{ttU}, let\footnote{With a minor abuse of notation and in order to lighten formulas, whenever we write $\mathtt U^h_{a,b,c}$ we actually mean $\mathtt U^{h'}_{a,b,c}$, with $h'$ the largest element of  $(1/N)(\mathbb Z+1/2)$ smaller than $h$.\label{ftn:minorabuse}} $\mathtt U^h_{a,b,c}$ denote the level line at height $h$ of the height function of a uniformly random tiling of $\mX_{a,b,c}$. As above, we assume that the hexagon $\mX_{a,b,c}$ is suitably translated, so that the level line at height $j/N$ of its limit shape coincides with $\mathcal U^{j/N}_{a,b,c}$. 
In the cell $\mathcal Y$, we have that $\widecheck U_{k-1}$ is below $\mathcal U^{k/N}_{a,b,c}$ (by definition of the Bernoulli path ensemble $\widecheck U$), and the latter is w.o.p. below
$\mathtt U^{k/N+N^{\delta-1}}_{a,b,c}$, by
\Cref{cor:conclines}\footnote{When applying \Cref{cor:conclines}, note that
$\mathtt d_z\asymp\md$ in $\mathcal Y$, since $\md\ge N^{-1/2}$ by assumption \eqref{condizd}.}.
Also by \Cref{cor:conclines} and by definition of $\widecheck U$, the endpoints of $\widecheck U_j,j\le k-1$, are w.o.p. below the intersections of $\mathtt U^{j/N+N^{\delta-1}}_{a,b,c}$ with $\bf L^\pm$. Applying \Cref{le:stochdomB}, we conclude that the Bernoulli path ensemble $\widecheck U$ can be coupled with the level lines $\mathtt U_{a,b,c}$ so that, w.o.p.,
$\widecheck U_j$ is below $\mathtt U^{j/N+N^{\delta-1}}_{a,b,c}$ for every $j\le k-1$,  in the cell $\mathcal Y$.
Applying once more \Cref{cor:conclines}, we conclude w.o.p. that
$\widecheck U_j$ is below $\mathcal U^{j/N+2N^{\delta-1}}_{a,b,c}$.
  It remains to prove that the latter is below
$ \mathcal U^{(j/N+N^{-1+(6A+R)\delta})}_{q'}$. To this purpose, we apply \eqref{heighthabcq} (with $h_0,\md_z$ there replaced by $k/N,\md$ here) and note that  $C q\mathfrak{d} kN^{-1} \asymp L^3 q \omega^{-3} \lesssim L^3 N^{6A\delta-1}$ by \eqref{domega}. In particular, $N^{\delta-1}+Cq\mathfrak{d} kN^{-1}\le N^{-1+(6A+R)\delta}$. Then, \eqref{heighthabcq} implies that 
$\mathcal U^{j/N+2N^{\delta-1}}_{a,b,c}$ is below $ \mathcal U^{(j/N+N^{-1+(6A+R)\delta})}_{q'}$, which concludes the proof.
\end{proof}

The next lemma indicates that the Bernoulli path ensemble $(\widehat{U}_j)$ is closely approximated by the level lines $(\mathcal{U}_{a,b,c}^{j/N})$ of $\mathcal{H}_{a,b,c}$; it quickly follows from first comparing $\widehat{U}$ to $\widecheck{U}$ through \Cref{lem:whpcoupling} and then using \Cref{cor:conclines}  (which follows from the concentration bound \Cref{volumeheight}) to approximate the latter. The subsequent lemma compares the level lines $(\mathcal{U}_{a,b,c}^{j/N})$ to $\mathcal{U}_{q'}^{j/N}$ at $x$-coordinate $\bar{x}$ (where we recall that we have fixed a point $\bar z=(\bar x,\bar y)\in \mathcal Y^-$, and that ${\bf L}^{\bar z}$ is the vertical line containing $\bar z$).

\begin{lem}
  \label{lem:fluttU}
 W.o.p., the following holds for every $1\le j<k$:
 the Bernoulli path  $\widehat U_j$ is between $\mathcal U^{j/N - N^{-1+(3/4)\delta}}_{a,b,c}$ and $\mathcal U^{j/N + N^{-1+(3/4)\delta}}_{a,b,c}$ if $j-N^{(3/4)\delta}\ge0$, or  below $\mathcal U^{j/N + N^{-1+(3/4)\delta}}_{a,b,c}$ if $j-N^{(3/4)\delta}<0$.
\end{lem}
\begin{proof}[Proof of \Cref{lem:fluttU}]
  By  \Cref{lem:whpcoupling} it is enough to prove the statement for $\widecheck{U}$. Then, the statement follows from
\Cref{cor:conclines}.
\end{proof}

\begin{lem}
\label{lem:UbarGbar}  
For every $1\le j<k$ the following holds:
$\mathcal U^{j/N}_{a,b,c}(\bar x)\ge \mathcal U^{j/N}_{q'}(\bar x)$.
\end{lem}
\begin{proof}[Proof of \Cref{lem:UbarGbar}  ]

Since \eqref{eq:tangent} implies that $\cU_{q'}^{k/N} (\bar x) = \cU_{a,b,c}^{k/N} (\bar x)$, the lemma follows from \eqref{e:primamail} (with the $h$ there equal to $j/N$ here).
\end{proof}

Finally we can prove \Cref{htqhtq02} (following the outline in the beginning of \Cref{Uq00h}) in the case $\bd=\mathfrak d$.

\begin{proof}[Proof of \Cref{htqhtq02} if $\bd=\mathfrak d$]
  As noted in 
  \Cref{rem:splits},
  we need only to prove that  w.o.p, for all times in the time interval $[q^{-1} N^{1+10B\delta}\omega, N^{5}]$ the event
  \begin{equation}
    \label{eq:theevento}
H_t(\bar z)\le \cH_{q'}(\bar z)+N^{\delta-1} \quad \text{for all}\quad \bar z\in\mathcal Y^-\cap\mathcal R^\downarrow
  \end{equation}
  occurs.  Since the event in \eqref{eq:theevento} is decreasing, by
  monotonicity (\Cref{prop:cancouple}) and by \Cref{rem:splits} we can
  and will impose that the height function $(H_t)_{t \ge 0}$ {in
    $\mathcal Y\cap \mathcal R$} (and in particular in
  $\mathcal Y\cap \mathcal R^\uparrow$, that is the region above the
  line $\mathcal U^{k/N}_{q'}$) is frozen to the maximal configuration
  lower than the ceiling $\cH_{q'}$, and that $(H_t)_{t \ge 0}$ in
  $\mathcal Y$ is constrained (in the sense of \Cref{def:constrained})
  to be above the floor $\cH_{q'}-N^{-1+(6A+R)\delta}$.

For the dynamic $(H_t)_{t \ge 0}$ in the cell $\mathcal Y$ with these floor/ceiling constraints, we can upper bound via \Cref{prop:tmixconstrained} its mixing time $T_{\rm mix}^{\mathcal Y}$ by a constant  times 
\begin{eqnarray}
  \label{eq:ctimes}
  N^2\times L^4 (\me_{w_0;q}\md_{w_0})N^2(\log N)^2  N^{(12A+3R)\delta-2}\lesssim L^5 N ^{1+2(B+(3/2)R+6A)\delta}\frac\omega q\lesssim N^{1+5 B\delta}\frac\omega q,
\end{eqnarray}
where we used also \eqref{e:dew0} and \eqref{eq:condizioni}.
In fact, in applying \Cref{prop:tmixconstrained} we use the fact that
the difference  $|\cH_q(z)-\cH_{q'}(z)|$  is (by \Cref{qlineq}) at most of order
\begin{equation}
	\label{q0q0} 
(q-q')(\me_{z;q}\md_{z})^{1/2}+(q-q')^{3/2}\md_z^{1/2}\lesssim (\me_{z;q}\md_{z})^{1/2}N^{R\delta-1}\omega+N^{-1+(3/2)R\delta}\lesssim N^{(3/2)R\delta-1}  
\end{equation}
 for $z\in\mathcal Y\setminus\mathcal R$, where the first bound follows from the assumption \eqref{eq:q'} on $q'$ and the upper bound \eqref{domega} on $\omega$, while the last bound follows from the fact that $z\in\mathcal Y\setminus \mathcal R$, where $(\me_{z;q}\md_{z})^{1/2}\lesssim \omega$ by \Cref{rem:splits} (recall that we are restricting our attention to points in $\mathcal R^\downarrow$). Because of \eqref{e:submult}, we  deduce that the law of the dynamic at time  $t= N^{1+10B\delta}\frac\omega q\ge  N^{B\delta} T_{\rm mix}^{\mathcal Y} $ is at  variation distance smaller than any inverse power of $N$  from its stationary  measure $\pi$. 

 We prove now that, w.o.p., if $H^{\eq}$ denotes the discrete height of the configuration sampled according to  $\pi$,
 \begin{equation}
   \label{e:eventopi}
H^{\eq}(\bar z)\le \cH_{q'}(\bar z)+N^{\delta-1} \quad \text{for all}\quad \bar z\in\mathcal Y^-\cap\mathcal R^\downarrow. 
\end{equation}

Once that is proven, it follows by \Cref{rem:forall} and \eqref{e:eventopi} that the event \eqref{eq:theevento} holds w.o.p. in the whole time interval $[q^{-1} N^{1+10B\delta}\omega, N^{5}]$, which concludes the proof of the proposition. 

 To prove \eqref{e:eventopi}, given $j\in\mathbb N,$ with $ j\le k$, call $W^j$ the  level line of $H^{\eq}$ at height $j/N-1/(2N)$ and 
    $(x^\pm,Y^{\pm,j})$ the intersection of $W^j$ with $\bf L^\pm$.
Because of assumption (i) of \Cref{htqhtq02} and  \Cref{lem:rRlL}, we have
  $Y^{\pm,j}\ge y^{\pm,j}_q\ge y^{\pm,j}_{a,b,c}$ for all $j\le k$.
 By assumption (ii) of \Cref{htqhtq02} and \Cref{rem:splits},
the level line $W^k$ is deterministically weakly above
$\mathcal U^{k/N}_{q'}$. Hence, by \Cref{le:stochdomB}, the Bernoulli path ensembles $W=\{W^j\}_{1\le j\le k-1}$ and $\widehat U$ (cf. \Cref{def:lineensembles}) can be coupled in such a way that $W^j$ is weakly above $\widehat U_j$ for every $j\le k-1$ (note that the constraint that the $j$-th line is weakly below $\mathcal U_{q'}^{j/N+N^{-1+(6A+R)\delta}}$ is present in both ensembles). Note also that the event in \eqref{e:eventopi} is decreasing in terms of height function, but increasing in terms of Bernoulli paths. Therefore, we can replace the ensemble $W$ by $\widehat U$.

Call $j\le k$ the index such that
$\lfloor N\cH_{q'}(\bar z)\rfloor=j$ (recall that, as noted in \Cref{rem:splits}, $\cH_{q'}(\bar z)\le k/N$).
The event \eqref{e:eventopi} is implied by  the event that the
level line $U_{j+N^\delta}$ stays above the point $\bar z$, that is, at the
intersection with ${\bf L}^{\bar z}$ the level line $\widehat U_{j+N^\delta}$ is above
$\mathcal U^{j}_{q'}$. We can assume that $j+N^\delta<k$, since (as noted above) the level line at height $k/N$ is above $\mathcal U_{q'}^{k/N}$, hence above $\bar z$. 
 On the other hand, for $j< k-N^\delta$ we have
that $\widehat U_{j+N^\delta}$ is w.o.p. above $\mathcal U_{a,b,c}^{j/N}$ (by
\Cref{lem:fluttU}), which at the intersection with ${\bf L}^{\bar z}$ is above $\mathcal U^{j/N}_{q'}$ (by
\Cref{lem:UbarGbar}), as desired.
\end{proof}

\subsection{Case II: \texorpdfstring{$\bd>\md$}{}}

\label{Proofd}

  In this section we establish \Cref{htqhtq02} in the case when $ \bm{d}>\md$. First observe since $\bm{d} > \mathfrak{d}$ and since $\bm{d} \le N^{(B+1)\delta} (\omega/qN)^{1/2} \le N^{(1-B/2)\delta}$ (by \eqref{domega} and \eqref{e0}) that we have 
  \begin{flalign}
    \label{eq:mdsmall}
	\mathfrak{d} < N^{-B\delta/4}.
\end{flalign} 

\subsubsection{Inductive setup}

The proof of \Cref{htqhtq02} if $\bm{d} > \mathfrak{d}$ will have a number of similarities to that when $\mathfrak{d} = \bm{d}$. However, one relevant difference relates to the estimates \eqref{eq:meccia2} and \eqref{eq:primamail'}. When $\mathfrak{d} = \bm{d}$, we used that the former ``increase in curvature'' dominates the latter ``error,'' when showing that $y_{a,b,c}^{+,j} \le y_{q}^{+,j}$ (recall Definition \ref{def:lLrR}, where the right endpoints will be the main relevant ones when $\bm{d} > \mathfrak{d}$) in the proof of Lemma \ref{lem:rRlL}. Now, when $ \bm{d} > \mathfrak{d}$, a direct attempt to follow the same reasoning could have us take $s_0 = \mathfrak{d}$ (from reducing $H_t(z)$ where $\mathfrak{d}_z = \mathfrak{d}$) and $s_+ = \bm{d}$ (from the length of the cell $\mathcal{Y}$) in Proposition \ref{prop:improvedconv}. However, if $\mathfrak{d} \ll \bm{d}$, then the error in \eqref{eq:primamail'} could dominate the increase in curvature from \eqref{eq:meccia2}, which would prevent us from ordering $\mathcal{H}_q$ with respect to $\mathcal{H}_{a,b,c}$. To circumvent this, one would like to reduce the value of $s_+$ in Proposition \ref{prop:improvedconv}, by applying that proposition on a smaller length scale now only slightly larger than $\mathfrak{d}$. This would enable us to show that $y_{a,b,c}^{\pm,j} \le y_{q'}^{\pm,j}$. However, if the length scale is too small, then we would not be able to deduce that $y_{a,b,c}^{\pm,j} \le y_q^{\pm,j}$. Indeed, the difference $y_q^{+,j} - y_{a,b,c}^{+,j}$ would be small (as the increase in curvature \eqref{eq:meccia2} would be used on a fairly short interval), and insufficient to overcome the difference  $y_{q'}^{+,j} - y_q^{+,j}$. 

To address this, we will effectively ``induct on $\mathfrak{d}$''; in particular, we will first show that $H_t (z)$ reduces when $\mathtt{d}_z \sim \bm{d}$, and then decrement $\mathtt{d}_z$ by factors of $N^{\delta}$ until it is close to $\mathfrak{d}$. For a given value of $\bar{\mathtt{d}} = \mathtt{d}_z$, we will apply Proposition \ref{prop:improvedconv} with the $(s_0, s_+)$ approximately equal to $(\bar{\mathtt{d}}, N^{B\delta/9} \bar{\mathtt{d}})$ here. After the first several decrements (where $\bar{\mathtt{d}} \approx \bm{d}$, so the arguments are similar to in the $\mathfrak{d} = \bm{d}$ case), the right endpoints of these cells will by the inductive hypothesis be in regions where we have already reduced the height of $(H_t)$ to $\mathcal{H}_{q'}$ (up a small additive error). In this way, we will only have to compare $\mathcal{H}_{a,b,c}$ to $\mathcal{H}_{q'}$ (and not to $\mathcal{H}_{q}$), for which the smaller length scale of $N^{B\delta/9} \bar{\mathtt{d}}$ of the interval will suffice.

The next proposition, which quickly implies Proposition \ref{htqhtq02}, provides this inductive statement. It indicates how long it takes for $H_t (z)$ to decrease to nearly $\mathcal{H}_{q'} (z)$, depending on the value of $\mathtt{d}_z$. 

\begin{prop} 
	
	\label{dhtqd} 
	
	For any integer $\mathfrak{K} \in [1, \delta^{-1}]$ and time $t \in [q^{-1} \delta N^{1+10B\delta} \omega \mathfrak{K}, N^5]$, w.o.p. the following holds. We have that $H_t (z) \le \mathcal{H}_{q'} (z) + \mathfrak{K} N^{\delta/2-1}$, for each $z \in \mathcal{Y}^-$ satisfying $\mathtt{d}_z \ge L^{-2} N^{1-\mathfrak{K}\delta} \mathfrak{d}$.
	
\end{prop} 

\begin{proof}[Proof of \Cref{htqhtq02} if $\mathfrak{d} < \bm{d}$]
	This follows from taking $\mathfrak{K} = \delta^{-1}$ in the second statement of Proposition \ref{dhtqd}.
\end{proof}

To show Proposition \ref{dhtqd}, we induct on $\mathfrak{K}$, the statement being true at $\mathfrak{K} = 1$, 
as then there does not exist any $z \in \mathfrak{X}$ for which $\mathtt{d}_z \ge L^{-2} N^{1 -\mathfrak{K}\delta} \mathfrak{d}$ (as $L^{-2} N^{1 - \mathfrak{K}\delta} \mathfrak{d} \gg 1$, since $\mathfrak{d} \ge N^{-1/2}$). So, assume that Proposition \ref{dhtqd} holds at $\mathfrak{K}-1$, and we will show it at $\mathfrak{K}$. 

Throughout, we fix 
a point $\bar z=(\bar x,\bar y) \in \mathcal Y^-$ satisfying $\mathtt{d}_{\bar{z}} \ge L^{-2} N^{1 - \mathfrak{K}\delta} \mathfrak{d}$. We further set $\mathsf{T} = q^{-1} \delta N^{1+10B\delta} \omega \mathfrak{K}$, which is the lower bound for the time interval considered in Proposition \ref{dhtqd}. Then we must show that, w.o.p., for all $t \in [\mathsf{T} , N^5]$ we have 
\begin{flalign} 
	\label{htkq} 
		H_t (\bar{z}) \le \mathcal{H}_{q'} (\bar{z}) + \mathfrak{K} N^{\delta/2-1}.
\end{flalign}

\noindent  To that end, we may assume that $\mathtt{d}_{\bar{z}} < L^{-2} N^{1-(\mathfrak{K}-1)\delta} \mathfrak{d}$, as otherwise these results follow from the inductive hypothesis given by the second part of Proposition \ref{dhtqd} (at $\mathfrak{K}-1$). Thus, we have
\begin{flalign}
	\label{0d0} L^{-2} N^{1-\mathfrak{K}\delta} \mathfrak{d} \le \mathtt{d}_{\bar{z}} \le L^{-2} N^{1-(\mathfrak{K}-1)\delta} \mathfrak{d}.
\end{flalign}

\noindent We can also assume by symmetry that 
\begin{flalign} 
	\label{dz00} 
	\bar{x} \ge x_q^{SW}. 
\end{flalign}

\subsubsection{Properties of cells and the $k$-th level line} 

To show \eqref{htkq}, we follow the reasoning in Section \ref{ProofEstimate0}, with the $\mathfrak{d}$ there replaced by $\td_{\bar z}$ here. We first define the analog of $k$ from \Cref{def:Lpm} (it will be larger here than there, which will eventually be a convenience to compare level lines). As in \Cref{ulevelnew}, let $\mathcal{U}_{q'}^{k/N}$ denote the $(k/N)$-level line of $\cH_{q'}$, where  
\begin{equation}
	\label{eq:h02}
	\frac kN = \frac{N^{5R\delta}}{\omega^3 \td_{\bar z}}.
\end{equation}	

\noindent Observe by \eqref{domega} (and the fact that $\td_{\bar z} \ge L^{-2} \mathfrak{d}$) that 
\begin{flalign}
	\label{knd} 
	\displaystyle\frac{k}{N} \le N^{(5R+6A)\delta} \cdot \displaystyle\frac{\mathfrak{d}^3}{\td_{\bar z}} \le L^2 N^{11R\delta} \td_{\bar z}^2,
\end{flalign}

\noindent meaning that $k/N$ is not too much larger than $\td_{\bar z}^2$. Also observe that 
\begin{flalign}
	\label{adeltad}
	N^{(A+4R)\delta-1} \ll kN^{-1} \le L^2 N^{5R\delta} (\mathfrak{d} \omega^3)^{-1} \ll N^{B\delta} (\mathfrak{d} \omega^3)^{-1} \ll \bm{d}^2,
\end{flalign} 

\noindent where the first and second bounds hold since $\td_{\bar z} \ge L^{-2} \mathfrak{d} \ge L^{-2} N^{-1/2}$; the third holds since $B \ge 10R$; and the fourth holds by the fact that $\bm{d} = L^{10} N^{B\delta} \omega^{1/2} (qN)^{-1/2}$ (by \eqref{e0}) and the two lower bounds in \eqref{domega} (with the fact that $B \ge 10A$). 

Hence, the curve $\mathcal U^{k/N}_{q'}$ intersects the arctic boundary $\mA_{q'}$ within $\mathcal{Y}^-$, at a point $\mathfrak{w}_0=(\mathfrak{l}_{q'} (k/N),k/N)$ with $\mathfrak{l}_{q'} (k/N) <x_{q'}^{SW}$, such that $\mathcal U^{k/N}_{q'}(x)=k/N$ for $x\le \mathfrak{l}_{q'} (k/N) $ (recall \Cref{ul}). Note that $\frak w_0\in \mF^{SW}_{q'}$, the frozen region defined in \Cref{def:pa}.

We next define the analog of the coordinates $x^{\pm}$ from \Cref{def:Lpm}. 

\begin{definition}
	\label{def:L=-caseII}
	
	Let ${\bf L}^{\pm}$ denote the line parallel to the $y$-axis and with $x$-coordinate $x^\pm$, defined by
	\begin{flalign*}
		 x^-= \mathfrak{l}_{q'} (kN^{-1});\qquad     x^+= \bar{x} + N^{B\delta/9} \td_{\bar z}.  
	\end{flalign*} 
\end{definition}
	
	\noindent Observe by Definition \ref{yy}, with the fact that $\bar{z} = (\bar{x}, \bar{y}) \in \mathcal{Y}^-$, that $\mathbf{L}^-$ and $\mathbf{L}^+$ intersect the cell $\mathcal{Y}$. Let $\bf S$ denote the strip between $\bf L^-$ and $\bf L^+$. Given $z\in \mathcal Y$, also let ${\bf L}^z$ denote the line parallel to the $y$-axis and containing $z$.

Let us make two reductions in showing the two statements above at $\bar{z}$. 
\begin{lem}

\label{zrdelta} 

Proposition \ref{dhtqd} holds under either of the following two conditions. 

\begin{enumerate} 
	\item The point $\bar{z}$ is above the level line $\mathcal{U}_{q'}^{k N^{-1-4R\delta}}$ (namely, $\bar{y} \ge \mathcal{U}_{q'}^{k N^{-1-4R\delta}} (\bar{x})$). 
	\item We have $(\mathtt{d}_{\bar{z}} \mathtt{e}_{\bar{z};q'})^{1/2} < N^{-R\delta-1} (q-q')^{-1}$ (which holds if $(\mathtt{d}_{\bar{z}} \mathtt{e}_{\bar{z};q'})^{1/2} < N^{-2R\delta} \omega^{-1}$).
\end{enumerate} 
\end{lem} 

\begin{proof} 

To show the first statement, suppose that $\bar{z}$ is above the level line $\mathcal{U}_{q'}^{k N^{-1-4R\delta}}$. Then we would have that $\mathcal{H}_{q'} (\bar{z}) \ge kN^{-1-4R\delta} \gg N^{A\delta-1}$ by \eqref{adeltad} and that $(\mathtt{d}_{\bar{z}} \mathtt{e}_{\bar{z};q'})^{1/2} \gtrsim (\mathtt{d}_{\bar{z}} \mathcal{H}_{q'} (\bar{z}))^{1/3} \ge (\td_{\bar z} kN^{-1-4R\delta})^{1/3} = (N^{R\delta} \omega^{-3})^{1/3} \gg \omega^{-1}$, by Lemma \ref{hde0} and \eqref{eq:h02}. As such, the ceiling constraint on the Glauber dynamics $(H_s)$ would then imply that $H_t (\bar{z}) \le \mathcal{H}_{q'} (\bar{z})$, thereby showing Proposition \ref{dhtqd}. 

To show the second, observe that $(\mathtt{d}_{\bar{z}} \mathtt{e}_{\bar{z};q'})^{1/2} < N^{-R\delta-1} (q-q')^{-1}$ holds if $(\mathtt{d}_{\bar{z}} \mathtt{e}_{\bar{z};q'})^{1/2} < N^{-2R\delta} \omega^{-1}$ does, by \eqref{eq:q'}. Now suppose that $(\mathtt{d}_{\bar{z}} \mathtt{e}_{\bar{z};q'})^{1/2} < N^{-R\delta-1} (q-q')^{-1}$. If $q-q' \le N^{R\delta} \mathtt{e}_{\bar{z};q}$, then we would have that 
\begin{flalign*}
	H_t (\bar{z}) - \mathcal{H}_{q'} (\bar{z}) \le \mathcal{H}_q (\bar{z}) - \mathcal{H}_{q'} (\bar{z}) \lesssim N^{R\delta} (q-q') (\mathtt{d}_{\bar{z}} \mathtt{e}_{\bar{z};q'})^{1/2} \le N^{-1},
\end{flalign*}

\noindent by the ceiling constraint on $(H_s)$ and Lemma \ref{qlineq}. This then yields Proposition \ref{dhtqd}. 

So, assume instead that $q-q' > N^{R\delta} \mathtt{e}_{\bar{z};q}$ (and $(\mathtt{d}_{\bar{z}} \mathtt{e}_{\bar{z};q'})^{1/2} < N^{-R\delta-1} (q-q')^{-1}$), so that $q-q' > N^{R\delta} \mathtt{e}_{\bar{z};q'}$ (since $\mathtt{e}_{\bar{z};q} \ge \mathtt{e}_{\bar{z};q'}$ for $\bar{x} \ge x_q^{SW}$ and $q \ge q'$). In this case, let
  $z' = (x',y')$ be above $\bar z$ (meaning that $x'=\bar{x}$ and $y'\ge \bar{y}$) so that $(\mathtt{d}_{z'} \mathtt{e}_{z';q'})^{1/2} = \max \{  \omega^{-1}, N^{R\delta} (\mathtt{d}_{\bar{z}} \mathtt{e}_{\bar{z};q'})^{1/2} \}$. Then, $H(z') \ge H(\bar z)$ for any admissible height function $H$. We also have that $\mathtt{d}_{z'} \ge \mathtt{d}_{\bar{z}}$ and $\mathtt{e}_{z';q'} \ge \mathtt{e}_{\bar{z};q'}$. 
 Hence, $\mathtt{d}_{z'}^{1/2} \mathtt{e}_{z';q'}^{3/2} \ge N^{R\delta} \mathtt{d}_{\bar{z}}^{1/2} \mathtt{e}_{\bar{z};q'}^{3/2} \ge N^{R\delta-1} \ge N^{A\delta-1}$. As such, the ceiling constraint of $\mathcal{H}_{q'}$ applies at $z'$, that is, $H_t (z') \le \mathcal{H}_{q'} (z')$.

We also have that 
\begin{flalign*} 
	\mathcal{H}_{q'} (z') - \mathcal{H}_{q'} (\bar{z}) \lesssim \mathtt{d}_{z'}^{1/2} \mathtt{e}_{z';q'}^{3/2} \le  N^{2R\delta-3} (q-q')^{-3} \mathtt{d}_{z'}^{-1} < N^{-3} (\mathtt{d}_{\bar{z}}^{1/2} \mathtt{e}_{\bar{z};q'}^{3/2})^{-2} \lesssim N^{-1},
\end{flalign*}
\noindent where the first statement holds by Lemma \ref{hde0}; the second by the fact that $(\mathtt{d}_{z'} \mathtt{e}_{z';q'})^{1/2} = \max \{ \omega^{-1}, N^{R\delta} (\mathtt{d}_{\bar{z}} \mathtt{e}_{\bar{z};q'})^{1/2} \} \le N^{R\delta-1} (q-q')^{-1}$; the third by the fact that $q-q' > N^{R\delta} \mathtt{e}_{\bar{z};q'}$ (and that $\mathtt{d}_{\bar{z}} \le \mathtt{d}_{z'}$); and the fourth from the fact that $\mathtt{d}_{\bar{z}}^{1/2} \mathtt{e}_{\bar{z};q'}^{3/2} \ge N^{-1}$. Hence, $H_t (\bar{z}) \le H_t (z') \le \mathcal{H}_{q'} (z') \le \mathcal{H}_{q'} (\bar{z}) + N^{\delta/2-1}$, which again yields Proposition \ref{dhtqd}. 
\end{proof}

\subsubsection{Comparison of level lines} 

In view of Lemma \ref{zrdelta}, for the remainder of this proof, we will assume that $\bar{z}$ lies below $\mathcal{U}_{q'}^{k/N^{1+4R\delta}}$, and that $(\mathtt{d}_{\bar{z}} \mathtt{e}_{\bar{z};q'})^{1/2} \ge N^{-R\delta-1} (q-q')^{-1} \ge N^{-2R\delta-1} \omega^{-1}$. Next, given $x\in[x^-,x^+]$, let $s_x:=x-\mathfrak l_{q'}(k/N)$, and denote 
\begin{flalign*}
	s_0 = s_{\bar{x}}; \qquad s_+ = s_{x^+}; \qquad s_- = \max \{ s_{x^-}, 0 \}. 
\end{flalign*}

 The following is analogous to \Cref{lem:hereandthere} and further provides an estimate on $s_0$. 

\begin{lem}
	
	The assumptions of \Cref{prop:improvedconv} are verified, with $q,\varepsilon,h_0$ there equal to $q',\varepsilon_0,(k/N)$ here. Moreover, 
	\begin{flalign}
		\label{dsx} 
		\td_{\bar z} \lesssim s_0 \lesssim N^{8R\delta} \td_{\bar z}.
	\end{flalign} 

  \label{lem:hereandthere2}
\end{lem}
\begin{proof}
	
	To show \eqref{dsx}, let $\bar{z}' \in \mathcal{U}_{q'}^{k/N}$ denote the point with $x$-coordinate $\bar{x}$, so that $\mathtt{d}_{\bar{z}'} \asymp s_{\bar{x}}$. Then, by Lemma \ref{zrdelta}, $\bar{z}$ is below $\bar{z}'$, so $\mathtt{d}_{\bar{z}} \le \mathtt{d}_{\bar{z}'} \asymp s_{\bar{x}}$. This shows the first bound in \eqref{dsx}; to confirm the second observe that, if $\mathtt{d}_{\bar{z}'} \ge 2 \mathtt{d}_{\bar{z}}$, then $\mathfrak{e}_{\bar{z}';q'} \gtrsim \mathtt{d}_{\bar{z}}$ (by the convexity of the arctic boundary from Lemma \ref{convex}). Therefore, either $s_{\bar{x}} \asymp \mathtt{d}_{\bar{z}'} < 2 \td_{\bar z}$ or $\mathfrak{e}_{\bar{z}';q'} \gtrsim \mathtt{d}_{\bar{z}'}$. The second bound of \eqref{dsx} holds in the first case. In the latter, we have $\mathfrak{e}_{\bar{z}';q'} \gtrsim \mathtt{d}_{\bar{z}'} \asymp s_{\bar{x}}$, so 
	\begin{flalign*}
		 \mathtt{d}_{\bar{z}}^{1/2} s_{\bar{x}}^{3/2}\lesssim \mathtt{d}_{\bar{z}'}^{1/2} \mathfrak{e}_{\bar{z}';q'}^{3/2} \asymp  \mathcal{H}_{q'} (\bar{z}') = kN^{-1} \lesssim L^2 N^{11R\delta} \td_{\bar z}^2,
	\end{flalign*}

	\noindent by Lemma \ref{hde0} and \eqref{knd}. This again verifies the second bound in \eqref{dsx}. 

  The proof of the first statement of the lemma is similar to that of \Cref{lem:hereandthere}, and it differs only in the way the upper bound in \eqref{eq:condiz1} is checked. In the present case, note by \Cref{def:L=-caseII} and \eqref{dsx} that $s^+-s^-=x^+-x^- = 2N^{B\delta/9}\td_{\bar z} \lesssim N^{B\delta/9} s_0$. 
  	 On the other hand, $s_0 \lesssim N^{-B\delta/4}$ (by Definition \ref{yy} for the cell $\mathcal{Y}^-$ and \eqref{eq:mdsmall}), and this implies $N^{B\delta/9}s_0\ll s_0^{1/2}$ as desired.
\end{proof}

By \Cref{prop:improvedconv},
we can then find $a,b,c\in [1-C\varepsilon_0,1+C\varepsilon_0]$
and $\mathfrak r,\mathfrak r'$ such that the statements of the proposition hold.
In analogy with the case $\md=\bd$, we will assume to simplify formulas that $\frak r=\frak r'=0$; see \Cref{rem:nor}.
The following is analogous to \Cref{def:lLrR}. The index shift by $\sigma$ below is to account for the fact that $H_t (w) \le \mathcal{H}_{q'} (w) + (\mathfrak{K}-1) N^{\delta/2-1}$ when $\mathtt{d}_w \ge L^{-2} N^{1 - \mathfrak{K}\delta} \mathfrak{d}$, by the inductive hypothesis. The vertical shift by $\mathtt{d}_{\bar z}^{1/2} \mathtt{e}_{\bar z;q'}^{-1/2} N^{\delta/4 - 1}$ is a convenience to compare level lines (particularly in the third part of Lemma \ref{lem:rRlL2}); observe that such a shift only affects the height function $\mathcal{H}_{q'}$ at $\bar{z}$ by a quantity of order $N^{\delta/4-1} \ll N^{\delta/2-1}$.

\begin{definition}
	
	\label{ujrho}

	Set
	\begin{flalign*} 
		\sigma = \sigma_{\mathfrak{K}} = \lfloor (\mathfrak{K}-1) N^{\delta/2}\rfloor.
	\end{flalign*} 

	\noindent In the following, to lighten notation, we drop the floors. For all $x\in[x^-,x^+]$ and $j \in [1, k+\sigma]$, let 
	\begin{flalign*} 
		\widetilde{U}_j(x):= \mathcal U^{(j-\sigma)/N}_{a,b,c}(x) - \mathtt{d}_{\bar z}^{1/2} \mathtt{e}_{\bar z;q'}^{-1/2} N^{\delta/4 - 1},
	\end{flalign*} 
	
	\noindent where $\mathcal{U}_{a,b,c}^r = -\infty$ if $r < 0$. Further let $y_q^{\pm,j}=\mathcal U_q^{j/N}(x^\pm)$ and $y^{\pm,j}_{a,b,c}=\widetilde{U}_j(x^\pm)$.
\end{definition}

Now, in addition to the inductive procedure explained in the beginning of Section \ref{Proofd}, a second difference between the cases $\mathfrak{d} = \bm{d}$ and $\mathfrak{d} < \bm{d}$ is as follows. When $\mathfrak{d} = \bm{d}$, we had $H_t \le \mathcal{H}_{q'}$ along the $k/N$-level line $\mathcal{U}_{q'}^{k/N}$ (see Remark \ref{rem:splits}). Now when $\mathfrak{d} < \bm{d}$, this will not necessarily be the case; it may occur $(\mathtt{d}_w \mathtt{e}_{w;q'})^{1/2} < \omega$ for some points $w \in \mathcal{U}_{q'}^{k/N}$ (as, unlike in Lemma \ref{l:inthecell}, $\mathtt{d}_w$ is now not of constant order for $w \in \mathcal{Y}^-$). At such points the ceiling constraint $H_t (w) \le \mathcal{H}_{q'} (w)$ must be replaced by the weaker one $H_t (w) \le \mathcal{H}_q (w)$. 

To address this, the following definition gives a ``transition point'' $\mathfrak{z} \in \mathcal{U}_{q'}^{k/N}$, where the ceiling constraint changes from $\mathcal{H}_q$ to $\mathcal{H}_{q'}$.

\begin{definition}
	
	\label{z00} 
	
	Let $\mathfrak{z}=(\mathfrak p,\mathfrak y)$ denote the rightmost point on $\mathcal U^{k/N}_{q'}\cap \mathcal Y$ with 
	$(\mathtt{d}_{\mathfrak{z}} \mathtt{e}_{\mathfrak{z};q'})^{1/2} \le \omega^{-1}$, if such a point exists. Otherwise, set $\mathfrak z:= \mathfrak{w}_0$ (defined above Definition \ref{def:L=-caseII}).
\end{definition}

The next lemma provides some properties of $\mathfrak{z}$. The relevance of the case $\mathtt{d}_{\bar{z}} < N^{2R\delta} \mathtt{e}_{\bar{z};q'}$ considered in its third statement will be explained above Lemma \ref{lem:rRlL2}.
\begin{lem} 
	
	\label{00z} 
	
	The following statements hold. 
	
	\begin{enumerate}
		\item If $\mathfrak{z}\ne \mathfrak{w}_0$, then $\mathfrak{p} < \bar{x}$ (that is, $\mathfrak{z}$ is to the left of $\bar{z}$). 
		\item If $\mathfrak{z} \ne \mathfrak{w}_0$ and $\mathfrak{p} \ge x^{SW}_q$ (that is, $\mathfrak{z}$ is to the right of the southwest tangency location), then $\mathtt{d}_{\mathfrak{z}} \le N^{-4R\delta} \mathtt{d}_{\bar{z}}$, 
		\item If $\mathtt{d}_{\bar{z}} < N^{2R\delta} \mathtt{e}_{\bar{z};q'}$, then $\mathtt{d}_{w} \ge N^{R\delta/4} \td_{\bar z}$ holds for all $w = (x,y) \in \mathcal{U}_{q'}^{k/N}$ with $x \in [x^-, \mathfrak{p}]$ (that is, to the left of $\mathfrak{z}$).	
	\end{enumerate} 
	
\end{lem} 

\begin{proof} 
	
	To show the first statement of the lemma, let $z' = (x',y') \in \mathcal{U}_{q'}^{k/N}$ denote any point on the $kN^{-1}$-level line with $x' \ge \bar{x}$ (namely, to the right of $\bar{z}$). Then, 
	\begin{flalign*}
		(\mathtt{d}_{z'} \mathtt{e}_{z';q'})^{1/2} \asymp (\mathtt{d}_{z'} \mathcal{H}_{q'} (z'))^{1/3} \ge (\mathtt{d}_{\bar z} kN^{-1})^{1/3} \ge  N^{5R\delta/3} \omega^{-1},
	\end{flalign*}
	
	\noindent where the first statement holds by Lemma \ref{hde0}; the second since $\mathtt{d}_{z'} \ge \mathtt{d}_{\bar z}$ (as $z'$ is to the right of $\bar{z}$, and $\bar{z}$ is below $\mathcal{U}_{q'}^{k/N}$ by Lemma \ref{zrdelta}); and the third by \eqref{eq:h02}. Hence, $(\mathtt{d}_{z'} \mathtt{e}_{z';q'})^{1/2} > \omega^{-1}$, so $z' \ne \mathfrak{z}$, meaning that $\mathfrak{z}$ must lie to the left of $\bar{z}$. 
	
	To show the second statement of the lemma, observe that if $\mathfrak{z}$ exists and satisfies $\mathfrak{p} \ge x^{SW}_q$, then by \Cref{hde0} and \eqref{eq:h02} we have 
	\begin{flalign}
		\label{dzd} 
		\mathtt{d}_{\mathfrak{z}} \asymp (\mathtt{d}_{\mathfrak{z}} \mathtt{e}_{\mathfrak{z};q'})^{3/2} \cdot \mathcal{H}_{q'} (\mathfrak{z})^{-1} = \omega^{-3} k^{-1} N = N^{-5R\delta}  \td_{\bar z}. 
	\end{flalign}
	
	To show the third statement of the lemma, suppose that $\mathtt{d}_{\bar{z}} < N^{2R\delta} \mathtt{e}_{\bar{z};q'}$. 
	Observe by Lemma \ref{hde0} that we either have $kN^{-1} = \mathcal{H}_{q'} (w) \asymp \mathtt{d}_{w}^{1/2} \mathtt{e}_{w;q'}^{3/2}$ (if $x \ge x_{q'}^{SW}$) or that $kN^{-1} = \mathcal{H}_{q'} (w) \asymp \mathtt{d}_{w}^2$ (if $x < x_{q'}^{SW}$). In either case, we have $\mathtt{d}_{w}^2 \gtrsim kN^{-1}$. So, $\mathtt{d}_{w}^2 \gtrsim kN^{-1} \ge N^{4R\delta} \mathcal{H}_{q'} (\bar{z}) \gtrsim N^{R\delta} \mathtt{d}_{\bar{z}}^2$, where the second inequality holds by Lemma \ref{zrdelta} and the third holds by Lemma \ref{hde0}, with the fact that $\mathtt{d}_{\bar{z}}^{3/2} \le N^{3R\delta} \mathtt{e}_{\bar{z};q'}^{3/2}$. This shows the lemma.
\end{proof}

The next lemma is an analog of Lemma \ref{lem:rRlL}. Observe that its third part compares $\widetilde{U}_k (x)$ to $\mathcal{U}_{q'}^{k/N} (x)$ for all $x \in [x^-, x^+]$, while its fourth compares $\widetilde{U}_k (x)$ to $\mathcal{U}_{q}^{k/N} (x)$ to the left of $\mathfrak{z}$, where the ceiling might be $\mathcal{H}_q$, when $\mathtt{d}_{\bar{z}} \ge N^{2R\delta} \mathtt{e}_{\bar{z};q'}$. If instead $\mathtt{d}_{\bar{z}} < N^{2R\delta} \mathtt{e}_{\bar{z};q'}$, we might no longer have that $\mathcal{U}_q^{k/N} (x) \ge \widetilde{U}_{k+\sigma} (x)$ for $x \in [x^-, \mathfrak{p}]$. However, in this case, the third part of Lemma \ref{00z} gives $\mathtt{d}_{\mathfrak{z}} \ge N^{R\delta/4} \td_{\bar z}\gg N^{\delta} \td_{\bar z}$. So we will instead use the inductive hypothesis, \eqref{htkq} at $\mathfrak{K}-1$, to obtain the improved ceiling of $\mathcal{H}_{q'}$ (up to a small additive error) also to the left of $\mathfrak{z}$ (see the proof of Lemma \ref{utestimate0} below).  

\begin{lem}
\label{lem:rRlL2}	
The following four statements hold, for some constant $c>0$.
\begin{enumerate} 
	\item For each $1-\sigma\le j \le k$, we have $y^{+,j}_{q'}  \ge y^{+,j+\sigma}_{a,b,c} + c q \td_{\bar z} (x^+ - \bar{x})^2$. 
	\item For each $1-\sigma \le j \le k$, we have $y^{-,j}_{q}  \ge y^{-,j+\sigma}_{a,b,c}$. 
	\item We have $\mathcal U^{k/N}_{q'}(x)  \ge \widetilde{U}_{k+\sigma} (x) $ for all $x\in[x^-,x^+]$.
	\item If $\mathtt{d}_{\bar{z}} \ge N^{2R\delta} \mathtt{e}_{\bar{z};q'}$, then $\mathcal U^{k/N}_{q}(x) \ge \widetilde{U}_{k+\sigma} (x)$ for $x\in[x^-,\mathfrak p]$ (to the left of $\mathfrak{z}$).
        \end{enumerate}
        
\end{lem}

\begin{proof}
	
	Let us first show $y^{+,j}_q \ge y^{+,j}_{a,b,c} + c q \td_{\bar z} (x^+ - \bar{x})^2$. To that end, we have that 
\begin{flalign}
	\label{uqkn1} 
\mathcal{U}_{q'}^{k/N} (x^+) - \mathcal{U}_{a,b,c}^{k/N} (x^+) \ge C^{-1} q \td_{\bar z} (x^+ - \bar{x})^2,
\end{flalign} 

\noindent where we used the facts that $\mathcal{U}_{q'}^{k/N} (x) - \mathcal{U}_{a,b,c}^{k/N} (x) \gtrsim q \mathtt{d}_{\bar{z}} (x - \bar{x})^2$ for $x \in [\bar{x}, x^+]$ (by \eqref{eq:meccia2}); that $\mathcal{U}_{q'}^{k/N} (\bar{x}) = \mathcal{U}_{a,b,c}^{k/N} (\bar{x})$ (by \eqref{eq:r'}); that $\partial_x \mathcal{U}_{q'}^{k/N} (\bar{x}) = \partial_x \mathcal{U}_{a,b,c}^{k/N} (\bar{x})$ (by \eqref{eq:meccia}); and that $s_0 \gtrsim \td_{\bar z}$ by \eqref{dsx}. Since $s_{x^+} = s_{\bar{x}} + N^{B\delta/9} \td_{\bar z}$ and $s_{\bar{x}} \lesssim N^{8R\delta} \mathtt{d}_{\bar{z}} \ll N^{B\delta/9} \mathtt{d}_{\bar{z}}$ by \eqref{dsx}, we have that $s_{x^+} \asymp x^+ - \bar{x}$. Together with \eqref{eq:primamail'}, this gives 
\begin{flalign}
	\label{uqkjn} 
\big| \big(\mathcal{U}_{q'}^{k/N} (x^+) - \mathcal{U}_{a,b,c}^{k/N} (x^+) \big) - \big(\mathcal{U}_{q'}^{j/N} (x^+) - \mathcal{U}_{a,b,c}^{j/N} (x^+)\big) \big| \le C q (x^+ - \bar{x})^{5/3} (kN^{-1})^{2/3}.
\end{flalign} 

\noindent We moreover have $kN^{-1} \le L^2 N^{11R\delta} \td_{\bar z}^2$ by \eqref{knd} and that $x^+ - \bar{x} = N^{B\delta/9} \td_{\bar z}$. Combining these estimates yields $\td_{\bar z} (x^+ - \bar{x})^{1/3} \gg N^{(B/28-22R)\delta} (k/N)^{2/3}$. As such, for $B \ge 650R$, it follows that $qs_0 (x^+ - \bar{x})^2 \gg q (x^+ - \bar{x})^{5/3} (k/N)^{2/3}$. Hence, it follows that 
\begin{flalign*} 
	y_{q'}^{+, j} = \mathcal{U}_{q'}^{j/N} (x^+) > \mathcal{U}_{a,b,c}^{j/N} (x^+) +  N^{B\delta/10} q  \td_{\bar z}^3 \ge y_{a,b,c}^{+,j} + (2C)^{-1} q  \td_{\bar z} (x^+ - \bar{x})^2,
\end{flalign*} 

\noindent which confirms the first part of the lemma (upon recalling Definition \ref{ujrho}).

We next show the second part of the lemma. From \Cref{314item1} and
  \Cref{314item1bis} of \Cref{prop:improvedconv}, the curve
  $\mathcal U^{k/N}_{q'}$ is above $\mathcal U^{k/N}_{a,b,c}$. 
  Furthermore, both $\mathcal U^{k/N}_{q'}$ and $\mathcal U^{k/N}_{q}$  intersect the vertical line $\mathbf{L}^-$ at the same point 
  $(x^-,k/N)\in\mF^{SW}_q$. Thus,
  $y^{-,k}_q = kN^{-1} \ge y^{-,k}_{a,b,c}$.
  Now, again by the fact that the left
  endpoints of $\mathcal U^{j/N}_{q}$, for all $j\le k$, are in $\mF^{SW}_q$, we have $y^{-,k}_q-y^{-,j}_q = (k-j)N^{-1}$.
  On the other hand, $y^{-,k}_{a,b,c} - y^{-,j}_{a,b,c} \ge (k-j)N^{-1}$ (as the gradient of the height function is in the triangle $\mathcal T$ in \eqref{e:cT}), which implies that $y^{-,j}_q \ge y^{-,j}_{a,b,c}$
for every $j\le k$.

To show the third part of the lemma observe that, since $\mathcal{U}_{q'}^{k/N}$ is above $\mathcal U_{a,b,c}^{k/N} $, by Definition \ref{ujrho} we have for all $x \in [x^-, x^+]$ that
\begin{flalign} 
	\label{uquq} 
\mathcal{U}_{q'}^{k/N}(x)- \widetilde{U}_{k+\sigma} (x)  \ge \mathtt{d}_{\bar z}^{1/2} \mathtt{e}_{\bar z;q'}^{-1/2} N^{\delta/4 - 1} \ge 0.
\end{flalign} 
\noindent To show the fourth part of the lemma, we may assume $\mathfrak z \ne \frak w_0$. By \eqref{uquq}, we must then verify that 
$  \mathcal{U}_{q}^{k/N}(x)\ge \mathcal{U}_{q'}^{k/N}(x)-\mathtt{d}_{\bar z}^{1/2} \mathtt{e}_{\bar z;q'}^{-1/2} N^{\delta/4-1}$. 
Let $w=(x,y)$ and $w'=(x,y')$ be points on
  $\mathcal{U}_{q}^{k/N}$ and $\mathcal{U}_{q'}^{k/N}$, respectively, with the same $x$-coordinate. We must then
show 
\begin{eqnarray}
  \label{eq:diff?}
  |y-y'| \le \mathtt{d}_{\bar z}^{1/2} \mathtt{e}_{\bar z;q'}^{-1/2} N^{\delta/4-1}.
\end{eqnarray}

Consider first the case when $x\in[ x^{SW}_q,\mathfrak z]$, if this interval is nonempty. By \Cref{qlineq} and the second part of Lemma \ref{hde0},  
  $|y-y'|\asymp\mathtt{d}_{w} (q-q')$. On the other hand, $\td_w$ is increasing along the curve
  $\mathcal{U}_{q}^{k/N}$: this follows from the fact that the gradient of
  height functions belongs to the triangle $\mathcal T$ defined in
  \eqref{e:cT}, so that the vertical coordinate of any level line is
  increasing as a function of the horizontal coordinate. This, together with the second statement of Lemma \ref{00z} and Lemma \ref{zrdelta}, yields
  	\begin{flalign*}
|y-y'| \lesssim (q-q') \mathtt{d}_{\mathfrak z}	 \lesssim   N^{-4R\delta} (q-q')  \mathtt{d}_{\bar z} \ll N^{-1} \mathtt{d}_{\bar{z}}^{1/2} \mathtt{e}_{\bar{z};q'}^{-1/2},
	\end{flalign*}

	 \noindent which shows \eqref{eq:diff?}. 
        
If instead $x\in[x^-, \min(x^{SW}_q,\mathfrak p)]$, then $|\mathcal{H}_q (w) - \mathcal{H}_{q'} (w)| \le |q-q'| (\mathfrak{d}_w \mathtt{e}_{w;q'})^{1/2}$ (by \Cref{qlineq}). Together with the fact that $\partial_y \mathcal{H}_{q'} (z) \asymp 1$ for $z$ to the left of $x^{SW}_q$ (see the fourth part of \Cref{hde0}), this yields $|y-y'|\lesssim|q-q'| (\mathtt{d}_w \mathtt{e}_{w;q'})^{1/2}\lesssim |q-q'| (\mathtt{d}_{\mathfrak{z}} \mathtt{e}_{\mathfrak{z};q'})^{1/2} = |q-q'|\omega^{-1}$
          where the last bound holds because $w$ is to the left of $\mathfrak{z}$ (and $\mathtt{d}_w$ and $\mathtt{e}_{w;q'}$ are increasing in $x$ near the left edge of the arctic boundary, due to its convexity). 
          By \eqref{eq:q'}, we obtain that $|y-y'|\lesssim N^{R\delta-1}$, so \eqref{eq:diff?} follows if  $\mathtt{d}_{\bar{z}} \ge N^{2R\delta} \mathtt{e}_{\bar{z};q'}$.   
         \end{proof}
  
  \subsubsection{Proof of \Cref{htqhtq02} if $\mathfrak{d} < \bm{d}$} 
  
  The next definition prescribes an analog of Definition \ref{ujrho}, now with respect to the discrete (random) height function $H_t$. Below, $\mathsf{T}' = \mathsf{T} - q^{-1} \delta N^{1+10B\delta} \omega$ is chosen so that the inductive hypothesis, given by the $\mathfrak{K}-1$ case of \eqref{htkq}, holds for all $t \in [\mathsf{T}',N^5]$.
  
  \begin{definition} 
  	
  	\label{ytut} 
  	
  Set $\mathsf{T}' = q^{-1} \delta N^{1+10B\delta} \omega (\mathfrak{K}-1)$. For each integer $j \in [1, k]$ and real number $t \in [\mathsf{T}', N^5]$, let $\mathtt{U}_t^{j/N} (x)$ denote the  $jN^{-1}$-level line (as in Definition \ref{ttU}, recall also \Cref{ftn:minorabuse}) of $H_{t} - \mathtt{d}_{\bar z}^{1/2} \mathtt{e}_{\bar z;q'}^{-1/2} N^{\delta/4 - 1}$, and set $ \mathtt{y}_t^{\pm, j} = \mathtt{U}_t^{j/N} (x^{\pm})$.

	\end{definition}

  The following lemma indicates that the endpoints and curves from Definition \ref{ytut} lie above the ones associated with $\mathcal{H}_{a,b,c}$ from Definition \ref{ujrho}. In the case when $\mathfrak{d} = \bm{d}$, it was a quick consequence of Lemma \ref{lem:rRlL}, with the ceiling constraint on $(H_s)$. 
  
  \begin{lem} 
  	
  \label{utestimate0} 
  
  We have w.o.p. that the following statements hold for all $t \in [\mathsf{T}', N^5]$. 
  
  \begin{enumerate} 
  	\item We have that $y_{a,b,c}^{\pm,j} \le \mathsf{y}_t^{\pm, j}$ for each $j \in [1,k]$.  
  	\item We have that $\widetilde{U}_k (x) \le \mathtt{U}_t^{k/N} (x)$ for each $x \in [x^-, x^+]$. 
  \end{enumerate} 
   
   \end{lem} 

\begin{proof} 
	
	By the ceiling constraint of $\mathcal{H}_q$ on $(H_s)_{s \ge 0}$, we have that $\mathsf{y}_t^{-,j} = \mathtt{U}_t^{j/N} (x^-) \ge y_q^{-,j}$. This, together with the second part of Lemma \ref{lem:rRlL2}, implies that $\mathsf{y}_t^{-,j} \ge y_{a,b,c}^{-,j}$. Thus, to verify the first part of the lemma, it remains to show that $\mathsf{y}_t^{+,j} \ge y_{a,b,c}^{+,j}$. 
	
	To do so, we separately consider two cases. The first is if $y_{a,b,c}^{+,j} \in \mathcal{Y}^-$ for each $j \in [1,k]$. By Definition \ref{yy}, this is equivalent to imposing that $y_{a,b,c}^{+,k} \in \mathcal{Y}^-$. Then, since $x^+ = \bar{x} +  N^{B\delta/9} \mathtt{d}_{\bar{z}} \ge x^{SW}_q + N^{B\delta/9} \td_{\bar z}$ by Definition \ref{def:L=-caseII}, we have that $\mathtt{d}_w \gtrsim N^{B\delta/9} \td_{\bar z} \gg N^{\delta} \td_{\bar z} \ge L^{-2} N^{1-(\mathfrak{K}-1)\delta} \mathfrak{d}$ for any $w = (x^+, y) \in \mathfrak{L}_{q'}$. Hence, by the inductive hypothesis \eqref{htkq} applied at a time $t \in [\mathsf{T}', N^5]$, we find that $H_t (x^+, y) \le \mathcal{H}_{q'} (x^+, y) + (\mathfrak{K}-1) N^{\delta/2-1}$ for any $(x^+, y) \in \mathfrak{L}_{q'}$. Recalling that $\sigma = (\mathfrak{K}-1) N^{\delta/2}$, this yields $\mathsf{y}_t^{+,j} \ge y_{q'}^{+,j-\sigma}$, which together with the first part of Lemma \ref{lem:rRlL2} yields $\mathsf{y}_t^{+,j} \ge y_{a,b,c}^{+,j}$. 
	
	Now suppose that $y_{a,b,c}^{+,k} \notin \mathcal{Y}^-$, in which case the proof will follow that of Lemma \ref{lem:rRlL}. First, we lower bound $\td_{\bar z}$. Since $y_{a,b,c}^{+,k} \notin \mathcal{Y}^-$, we have by Definition \ref{yy} and \eqref{e:dew0} that 
	\begin{flalign} 
	\label{x12} 
	x^+ \ge x^{SW}_q + N^{B\delta} \omega^{1/2} (qN)^{-1/2}.
	\end{flalign} 

	\noindent Since $\bar{x} - x^{SW}_1 \le s_{\bar{x}} \lesssim N^{8R\delta} \td_{\bar z}$ by \eqref{dsx} and $x^+ - \bar{x} = s_{x^+} - s_{\bar{x}} = N^{B\delta/9} \td_{\bar z} \gg N^{8R\delta} \td_{\bar z}$, we have that $x^+ - x^{SW} \asymp N^{B\delta/9} \td_{\bar z}$. Together with \eqref{x12}, this gives $s_+ \asymp x^+ - \bar{x}$ and that $\td_{\bar z} \ge N^{8B\delta/9} \omega^{1/2} (qN)^{-1/2}$. By the first part of Lemma \ref{lem:rRlL2}, there thus exists a constant $c>0$ such that
	\begin{flalign}
		\label{x13} 
		y_{q'}^{+,j} - y_{a,b,c}^{+,j} \ge y_{q'}^{+,j} - y_{a,b,c}^{+,j+\sigma} \ge cq s_+ \td_{\bar z} (x^+ - \bar{x}) \ge c s_+ N^{B\delta-1} \omega.
   \end{flalign}  

	\noindent We also have that $y_{q'}^{+,j} - y_q^{+,j} \lesssim (q-q') s_+ \lesssim N^{R\delta-1} s_+ \omega$, where the first inequality holds by Lemma \ref{qlineq} and Lemma \ref{hde0}, and the second holds by \eqref{eq:q'}. Since $B>R$, this with \eqref{x13} implies that $y_q^{+,j} \ge y_{a,b,c}^{+,j}$. 	By the ceiling of $\mathcal{H}_q$ on $(H_s)$, we then have that $H_t (x^+, y) \le \mathcal{H}_q (x^+, y)$. Thus, $\mathsf{y}_t^{+,j} \ge y_q^{+,j} \ge y_{a,b,c}^{+,j}$, confirming the first part of the lemma. 
	
	To show the second, first observe that for $x \in [\mathfrak{p}, x^+]$ we have $\mathtt{U}_t^{k/N} (x) \ge \mathcal{U}_{q'}^{k/N} (x) \ge \widetilde{U}_k (x)$, where the first inequality holds by the ceiling constraint of $\mathcal{H}_{q'}$ on $(H_s)$ for points $(x,y) \in \mathcal{U}_{q'}^{k/N}$ to the right of $\mathfrak{z}$ (that is, with $x > \mathfrak{p}$), and the second holds by the third part of Lemma \ref{lem:rRlL2} (and the ordering of the curves in $\widetilde{U}$). Moreover, if $\td_{\bar z} \ge N^{2R\delta} \mathtt{e}_{\bar{z};q'}$, we have for $x \in [x^-, \mathfrak{p}]$ that $\mathtt{U}_t^{k/N} (x) \ge \mathcal{U}_q^{k/N} (x) \ge \widetilde{U}_k (x)$, where the first statement holds by the ceiling constraint of $\mathcal{H}_q$ on $(H_s)$, and the second holds by the fourth part of Lemma \ref{lem:rRlL2}. 
	
	Hence, it remains to verify $\mathtt{U}_t^{k/N} (x) \ge \widetilde{U}_k (x)$ for $x \in [x^-, \mathfrak{p}]$, when $\td_{\bar z} < N^{2R\delta} \mathtt{e}_{\bar{z};q'}$. In this case, fix $w = (x,y) \in \mathcal{U}_{q'}^{k/N}$ with $x \in [x^-, \mathfrak{p}]$. By the third part of Lemma \ref{00z}, we have $\mathtt{d}_w \ge N^{R\delta/4} \td_{\bar z} \ge N^{\delta} \td_{\bar z} \ge L^{-2} N^{1-(\mathfrak{K}-1)\delta} \mathfrak{d}$. Hence, the inductive hypothesis \eqref{htkq} at $\mathfrak{K}-1$ (with the $\bar{z}$ there equal to $w$ here) applies and yields that we w.o.p. have $H_t (w) \le \mathcal{H}_{q'} (w) + (\mathfrak{K}-1) N^{\delta/2-1}$. Recalling that $\sigma = (\mathfrak{K}-1) N^{\delta/2}$, it follows that $\mathtt{U}_t^{k/N} (x) \ge \mathcal{U}_{q'}^{(k-\sigma)/N} (x) \ge \widetilde{U}_k (x)$, by the third part of Lemma \ref{lem:rRlL2}. This establishes the lemma.
   \end{proof} 
 
The following is analogous to \Cref{def:lineensembles}. 
\begin{definition}
  \label{def:lineensembles2}
  Let $\widehat U=\{\widehat{U}_j\}_{j<k}$ denote the Bernoulli path ensemble with endpoints $(x^\pm,y^{\pm,j}_{a,b,c})$, conditioned to the event that $\widehat{U}_{k-1}$ is weakly below $\widetilde{U}_k$ and that each $\widehat{U}_j$ for $j\le k-1$ is weakly below {$\mathcal U^{(j/N+N^{-1+(B/2)\delta})}_{q'}$}. Similarly, let $\widecheck{U} =\{\widecheck{U}_j\}_{j<k}$ denote the Bernoulli path ensemble with endpoints $(x^\pm,y^{\pm,j}_{a,b,c})$, conditioned to the event that $\widecheck{U}_{k-1}$ is weakly below $\widetilde{U}_k$. 
      \end{definition}

     The following three results are analogs of  Lemma \ref{lem:whpcoupling}, Lemma \ref{lem:fluttU}, and Lemma \ref{lem:UbarGbar} (with similar proofs).

      \begin{lem}
        \label{lem:ennesimo}
	 The Bernoulli path ensemble  $\widehat U$ can be coupled  with $\widecheck U$, so that w.o.p. we have $\widehat U=\widecheck U$.
       \end{lem}
       \begin{proof}
         [Proof of \Cref{lem:ennesimo}]
         The proof is similar to that of \Cref{lem:whpcoupling}.
           The claim follows once we prove
  \begin{eqnarray}
    \label{eq:lpa1bis}
\text{w.o.p., } \widecheck U_j \text{  is weakly below }
 \mathcal U^{j/N+N^{-1+(B/2)\delta}}_{q'}, \text{ for every } j\le k-1.
  \end{eqnarray}
  As in the proof of \Cref{lem:whpcoupling}, an application of \Cref{cor:conclines} implies that, w.o.p., $\widecheck U_j$ is below $\mathcal U^{j/N+2N^{\delta-1}}_{a,b,c}$ and it remains to prove that the latter is below $\mathcal U^{j/N+N^{(B/2)\delta-1}}_{q'}$. To this purpose, we apply \eqref{heighthabcq} with $h_0,\md_z$ replaced by $k/N,s$, where $s$ is the maximal value of $\mathfrak{d}_w$ in the part of the strip $\mathbf{S}$ below $\mathcal{U}_{q'}^{k/N}$.  We have
$ s \lesssim N^{B\delta/9} s_0 \lesssim N^{(B/9+8R)\delta} \td_{\bar z}$, where the {first inequality follows from \Cref{def:L=-caseII} defining $x^+$ and the last follows from \eqref{dsx}}. Therefore,  $C q \mathfrak{d} kN^{-1} \lesssim N^{(13R+B/9)\delta} q \omega^{-3} \lesssim N^{(6A+13R+B/9)\delta-1}\lesssim N^{(B/2)\delta-1}$ by \eqref{domega} and \eqref{eq:condizioni}. Then, \eqref{heighthabcq} implies that $\mathcal U^{j/N+2N^{\delta-1}}_{a,b,c}$ is below $\mathcal U^{j/N+N^{(B/2)\delta-1}}_{q'}$, which concludes the proof.
       \end{proof}
	
\begin{lem}
	\label{lem:fluttU2}
 W.o.p., the following holds for every $1\le j<k$:
 the Bernoulli path  $\widehat U_j$ is between $\widetilde{U}_{j-N^{\delta/4}}$ and $\widetilde{U}_{j+N^{\delta/4}}$ if $j-N^{\delta/4}\ge0$, or  below $\widetilde{U}_{j+N^{\delta/4}}$ if $j-N^{\delta/4}<0$.
\end{lem}
\begin{proof}[Proof of \Cref{lem:fluttU2}]
	
	As in the proof of Lemma \ref{lem:fluttU}, it suffices by \Cref{lem:ennesimo} to prove the statement for $\widecheck U$; this follows from \Cref{cor:conclines}.
\end{proof}

\begin{lem}
	\label{lem:attheint2}
For every $1-\sigma\le j<k-\sigma$	we have $\widetilde{U}_{j+\sigma} (\bar x) \ge \mathcal U^{j/N}_{q'}(\bar x) - \mathtt{d}_{\bar{z}}^{1/2} \mathtt{e}_{\bar{z};q'}^{-1/2} N^{\delta/4-1}$.
	
\end{lem}

\begin{proof}[Proof of \Cref{lem:attheint2}]
	As in the proof of \Cref{lem:UbarGbar}, we have $\mathcal U^{j/N}_{a,b,c} (\bar x)\ge  \mathcal U^{j/N}_{q'}(\bar x)$ by \eqref{eq:r'} and \eqref{e:primamail}. This, together with Definition \ref{ujrho} for $\widetilde{U}_j$, implies the lemma. 
\end{proof}

Finally we can prove \Cref{dhtqd}.
We need a few preliminary definitions.

  We let $\mathcal{Y}' = \mathbf{S} \cap \mathcal{Y}.$
  Let $G^+$ denote the ceiling for the dynamics (as imposed by \Cref{htqhtq02}), namely, $G^+$ equals  the maximal admissible height function that is at most equal to $\mathcal{H}_q$ on $\mathcal{Y}'$ and at most equal to $\mathcal{H}_{q'}$ on $w \in \mathcal{Y}'$ for which $(\mathtt{d}_w \mathtt{e}_{w;q'})^{1/2} \ge \omega^{-1}$ and $\mathtt{d}_w^{1/2} \mathtt{e}_{w;q'}^{3/2} \ge N^{A\delta-1}$.
Let $H^+$ denote  maximal admissible height function that is at most equal to $G^+$ on all of $\mathcal{Y}'$; satisfies $H^+ (x^{\pm}, y_{a,b,c}^{\pm,j}) \le j/N$ for $j\in[1,k]$; and satisfies $H^+ (z) \le k/N$ for $z=(x,y)$ with $x \in [x^-,x^+]$ and $y = \widetilde{U}_k (x)$.
Finally, let $H^- = \mathcal{H}_{q'} - N^{B\delta/2-1}$. Note that the additional constraints (with respect to $G^+$) are those imposed by \Cref{utestimate0}.

    We have:
    \begin{lem}
      For every $z\in\mathcal Y'$ it holds
\begin{eqnarray}
  \label{eq:G+G-}
  H^+(z)-H^-(z)\le N^{B\delta-1}.
\end{eqnarray}
\label{lem:G+G-}
\end{lem}

\begin{proof}[Proof of \Cref{lem:G+G-}]
Because of the definition of $H^-$, it is enough to show $H^+ (z) \le \mathcal{H}_{q'} (z) + 2N^{6R\delta-1}$. 
As explained in the beginning of the proof of Lemma \ref{zrdelta}, this holds if $z$ is above the level line $\mathcal{U}_{q'}^{k/N^{1+4R\delta}}$. So, assume that $z$ is below this level line. We may also assume that
$\md_z^{1/2}\me_{z;q'}^{3/2}\ge N^{A\delta-1}$, as otherwise we could replace $z$ by a point $z'$ vertically above it so that $\md_{z'}^{1/2}\me_{z';q'}^{3/2}= N^{A\delta-1}$; if there we have $H^+ (z') \le \mathcal{H}_{q'} (z') + N^{6R\delta-1}$, then it would follow that $H^+ (z) \le H^+ (z') \le \mathcal{H}_{q'} (z') + N^{6R\delta-1} \le C N^{A\delta-1} + N^{6R\delta-1} \le 2N^{6R\delta-1}$ by \eqref{eq:condizioni} (here, $C$ is a positive constant and the third inequality uses \Cref{hde0}). Therefore, we assume $\md_z^{1/2}\me_{z;q'}^{3/2}\ge N^{A\delta-1}$.  If $(\mathtt{d}_{z} \mathtt{e}_{z;q'})^{1/2} \ge \omega^{-1}$ then $H^+ (z) \le  \mathcal{H}_{q'} (z)$ (by definition of $H^+$). Therefore,  assume in what follows that $(\mathtt{d}_{z} \mathtt{e}_{z;q'})^{1/2} \le \omega^{-1} < N^{R\delta-1} (q-q')^{-1}$ (by \eqref{eq:q'}). If $N^{R\delta} (q-q') \le \mathtt{e}_{z;q}$, then we would have that \begin{flalign*}
	H^+ (z) - \mathcal{H}_{q'} (z) \le \mathcal{H}_q (z) - \mathcal{H}_{q'} (z) \lesssim N^{R\delta} (q-q') (\mathtt{d}_{z} \mathtt{e}_{z;q'})^{1/2} \le N^{2R\delta-1},
\end{flalign*}

\noindent by the fact that $H^+\le \cH_q$  and Lemma \ref{qlineq}, which establishes the above estimate.

So, assume instead that $N^{R\delta} (q-q') > \mathtt{e}_{z;q}$ (and $(\mathtt{d}_{z} \mathtt{e}_{z;q'})^{1/2} < N^{R\delta-1} (q-q')^{-1}$). 
As in the proof of Lemma \ref{zrdelta}, let
  $z'$ be above $z$ so that $(\mathtt{d}_{z'} \mathtt{e}_{z';q'})^{1/2} = \max \{ \omega^{-1}, N^{R\delta} (\mathtt{d}_{z} \mathtt{e}_{z;q'})^{1/2}  \}$. Then, $H(z') \ge H(z)$ for any admissible height function $H$. We also have that $\mathtt{d}_{z'} \ge \mathtt{d}_{{z}}$ and $\mathtt{e}_{z';q'} \ge \mathtt{e}_{{z};q'}$. 
 Hence, $\mathtt{d}_{z'}^{1/2} \mathtt{e}_{z';q'}^{3/2} \ge N^{R\delta} \mathtt{d}_{{z}}^{1/2} \mathtt{e}_{{z};q'}^{3/2} \ge N^{R\delta-1} \ge N^{A\delta-1}$. As such, the ceiling constraint of $\mathcal{H}_{q'}$ applies at $z'$, that is, $H (z') \le \mathcal{H}_{q'} (z')$. 
We also have that 
\begin{flalign*} 
	\mathcal{H}_{q'} (z') - \mathcal{H}_{q'} (z) \lesssim \mathtt{d}_{z'}^{1/2} \mathtt{e}_{z';q'}^{3/2} \le N^{3R\delta-3} (q-q')^{-3} \mathtt{d}_{z'}^{-1} < N^{6R\delta-3} (\mathtt{d}_{z}^{1/2} \mathtt{e}_{z;q'}^{3/2})^{-2} \lesssim N^{6R\delta-1},
\end{flalign*}

\noindent where the first statement holds by Lemma \ref{hde0}; the second by the fact that $(\mathtt{d}_{z'} \mathtt{e}_{z';q'})^{1/2} = \max \{ \omega^{-1}, N^{R\delta} (\mathtt{d}_{z} \mathtt{e}_{z;q'})^{1/2}  \} \le N^{R\delta-1} (q-q')^{-1}$; the third by the fact that $N^{R\delta} (q-q') > \mathtt{e}_{z;q'}$ (and that $\mathtt{d}_{z} \le \mathtt{d}_{z'}$); and the fourth from the fact that $\mathtt{d}_{z}^{1/2} \mathtt{e}_{z;q'}^{3/2} \ge N^{-1}$. Hence, $H (z) \le H (z') \le \mathcal{H}_{q'} (z') \le \mathcal{H}_{q'} (z) + N^{6R\delta-1}$, which again (recalling \eqref{eq:condizioni}) yields the above bound.
\end{proof}

\begin{proof}
  [Proof of \Cref{dhtqd}] This proof is quite close to that of
    \Cref{htqhtq02} in the case $\md=\bd$. 
 To show \Cref{dhtqd}, we may impose that the height function is above the floor $H^-$ by monotonicity. Moreover, by \Cref{utestimate0}, w.o.p. we may also impose the ceiling of $H^+$ on the time interval $[\mathsf{T}',N^5]$.
 By \Cref{lem:G+G-}, \Cref{prop:tmixconstrained} and the definition \eqref{otherwise} of the cell  $\mathcal Y$, the mixing time $T_{\rm mix}^{\mathcal Y'}$ for the Glauber dynamics on $\mathcal{Y}'$ with this ceiling and floor is at most of order $L^4 N^{2+5B\delta} \bm{d}^2(\log N)^2 \le \delta q^{-1} N^{8B\delta+1} \omega$ for $N$ large enough,
 where in the last bound we used \eqref{e0} in the case $\bd>\md$.
 Note that  $\mathsf{T} - \mathsf{T}'\gtrsim N^{2B\delta}T_{\rm mix}^{\mathcal Y'}$. Therefore, by \eqref{e:submult}, with overwhelming probability, the dynamics at times in $[\mathsf{T}, N^5]$ can be coupled to
 the stationary process constrained  between $H^-$ and $H^+$.
 Thanks to \Cref{rem:forall}, the claim of \Cref{dhtqd} follows if we prove that the following holds w.o.p. under the uniform measure $\pi$ constrained  between $H^-$ and $H^+$:
 \begin{eqnarray}
   \label{eq:pi}
   H(\bar z)\le \cH_{q'}(\bar z)+\frak K N^{\delta/2-1}
 \end{eqnarray}
(recall that the point $\bar z$ has been fixed and  satisfies the conditions of \Cref{dhtqd}).
Since the event in \eqref{eq:pi} is decreasing, and $H^+$ was defined as the minimum of $G^+$ with the additional (decreasing) constraints in \Cref{utestimate0},  by monotonicity we can replace $H^+$ with the maximal height function subject just to the upper constraints in \Cref{utestimate0}. That is, we can replace the level lines (with indices $j\le k-1$) of the height function under $\pi$ by the Bernoulli path ensemble on $[x^-, x^+]$, with left endpoints   $(x^-,y_{a,b,c}^{-,j})$, right endpoints at  $(x^+,y_{a,b,c}^{+,j})$, conditioned to the event that the path with index $k-1$ is weakly  below $\widetilde{U}_k$, and the path with index $j$ is weakly below $\mathcal{U}_{q'}^{j/N+N^{B\delta/2-1}}$. Note that this is just the ensemble $\widehat{U}$ from \Cref{def:lineensembles2}.
 By \Cref{lem:ennesimo}, we may further replace $\widehat{U}$ with $\widecheck{U}$, and recall  by \Cref{lem:fluttU2} that $\widecheck U_j$ is above $\widetilde{U}_{j-N^{\delta/4}}$. By \Cref{lem:attheint2}, the latter is above $\mathcal{U}_{q'}^{(j-\sigma)/N-N^{\delta/4-1}} - \mathtt{d}_{\bar{z}}^{1/2} \mathtt{e}_{\bar{z};q'}^{-1/2} N^{\delta/4-1}$. 
(Recall that $\sigma=(\frak K-1)N^{\delta/2}$.)
 Wrapping up, we have that w.o.p.
$H (\bar{z}) \le \mathcal{H}_{q'} (\bar{x}, \bar{y} + \mathtt{d}_{\bar{z}}^{1/2} \mathtt{e}_{\bar{z};q'}^{-1/2} N^{\delta/4-1}) + \sigma/N + N^{\delta/4-1}
\le \mathcal{H}_{q'} (\bar{x}, \bar{y}) + \sigma/N + N^{\delta/2-1} 
= \mathcal{H}_{q'} (\bar{z}) + \mathfrak{K} N^{\delta/2-1},$
(where the second inequality follows from the second item in \Cref{hde0}, since we are under the assumption that $\bar{x} \ge x^{SW}_q\ge x^{SW}_{q'}$), thereby confirming \eqref{eq:pi}.
  \end{proof}

\section{Proof of \Cref{htq3}}

\label{ProofEstimate1} 

In this section we show \Cref{htq3}. At a first glance, this result may appear qualitatively different from \Cref{htqhtq02}, since it states an improved bound (reducing the errors of $N^{R\delta}$ to $N^{\delta}$) without an improved assumption. What we will use is that the distance between the $k$-th and $(k+1)$-th level lines of $\mathcal{H}_q$ decreases as $k$ increases (by, for example, the second part of \Cref{cor:conclines}). So, we will choose our $k$ here to be slightly larger (by a factor of about $N^{8R\delta}$; see \eqref{eq:h0new}) than the value \eqref{e:h0} of $k$ in the proof of Proposition \ref{htqhtq02}. In this way, the error of $N^{R\delta-1}$ in the height function around the $k$-th curve will correspond to a fairly small error in the vertical (spatial) direction, of order less than $\mathfrak{d}^{2/3} N^{\delta/2-2/3}$.

This error is negligible from the perspective of a point $\bar{z} \in \mathfrak{X}$, with $(\mathtt{d}_{\bar{z}} \mathtt{e}_{\bar{z};q})^{1/2} \le N^{R\delta-1} \omega^{-1}$, at which we seek to decrease $(H_t)$ (see \eqref{eq:itrem}). So, we will be able to vertically shift level lines under consideration by $\mathfrak{d}^{2/3} N^{\delta/2-2/3}$ (see \eqref{eq:widetildeU} and Definition \ref{u20}) while minimally affecting the height function around $\bar{z}$, in order to account for the height difference of $N^{R\delta-1}$ along the $k$-th level line (see Lemmas \ref{lem:isabove} and \ref{lem:alsoright}).

Once this shift is taken into account, the proof of \Cref{htq3} will be quite similar to that of \Cref{htqhtq02}, and so some arguments in this section will be only outlined. In fact, several aspects in the proof of \Cref{htqhtq02} will simplify when showing \Cref{htq3}, since the number $k$ of curves under consideration here will be fairly small. Indeed, we will have $k \ll N^{C\delta}$ for some constant $C>1$, so the height function $(H_t)$ of the Glauber dynamics below the $k$-th curve is at most $N^{C\delta-1}$. This directly enables us to take the essentially minimal value of $H_{\max} \sim N^{C\delta-1}$ in \Cref{prop:tmixconstrained}, without having to implement an auxiliary ceiling constraint (as in the $\widehat{U}$ Bernoulli path ensemble from Definitions \ref{def:lineensembles} and \ref{def:lineensembles2}) and then showing that this constraint is irrelevant (as in Lemmas \ref{lem:whpcoupling} and \ref{lem:ennesimo}). 

Throughout this section, we recall from Remark \ref{rem:wlog} that we may assume that the tangency point closest to $w_0$ is $p^{SW}_q$, and that the $x$-coordinate of $w_0$ is at least $x^{SW}_q$ (that is, $w_0$ is to the right of $p_q^{SW}$). We further fix a point $\bar z=(\bar x,\bar y)\in\mathcal Y^-
$ and let ${\bf L}^{\bar z}$ be
the vertical line through $\bar z$. We then define $k$ via
\begin{eqnarray}
  \label{eq:h0new}
  \frac kN=\frac{N^{30R \delta}}{\omega^3 \md}.
\end{eqnarray}

\noindent Observe by \eqref{eq:q'2} that $k \le N^{33R\delta}$. 

\subsection{Case I: \texorpdfstring{$\bd=\md$}{}}
\label{sec:77caseI}
If $\bm{d} = \mathfrak{d}$,
we adopt \Cref{def:Lpm} for $x^\pm,\bf L^\pm$ and  $\bf S$. 
Call
\begin{equation}
  \label{eq:jz}
j_{\bar z}=\lfloor N\cH_q(\bar z)\rfloor+1.
\end{equation}

\begin{rem}
 In the proof of  \Cref{66.1} of \Cref{htq3}, we can assume that  $\bar z\in\mathcal Y^-\cap\mathfrak L_q$ and 
 $(\mathtt d_{\bar z}\mathtt e_{\bar z;q})^{1/2}\le N^{R\delta}\omega^{-1}$.
 The assumption  $\bar z\in\mL_q$ can be made since the height
is non-decreasing in the vertical direction, and $\cH_q(z)$ in the
r.h.s. of \eqref{eq:statmesoe3} is zero in $\mathfrak{F}_q^S$. 
In this case, note that (since $\md_{\bar z}\asymp \md$; see \Cref{l:inthecell})
\begin{eqnarray}
  \label{eq:12or30}
  {j_{\bar z}}\asymp N\frac{(\me_{\bar z;q}\md_{\bar z})^{3/2}}\md\le N\frac{(\mathtt e_{\bar z;q}\mathtt d_{\bar z})^{3/2}}\md\lesssim
  N^{3R\delta}\frac N{\omega^3\md}\le N^{-9R\delta}k. 
\end{eqnarray}
 To prove  \Cref{66.2} of \Cref{htq3}, we will have $\bar z$ outside $\mathfrak L_q$: in this case,
\[k=\frac{N^{1+30R\delta}}{\omega^3\md}\ge N^{30R\delta}\ge N^{30R\delta}j_{\bar z} 
\]
using $j_{\bar z}=1$ and the upper bound on $\omega$ in \eqref{domega}.
In either case, we have then
\begin{eqnarray}
  \label{eq:kj}
  k\gtrsim N^{8R\delta}j_{\bar z} .
\end{eqnarray}
\label{rem:kj}
\end{rem}

\subsubsection{Comparison of level lines}

The following lemma is analogous to \Cref{prop:Ucontained}. Its proof is also very similar to that of \Cref{prop:Ucontained} (essentially upon replacing $L^3$ there with $N^{30R\delta}$ here), so we omit it.

\begin{lem}
	\label{prop:Ucontained2}
	The intersection $\mathcal U_{q}^{k/N}\cap \bf S$ is contained in $\mathcal Y$.
\end{lem}

As in the proof of \Cref{htqhtq02}, we may apply \Cref{prop:improvedconv} (whose applicability can be proven entirely analogously to as in  \Cref{lem:hereandthere}), which implies the existence of  $a,b,c,\mathfrak r,\mathfrak r'$  (here, as in the proof of \Cref{htqhtq02}, we drop $\mathfrak r,\mathfrak r'$, to lighten notation; see \Cref{rem:nor}) such that \eqref{eq:tangent} and \eqref{eq:tangent2} hold, and in particular $\mathcal U^{k/N}_q(x)\ge \mathcal U^{k/N}_{a,b,c}(x),x\in[x^-,x^+].$ For 
$x\in[x^-,x^+]$ define 
\begin{eqnarray}
  \label{eq:widetildeU}
   \widetilde U_j(x)=\mathcal U^{j/N}_{a,b,c}(x)-
\md^{2/3}N^{-2/3+\delta/2}j_{\bar z}^{-1/3}.
\end{eqnarray}

Observe that $\widetilde U_{k}$ is below $\mathcal U^{k/N}_q$, since  $\mathcal U^{k/N}_{a,b,c}$ is. As mentioned in the beginning of Section \ref{ProofEstimate1}, the vertical shift in the definition of $\widetilde{U}$ takes into account the height error of $N^{R\delta-1}$ in the ceiling constraint on $(H_t)$. Specifically, we have the following, which is the analog of the second statement of Lemma \ref{lem:rRlL}.

\begin{lem}
  \label{lem:isabove}
We have $\mathcal U_q^{k/N-N^{R\delta-1}}(x)\ge \widetilde U_k(x)$ for all $x\in[x^-,x^+]$.
\end{lem}
\begin{proof}
  By \Cref{hde0} we have
  \[
\mathcal U_q^{j/N}(x)-\mathcal U_q^{(j-1)/N}(x)\asymp
    \md^{2/3}j^{-1/3}N^{-2/3},
\]uniformly in $x\in[x^-,x^+]$ and in $j$.  Since $k \gg N^{R\delta}$ by \eqref{eq:kj}, the level line $\mathcal U^{k/N-N^{R\delta-1}}_q$ is at a distance at most of order
    \begin{equation*}
 \sum_{j=k-N^{R\delta}}^k \md^{2/3}j^{-1/3}N^{-2/3}\asymp N^{R\delta-2/3}\md^{2/3}k^{-1/3}
    \end{equation*}
    below $\mathcal U^{k/N}_{a,b,c}$. On the other hand, again by \eqref{eq:kj},
    \begin{equation*}
      N^{R\delta-2/3}\md^{2/3}k^{-1/3}\le N^{-R\delta/3-2/3}\md^{2/3}j_{\bar z}^{-1/3}\ll \md^{2/3}N^{-2/3+\delta/2}j_{\bar z}^{-1/3}.
    \end{equation*}
 Therefore, the claim follows.
\end{proof}
The following is analogous to \Cref{def:lLrR} and \Cref{def:lineensembles} (except, as mentioned in the beginning of Section \ref{ProofEstimate1}, we no longer require a Bernoulli path ensemble $\widecheck U$ with its paths bounded above).

\begin{definition}
  \label{def:lineensemblesnew}
  For $j\le k-1$, let $y^{\pm,j}_q=\mathcal U^{j/N}_q(x^\pm)$ and $y^{\pm,j}_{a,b,c}=\widetilde U_j(x^\pm)$.
  Let $\widehat U=\{\widehat U_j\}_{1\le j\le k-1}$ denote the Bernoulli path ensemble with  endpoints given by\footnote{As in \Cref{def:lineensembles}, the endpoints are points in $\mathbb T_N+(0,1/(2N))$ at minimal distance from   $\{(x^\pm,y^{j,\pm}_{a,b,c})\}_{1\le j\le k-1}$} $\{(x^\pm,y^{j,\pm}_{a,b,c})\}_{1\le j\le k-1}$, conditioned to the event that $\widehat U_{k-1}$ is weakly below $\widetilde U_k$.
\end{definition}
The following is analogous to the first statement of \Cref{lem:rRlL}.
\begin{lem}
  \label{lem:withmargin}
  For $N^{R\delta}<j\le k$, it holds  $y^{\pm,j-N^{R\delta}}_q-y^{\pm,j}_{a,b,c}\gtrsim  N^{2B\delta-1}{\omega\md}$.
\end{lem}

\begin{proof}
  The proof is very similar to that of \Cref{lem:rRlL}, so we just sketch it.
  For $j=k$ the claim {readily follows from \eqref{eq:tangent2} and $(x^+-x^-)^2\gtrsim L^3 N^{2B\delta}\omega(qN)^{-1}$}.  To prove the claim for $j<k$, it is enough to show that 
   \begin{eqnarray}\label{eq:tdb}
  |(y^{\pm,k-N^{R\delta}}_{a,b,c}-y^{\pm,j-N^{R\delta}}_{a,b,c})-(y^{\pm,k}_{q}-y^{\pm,j}_{q})|\ll   N^{2B\delta-1}{\omega\md}.
   \end{eqnarray}
   
   \noindent Arguing like in \eqref{eq:imaa},  
   we have thanks to \eqref{eq:primamail'} that
\begin{equation*}
  |(y^{\pm,k-N^{R\delta}}_{a,b,c}-y^{\pm,j-N^{R\delta}}_{a,b,c})-(y^{\pm,k}_{q}-y^{\pm,j}_{q})|\lesssim q \md\left(\md \frac
    kN\right)^{2/3}= \md \omega \frac
  q{\omega^3}N^{20R\delta}\le {\omega\md} N^{-1+(20R+6A)\delta}
\end{equation*}
using \eqref{eq:h0new} and the lower bound \eqref{domega} on $\omega$. Because of \eqref{eq:condizioni2}, \eqref{eq:tdb} follows.
\end{proof}

\subsubsection{Proof of \Cref{htq3} if $\mathfrak{d} = \bm{d}$} 

The following lemma is analogous to  \Cref{lem:fluttU}.

\begin{lem}
  \label{lem:fluttUancora}
  W.o.p., the following holds:  $\widehat U_j$ is between $\widetilde U_{j- N^{3\delta/4}}$ and $\widetilde U_{j+ N^{3\delta/4}}$ if $j-N^{3\delta/4}>0$ and between  $\widetilde U_{j+N^{3\delta/4}}$ and $\mathcal U_{a,b,c}^{1/N}-\md^{2/3}N^{-2/3+3\delta/4}$ if $j$ is smaller.
\end{lem}
\begin{proof}
   The proof (based on \Cref{cor:conclines} and on the second claim of \Cref{volumeheight} to address the case when $j \le N^{3\delta/4}$) is almost identical to that of \Cref{lem:fluttU} and \Cref{lem:fluttU2}, so we skip it.
\end{proof}

\begin{proof}[Proof of \Cref{htq3} if $\md=\bd$]
   Since the events in both claims of the proposition are decreasing, we can by monotonicity impose that above $\mathcal U^{k/N}_q$,  the height function is frozen to the maximal configuration lower than the ceiling $\cH_q+N^{R\delta-1}$. Since in the region  below $\mathcal U^{k/N}_q$ the height function is at least zero and at most $k/N\le  N^{30R\delta} /(\omega^3 \mathfrak{d}) \le N^{33R\delta-1}$  (where in the last inequality we used \eqref{eq:q'2}), \Cref{prop:tmixconstrained} implies that the mixing time satisfies $T_{\rm mix}^{\mathcal Y}\lesssim N^{1+5B\delta}\frac\omega q$ (the argument is very similar to that in \eqref{eq:ctimes}).
  Since $q^{-1}N^{1+10B\delta}\ge N^{\delta}T_{\rm mix}^{\mathcal Y}$, applying \eqref{e:submult} and \Cref{rem:forall}, it is sufficient to prove that the two claims of the proposition hold, w.o.p., under the uniform measure $\pi$ subject to the constraints that the height is  lower than $\cH_q+N^{R\delta-1}$, lower than $\cH_q$ on $\mathcal Y\setminus \mathfrak L_q^+(R\delta)$ {and frozen in the region above $\mathcal U^{k/N}_q$ to the maximal configuration smaller than $\cH_q+N^{R\delta-1}$.}

  Call $W_j$ the level line of the height function $H$ in the cell
  $\mathcal Y$.
We claim that the Bernoulli path ensemble $\{W_j\}_{1\le j\le k-1}$ can be coupled with $\widehat U$ in such a way that $W_j$ is weakly above $\widehat U_j$ for all $j$, almost surely. To see this, note that, by the ceiling constraint, $W_{j+N^{R\delta}}$ is
  deterministically above $\mathcal U_q^{j/N}$; in particular,
$W_{k}$ is above $\mathcal U_q^{k/N-N^{-1+R\delta}}$ that is above
$\widetilde U_k$ by \Cref{lem:isabove}. Denoting by $(x^\pm,Y^{\pm,j})$  the endpoints of the lines $W_j,j\le k$, we have for some positive constant $\kappa$ that $Y^{\pm,j}\ge y^{\pm,j-N^{R\delta}}_q\ge y^{\pm,j}_{a,b,c}+\kappa N^{2B\delta-1}{\omega\md}$
 for
$j> N^{R\delta}$, by \Cref{lem:withmargin}. 
On the other hand,  again by the ceiling constraint, $W_j$ is inside $\mathfrak L^+_q(R\delta)$
for $j\le N^{R\delta}$, that is, $W_j$ is above
$\mathcal U^{1/N}_q-\md^{2/3}N^{-2/3+R\delta}$, and therefore
\begin{eqnarray}
Y^{\pm,j}\ge y^{\pm,1}_{q} -\md^{2/3}N^{-2/3+R\delta}
  \ge y^{\pm,N^{R\delta}}_{a,b,c}+\kappa N^{2B\delta-1}\omega\md-\md^{2/3}N^{-2/3+R\delta}\ge  y^{\pm,N^{R\delta}}_{a,b,c}\ge y^{\pm,j}_{a,b,c}
\end{eqnarray}
where we first applied \Cref{lem:withmargin}, then  the
bound \eqref{eq:q'2} on $\omega$ and \eqref{eq:condizioni2}. Therefore, the endpoints of $W_j$ are higher than those of $\widehat U_j$ for all $j\le k-1$.
The claimed domination of the Bernoulli ensemble $\widehat U$ by $W$ then follows from \Cref{le:stochdomB}.

To prove the first claim of the proposition it is then enough to show that,
w.o.p., $\widehat U_{j_{\bar z}+N^\delta}(\bar x)\ge \bar y$. By \Cref{lem:fluttUancora}, $\widehat U_{j_{\bar z}+N^\delta}$
  is  w.o.p. above
  $\widetilde U_{j_{\bar z}+N^\delta-N^{3\delta/4}}
  $, therefore above
  $\widetilde U_{j_{\bar z}+N^{\delta}/2}=\mathcal U_{a,b,c}^{j_{\bar z}/N+N^{\delta-1}/2}-\md^{2/3}N^{-2/3+\delta/2}j_{\bar z}^{-1/3}$. 
On the other hand, by
\Cref{lem:UbarGbar},  $ \mathcal U^{j/N}_{a,b,c}(\bar x)\ge \mathcal U_q^{j/N}(\bar x)$ for $j\le k$.  
Wrapping up, we have $\widehat U_{j_{\bar z}+N^\delta}(\bar x)\ge \mathcal U_{q}^{j_{\bar z}/N+N^{\delta-1}/2}(\bar x)-\md^{2/3}N^{-2/3+\delta/2}j_{\bar z}^{-1/3}$
and it remains to show that
\begin{eqnarray}
  \label{eq:itrem}
  \mathcal U_{q}^{j_{\bar z}/N+N^{\delta-1}/2}(\bar x)-\md^{2/3}N^{-2/3+\delta/2}j_{\bar z}^{-1/3}\ge \bar y= \mathcal U_{q}^{j_{\bar z}/N}(\bar x).
\end{eqnarray}
 On the other hand, the third statement of \Cref{hde0} (see also \eqref{eq:lpa2}) implies
 \begin{multline}
     \mathcal U_{q}^{j_{\bar z}/N+N^{\delta-1}/2}(\bar x)-\mathcal U_{q}^{j_{\bar z}/N}(\bar x)\gtrsim
     \sum_{j=j_{\bar z}}^{j_{\bar z}+N^{\delta}/2}\frac{\md^{2/3}N^{-2/3}}{j^{1/3}}\asymp
     \md^{2/3}N^{-2/3}(N^{2\delta/3}{\bf 1}_{j_{\bar z}\le N^\delta}+\frac {N^\delta}{j_{\bar z}^{1/3}}{\bf 1}_{j_{\bar z}\ge N^\delta})\\\gtrsim
\frac{\md^{2/3}N^{-2/3}}{j_{\bar z}^{1/3}}N^{(2/3)\delta}
\gg 
\frac{\md^{2/3}N^{-2/3+\delta/2}}{j_{\bar z}^{1/3}}.
\end{multline}
\Cref{66.1} is proven.

Next, we prove \Cref{66.2} of \Cref{htq3} and we take
$\bar z\in\mathcal Y^- \setminus\mathfrak L_q^+(\delta)$ (the
condition $\td_{\bar z}\in(L^{-2}\md,L^2\md)$ is automatically
satisfied in the case $\md=\bd$ if $L$ is large, thanks to
\Cref{l:inthecell} and to \eqref{condizd}, that guarantees that
$\td_{\bar z}= \md_{\bar z}$). The goal is to prove that, w.o.p. under the measure $\pi$,
$W_{1}$ is above $\bar z$ (that is,  $H(\bar z)
=\cH_q(\bar z)=0$). We can assume that
$\bar z\in\mathfrak L_q^+(R\delta)$ since outside of
$\mathfrak L_q^+(R\delta)$ the height $H$ is zero, because of
the ceiling constraint $H\le \cH_q=0$.
Since $W$ dominates $\widehat U$, it suffices to show that $\widehat U_1(\bar x)\ge y$, w.o.p.; by \Cref{lem:fluttUancora} and \Cref{lem:UbarGbar}, $\widehat U_1(\bar x)\ge \mathcal U^{1/N}_{a,b,c}(\bar x)-\md^{2/3}N^{-2/3+3\delta/4}\ge \mathcal U_q^{1/N}(\bar x)-\md^{2/3}N^{-2/3+3\delta/4}$. The
claim then follows since  $\bar z\not\in\mL^+_q(\delta)$ implies that $\bar y\le 
\mathcal U_q^{1/N}(\bar x)-\md^{2/3}N^{-2/3+4\delta/5}$.
\end{proof}

\subsection{Case II: \texorpdfstring{$ \bm{d}>\mathfrak{d}$}{}}

\label{Proofh000}

In this section we prove \Cref{htq3} when $\bm{d}> \mathfrak{d} $. The
proof will be quite similar to that of \Cref{htq3} when
$\mathfrak{d} = \bm{d}$, and to that of \Cref{htqhtq02} when $\bm{d}> \mathfrak{d} $. In fact, it will admit some simplifications
with respect to the latter; specifically, since the number of curves
$k \le N^{33R\delta}$ (as noted after \eqref{eq:h0new}) is reasonably small, we will be able to
circumvent the inductive procedure implemented there, by instead
shifting the curve $\widetilde{U}$ slightly to the right (see
Definition \ref{u20} below), in addition to down (as in
\eqref{eq:widetildeU}). We will also not require the point
$\mathfrak{z}$ from \Cref{z00}, as we will not need to make use of a
``transition'' between the two different ceiling constraints.

\subsubsection{Properties of cells and the $k$-th level line} 

Recall from the beginning of Section \ref{ProofEstimate1} that we have fixed $\bar{z} = (\bar{x}, \bar{y}) \in \mathcal{Y}^-$. By symmetry (see also \Cref{rem:wlog}), we may assume that $\bar{x} \ge x_q^{SW}$ (that is, $\bar{z}$ is to the right of $p^{SW}_q$). In analogy with Definition \ref{def:L=-caseII},  we define $x^\pm$ as follows.

\begin{definition}
	\label{def:L=-caseIII}
	
	Recalling the definition of $k$ from \eqref{eq:h0new}, set 
	\begin{flalign*}
		x^{\pm} =\bar{x} \pm N^{B\delta/9} \td_{\bar{z}}
	\end{flalign*} 
	
	\noindent and let ${\bf L}^{\pm}$ denote the line parallel to the $y$-axis and with $x$-coordinate $x^\pm$.
	
\end{definition}
\begin{rem}
	Since $\bar{z} \in \mathcal{Y}^-$, it is quickly seen from \eqref{otherwise} and the fact that $\bm{d} \ge L^{-2} \td_{\bar{z}}$,  together with the definition \eqref{otherwise2}  of the cell $\mathcal Y$, that the vertical lines $\mathbf{L}^\pm$  intersect the cell $\mathcal{Y}.$ Also observe by \eqref{eq:h0new}, \eqref{domega} and \eqref{eq:condizioni2} (and that $\td_{\bar z} \in (L^{-2} \mathfrak{d}, L^2 \mathfrak{d})$) that 
	\begin{flalign}
		\label{kndnantra} 
		\displaystyle\frac{k}{N} \lesssim N^{(30R+6A)\delta} \md^2\ll N^{36R\delta}\md^2.
	\end{flalign}
\end{rem}

\begin{rem}\label{rem:wcat}
	We can assume that
	\begin{eqnarray}
		\label{eq:wecan}
		(\mathtt d_{\bar z}\mathtt e_{\bar z;q})^{1/2}\le N^{R\delta}\omega^{-1}.  
	\end{eqnarray}
	Indeed, for \Cref{66.1} this is an assumption of the proposition, and for \Cref{66.2}, $(\mathtt d_{\bar z}\mathtt e_{\bar z;q})^{1/2}\asymp \md^{1/3}N^{-1/3}$ (because $\td_{\bar z}=\md_{\bar z}\asymp \md$, the latter thanks to the assumption $\td_{\bar z}\in(L^{-2}\md,L^2\md)$, and $\te_{\bar z;q}=\td_{\bar z}^{-1/3}N^{-2/3}\asymp \md^{-1/3}N^{-2/3}$, the first equality in the latter formula being true by \Cref{lz} because $\bar z$ is outside the liquid region), while $N^{R\delta}\omega^{-1}\ge \md^{1/3}N^{(R+A)\delta-1/3}$ by \eqref{domega}. 
	Then, as in \Cref{rem:kj}, it follows that
	\begin{eqnarray}
		\label{eq:jzg}
		k\ge N^{8R \delta}j_{\bar z}
	\end{eqnarray}
with $j_{\bar z}$ defined in \eqref{eq:jz}. \end{rem}
The following is the analog of  \Cref{prop:Ucontained2}:
\begin{lem}
	\label{prop:Ucontained3}
	The graph of the curve $x\in[x^-,x^+]\mapsto \mathcal U_q^{k/N}(x)$ is contained in  $\mathcal Y$.
\end{lem}
\begin{proof} 
	Letting $z=(x,y) \in \mathcal{Y}$ denote any point on the top boundary of $\mathcal{Y}$ with $x\in[x^-,x^+]$, it suffices to show that $\mathcal{H}_{q} (z) \ge kN^{-1}$.
	Since the top boundary of $\mathcal Y$ is horizontal and has vertical coordinate  $L^{5}N^{2B\delta}\bd^2$ (in the case $\md<\bd$) and the height is decreasing in the horizontal direction, it is enough to show the inequality for $x=x^+,y=L^{5}N^{2B\delta}\bd^2$.
	Observe from \Cref{def:L=-caseIII} that $x^+ - x_q^{SW} \lesssim N^{B\delta/2} \mathfrak{d}_{\bar{z}}$ and $y \gg (x^+ - x_q^{SW})^2$ (as $\mathfrak{d}_{\bar{z}} \le L^2 \mathfrak{d} $). Hence, since there exist constants $c>0$ and $C>1$ for which $\mathfrak{e}_{z;q} \ge c y^{1/2} - C(x^+ - x_q^{SW})$ (because the arctic boundary has positive and finite curvature, by \Cref{convex}), it follows that $\mathfrak{e}_{z;q} \gtrsim y^{1/2} = \mathfrak{d}_{z}$. Thus, by \Cref{hde0}, 
	$\mathcal{H}_q (z) \gtrsim \mathfrak{d}_{z}^2 = y= L^{5}N^{2B\delta}\bd^2$.
	Since $k/N=   N^{30R\delta} (\omega^3 \mathfrak{d})^{-1}\le N^{36R\delta}\md^2\le N^{36R\delta}\bd^2$ (by \eqref{domega} and the assumption $\md<\bd$), the claim $\mathcal{H}_q (z)\ge k/N$ follows from \eqref{eq:condizioni2}.
\end{proof}

\begin{lem}
	If $z_x=(x,\mathcal U_q^{k/N}(x))$  and $x\in
	[x^-,x^+]$, then $\mathfrak d_{z_x}\le N^{B\delta/2}\md$.
	Moreover, $\mathcal U_q^{k/N}$ intersects the arctic boundary $\mA_q$ at a point $\mathfrak w_0=(x_0,k/N)$ with $x_0\in[x^-,x^+]$, and $\mathcal U_q^{k/N}(x)=k/N$ for $x\le x_0$.
\end{lem}
\begin{proof} The implicit constants in the bounds $\lesssim,\asymp$ below can depend on $L$.
	Note that {$x\mapsto
U^{k/N}(x)$}  is increasing, so it is enough to prove the first statement of the lemma for
	$x=x^+$.  Observe that $\mathcal{U}_q^0 (x_+) \lesssim (x^+ -
	x_q^{SW})^2 \lesssim N^{B\delta}
	\mathfrak{d}^2$ by the fact that $x^+ - \bar{x} \lesssim N^{B\delta/2}
	\mathfrak{d}$, that $\bar{x} - x_q^{SW} \lesssim
	\mathfrak{d}$ {and that the curvature of the arctic boundary is bounded}. Hence, for some constant $C>1$, we have
	\begin{flalign*}
		\mathfrak{d}_{z_{x^+}}^2 = \mathcal{U}_q^{k/N} (x_+) \le \mathcal{U}_q^{k/N} (x^+) - \mathcal{U}_q^0 (x^+) + C N^{B\delta} \mathfrak{d}^2
		\lesssim \mathfrak{d}_{z_{x^+}}^{2/3} (kN^{-1})^{2/3} + N^{B\delta} \mathfrak{d}^2
		\\\lesssim \mathfrak{d}_{z_{x^+}}^{2/3} N^{24R\delta-2/3} + N^{B\delta} \mathfrak{d}^2,      
	\end{flalign*}
	where the first statement follows from the definition of $z_{x_+}$; the second from the above bound on $\mathcal{U}_q^0 (x^+)$; the third from \Cref{hde0}; and the fourth from the definition of $k$, with \eqref{eq:q'2}. This implies that either $\mathfrak{d}_{z_{x^+}} \lesssim N^{B\delta/2} \mathfrak{d}$ or that $\mathfrak{d}_{z_{x^+}} \lesssim N^{18R\delta-1/2} \lesssim N^{18R\delta} \mathfrak{d}$ by \eqref{condizd}. So, by \eqref{eq:condizioni2} we have either way that $\mathfrak{d}_{z_{x^+}} \le N^{B\delta/2} \mathfrak{d}$.
	
	Since the arctic boundary is uniformly convex, the point $\frak w_0=(x_0,k/N)$ satisfies $x^{SW}_q-x_0\asymp \sqrt{k/N}\lesssim N^{18 R\delta}
	\md$, by \eqref{kndnantra}. Since $\bar x-x^{SW}_q\lesssim \md$ and $\bar x-x^-\asymp N^{B\delta/9}\md\gg N^{18R\delta}
	\md$ by \eqref{eq:condizioni2}, it follows that $x_0\in[x^-,x^+]$.
	Finally, the fact that $\mathcal U^{k/N}_q(x)=k/N$ for $x\in[x^-,x_0]$ follows from \Cref{ul}.
\end{proof}

\subsubsection{Comparison of level lines}

Now let $s_x = x - \mathfrak{l}_q (k/N)$, and denote
\begin{flalign*}
	& \breve{x} = \bar{x} - \mathfrak{d}^{-1/3} j_{\bar{z}}^{-1/3} N^{-2/3-\delta}; \qquad s_0 = s_{\breve{x}} = s_{\bar{x}} - \mathfrak{d}^{-1/3} j_{\bar{z}}^{-1/3} N^{-2/3-\delta}; \\
	& s_+ = s_{x^+}; \qquad s_- = 0. 
\end{flalign*} 

\noindent By \Cref{prop:improvedconv} (as in \Cref{lem:hereandthere2}, one checks that the assumptions of the proposition are satisfied with the above {choice of $s_0,s_\pm$}),
we can find $a,b,c\in [1-C\varepsilon_0,1+C\varepsilon_0]$ (once again, set $\frak r,\frak r'$ from those formulas to $0$, for lightness of notation) such that
\begin{flalign}
	\label{eq:tangent3}
	\begin{aligned} 
		& \mathcal U^{k/N}_{q}(\breve{x})=\mathcal U^{k/N}_{a,b,c}(\breve{x}), \quad  \partial_x\mathcal U^{k/N}_{q}(\breve{x})=\partial_x\mathcal U^{k/N}_{a,b,c}(\breve{x}) ; \\
		&
\mathcal U_{q}^{k/N}(x)-\mathcal U_{q}^{k/N}(\breve x)\ge \mathcal U_{a,b,c}^{k/N}(x)-\mathcal U_{a,b,c}^{k/N}(\breve x) + C^{-1} q\md(x-\breve x)^2 \quad \forall x\in[\breve{x},x^+]; \\
		& \mathcal U_q^{k/N}(x) \ge \mathcal U_{a,b,c}^{k/N}(x) \quad \forall x\in[x^-,x^+].
	\end{aligned}  
\end{flalign}

\begin{definition} 
	
	\label{u20} 
	
	Define $\widetilde{U}_j$ by shifting $\mathcal{U}_{a,b,c}^{j/N}$ (vertically) down by $\mathfrak{d}^{2/3} j_{\bar{z}}^{-1/3} N^{\delta/2-2/3}$ and (horizontally) to the right by $\mathfrak{d}^{-1/3} j_{\bar{z}}^{-1/3} N^{-2/3-\delta}$. Specifically, for $j\le k$ and $x\in[x^-,x^+]$, set 
	\begin{flalign*} 
		\widetilde U_j(x) = \mathcal{U}_{a,b,c}^{j/N} (x - \mathfrak{d}^{-1/3} j_{\bar{z}}^{-1/3} N^{-2/3-\delta}) - \mathfrak{d}^{2/3} j_{\bar{z}}^{-1/3} N^{\delta/2-2/3}.
	\end{flalign*} 
	
	\noindent Furthermore, for each $t \in [0, N^5]$, let $\mathtt{U}_t^{j/N} (x)$ denote the $jN^{-1}$-level line (as in Definition \ref{ttU}, recalling also \Cref{ftn:minorabuse}) of $H_t$, and set $\mathtt{y}_t^{\pm,j} = \mathtt{U}_t^{j/N} (x^{\pm})$.
\end{definition} 

The following two lemmas provide the analog of \Cref{utestimate0} (which will also replace \Cref{lem:isabove} and \Cref{lem:withmargin}) in our setting. The vertical shift of $\widetilde{U}_j$ in \Cref{u20} will facilitate in ensuring in Lemma \ref{lem:alsoright} that $\mathtt{U}_t^{k/N} (x) \ge \widetilde{U}_k (x)$ for $x \le \breve{x}$; the horizontal shift will facilitate in ensuring that $\mathtt{U}_t^{1/N} (x^+) \ge \widetilde{U}_1 (x^+)$ in Lemma \ref{uxux}. 

\begin{lem}
	\label{lem:alsoright}
	
	We have w.o.p. that $\mathtt{U}_t^{k/N} (x)\ge \widetilde U_k(x)$, for all $x\in[x^-,x^+]$ and $t \in [0, N^5]$.
\end{lem}
\begin{proof}
	
  By Definition \ref{u20} and the last bound in \eqref{eq:tangent3}, it suffices to show w.o.p. that $\mathtt{U}_t^{k/N} (x) \ge \mathcal{U}_q^{k/N} (x - \mathfrak{d}^{-1/3} j_{\bar{z}}^{-1/3} N^{-2/3-\delta}) - \mathfrak{d}^{2/3} j_{\bar{z}}^{-1/3} N^{\delta/2-2/3}$ for all $x \in [x^-, x^+]$ and $t \in [0, N^5]$. We may assume below that
  \begin{eqnarray}
    \label{eq:nominimal}
  (x - \mathfrak{d}^{-1/3} j_{\bar{z}}^{-1/3} N^{-2/3-\delta}, \mathcal{U}_q^{k/N} (x - \mathfrak{d}^{-1/3} j_{\bar{z}}^{-1/3} N^{-2/3-\delta})) \notin \mathfrak{F}_q^{SW},  
  \end{eqnarray}
   for otherwise $\mathcal{U}_q^{k/N} (x - \mathfrak{d}^{-1/3} j_{\bar{z}}^{-1/3} N^{-2/3-\delta})$ {takes the  minimal value $k/N$} and therefore $\mathtt{U}_t^{k/N} (x) \ge \mathcal{U}_q^{k/N} (x - \mathfrak{d}^{-1/3} j_{\bar{z}}^{-1/3} N^{-2/3-\delta})$ holds.
	
	Next, by the ceiling constraint on $H_t$, we have $\mathtt{U}_t^{k/N} (x) \ge \mathcal U_q^{k/N-N^{R\delta-1}}(x)$, so it remains to show 
	\begin{flalign} 
		\label{uquq2} 
		\mathcal{U}_q^{k/N - N^{R\delta-1}} (x) \ge \mathcal{U}_q^{k/N} (x - \mathfrak{d}^{-1/3} j_{\bar{z}}^{-1/3} N^{-2/3-\delta}) - \mathfrak{d}^{2/3} j_{\bar{z}}^{-1/3} N^{\delta/2-2/3}.
	\end{flalign} 
	
	To that end, first observe for any $h \ge 0$ that by \Cref{ulevelnew} we have $ \partial_x \mathcal{U}_q^h (x) = - \partial_x \mathcal{H}_q (z) \cdot \partial_y \mathcal{H}_q (z)^{-1}$, where $z = (x, \mathcal{U}_q^h (x))$. Thus, 
	\begin{flalign}
		\label{derivativexuqh} 
		\partial_x \mathcal{U}_q^h (x) \asymp \mathfrak{d}_z, \quad \text{if $x \ge x_q^{SW}$}; \qquad \partial_x \mathcal{U}_q^h (x) \asymp (\mathfrak{d}_z \mathfrak{e}_{z;q})^{1/2}, \quad \text{if $x \le x_q^{SW}$}.
	\end{flalign}
	
	\noindent where the first statement follows from \Cref{hde0}, and the second follows from the facts (see the fourth part of \Cref{hde0}) that $-\partial_x \mathcal{H}_q (z) \asymp (\mathfrak{d}_z \mathfrak{e}_{z;q})^{1/2}$ and $\partial_y \mathcal{H}_q (z) \asymp 1$ for $x \le x_q^{SW}$. We additionally have for any $x \in [x_-, x_+]$ that, denoting $z = (x, \mathcal{U}_q^h (x))$,
	\begin{flalign}
		\label{dz} 
		\mathfrak{d}_z \asymp |x - x_q^{SW}| + h^{1/2},
	\end{flalign} 
	
	\noindent where this holds for $x \le x_q^{SW}$ since then $h = \mathcal{H}_q (z)  \asymp \mathfrak{d}_z^2$ (and $|x - x_q^{SW}| = x_q^{SW} - x \lesssim h^{1/2}$ in this case) and it holds for $x \ge x_q^{SW}$ by \eqref{derivativexuqh} (with the fact that $\mathcal{U}_q^h (x) = \mathfrak{d}_z^2$). 
	
	We next use \eqref{derivativexuqh} and \eqref{dz} to show \eqref{uquq2}; denote $z_0 = (x, \mathcal{U}_q^{k/N} (x)) \in \mathcal{U}_q^{k/N}$. First assume that $x \ge x_q^{SW}$. Then, for a constant $c>0$, we have 
	\begin{flalign}
		\label{uqknx2} 
		\mathcal{U}_q^{k/N} (x - \mathfrak{d}^{-1/3} j_{\bar{z}}^{-1/3} N^{-2/3-\delta}) \le \mathcal{U}_q^{k/N} (x) -c  \mathfrak{d}_{z_0} \mathfrak{d}^{-1/3} j_{\bar{z}}^{-1/3} N^{-2/3-\delta}.
	\end{flalign}
	
	\noindent since $\partial_x \mathcal{U}_q^{k/N} (x') \gtrsim \mathfrak{d}_{z_0}$ for $x' \in [x-\mathfrak{d}^{-1/3} j_{\bar{z}} N^{-2/3-\delta}, x]$. To verify the latter bound, observe it holds at $x' = x$, by \eqref{derivativexuqh}; for other $x'$, we have denoting $z' = (x', \mathcal{U}_q^{k/N} (x'))$ that $\partial_x \mathcal{U}_q^{k/N} (x') \gtrsim \mathfrak{d}_{z'}$, by \eqref{derivativexuqh}. So, we must verify that $\mathfrak{d}_{z'} \gtrsim \mathfrak{d}_{z_0}$. If this were false then, since $x-x' \le \mathfrak{d}^{-1/3} j_{\bar{z}}^{-1/3} N^{-2/3-\delta} \ll N^{-\delta} \mathfrak{d}_{z_0}$ (as $\mathfrak{d}_{z_0} \ge N^{-1/2}),$ $\mathcal{U}_q^{k/N}$ would have decreased too quickly at some point in  $[x-\mathfrak{d}^{-1/3}j_{\bar{z}}^{-1/3} N^{-2/3-\delta}, x]$. Specifically, we would have that $\partial_x \mathcal{U}_{k/N} (x'') \gg \mathfrak{d}_{z_0}$ for some $x'' \in  [x-\mathfrak{d}^{-1/3} j_{\bar{z}}^{-1/3}N^{-2/3-\delta}, x]$, which contradicts \eqref{derivativexuqh}. Thus $\partial_x \mathcal{U}_q^{k/N} (x') \gtrsim \mathfrak{d}_{z_0}$, and hence \eqref{uqknx2}, holds.  Moreover, for any  $j\in[k-N^{R\delta},k]$, we have by \Cref{hde0}, \eqref{eq:jzg} and \eqref{eq:lpa2} that
	\begin{eqnarray}
		\label{uqjnx} 
		0\le  \mathcal U_q^{j/N}(x)-\mathcal U_q^{(j-1)/N}(x)\lesssim \md_{z_0}^{2/3}j^{-1/3} N^{-2/3}.
	\end{eqnarray}
	
	\noindent Therefore, since $j \ge k/2$,
	\begin{eqnarray}
		\label{uqknx1} 
		\mathcal U_q^{k/N}(x) -  \mathcal U_q^{k/N-N^{R\delta-1}}(x) \lesssim \md_{z_0}^{2/3} k^{-1/3} N^{R\delta-2/3}.
	\end{eqnarray}
	
	\noindent By \eqref{uqknx2} and \eqref{uqknx1}, to confirm \eqref{uquq2} when $x \ge x_q^{SW}$, it suffices to verify that 
	\begin{flalign*}
		\md_{z_0}^{2/3} k^{-1/3} N^{R\delta-2/3} \ll \mathfrak{d}_{z_0} \mathfrak{d}^{-1/3} j_{\bar{z}}^{-1/3} N^{-2/3-\delta} +  \mathfrak{d}^{2/3} j_{\bar{z}}^{-1/3} N^{\delta/2-2/3},
	\end{flalign*}
	
	\noindent which holds since $k \ge N^{4R\delta} j_{\bar{z}}$ by \eqref{eq:kj} (together with the fact that $a+b \ge (a^2 b/4)^{1/3}$ for any real numbers $a, b \ge 0$).
	
	Now assume instead that $x \in [x^-, x_q^{SW}]$ is to the left of the southwest tangency location. Note that \eqref{derivativexuqh} and \eqref{dz} (together with $\md_z\gtrsim \md_z^{1/2}\me_{z;q}^{1/2}$) give for a constant $c>0$ that 
	\begin{flalign*}
		\mathcal{U}_q^{k/N} (x - \mathfrak{d}^{-1/3} j_{\bar{z}}^{-1/3} N^{-2/3-\delta}) \le \mathcal{U}_q^{k/N} (x) - c (kN^{-1})^{1/4} \mathfrak{e}_{z_0;q}^{1/2} \mathfrak{d}^{-1/3} j_{\bar{z}}^{-1/3} N^{-2/3-\delta}.
	\end{flalign*}
        
	\noindent Moreover, since $\partial_y \mathcal{H}_q (x,\mathcal{U}_q^{j/N} (x)) \asymp 1$ for $x \le x_q^{SW}$, we have 
	\begin{flalign*}
		\mathcal{U}_q^{k/N} (x) - \mathcal{U}_q^{k/N-N^{R\delta-1}} (x) \asymp N^{R\delta-1}.
	\end{flalign*}
	
	\noindent Hence, to confirm \eqref{uquq2} when $x \le x_q^{SW}$, it suffices to verify that 
	\begin{flalign}
		\label{0kn0} 
		N^{R\delta-1} \ll (kN^{-1})^{1/4} \mathfrak{e}_{z_0;q}^{1/2} \mathfrak{d}^{-1/3} j_{\bar{z}}^{-1/3} N^{-2/3-\delta} +  \mathfrak{d}^{2/3} j_{\bar{z}}^{-1/3} N^{\delta/2-2/3}.
	\end{flalign}
	
	\noindent Due to the second term on the right side of \eqref{0kn0}, this bound holds if $\mathfrak{d} \ge j_{\bar{z}}^{1/2} N^{3R\delta/2-1/2}$. So, assume instead that $\mathfrak{d} < j_{\bar{z}}^{1/2} N^{3R\delta/2-1/2}$. In this case, we use the fact that $\mathfrak{e}_{z_0;q} \ge  \mathfrak{d}^{-1/3} N^{-2/3-\delta} \ge j_{\bar{z}}^{-1/6} N^{-(R/2+1)\delta-1/2}$, where the first bound follows from \eqref{eq:nominimal}, and so \eqref{0kn0} becomes equivalent to
	\begin{flalign*} 
		N^{R\delta-1} \ll k^{1/4}  \cdot j_{\bar{z}}^{-7/12} N^{-1-3R\delta/4-3\delta/2}.
	\end{flalign*} 
	 Therefore, it suffices to verify that $k^{1/4} j_{\bar{z}}^{-7/12} \gg N^{2R\delta}$.
         By  \eqref{eq:12or30} and \eqref{eq:q'2}, we have that $j_{\bar{z}} \lesssim N^{3R\delta} \cdot N/(\omega^3 \md)\le N^{6R\delta}$, so that $j_{\bar{z}}^{7/12} \le N^{7R\delta/2}$. Thus, it suffices to verify that $k \gg N^{22R\delta}$, which holds by \eqref{eq:h0new} and \eqref{domega}.
This establishes the lemma.
\end{proof}

\begin{lem} 
	
	\label{uxux} 
	
	For both $x \in \{ x^-, x^+ \}$, we have w.o.p. that  $\mathtt{U}_t^{j/N} (x) \ge \widetilde{U}_j (x)$, for all $j \in [1,k-1]$ and $t \in [0,N^5]$. 
	
\end{lem} 

\begin{proof} 
	
	Observe for all $j \in [1,k]$ that $\mathtt{U}_t^{j/N} (x^-) \ge \widetilde{U}_j (x^-)$, since $\mathtt{U}_t^{k/N} (x^-) \ge \widetilde{U}_k (x^-)$ (by Lemma \ref{lem:alsoright}) and $\widetilde{U}_k (x^-) - \widetilde{U}_j (x^-) = (k-j)N^{-1} \le \mathtt{U}_t^{k/N} (x^-) - \mathtt{U}_t^{j/N} (x^-)$ (the first since $(x^-, \widetilde{U}_k (x^-)) \in \mathfrak{F}_{a,b,c}^{SW}$, and the second since the $y-$derivative of the height is at most $1$). It therefore remains to show that $\mathtt{U}_t^{j/N} (x^+) \ge \widetilde{U}_j (x^+)$ holds w.o.p. for all $j \in [1,k]$.
	
	To that end, set $z_0 = (x^+, \mathcal{U}_q^{k/N} (x^+)) \in \mathcal{U}_q^{k/N}$ and let $j \in [1,k]$. Then, w.o.p., we have for some constant $C>1$ that 
	\begin{flalign*}
		\mathtt{U}_t^{j/N} (x^+) - \widetilde{U}_j (x^+) & \ge \max \{ \mathcal{U}_q^{j/N-N^{R\delta-1}} (x^+), \mathcal{U}_q^{1/N} (x^+) - \mathfrak{d}_{z_0}^{2/3} N^{R\delta-2/3} \} - \widetilde{U}_k (x^+) \\ 
		& \ge   \mathcal{U}_q^{k/N} (x^+) - C k^{2/3} \mathfrak{d}_{z_0}^{2/3} N^{R\delta-2/3}  - \widetilde{U}_k (x^+)  \\
		& \ge    \mathcal{U}_q^{k/N} (x^+) - \mathcal{U}_q^{k/N} (x^+ - \mathfrak{d}^{-1/3} j_{\bar{z}}^{-1/3} N^{-2/3-\delta}) - C k^{2/3} \mathfrak{d}_{z_0}^{2/3} N^{R\delta-2/3} ,
	\end{flalign*} 

	\noindent where the first statement follows from the ceiling constraints \eqref{eq:cond12} and \eqref{eq:cond32} on $(H_t)$; the second follows from summing \eqref{uqjnx} over $j \in [1,k]$; and the third follows from Definition \ref{u20} with the third statement of \eqref{eq:tangent3}. Thus, by \eqref{uqknx2}, we have w.o.p. for some constant $c>0$ that
	\begin{flalign}
		\label{utjn} 
		\mathtt{U}_t^{j/N} (x^+) - \widetilde{U}_j (x^+) & \ge c \mathfrak{d}_{z_0} \mathfrak{d}^{-1/3} j_{\bar{z}}^{-1/3} N^{-2/3-\delta} - C k^{2/3} \mathfrak{d}_{z_0}^{2/3} N^{R\delta-2/3}.
	\end{flalign} 
	
	\noindent Further observe that 
	\begin{flalign*}
		\mathfrak{d}_{z_0} \mathfrak{d}^{-1/3} j_{\bar{z}}^{-1/3} N^{-2/3-\delta} \cdot (k^{2/3} \mathfrak{d}_{z_0}^{2/3} N^{R\delta-2/3})^{-1} & =  \mathfrak{d}_{z_0}^{1/3} \mathfrak{d}^{-1/3} \cdot j_{\bar{z}}^{-1/3} k^{-2/3} N^{-(R+1)\delta} \\ 
		& \ge N^{B\delta/30} N^{-33R\delta}  \gg 1, 
	\end{flalign*}

	\noindent where the second statement follows from \eqref{dz}, \eqref{eq:kj}, the fact that $x^+ - x_q^{SW} \ge x^+ - \bar{x} = N^{B\delta/9} \mathtt{d}_{\bar{z}} \ge N^{B\delta/10} \mathfrak{d}$ (by Definition \ref{def:L=-caseIII} and the fact that $\mathtt{d}_{\bar{z}} \ge L^{-2} \mathfrak{d}$), \eqref{eq:kj}, and the fact from \eqref{eq:h0new} and \eqref{eq:q'2} that $k \le N^{33R\delta}$; the third follows from \eqref{eq:condizioni2}. This, together with \eqref{utjn} yields w.o.p. $\mathtt{U}_t^{j/N} (x^+) \ge \widetilde{U}_j (x^+)$, which shows the lemma. 	
\end{proof} 

\subsubsection{Proof of \Cref{htq3} when $\mathfrak{d} < \bm{d}$} 

The following is the analog of \Cref{def:lineensemblesnew}.
\begin{definition}
	
	\label{3u} 
	
	For $j\le k-1$, set $y^{\pm,j}_{a,b,c}=\widetilde U_j(x^\pm)$.  Let $\widehat U=\{\widehat U_j\}_{1\le j\le k-1}$ denote the Bernoulli path ensemble with  endpoints given by $\{(x^\pm,y^{j,\pm}_{a,b,c})\}_{1\le j\le k-1}$, conditioned to the event that $\widehat U_{k-1}$ is weakly below $\widetilde U_k$. 
\end{definition}

The following is the analog of Lemma \ref{lem:fluttUancora} in our setting, with a similar proof.

\begin{lem}
	\label{lem:fluttUancora2}
	W.o.p., the following holds for all $x \in [x^-, x^+]$. For $j \in (N^{3\delta/4}, k-1]$, we have $\widetilde U_{j- N^{3\delta/4}} (x) \le \widehat{U}_j (x) \le \widetilde U_{j+ N^{3\delta/4}} (x)$. Moreover, for all $j \in [1, N^{3\delta/4}]$, we have that $\mathcal U_{a,b,c}^{1/N} (x - \mathfrak{d}^{-1/3} j_{\bar{z}}^{-1/3} N^{-2/3-\delta})-\md^{2/3} j^{-1/3} N^{-2/3+3\delta/4} \le  \widehat U_j (x) \le \widetilde U_{j+N^{3\delta/4}} (x)$. 
\end{lem}
\begin{proof}
	As in the proof of Lemma \ref{lem:fluttUancora}, the first statement of the lemma follows from \Cref{cor:conclines}, and the second follows from the second claim of \Cref{volumeheight}. 
\end{proof}

\begin{proof}[Proof of \Cref{htq3} if $\md<\bd$]
	
	Given the above estimates, the proof of this proposition when $\bm{d} < \mathfrak{d}$ is similar to that in the case $\mathfrak{d} = \bm{d}$ provided in Section \ref{sec:77caseI}; so, we only outline it here. 
	
	First, by  \Cref{rem:wcat} (in particular, \eqref{eq:jzg}), we can restrict our attention to points $\bar z$ satisfying $\cH_q(\bar z)\le k/N$. Then, by monotonicity (\Cref{prop:PW}), we can freeze the height function above the line $\mathcal U_q^{k/N}$ to the maximal configuration compatible with the ceiling constraints \eqref{eq:cond12} and \eqref{eq:cond32}.
	
	As in the proof of \Cref{htq3} when $\md=\bd$, using \Cref{prop:tmixconstrained}, we can quickly verify that the mixing time $T_{\rm mix}^{\mathcal Y}$ in the cell $\mathcal{Y}$ is at most $q^{-1} N^{1+5B\delta}\omega$. Therefore, in the time interval $[q^{-1}N^{1+10B\delta}\omega,N^5]$, we can couple height $(H_t)$ to coincide w.o.p. with the uniform measure $\pi$, subject to the constraints that the height is  lower than $\cH_q+N^{R\delta-1}$, lower than $\cH_q$ on $\mathcal Y\setminus \mathfrak L_q^+(R\delta)$, and frozen in the region above $\mathcal U^{k/N}_q$ to the maximal configuration smaller than $\cH_q+N^{R\delta-1}$. Hence, by \Cref{rem:forall}, it suffices to show that both claims in \Cref{htq3} hold w.o.p. under this measure $\pi$.
	
	By Lemma \ref{lem:alsoright} and Lemma \ref{uxux}, it suffices by monotonicity (\Cref{prop:PW}) to verify that both claims in \Cref{htq3} hold w.o.p. for the Bernoulli path ensemble $\widehat{U}$ from Definition \ref{3u}. By Lemma \ref{lem:fluttUancora2} and Definition \ref{u20}, we have that
	\begin{flalign}
		\label{ujx} 
		\widehat{U}_j (\bar{x}) \ge \max \{ \mathcal{U}_{a,b,c}^{j/N-N^{3\delta/4-1}} (\breve{x}), \mathcal{U}_{a,b,c}^{1/N} (\breve{x})\} - \mathfrak{d}^{2/3} (j^{-1/3} + j_{\bar{z}}^{-1/3}) N^{3\delta/4-2/3}.
	\end{flalign} 
	
	\noindent By \eqref{e:primamail} in Proposition \ref{prop:improvedconv}, together with the first statement in \eqref{eq:tangent3}, we have  $ \mathcal{U}_{a,b,c}^{j/N-N^{3\delta/4-1}} (\breve{x}) \ge \mathcal{U}_q^{j/N-N^{3\delta/4-1}} (\breve{x})$ and $\mathcal{U}_{a,b,c}^{1/N} (\breve{x}) \ge \mathcal{U}_q^{1/N} (\breve{x})$. Setting $j = \max \{ j_{\bar{z}}-1, 1 \}$, by these bounds with \eqref{ujx} (and recalling that $\breve{x} = \bar{x} - \mathfrak{d}^{-1/3} j_{\bar{z}}^{-1/3} N^{-2/3-\delta}$), it remains to show that 
	\begin{flalign*}
		\max \{&  \mathcal{U}_q^{j/N-N^{3\delta/4-1}} (\breve{x}), \mathcal{U}_q^{1/N} (\breve{x})\} \\
		& \ge 	\max \{ \mathcal{U}_q^{j/N - N^{\delta-1}} (\bar{x}), \mathcal{U}_q^{1/N} (\bar{x}) - \mathfrak{d}^{2/3} N^{\delta - 2/3} \} + 3 \mathfrak{d}^{2/3} j_{\bar{z}}^{-1/3} N^{3\delta/4-2/3},
	\end{flalign*} 
	\noindent (where $\mathcal{U}_q^{h} (\breve{x}):=-\infty$ if $h<0$) which follows quickly from \eqref{derivativexuqh} and summing the bound $\mathcal{U}_q^{j/N} (\bar{x}) - \mathcal{U}_q^{(j-1)/N} (\bar{x}) \gtrsim \mathfrak{d}_{\bar{z}}^{2/3} j_{\bar{z}}^{-1/3} N^{-2/3}$ (as a consequence of \eqref{eq:lpa2}) over $j \in [j_{\bar{z}} - N^{\delta}, j_{\bar{z}} - N^{3\delta/4}]$. This confirms the proposition.
\end{proof}

	\appendix
	
	\section{Proof of height concentration estimate}

	\subsection{Proof of \Cref{volumeheight}} 
	
	\label{ProofHeightH} 
	
	In this section we establish \Cref{volumeheight}. Throughout this section, fix real numbers $R > 1 > \delta > 0$ and $a, b, c \in [R^{-1}, R]$, as in \Cref{volumeheight}, and further let $H$ denote the height function associated with a uniformly random tiling of ${\mathfrak{X}}_{a, b, c}$. We begin with the following lemma, originally due to
        \cite[Lemma 5.6]{AURT} (but that appears in the form stated
        below as \cite[Equation (13)]{LTDMLS}). It indicates the
        first statement of \Cref{volumeheight} holds, which provides a
        concentration bound for $H$ everywhere in
        ${\mathfrak{X}}_{a, b, c}$.

	\begin{lem}[{\cite[Equation (13)]{LTDMLS}}]
	
	\label{1estimateabc}
	
	Item \ref{heightabc0} of \Cref{volumeheight} holds. 
	 
	\end{lem} 

	It remains to verify the second statement of \Cref{volumeheight}, indicating that $H$ is w.o.p. frozen outside of the augmented liquid region $\mathfrak{L}_{a,b,c}^+ (\delta)$. To that end, set $K = \lfloor N^{\delta/3} \rfloor$ and for each integer $i \in [1, bN]$ (recall $b$ is the vertical side-length of the hexagon) define the (random) level line $U_i (x) = \mathtt{U}_{a,b,c}^{(i-1/2)/N} (x)$ of the tiling height function $H$; see \Cref{ttU}. Then, applying \Cref{1estimateabc} (with the $\delta$ in \Cref{volumeheight} equal to $\delta/4$ here), we obtain w.o.p. that $U_K (x) \ge \mathcal{U}_{a,b,c}^0 (x)$ for all $x \in [0, a]$. 
	
	To show Item (2) of \Cref{volumeheight}, we must show that $U_1 (x) \ge \mathcal{U}_{a,b,c}^0 (x) - \mathtt{d}_x^{2/3} N^{\delta-2/3}$. To that end, abbreviate $\mathfrak{l} = \mathfrak{l}_{a,b,c} (0)$ and, for any $x \ge \mathfrak{l}$, set $d_x = \max \{ x - \mathfrak{l}, N^{-1/2} \} \lesssim \mathtt{d}_x$. For each $j \in [1,K]$, define the event 
	\begin{flalign*} 
		\Omega_j = \bigcap_{x \in [\mathfrak{l},a]} \{ U_{K-j+1} (x) \ge \mathcal{U}_{a,b,c}^0 (x) - (j-1) d_x^{2/3} N^{\delta/3-2/3} \},
	\end{flalign*} 

	\noindent which provides the above type of lower bound on $U_{K-j+1}$. As indicated above, $\Omega_1$ holds w.o.p.; the following lemma with induction will indicate that $\Omega_K$ holds w.o.p. as well.

	\begin{lem} 
		
		\label{habcj2} 
		
		Let $ j \in [1,K-1]$ be an integer. If $\Omega_j$ holds w.o.p., then $\Omega_{j+1}$ holds w.o.p., as well.
	\end{lem}

	\begin{proof}[Proof of \Cref{volumeheight}, Item (2)]
		
		By rotating $\mathfrak{X}_{a,b,c}$ and permuting $(a, b, c)$ if necessary, we may assume that $z = (x,y)$ lies in the connected component of $\mathfrak{X} \setminus \mathfrak{L}_{a,b,c}^+(\delta)$  containing the south corner $(a, 0)$ of ${\mathfrak{X}}_{a, b, c}$. By symmetry, we may further assume that the closest side of $\partial \mathfrak{X}_{a,b,c}$ to $z$ is its SW one, connecting $(0, 0)$ to $(a, 0)$; hence, $x \in [\mathfrak{l}, a]$. It suffices to show that $H(z) = 0$, or equivalently that $U_1 (x) > y$. Since $y \le \mathcal{U}_{a,b,c}^0 (x) - \mathfrak{d}_x^{2/3} N^{\delta-2/3}$, it remains to show (as $\mathtt{d}_x \lesssim \mathfrak{d}_x$) that $U_1 (x) \ge \mathcal{U}_{a,b,c}^0 (x) - K \mathtt{d}_x^{2/3} N^{\delta/3-2/3}$. This is equivalent to verifying that w.o.p. $\Omega_K$ holds, which is true by the fact that $\Omega_1$ holds w.o.p., together with Lemma \ref{habcj2} and induction on $j$. 
	\end{proof}

           To show \Cref{habcj2}, we fix $j \in [1, K-1]$ and restrict to $\Omega = \Omega_j$. Then, conditional on $U_{K-j+1}$, the law of $\bm{U} = (U_1, U_2, \ldots , U_{K-j})$ is given the uniform measure on families of $K-j$ non-intersecting Bernoulli paths on $[0, a+c]$, conditioned to satisfy, for each integer $i \in [1, K-j]$,
		\begin{multline*} 
                  U_i (0) = \frac {i-(1/2)}N
                  ; \qquad U_i (a+c) = c+\frac{i-(1/2)}N\\ U_{K-j} (x) \le U_{K-j+1} (x), \quad \text{for all $x \in [0, a+c]$}.
		\end{multline*}
		Next let $V : [ \mathfrak{l}
                , a+c] \rightarrow \mathbb{R}$ denote a uniformly random Bernoulli path conditioned to start at $V( \mathfrak{l}
) = 0$, end at $V(a+c) = c$, and satisfy $V(x) < \mathcal{U}_{a,b,c}^0 (x)$ for each $x \in ( \mathfrak{l}
, a+c)$.

                Since we have restricted to the event $\Omega_j$, monotonicity (Lemma \ref{le:stochdomB}) yields a coupling between $V$ and $\bm{U}$ such that 
		\begin{flalign} 
			\label{ujvs}
			U_{K-j} (x) \ge V(x) - (j-1) d_x^{2/3} N^{\delta/3-2/3} , \qquad \text{for each $x \in [0, a+c]$}.
		\end{flalign}
		\begin{lem}
		
		\label{vs} 
		
		We have w.o.p. that $V(x) \ge \mathcal{U}_{a,b,c}^0 (x) - d_x^{2/3} N^{\delta/3-2/3}$, for all $x \in [ \mathfrak{l}, a+c]$.
		\end{lem} 
	
		\begin{proof}[Proof of \Cref{habcj2}]
			
			By \Cref{vs}, \eqref{ujvs} yields $U_{K-j} \ge \mathcal{U}_{a,b,c}^0 (x) - j d_x^{2/3} N^{\delta/3-2/3}$ holds w.o.p., meaning that w.o.p. $\Omega_{j+1}$ holds.	
		\end{proof}

	To establish \Cref{vs}, for any real numbers $x \in [\mathfrak{l}, a+c]$ and $r > 0$, define the events
	\begin{flalign}
		\label{omegasr}
		\Omega(x; r) = \big\{ V(x) \ge \mathcal{U}_{a,b,c}^0 (x) - r d_x^{2/3} N^{-2/3} \big\}; \qquad \Omega (r) = \bigcap_{x \in [\mathfrak{l}, a+c]} \Omega (x;r). 
	\end{flalign}
	
	\noindent Lemma \ref{vs} indicates that $\Omega (x; r)$ very likely holds for $r \ge N^{\delta/3}$. The following lemma, to be established in \Cref{ProofOmega0} below, indicates that if $\Omega (r)$ holds for $r \ge N^{\delta/3}$, then $\Omega ( (1-N^{-\delta}) r)$ holds w.o.p.; repeatedly applying this  statement yields \Cref{vs}.

	\begin{lem} 
		
		\label{omegasr0} 
		
			If $\mathfrak{R} \ge N^{\delta/3}$, then w.o.p. $\Omega (\mathfrak{R})^{\complement} \cup \Omega ( (1- N^{-\delta}) \mathfrak{R})$ holds.

	\end{lem}

	\begin{proof}[Proof of \Cref{vs}]
		
		Observe that $\Omega (N^2)$ holds deterministically, as since $d_x \ge N^{-1/2}$ we have $r d_x^{2/3} N^{-2/3} \gg 1$ and thus $\mathcal{U}_{a,b,c}^0 (x) - rd_x^{2/3} N^{-2/3} < 0$ for $r \gg N$. Thus, applying \Cref{omegasr0} $N^{2\delta}$ times, since $N^{\delta/3} \ge (1-N^{-\delta})^{N^{2\delta}} \cdot N^2$, we deduce that $\Omega (N^{\delta/3})$ holds w.o.p.; this implies the lemma.		
	\end{proof} 
	
	\subsection{Proof of \Cref{omegasr0}}
	
	\label{ProofOmega0}
	
	In this section we establish \Cref{omegasr0}. The proof uses ideas reminiscent of the tangent method introduced in \cite{ACSVMGD}; see also \cite[Section 5.2]{ABIM}.

	\begin{proof} 
		
		It suffices to show for any fixed $x_0 \in [\mathfrak{l},a+c]$ that w.o.p. $(\Omega(\mathfrak{R}) \cap \Omega(x_0; \mathfrak{R} - \mathfrak{R} N^{-\delta})^{\complement})^{\complement}$ holds. We will assume in what follows that $x_0 \le a$ (as the case when $x_0 \ge a$ follows similarly, by rotating the hexagon clockwise).  We also assume $\mathcal{U}_{a,b,c}^0 (x_0) \ge (1 - N^{-\delta}) \mathfrak{R} d_{x_0}^{2/3} N^{-2/3}$, as otherwise $\mathcal{U}_{a,b,c}^0 (x_0) - (1 - N^{-\delta}) \mathfrak{R} d_{x_0}^{2/3} N^{-2/3} \le 0 \le V (x_0)$ and so $\Omega (\mathfrak{R} - \mathfrak{R} N^{-\delta})$ (and thus $(\Omega(\mathfrak{R}) \cap \Omega(x_0; \mathfrak{R} - \mathfrak{R} N^{-\delta})^{\complement})^{\complement}$) holds deterministically. 
		
		Then, since the strict convexity of $\mathfrak{A}_{a,b,c}$ gives $\mathcal{U}_{a,b,c}^0 (x_0) \asymp d_{x_0}^2$, it follows that $d_{x_0}^2 \gtrsim \mathfrak{R} d_{x_0}^{2/3} N^{-2/3}$. This yields  
		\begin{flalign}
			\label{dx0r} 
			d_{x_0} \gtrsim \mathfrak{R}^{1/2} d_{x_0}^{1/3} N^{-1/3},
		\end{flalign} 
		
		\noindent which in particular implies that $d_{x_0} \ge N^{\delta/4-1/2}$ (as $\mathfrak{R} \ge N^{\delta/3}$) and hence $d_{x_0} = x_0 - \mathfrak{l}$.

		Now, let $u > 0$ be a small constant and $U > 1$ be a large constant, both to be fixed later, and set 
		\begin{flalign}
			\label{x1x2} 
			x_1 =  x_0 - u \mathfrak{R}^{1/2} d_{x_0}^{1/3} N^{-1/3}; \qquad x_2 = \min \{ x_0 + U \mathfrak{R}^{1/2} d_{x_0}^{1/3} N^{-1/3}, a+c \}.
		\end{flalign} 
		
		\noindent Observe for sufficiently small $u > 0$ that $x_1 > \mathfrak{l}$, as $x_1 - \mathfrak{l} = x_0 - \mathfrak{l} - u \mathfrak{R}^{1/2} d_{x_0}^{1/3} N^{-1/3} \gtrsim d_{x_0}$, by \eqref{dx0r}. By similar reasoning, we have $x_2 - x_0 \lesssim U d_{x_0}$.  Moreover, we have 
		\begin{flalign}
			\label{dx012}
			d_{x_1} < d_{x_0} < d_{x_2} \lesssim  U d_{x_0},
		\end{flalign}
		
		\noindent where the first two statements follow from \eqref{x1x2}, and the third from \eqref{dx0r} (with \eqref{x1x2}). In what follows, we restrict to the event $\Omega (x_1; \mathfrak{R}) \cap \Omega (x_2; \mathfrak{R})$. We must then show that w.o.p. that 
		\begin{flalign} 
			\label{vx0} 
			V(x_0) \ge \mathcal{U}_{a,b,c}^0 (x_0) - (1-N^{-\delta}) \mathfrak{R} d_{x_0}^{2/3} N^{-2/3}.
		\end{flalign} 
	
		\noindent Now, let $\ell$ denote the line tangent to $\mathfrak{A}_{a,b,c}$ at the point $(x_0, \mathcal{U}_{a,b,c}^0 (x_0))$. For any $x \in \mathbb{R}$, let $\ell (x)$ be such that $(x, \ell(x)) \in \ell$; set $\rho = \ell' (x_0)$, which is the slope of $\ell$ at $x_0$. By the strict convexity of $\mathfrak{A}_{a,b,c}$, we have 
		\begin{flalign}
			\label{rhodx0} 
			\rho \asymp d_{x_0}.
		\end{flalign}
		
		\noindent Again by the strict convexity of $\mathfrak{A}_{a,b,c}$ there exists a constant $v > 0$ such that $\ell (x_1) \le \mathcal{U}_{a,b,c}^0 (x_1) - 2 v (x_1 - x_0)^2$. By \eqref{x1x2}, it follows that 
		\begin{flalign} 
			\label{12x1} 
			\ell (x_1) \le \mathcal{U}_{a,b,c}^0 (x_1) - 2vu^2 \mathfrak{R} d_{x_0}^{2/3} N^{-2/3}.
		\end{flalign}
		
		\noindent By similar reasoning, we have
		\begin{flalign} 
			\label{12x2} 
			\begin{aligned} 
				& \ell (x_2) \le \mathcal{U}_{a,b,c}^0 (x_2) - 2vU^2 \mathfrak{R} d_{x_0}^{2/3} N^{-2/3}, \qquad \text{if $x_2 = x_0 + U \mathfrak{R}^{1/2} d_{x_0}^{1/3} N^{-1/3}$}.
			\end{aligned}
		\end{flalign}

		Next observe that, conditional on $V(x)$ for $x \notin (x_1, x_2)$, the law of $V(x)$ for $x \in [x_1, x_2]$ is that of a uniformly random Bernoulli path on $[x_1, x_2]$, with endpoints given by\footnote{Here and in the rest of the proof, endpoints of Bernoulli paths should be taken as points in $\mathbb T_N+(0,1/(2N))$ at minimal distance from the points mentioned in the text.} $(x_1,V (x_1))$ and $(x_2,V (x_2))$, conditioned to satisfy $V(x) \le \mathcal{U}_{a,b,c}^0 (x)$ for each $x \in [x_1, x_2]$. We will compare $V$ to a uniformly random Bernoulli path without any such conditioning. To set this up, let $\ell_1$ denote the vertical shift of $\ell$ down by $(1 - vu^2) \mathfrak{R} d_{x_0}^{2/3} N^{-2/3}$, that is, set 
		\begin{flalign}
			\label{10} 
			\ell_1 (x) = \ell (x) - (1 - vu^2) \mathfrak{R} d_{x_0}^{2/3} N^{-2/3},
		\end{flalign} 
		
		\noindent for each $x \in \mathbb{R}$. We then have 
		\begin{flalign*} 
			\ell_1 (x_1) \le \mathcal{U}_{a,b,c}^0 (x_1) - \mathfrak{R} d_{x_0}^{2/3} N^{-2/3} \le \mathcal{U}_{a,b,c}^0 (x_1) - \mathfrak{R} d_{x_1}^{2/3} N^{-2/3} \le V(x_1),
		\end{flalign*} 
		
		\noindent where the first statement follows from \eqref{10} and \eqref{12x1}; the second from \eqref{dx012}; and the third from our restriction to $\Omega (x_1; \mathfrak{R})$. We further have for sufficiently large $U$ that 
		\begin{flalign*}
			\ell_1 (x_2) \le V(x_2).
		\end{flalign*}
		
		\noindent Indeed, if $x_2 = a+c$, then this follows from the fact $V(a+c) = \mathcal{U}_{a,b,c}^0 (a+c) \ge \ell_1 (a+c)$ (the first by the boundary conditions of $V$, imposed by those of the tiling, and the second by the convexity of $\mathfrak{A}_{a,b,c}$). If instead $x_2 = x_0 + U \mathfrak{R}^{1/2} d_{x_0}^{1/3} N^{-1/3}$, then this holds, as for sufficiently large $U$ we have  for some constant $C>1$ that 
		\begin{flalign*}
			\ell_1 (x_2) & \le \mathcal{U}_{a,b,c}^0 (x_2) - (1 -vu^2 + 2vU^2) \mathfrak{R} d_{x_0}^{2/3} N^{-2/3} \\
			& \le \mathcal{U}_{a,b,c}^0 (x_2) - (1 + vU^2) (CU)^{-1} \mathfrak{R} d_{x_2}^{2/3} N^{-2/3}  \le \mathcal{U}_{a,b,c}^0 (x_2) - \mathfrak{R} d_{x_2}^{2/3} N^{-2/3} \le V(x_2),
		\end{flalign*} 
		
		\noindent where the first statement follows from \eqref{10} and \eqref{12x2}; the second from \eqref{dx012}; the third from the fact that $U$ is sufficiently large; and the fourth from our restriction to $\Omega (x_2; \mathfrak{R})$.
		
		Let $W_1: [x_1, x_2] \rightarrow \mathbb{R}$ denote a uniformly random Bernoulli path starting at $(x_1,W_1(x_1))$ with $W_1 (x_1) = \ell_1 (x_1)$, ending at $(x_2,W_1(x_2))$ with  $W_1 (x_2) = \ell_1 (x_2)$, and satisfying $W_1 (x) \le \mathcal{U}_{a,b,c}^0 (x)$. By monotonicity (\Cref{le:stochdomB}), we may couple $(V, W_1)$ such that $W_1 (x) \le V(x)$ for all $x \in [x_1, x_2]$. Also let $W_2: [x_1, x_2] \rightarrow \mathbb{R}$ denote a uniformly random Bernoulli walk satisfying $W_2 (x_1) = \ell_1 (x_1)$ and $W_2 (x_2) = \ell_1 (x_2)$. We would like to couple $W_1$ and $W_2$ to coincide w.o.p., to which end it suffices to show that w.o.p. $W_2$ remains below $\mathcal{U}_{a,b,c}^0$. 
		
		To that end, observe that, w.o.p., 
		\begin{flalign}
			\label{w20} 
				\sup_{x \in [x_1, x_2]} |W_2 (x) - \ell_1 (x)| \le  N^{\delta/200} \cdot \rho^{1/2} \cdot \Big( \displaystyle\frac{x_2-x_1}{N} \Big)^{1/2} \le N^{\delta/100} \mathfrak{R}^{1/2} d_{x_0}^{2/3} N^{-2/3}.
		\end{flalign}
		
		\noindent where the first statement follows from a concentration bound for uniformly random Bernoulli walks (with the fact that the slope of $\ell_1$ is $\rho$) and the second follows from \eqref{x1x2} and \eqref{rhodx0}. Since $\dist (\ell_1, \mathfrak{A}_{a,b,c}) \gtrsim \mathfrak{R} d_{x_0}^{2/3} N^{-2/3}$ and $\mathfrak{R} \ge N^{\delta/3}$, it follows that $W_2$ remains w.o.p. below $\mathcal{U}_{a,b,c}^0$. Hence, we may couple $W_1 = W_2$ w.o.p., meaning that w.o.p. $V(x_0) \ge W_2 (x_0)$. 
		
		Thus, we have w.o.p. that 
		\begin{flalign*} 
			V(x_0) \ge W_2 (x_0) & \ge \ell_1 (x_0) -  N^{\delta/100} \mathfrak{R}^{1/2} d_{x_0}^{2/3} N^{-2/3} \\
			& \ge \mathcal{U}_{a,b,c}^0 (x_0) - (1 - vu^2 + N^{\delta/100} \mathfrak{R}^{-1/2})  \cdot \mathfrak{R} d_{x_0}^{2/3} N^{-2/3} \\
			& \ge \mathcal{U}_{a,b,c}^0 (x_0) - (1 - N^{-\delta}) \mathfrak{R} d_{x_0}^{2/3} N^{-2/3}, 
		\end{flalign*} 
		where the second statement follows from \eqref{w20}; the third from the definition of $\ell_1$ (and the fact that $\ell (x_0) = \mathcal{U}_{a,b,c}^0 (x_0)$); the fourth from the fact that $\mathfrak{R} \ge N^{\delta/3}$ and that $N$ is sufficiently large. This establishes \eqref{vx0} and thus the lemma.
	\end{proof}

	\section{Properties of limit shapes} 
	
	\label{Shape}
	
	Recall from \Cref{sec:Vtilted} that $\mathcal{H}_{q;a,b,c}$ denotes the limit shape for the volume-tilted random tiling on $\mathfrak{X}_{a,b,c}$. In this section we prove the properties stated in \Cref{sec:Vtilted}. Throughout, we fix $q \in \mathbb{R}$ and $a,b,c>0$. IWe also set $\mathfrak{q} = e^q$ and let $\mathbb{H} = \{ z \in \mathbb{C}: \Imaginary z > 0 \}$ denote the upper half plane.

	\subsection{Formula for the limit shape}
        
	\label{CurvatureDensity}
	
	In this section we recall from \cite{DBPP} (though we will use the conventions from \cite{DLE}) a formula for $\mathcal{H}_{q;a,b,c}$. To that end, the boundary values of $\mathcal{H}_{q;a,b,c}$ along $\partial \mathfrak{X}_{a,b,c}$ are fixed, as depicted in \Cref{fig:esagono}. Specifically, 
	\begin{eqnarray}
		\label{hqabcboundary}
	\mathcal{H}_{q;a,b,c} (x,y)=
          \begin{cases}
            0  & \text{ if } y = 0 \text{ or if }  x-y=a\\
            b & \text{ if }  y-x=b \text{ or if }  y=b+c\\
            y & \text{ if }   x=0 \\
            y-c & \text{ if }   x=a+c
          \end{cases}.
	\end{eqnarray}

        \noindent It  remains to determine its gradient $\nabla \mathcal{H}_{q;a,b,c}$. Following \cite{DLE}, we do this through its complex slope, defined below. Here, \eqref{f0uv}  is \cite[Equations (14) and (15)]{DLE}, and \eqref{fxy} and \eqref{equationu} are \cite[Theorem 1.6]{DLE}. Moreover, the first part of \Cref{fqxy} is \cite[Equation (17)]{DLE}; its second part is \cite[Proposition 5.19]{DLE}; and its third part is explained in the proof of \cite[Theorem 1.3]{DLE}. The parameters $(a_1, b_1, \mathsf{N}, \mathsf{T}; \kappa; \mathfrak{q})$ from \cite{DLE} are taken to be $(-c,b-c,b,a+c; 0; \mathfrak{q}^{-1})$ here, the coordinates $(t, x)$ from \cite{DLE} are taken to be $(x,y)$ here, and the parameter $u$ from \cite{DLE} is taken to be $u-c$ here.

	\begin{definition}[Complex slope]
		
		For any $u \in \mathbb{C}$, define 
		\begin{flalign}
                  \label{f0uv}
                  	f_0 (u) = \displaystyle\frac{(\mathfrak{q}^{u-b - c} - 1) (\mathfrak{q}^{u} - 1)}{(\mathfrak{q}^{u+a} - 1) (\mathfrak{q}^{u - b} - 1)}; \qquad U(u) = \mathfrak{q}^{u} \cdot \displaystyle\frac{f_0 (u)}{1 - f_0 (u)}; \qquad V(u) = \displaystyle\frac{\mathfrak{q}^{u}}{1 - f_0 (u)}. 
		\end{flalign} 
		
		\noindent For any $(x, y) \in \mathbb{R}^2$, define $u = u(x,y) = u_q (x, y)$ to be the unique solution in $\overline{\mathbb{H}}$ to the equation
		\begin{flalign}
                  \label{equationu}
                  	\mathfrak{q}^y =  V(u) - \mathfrak{q}^x \cdot U(u),
		\end{flalign}
		whenever \eqref{equationu}, which is quickly verified to be a quadratic equation in $\mathfrak q^u$, has at most one real root. In this case, define the \emph{complex slope} $f_x(y):=f (x, y) $ by
		\begin{flalign}
			\label{fxy} 
			f (x,y) = \displaystyle\frac{\mathfrak{q}^{u} - \mathfrak{q}^{y}}{\mathfrak{q}^{u} - \mathfrak{q}^{y-x}}.
		\end{flalign}
              \end{definition}

	\begin{lem}[{\cite{DLE}}]
		
		\label{fqxy} 
		
		The following statements hold for any $(x, y) \in \overline{\mathfrak{L}}_{q;a,b,c}$.
		
		\begin{enumerate} 
			
			\item We have that  
			\begin{flalign}
				\label{fxygradient} 
				\arg f (x,y) = \pi (1 - \partial_y \mathcal{H}_{q;a,b,c} (x,y)); \qquad \arg (f (x,y)-1) = \pi (1 + \partial_x \mathcal{H}_{q;a,b,c} (x,y)).
			\end{flalign} 
			
			\item We have that $(x, y) \notin \mathfrak{L}_{q;a,b,c}$ if and only if $f(x, y) \in \mathbb{R}$. 
			
			\item We have that $(x,y) \in \mathfrak{L}_{q;a,b,c}$ if and only if both roots of \eqref{equationu} are not real. Moreover,  $(x, y) \in \mathfrak{A}_{q;a,b,c}$ if and only if \eqref{equationu} has a double (real) root for $u$.
			
		\end{enumerate} 
		
	\end{lem}

	\begin{rem} 
		
		\label{derivativeh}
		
		Using \eqref{fxygradient}, if $f(x,y) \in \mathbb{R}$, we can quickly determine $\nabla \mathcal{H}_{q;a,b,c} (x, y)$ depending on the interval in which $f(x,y)$ resides. Indeed, we have 
		\begin{flalign*}
			& \nabla \mathcal{H}_{q;a,b,c} (x, y) = (0, 1), \qquad \text{if $f(x,y) \in (0, 1)$}; \\
			& \nabla \mathcal{H}_{q;a,b,c} (x, y) = (-1, 1), \qquad \text{if $f(x,y) \in (1, \infty)$}; \\
			& \nabla \mathcal{H}_{q;a,b,c} (x,y) = (0, 0), \qquad \text{if $f(x,y) \in (-\infty, 0)$}.
		\end{flalign*}
		
	\end{rem}

	Together, \Cref{fqxy} and \Cref{hlx} determine the gradient $\nabla \mathcal{H}_{q;a,b,c}$  everywhere in $\mathfrak{X}_{a,b,c}$. This fixes the function $\mathcal{H}_{q;a,b,c} : \mathfrak{X}_{a,b,c} \rightarrow \mathbb{R}$ by its boundary condition \eqref{hqabcboundary}.
	
	We next provide the Euler--Lagrange equation satisfied by $\mathcal{H}_{q;a,b,c}$. The first part of the below lemma is \cite[Proposition 9.4]{RT} (which appeared earlier as \cite[Equation (3)]{KO}); the second follows quickly from the explicit form of the $(\mathfrak{a}_{ij})$ below, provided by \cite[Equations (5), (6), and (7)]{LTDMLS}. Recall the definition of $\mathcal T$ in \eqref{e:cT}.
        
	\begin{lem}[{\cite{KO,RT,LTDMLS}}]
		
		\label{ahl} 
		
		For any indices $i, j \in \{ x, y \}$, there exists $\mathfrak{a}_{ij} : \mathcal{T} \rightarrow \mathbb{R}$
                satisfying $\mathfrak a_{xy}=\mathfrak a_{yx}$ as well as the following two statements:
		\begin{enumerate} 
			\item For any $z \in \mathfrak{L}_{q;a,b,c}$, we have 
			\begin{flalign}
                          \displaystyle\sum_{i,j \in \{ x, y \}} 	\mathfrak{a}_{ij} (\nabla \mathcal{H}_{q;a,b,c} (z)) \cdot \partial_{ij} \mathcal{H}_{q;a,b,c} (z) = -q.
                          \label{eq:pdeq}
			\end{flalign}
			
			\item For any $\omega > 0$, there exists a constant $C = C(\omega) > 1$ such that the following holds. Whenever $\gamma\in\mathcal T$ with  $\dist (\gamma, \partial \mathcal{T}) \ge \omega$, we have
			\begin{flalign*} 
				C^{-1} \le \mathfrak{a}_{xx} (\gamma) \le C; \qquad |\mathfrak{a}_{xy} (\gamma)| \le C; \qquad C^{-1} \le \mathfrak{a}_{yy} (\gamma) \le C.
			\end{flalign*} 
			
		\end{enumerate} 
		
	\end{lem}

	\subsection{Reformulation of the complex slope} 
	
	\label{Slope2}
	
	In this section we change variables to express the quadratic equation \eqref{equationu} defining $u$ in a slightly cleaner way. In what follows, for any $z \in \mathbb{C}$, define $[z] = [z]_{\mathfrak{q}}$ by
	\begin{flalign*}
		[z] = \displaystyle\frac{\mathfrak{q}^{-z}-1}{\mathfrak{q}^{-1}-1},
	\end{flalign*}
	
	\noindent where we let $\mathfrak{q}^z = e^{z \log \mathfrak{q}}=e^{q z}$. Note that
          $
\lim_{q\to 0}          [z]= z.$
           Then, define $v = v(x,y)$ by 
	\begin{flalign}
		\label{v} 
		v = [u(x,y)],
	\end{flalign}
	
	\noindent and denote 
	\begin{flalign}
		\label{abcxy} 
		A = -[-a]; \qquad B = [b]; \qquad C = [b+c] - [b]; \qquad X = -[-x]; \qquad Y = -[-y],
	\end{flalign}
	
	\noindent Observe, by \eqref{f0uv}, that we have under this notation 
	\begin{flalign}
		\label{f0uv2} 
		f_0 (u) = \displaystyle\frac{v(v-B-C)}{(v + A) (v - B)}.
	\end{flalign}
	\begin{lem} 
		
		\label{equationv} 
		
		We have that 
		\begin{flalign}
			\label{equationv2} 
			\begin{aligned} 
				& v^2 \big(A + C - X + (1 - \mathfrak{q}^{-1}) (A+C) Y \big) \\
				& \qquad + v \big( X(B+C) - Y(A+C) - AB + (\mathfrak{q}^{-1} - 1) ABY \big) + ABY = 0.
			\end{aligned} 
		\end{flalign}
		
	\end{lem} 
	
	\begin{proof} 
		
		Using the definitions \eqref{f0uv} of $U$ and $V$, and recalling \eqref{v} and \eqref{abcxy}, the equation \eqref{equationu} defining $u$ is equivalent to 
		\begin{flalign}
			\label{v0u2} 
			(1 - \mathfrak{q}^{-1}) Y v + v = \mathfrak{q}^y v = X \cdot \displaystyle\frac{f_0 (u)}{1-f_0 (u)} + Y.
		\end{flalign}
		
		\noindent Inserting \eqref{f0uv2} into \eqref{v0u2} then yields the lemma.
	\end{proof} 
	
	\begin{cor}
		
		\label{axy} 
		
		We have that $(x, y) \in \mathfrak{A}_{q;a,b,c}$ if and only if 
		\begin{flalign}
			\label{axy2} 
			\begin{aligned} 
				\big( X(B+C) - Y(& A+C) - AB + (\mathfrak{q}^{-1} - 1) ABY \big)^2 \\
				& - 4 AB Y \big( A + C - X + (1 - \mathfrak{q}^{-1}) (A+C) Y \big) = 0.
			\end{aligned} 
		\end{flalign}
	\end{cor}
	
	\begin{proof}
		
		This follows from the third part of \Cref{fqxy}, with \Cref{equationv} (and the fact that \eqref{equationu} having a double root for $u$ is equivalent to \eqref{equationv2} having a double root for $v$). 
	\end{proof}
	
	\begin{proof}[Proof of \Cref{axy0}]

		By \Cref{axy}, the arctic boundary $\mathfrak{A}_{q;a,b,c}$ defines a conic equation for $(X,Y)$; this arctic boundary therefore defines a smooth curve in the coordinates $(X, Y)$. Since (by \eqref{abcxy}) the map from $(X, Y)$ to $(x, y)$ and its inverse are smooth, it follows that $\mathfrak{A}_{q;a,b,c}$ is a smooth curve (in the coordinates $(x, y)$). 
		
		It remains to show that $\mathfrak{A}_{q;a,b,c}$ is tangent to each edge of $\mathfrak{X}_{a,b,c}$ at exactly one point. To that end, it suffices (since $\mathfrak{A}_{q;a,b,c}$ is smooth) to show that $\mathfrak{A}_{q;a,b,c}$ meets any edge of the hexagon $\mathfrak{X}_{a,b,c}$ at exactly one point. All edges of $\mathfrak{X}_{a,b,c}$ can be analyzed similarly, so let us only focus on its southwest edge $\{ y = 0 \}$. Since $y = 0$, we have by \eqref{abcxy} that $Y=0$, and so \eqref{axy2} implies that $\big( X(B+C) - AB \big)^2 = 0$, meaning that $X = AB/(B+C)$. Hence, $\mathfrak{A}_{q;a,b,c}$ only intersects the southwest edge of $\mathfrak{X}_{a,b,c}$ at one point $p^{SW}_{q;a,b,c}=(x,0)$, which is the solution to $X = AB/(B+C)$; this establishes the lemma.
	\end{proof} 
	
	\begin{rem}
		
		\label{xabc}
		
		Observe from the proof of \Cref{axy0} that the tangency point $p^{SW}_{q;a,b,c}$ is defined by the solution $x$ to the equation $X = AB/(B+C)$. 
		
	\end{rem} 
	
	\begin{proof}[Proof of \Cref{hlx}]
		
		The identification of $\nabla \mathcal{H}_{q;a,b,c} (z)$ is entirely analogous for each of the six frozen regions in which $z$ is contained. So, let us assume that $z\in\mF^{SW}_{q;a,b,c}$. First observe that it suffices to show that $\nabla \mathcal{H}_q (z) = (0, 1)$ for any $z\in\mathfrak{A}^{SW}_{q;a,b,c}$  (recall \Cref{def:pa}). Indeed, denoting $z = (x, y)$, integrating this equality along $\mathfrak{A}^{SW}_{q;a,b,c}$ (and using the fact that $\mathcal{H}_{q;a,b,c} (x,0) = 0$, by \eqref{hqabcboundary}), it would follow that $\mathcal{H}_{q;a,b,c} (x, y) = y$ for each $(x, y) \in\mathfrak{A}^{SW}_{q;a,b,c}$. Since $\partial_y \mathcal{H}_{q;a,b,c} (z) \in [0, 1]$ for all $z \in \mathfrak{X}_{a,b,c}$, it follows that we must have $\mathcal{H}_{q;a,b,c} (x, y) \le y$ for all $(x, y) \in \mathfrak{X}_{a,b,c}$. Since equality $\mathcal{H}_{q;a,b,c} (x, y) = y$ is attained for $(x, y)$ on the boundary of $\mF^{SW}_{q;a,b,c}$, it follows that $\partial_y \mathcal{H}_{q;a,b,c} (x,y) = 1$ (and $\mathcal{H}_{q;a,b,c} (z) = y$) for all $(x,y)\in\mF^{SW}_{q;a,b,c}$. Hence, $\nabla \mathcal{H}_{q;a,b,c} (z) = (0, 1)$ for all $z\in\mF^{SW}_{q;a,b,c}$, which yields the lemma. 
		
		 To establish $\nabla \mathcal{H}_{q;a,b,c} (z) = (0, 1)$ for any $z \in\mathfrak{A}^{SW}_{q;a,b,c}$, it suffices by \Cref{derivativeh} to verify that for such $z$ we have $f(z) \in [0, 1]$; by \eqref{f0uv2}, this is equivalent to $v(z) \in [0, AB/(A+C)]$. It is also quickly verified from \eqref{equationv2} that $v = 0$ only if $Y = 0$ (equivalently, $y = 0$), and $v = AB/(A+C)$ only if $X = 0$ (equivalently, $x = 0$). Therefore, it suffices to show that $v(z_0) = 0$ for $z_0=p^{SW}_{q;a,b,c}$, and that $v(z_1) = AB/(A+C)$ for $z_1=p^W_{q;a,b,c}$. The two statements are proven very similarly, so let us only address the former. By \Cref{xabc}, we have at $p^{SW}_{q;a,b,c}$ that $X = AB/(B+C)$. Inserting this value of $X$ (and $Y = 0$) into \eqref{equationv2} yields $(B+C)^{-1} C(A+B+C) \cdot v^2 = 0$, so $v = 0$; this establishes the lemma.		
	\end{proof} 
	
	\begin{proof}[Proof of \Cref{convex}]
		
		By \Cref{axy}, the arctic boundary $\mathfrak{A}_{q;a,b,c}$ defines a conic equation for $(X,Y)$, whose coefficients are smooth in the parameters $(q, a, b, c)$ (since $\mathfrak{q}$ and $(A, B, C)$ are). Setting $q = 0$ and $(a, b, c) = (1, 1, 1)$, this conic equation takes the form $(X+Y-2)^2 + 3(X-Y)^2 = 3$, which defines a smooth, convex curve, with curvature uniformly bounded above and below. Hence, for $\varepsilon > 0$ sufficiently small, it follows for any $q \in (-\varepsilon, \varepsilon)$ and $a,b,c \in (1-\varepsilon, 1+\varepsilon)$ that the conic equation from \eqref{axy2} defining $(X, Y)$ is a smooth, convex curve, with curvature uniformly bounded above and below. Since, by \eqref{abcxy}, we have that $X = x + O(\varepsilon)$, that $\partial_x X = 1 + O (\varepsilon)$, and that $\partial_x^2 X = O (\varepsilon)$ (and similarly for $Y$), it quickly follows for $\varepsilon$ sufficiently small that $\mathfrak{A}_{q;a,b,c}$ is a smooth, convex curve (in the coordinates $(x, y)$), with curvature uniformly bounded above and below; this establishes the claim.
	\end{proof}

	\begin{proof}[Proof of \Cref{ul}]
		Observe since $(\mathfrak{l}_{q;a,b,c}(h), h) \in \mathfrak{A}^{SW}_{q;a,b,c}$ that $\mathcal{H}_{q;a,b,c} (\mathfrak{l}_{q;a,b,c}(h),h) = h$, by integrating the second statement of \Cref{hlx} (from $z$ starting at $(x^{SW}_{q;a,b,c}, 0)$, and ending at $(\mathfrak{l}_{q;a,b,c}(h), h)$, along $\mathfrak{A}^{SW}_{q;a,b,c}$). Moreover, since $\mathfrak{A}^{SW}_{q;a,b,c}$ is convex by \Cref{convex}, we have that $(x, y)\in\mF^{SW}_{q;a,b,c}$ for any  $(x,y) \in [0, \mathfrak{l}_{q;a,b,c}(h)] \times [0, h]$. Therefore, it follows again from integrating the second part of \Cref{hlx} that $\mathcal{H}_{q;a,b,c} (x, y) =y$ for such $(x, y)$. Recalling the definition of $\mathcal{U}_{q;a,b,c}^h$ from \eqref{e:geninversa}, this implies that $\mathcal{U}_{q;a,b,c}^h (x) = h$ for any $x \in [0, \mathfrak{l}_{q;a,b,c}(h)]$; this establishes \Cref{ul}. 
	\end{proof}

	\subsection{Complex slope behavior away from the arctic boundary} 
	
	\label{SlopeBoundary00} 
	
	Throughout this section, we let $\varepsilon \in (0, 1)$ denote a sufficiently small constant, and we fix $q \in (-\varepsilon, \varepsilon)$ and $a, b, c \in (1-\varepsilon, 1 + \varepsilon)$. The proof of \Cref{prop:formerassumption}, \Cref{hde0} and \Cref{prop:last} require distinct arguments according to whether one considers  the region mesoscopically far from the arctic boundary (this is dealt with in the present section) or  the one close to it (this is the topic of  \Cref{SlopeScale}).
        
	\label{Estimatef}
	
	In what follows, we define $\xi, \eta, \zeta : \mathbb{R}^2 \rightarrow \mathbb{R}$ as
	\begin{flalign}
		\label{xietazeta} 
		\begin{aligned} 
			\xi (X, Y) & = \big( X(B+C) - Y(A+C) - AB + (\mathfrak{q}^{-1} -1 )ABY \big)^2 \\
			& \qquad - 4 ABY \big( A+C-X + (1-\mathfrak{q}^{-1}) (A+C) Y \big); \\	
			\eta (X, Y) & = Y(A+C) -X(B+C) + AB + (1-\mathfrak{q}^{-1}) ABY;  \\
			\zeta (X, Y) & = 2 \big( A+C  - X + (1-\mathfrak{q}^{-1}) (A+C) Y \big);
		\end{aligned} 
	\end{flalign}
	
	\noindent in particular, $\xi$ is defined as the left side of \eqref{axy2}. Then,
	\begin{flalign}
		\label{vxietazeta} 
		v = \displaystyle\frac{\eta (X, Y) + \xi(X,Y)^{1/2}}{\zeta (X,Y)}, 
	\end{flalign}
	
	\noindent whenever $\xi (X, Y) \le 0$ (if $\xi (X, Y) < 0$, then we take the root above with positive imaginary part). 
	
	We next require the following lemma. Its first statement is \cite[Theorem 4.4]{LTDMLS}, and its second follows quickly from \cite[Equations (28), (29), and (30)]{LTDMLS}.
	
	\begin{lem}[\cite{LTDMLS}]
		
		\label{wa} 
		Set $c=3-a-b$ and define the sets 
		\begin{flalign*}
			& \mathcal{A} = \{ (s, t, u, v) \in \mathbb{R}^4 : (s, t) \in \mathcal{T} \}; \\
			& \mathcal{W} = \{ (a, b, x, y) \in \mathbb{R}^4 : 0 < \min \{ a, b \} < a+b < 3, (x, y) \in \mathfrak{L}_{a,b,c} \},
		\end{flalign*}
		
		\noindent and define the map $\Phi :  \mathcal{W} \mapsto\mathcal{A}$ by setting 
		\begin{flalign*} 
			\Phi (a,b, x, y) = ( \partial_x \mathcal{H}_{a,b,c} (x,y), \partial_y \mathcal{H}_{a,b,c}(x,y), \partial_{xx} \mathcal{H}_{a,b,c} (x, y), \partial_{xy} \mathcal{H}_{a,b,c} (x,y)).
		\end{flalign*}
		
		\begin{enumerate} 
			\item The map $\Phi : \mathcal{W} \rightarrow \mathcal{A}$ is a diffeomorphism. 
			\item For any compact subset $\mathcal{W}_0 \subset \mathcal{W}$, there exists a constant $C = C(\mathcal{W}_0)$ such that the  Jacobian of $\Phi$ is bounded above by $C$ and below by $C^{-1}$ everywhere on $\mathcal{W}_0$. 
		\end{enumerate} 
		
	\end{lem} 
        
	\begin{lem} 
		
		\label{omegaf} 
		
		Fix $\omega > 0$. For $\varepsilon > 0$ sufficiently small (dependent on $\omega$), \Cref{prop:formerassumption}, \Cref{hde0}, and \Cref{prop:last} hold (with their constant $C$ dependent on $\omega$) whenever $\mathfrak{e}_{z_0;q} \ge \omega$.
		
	\end{lem} 

	\begin{proof} 
		
		Throughout this proof, we set $z_0 = (x_0, y_0)$. Let $\mathfrak{B} = \mathfrak{B}(\omega) > 1$ be a constant, to be fixed later, and let $z = (x, y) \in \mathfrak{L}_{q;a,b,c}$ be a point satisfying $|z-z_0| \le \mathfrak{B}^{-1}$. Denote the $(X, Y)$ associated with $(x_0, y_0)$ and $(x, y)$ by $(X_0, Y_0)$ and $(X, Y)$, respectively.
		
		We begin by verifying \Cref{prop:formerassumption}. Observe since $\mathfrak{e}_{z_0;q} \ge \omega \gtrsim 1$, we have $\mathfrak{d}_{z_0} \gtrsim \mathfrak{e}_{z_0;q} \gtrsim 1$, where the first statement is a consequence of \Cref{convex} (with the fact that $\mathfrak{A}_{q;a,b,c}$ is tangent to each edge of $\mathfrak{X}_{a,b,c}$). Therefore,
		\begin{flalign}
			\label{dz02} 
			\mathfrak{d}_{z_0} \asymp \mathfrak{e}_{z_0;q} \asymp 1.
		\end{flalign} 
		
		\noindent By \eqref{e:funzrescaled}, it therefore suffices to show that $|\partial_{\gamma} \mathcal{H}_{q;a,b,c} (z)| \lesssim 1$ for each differential operator $\partial_\gamma$ with $|\gamma| \in [1, B]$. One way of deducing this would be through Morrey's theorem, using the Euler--Lagrange equation given by \Cref{ahl}. However, let us do this more directly (as it will be anyways useful for confirming the second part of the lemma), which will proceed by showing that $v$ is uniformly smooth in $(a, b, c, x, y)$. 
		
		To that end, we will require several bounds on $(\xi, \eta, \zeta)$. Since $\mathfrak{e}_{z_0;q} \asymp 1$ (and $\mathfrak{A}_{q;a,b,c}$ is defined by the set of $(x, y)$ for which $\xi (X,Y) = 0$), it is quickly verified that $-\xi (X_0, Y_0) \asymp 1$; by the continuity of $\xi$, it follows that $-\xi (X, Y) \asymp 1$. One also clearly has $\eta (X, Y) \lesssim 1$. By \eqref{dz02}, we also have $b+c-x  \gtrsim 1$ (where the  bound follows from the fact that $b+c-x$ is the distance from $z$ to the northeast edge of $\mathfrak{X}_{a,b,c}$). Therefore, for $\varepsilon$ sufficiently small, $\zeta (X, Y) = 2(A+C-X + (1- \mathfrak{q}^{-1} ) (A+C) Y) \gtrsim 1$, as $\zeta (X, Y) \approx 2 - x \approx b+c-x$ for $\varepsilon$ sufficiently small; we also have that $\zeta (X, Y) \lesssim 1$, since $(A, B, C, X, Y, q)$ are bounded. Together, these give
		\begin{flalign}
			\label{xietazeta1} 
			-\xi (X, Y) \asymp 1; \qquad \eta (X, Y) \lesssim 1; \qquad \zeta (X, Y) \asymp 1. 
		\end{flalign}
		
		\noindent By the explicit forms of $(\xi, \eta, \zeta)$ from \eqref{xietazeta}, $\xi$, $\eta$, and $\zeta$ are all uniformly smooth (meaning that their derivatives of order at most $m$ are bounded by a constant depending on $m$) in the parameters $(A, B, C, X, Y, q)$, whenever $|z-z_0| \le \mathfrak{B}^{-1}$. Together with \eqref{vxietazeta} (and the fact that the map from $(A,B,C,X,Y)$ to $(a,b,c,x,y)$ is uniformly smooth for $\varepsilon$ sufficiently small), this implies that $v$ is uniformly smooth in $(a, b, c, x, y, q)$ for such $z = (x,y)$. 
		
		We also have from \eqref{xietazeta1} that $|v| \lesssim 1$ and that $\Imaginary v \asymp 1$. These two facts, with \eqref{fxygradient}, imply that $\partial_x \mathcal{H}_{q;a,b,c} (x,y) \asymp 1$ and $\partial_y \mathcal{H}_{q;a,b,c} (x,y) \asymp 1$ whenever $z = (x,y)$ satisfies $\mathfrak{e}_{z;q} \ge \omega$. This confirms the first two statements of \Cref{hde0} in this case; the third follows from integrating the first. The proof of the fourth is entirely analogous. 
		
		Moreover, by \eqref{f0uv2}, the map from $v = v(x,y)$ to $f = f (x, y)$ is (uniformly) smooth whenever $\Imaginary v \asymp 1$ and $|v| \lesssim 1$, in which case $\Imaginary f \asymp 1$ and $|f| \lesssim 1$. Furthermore, we have from \eqref{fxygradient} that the map from $f = f(x,y)$ to $\nabla \mathcal{H}_{q;a,b,c} (x,y)$ is (uniformly) smooth whenever $\Imaginary f \asymp 1$ and $|f| \lesssim 1$. Composing these maps, we deduce that $\nabla \mathcal{H}_{q;a,b,c} (x,y)$ is (uniformly) smooth in $(a, b, c, x, y, q)$, if $|z-z_0| \le \mathfrak{B}^{-1}$. The smoothness in $(x,y)$ in particular verifies \Cref{prop:formerassumption} whenever $\me_{z_0;q}\ge \omega$.
		
		To confirm \Cref{prop:last} in the same regime, we use \Cref{wa}. Let $\nu = \nu(\omega) > 0$ denote some sufficiently small constant to be fixed later, and set
		\begin{flalign}
                  (s_0, t_0, u_0, v_0) = (\partial_x \mathcal{H}_q (z_0), \partial_y \mathcal{H}_q (z_0), \partial_{xx} \mathcal{H}_q (z_0) + \nu q, \partial_{xy} \mathcal{H}_q (z_0)).
                  \label{eq:abqnew}
		\end{flalign} 
		\noindent Observe that $(s_0, t_0, u_0, v_0)$ is in a  compact subset of $\mathcal{A}$ (only dependent on $\omega$). Indeed, since $\Imaginary f(z_0) \asymp 1$ and $|f (z_0)| \lesssim 1$, we have from \eqref{fxygradient} that 
		\begin{flalign}
			\label{hz01}
			\dist (\nabla \mathcal{H}_q (z_0), \partial \mathcal{T}) \asymp 1.
		\end{flalign}
		
		\noindent We also have by the first part of the lemma $|\partial_{\gamma} \mathcal{H}_{q;a,b,c} (z_0)| \lesssim 1$, whenever $|\gamma| \le 2$.
                Recalling the function $\Phi$ from \Cref{wa}, set $(a, b, x', y') = \Phi^{-1} (s_0, t_0, u_0, v_0)$, set $c = 3 - a - b$ and $z_1 = (x', y')$; let us show that \Cref{prop:last} holds for this choice of parameters.
		
		To that end, set 
		\begin{flalign}
			\label{s1t1} 
			(s_1, t_1, u_1, v_1) = (\partial_x \mathcal{H} (z_0), \partial_y \mathcal{H} (z_0), \partial_{xx} \mathcal{H} (z_0), \partial_{xy} \mathcal{H} (z_0)),
		\end{flalign}
		
		\noindent which satisfies $\Phi (1,1,x_0, y_0) = (s_1, t_1, u_1, v_1)$. From the above fact that $\mathcal{H}_{q} (x,y)$ is uniformly smooth in $(x,y,q)$, we have 
		\begin{flalign}
			\label{s1s0} 
			|s_1 - s_0| + |t_1 - t_0| + |u_1 - u_0| + |v_1 - v_0| \lesssim q.
		\end{flalign}
		
		\noindent Observe that $(s_1, t_1, u_1, v_1)$ is also in a compact  subset of $\mathcal{A}$ that only depends on $\omega$ (as $(s_0, t_0, u_0, v_0)$ is, by \eqref{s1s0} and the fact that $|q| \le \varepsilon$ is sufficiently small). Therefore, by the second part of \Cref{wa}, the Jacobian of $\Phi^{-1}$ is bounded in a neighborhood of $(s_0, t_0, u_0, v_0)$ and $(s_1, t_1, u_1, v_1)$. Together with \eqref{s1s0} and the fact that $\Phi^{-1} (s_1, t_1, u_1, v_1) = (1, 1, x_0, y_0)$, it follows that 
		\begin{flalign}
			\label{abq} 
			|a-1| + |b-1| + |x'-x_0| + |y'-y_0| \lesssim q.
		\end{flalign} 
		
		\noindent Since $c = 3-a-b$, it also follows that $|c-1| \lesssim q$, which with \eqref{dz02} yields the claim in \Cref{prop:last} about $a,b,c$ being close to $1$. The first  statement in \eqref{e:thereexist} follows from the fact that $(s_0, t_0) = \nabla \mathcal{H}_q (z_0)$ by  \eqref{eq:abqnew}, with the definition of $(a,b,x',y') = \Phi^{-1} (s_0, t_0, u_0, v_0)$. 
		
		By \eqref{dz02} and \eqref{e:funzrescaled}, to establish the second statement in \eqref{e:thereexist}, we must show that $D^2{\mathcal{H}_{a,b,c} (z_1)} - D^2 {\mathcal{H}_q (z_0)} \gtrsim q \cdot \Id$. Since 
		\begin{flalign}
			\label{xy} 
			\partial_{xx} \mathcal{H}_{a,b,c} (z_1) = u_0 = \partial_{xx} \mathcal{H}_q (z_0) + \nu q; \qquad \partial_{xy}\mathcal{H}_{a,b,c} (z_1) = v_0 = \partial_{xy} \mathcal{H}_q (z_0),
		\end{flalign}
		it suffices to verify that 
		\begin{flalign} 
			\label{yy0} 
			\partial_{yy} \mathcal{H}_{a,b,c} (z_1) - \partial_{yy} \mathcal{H}_q (z_0) \gtrsim q, 
		\end{flalign} 
		if $\nu > 0$ is sufficiently small. To do this, observe from the first statement of \Cref{ahl} (with the fact that $\nabla \mathcal{H}_{a,b,c} (z_1) = (s_0, t_0) = \nabla \mathcal{H}_q (z_0)$) that 
		\begin{flalign*}
			& \mathfrak{a}_{xx} (s_0, t_0) \cdot \partial_{xx} \mathcal{H}_{a,b,c} (z_1) + 2 \mathfrak{a}_{xy} (s_0, t_0) \cdot \partial_{xy} \mathcal{H}_{a,b,c} (z_1) + \mathfrak{a}_{yy} (s_0, t_0) \cdot \partial_{yy} \mathcal{H}_{a,b,c} (z_1) = 0; \\
			& \mathfrak{a}_{xx} (s_0, t_0) \cdot \partial_{xx} \mathcal{H}_q (z_0) + 2 \mathfrak{a}_{xy} (s_0, t_0) \cdot \partial_{xy} \mathcal{H}_q (z_0) + \mathfrak{a}_{yy} (s_0, t_0) \cdot \partial_{yy} \mathcal{H}_q (z_0) = -q.
		\end{flalign*} 
	Subtracting these equations and using \eqref{xy} gives
		\begin{flalign}
			\label{yy2} 
			\mathfrak{a}_{yy} (s_0, t_0) \cdot \big( \partial_{yy} \mathcal{H}_{a,b,c} (z_1) - \partial_{yy}\mathcal{H}_q (z_0) \big) = q (1 - \mathfrak{a}_{xx} (s_0, t_0) \cdot \nu \big).
		\end{flalign}
		By the second part of \Cref{ahl}, together with \eqref{hz01} (and the fact that $\nabla \mathcal{H}_q (z_0) = (s_0, t_0)$), there exists a constant $C = C(\omega) > 1$ such that $C^{-1} \le \mathfrak{a}_{yy} (s_0, t_0) \le C$ and $\mathfrak{a}_{xx} (s_0, t_0) \le C$. Setting $\nu = (2C)^{-1}$ then yields \eqref{yy0}, by \eqref{yy2}.

		The last statement of \Cref{prop:last} follows from the above-mentioned fact that each derivative of $\mathcal{H}_{q;a,b,c} (z)$ is uniformly smooth in the parameters $(q, a, b, c)$, whenever $|z-z_0| \le \mathfrak{B}^{-1}$, together with the bound \eqref{abq} (and $|c-1| \lesssim q$).		
	\end{proof}

	\section{Behavior of limit shapes near the arctic boundary}	
\label{appC}

	Throughout this section, we let $\varepsilon \in (0, 1)$ denote a sufficiently small constant, and we fix $q \in (-\varepsilon, \varepsilon)$ and $a, b, c \in (1-\varepsilon, 1 + \varepsilon)$. 
	
	\subsection{Scaling of limit shapes}
	
	\label{SlopeScale}
	
        \begin{rem}
          \label{rem:rightoftangency}
	In view of \Cref{omegaf}, to establish \Cref{prop:formerassumption}, \Cref{hde0} and \Cref{prop:last} we may restrict our attention to the case when $z_0 = (x_0, y_0)$ is close to the arctic boundary $\mathfrak{A}_{q;a,b,c}$; in particular, $\dist (z_0, \mathfrak{A}_{q;a,b,c}) \le 1/100$. Let us abbreviate $(\mathfrak{d}, \mathfrak{e}) = (\mathfrak{d}_{z_0}, \mathfrak{e}_{z_0;q})$.
	
	We will assume in what follows that the closest edge $\ell_{z_0}$ of $\mathfrak{X}_{a,b,c}$ to $z_0$ is its south one $\{ y = 0 \}$, and that $x_0 \ge \mathfrak{l}_{q;a,b,c} (0) = x^{SW}_{q;a,b,c}$, as the other cases are entirely analogous. Then, it is quickly verified for $\varepsilon$ sufficiently small, which we assume throughout this section (hence, $\mq$ close to $0$ and $a,b,c$ close to $1$) that 
	\begin{flalign}
		\label{x0y0} 
		\displaystyle\frac{1}{2} - O(\varepsilon) \le x_0 \le 1 + O(\varepsilon); \qquad y_0 \le \displaystyle\frac{49}{100}, 
	\end{flalign} 
	
	\noindent as otherwise we would not have $\dist (z_0, \mathfrak{A}_{q;a,b,c}) \le 1/100$. This is the reason why some  statements in this section are given for $x_0\in[1/4,5/4]$.
        \end{rem}

         \subsubsection{Curvature and tangent approximations}

         \label{CurvatureDensity2}

         In this section, we obtain exact and asymptotic expressions for the  curvature and the slope of the arctic boundary of $\cH_{q;a,b,c}$. This involves asymptotic expressions for the 
         functions $f_0$ and $f_x$. These functions are in principle explicit, but a bit intricate, and so it will be more direct to study $f_0$ and $f_x$ through the method of characteristics \cite[Section 6.2]{DLE}.
         The computations here are similar to those in \cite[Sections 10.2 and 10.4]{SCLE}. 
         
         First, following \cite[Proposition 5.5]{DLE}, for any $x \in \mathbb{R}_{\ge 0}$ and $z \in \overline{\mathbb{H}}$, let
         \begin{flalign}
         	\label{fxz}
         	f_x (z) = \displaystyle\frac{\mathfrak{q}^u - \mathfrak{q}^z}{\mathfrak{q}^u - \mathfrak{q}^{z-x}},
         \end{flalign}
         
         \noindent where $u=u_x(z) \in \overline{\mathbb{H}}$ solves (see \cite[Remark 5.6]{DLE} for the choice of solution, if there is more than one) 
         \begin{flalign}\label{eq:exnonnum}
         	\mathfrak{q}^{z-x} = \mathfrak{q}^u \cdot \displaystyle\frac{f_0 (u)}{f_0 (u) - 1} + \displaystyle\frac{\mathfrak{q}^{u - x}}{1 - f_0 (u)}= \mathfrak q^{-x}V(u)-U(u).
         \end{flalign}
         If $z=y\in\mathbb R$, then $u_x(z)$ coincides with $u(x,y)$ defined by \eqref{equationu} and $f_x(z)$ coincides with $f(x,y)$ defined by \eqref{fxy}. Given this, the below definition \eqref{zxequation} coincides with \cite[Equation (214)]{DLE} (with the $(\mathsf{q},\kappa)$ there equal to $(\mathfrak{q}^{-1},0)$ here).
         
\begin{definition}

	\label{zxu}

	For any $x \in [0, a]$ and $u \in \mathbb{C}$, let $z_x = z_x (u) = z(x, u) \in \mathbb{C}$ solve the equation (with $f_x(z)$ as in \eqref{fxz})
	\begin{flalign}
		\label{zxequation} 
		\partial_x z_x = \displaystyle\frac{f_x(z_x)}{f_x (z_x) - 1}, \quad \text{where $z_0 = u$}.
	\end{flalign} 
\end{definition}

\begin{definition} 
	
	\label{gammaxi} 
	
	For any $x \in [0, a]$, let $\gamma_x = \gamma (x) \in \mathbb{R}_{\ge 0}$ denote the minimal real number such that $(x, \gamma_x) \in \mathfrak{A}_{q;a,b,c}$ is on the arctic boundary. In this way, we parameterize the (union of the southwest and south parts of) the arctic boundary by $(x, \gamma_x)_{x \in [0, a]}$. We further define $\xi_x = \xi(x)$ by  
	\begin{flalign*}
		\xi_x = u_x ( \gamma_x), \quad \text{so that} \quad  \gamma_x = z_x ( \xi_x)
	\end{flalign*}
	where the latter holds since we more generally have $z_x (u_x (y)) = y$, by \cite[Equation (222)]{DLE}.\footnote{Indeed, \cite[Equation (222)]{DLE} (with the $(\mathsf{q},\kappa)$ there equal to $(\mathfrak{q}^{-1},0)$ here) gives $q^{z_x (u_x(y))} = V(u_x(y)) - q^{x} U(u_x(y))$; by \eqref{equationu} (and \cite[Remark 5.6]{DLE}), this yields $z_x (u_x (y)) =y$.} We frequently abbreviate $\xi_x' = \xi' (x)$, $\gamma_x' = \gamma' (x)$, and $\gamma_x'' = \gamma'' (x)$. 
	
\end{definition}

\begin{rem} 
	
	\label{xiq1} 
	
	When $\mathfrak{q}=1$, it can be verified that $\xi_x$ denotes the $y$-coordinate of the point where the tangent line to $\mathfrak{A}_{a, b, c}$ at $(x, \gamma_x)$ meets the $y$-axis. Indeed, \eqref{fxf0} below indicates that $f_x (z_x)$ is constant in $x$ when $\mathfrak{q}=1$, and so the right side of the equation \eqref{zxequation} defining $z_x$ is constant in $x$. Moreover, by \Cref{gammatderivative} below, this right side is equal to $\gamma_x'$. Therefore, since $z_x = \gamma_x$, the path $(z_s)_{s \in \mathbb{R}}$ follows the line tangent to $\mathfrak{A}_{a,b,c}$ at $(x, \gamma_x)$. In particular, $\xi_x = z_0$ is where this line meets the $y$-axis.
	
\end{rem}

Recall the definition \eqref{fxy} of  $f_x (y)$. By \cite[Equation (215)]{DLE} (with the $(\mathsf{q},\kappa)$ there equal to $(\mathfrak{q}^{-1},0)$ here), we have  for any $(x, y) \in \mathfrak{L}_{q;a,b,c}$ that
\begin{flalign}
	\label{fxzxderivative} 
	\partial_x ( f_x (z_x) ) =  \log \mathfrak{q} \cdot f_x (z_x),
\end{flalign}

\noindent and integrating over $x$ yields
\begin{flalign}
	\label{fxf0} 
	f_x (z_x) = \mathfrak{q}^x \cdot f_0 (z_0) = \mathfrak{q}^x f_0 (u).
\end{flalign}

\noindent Then,
\begin{flalign}
	\label{zxuf0}
	\mathfrak{q}^{z_x} = \mathfrak{q}^u \cdot \displaystyle\frac{1 - \mathfrak{q}^x f_0 (u)}{1 - f_0 (u)} = \mathfrak{q}^{u} \cdot \displaystyle\frac{1 - f_x (z_x)}{1 - \mathfrak{q}^{-x} f_x (z_x)},
\end{flalign}

\noindent where the first equality follows from inserting \eqref{fxf0} into \eqref{fxzxderivative} and integrating over $x$, and the second follows from \eqref{fxf0}.

The following lemma expresses the slope $\gamma_x'$ of the tangent line to the arctic boundary at $(x, \gamma_x)$ in terms of $f_x (\gamma_x)$. 

\begin{lem} 
	
	\label{gammatderivative} 
	
	For any $x \in [0, a]$, we have  
	\begin{flalign*}
		\gamma_x' = \displaystyle\frac{f_x ( \gamma_x)}{f_x (\gamma_x) - 1}. 
	\end{flalign*}
	
\end{lem} 
      \begin{rem}
        \label{rem:gamma'}
        {If $(x,\gamma_x)\in\mA^S_{q;a,b,c}$ and is bounded away from $p^{SE}_{q;a,b,c}$,   then $\gamma'_x$ is positive and strictly smaller than $1$ (the value $1$ is the slope $\gamma'_x$ at $p^{SE}_{q;a,b,c}$). It also follows from the statement of \Cref{gammatderivative} that $f_x(\gamma_x)<0$.}        
      \end{rem}

\begin{proof} 
	Observe that the slope of the arctic boundary $\mathfrak A_{q;a,b,c}$ at $( x, \gamma_x)$ satisfies $\gamma_x' \partial_y \cH_{q;a,b,c} + \partial_x \cH_{q;a,b,c} = 0$. Near this point, we have by \eqref{fxygradient}
	\begin{flalign*}
		\displaystyle\frac{\Imaginary f}{f} \sim -\pi \partial_y \cH_{q;a,b,c}; \qquad \displaystyle\frac{\Imaginary f}{f-1} \sim \pi \partial_x \cH_{q;a,b,c},
	\end{flalign*}
	(where $\sim$ indicates that the ratio of the two sides tends to $1$ as the argument $z$ of $f$ converges to $(x,\gamma_x)$) which together with the equality $\gamma_x' = -\partial_x  \cH_{q;a,b,c} \cdot (\partial_y  \cH_{q;a,b,c})^{-1}$ yields the lemma.
\end{proof}

Now, we seek to understand the limit shape $\mathcal{H}_{q;a,b,c} (x,y)$ for $(x,y) \in \mathfrak{L}_{q;a,b,c}$ near the arctic boundary $\mathfrak{A}_{q;a,b,c}$. To do so, we must analyze $f_x (y)$ for $y \approx \gamma_x$. The exact form \eqref{fxy} of $f_x (y)$ involves the function $q^{u_x (y)}$ that solves the quadratic equation \eqref{equationu}. It will be more direct to approximate $f_x (y)$ without solving the latter equation, to which end we will instead use \eqref{fxf0} to express $f_x (y)$ in terms of $f_0 (u)$, where $u$ is such that $z_x (u) = y$. If $y \approx \gamma_x$, so $u \approx \xi_x$, by Definition \ref{gammaxi}. Therefore, we will be interested in Taylor expanding $f_0 (u)$ around $u \approx \xi_x$, which will involve the first two derivatives of $f_0$ at $\xi_x$. As these derivatives will be the underlying constants underlying our approximation for the limit shape, it will be convenient to express quantities of interest on the arctic boundary (such as its curvature) through these derivatives, as well.

The next lemma expresses the second derivative $\gamma_x''$ of the arctic boundary at $(x, \gamma_x)$ in $x$, $\xi_x$, $\gamma_x$, $f_0$, and $f_x (\gamma_x)$. The quantity $\Gamma_{2;x}$ will appear in \Cref{imaginaryftx0} below.

\begin{lem}
	
	\label{gammatderivative2} 
	
	For any $x \in [0,a]$, define $A = A_x$ by
	\begin{flalign} 
		\label{af0}
		A = f_0'' (\xi_x) - \displaystyle\frac{\mathfrak{q}^{x + \xi_x - \gamma_x} (\mathfrak{q}^{x + \xi_x - \gamma_x} + 1)}{(\mathfrak{q}^{x + \xi_x - \gamma_x}-1)^3} \cdot (\mathfrak{q}^{-x}-1) (\log \mathfrak{q})^2.
	\end{flalign}
	
	\noindent Further denoting
	\begin{flalign}
		\label{gammax12} 
		\Gamma_{x;1} = -\displaystyle\frac{f_x (\gamma_x)}{( f_x (\gamma_x) - 1 )^2} \cdot \log \mathfrak{q}; \qquad \Gamma_{x;2} = ( f_x (\gamma_x)-1 )^2 \cdot \displaystyle\frac{\mathfrak{q}^{3\xi_x - 3\gamma_x-x} (\log \mathfrak{q})^3}{A (1 - \mathfrak{q}^{-x})^3},
	\end{flalign}
	
	\noindent we have 
	\begin{flalign}
		\label{gammaderivative2} 
		\gamma_x'' = \Gamma_{x;1} + \Gamma_{x;2}.
	\end{flalign}

	\noindent In particular, 
	\begin{flalign}
		\label{gamma2gamma} 
		\Gamma_{x;2} = \gamma_x'' - \gamma_x' (1 - \gamma_x') \log \mathfrak{q}. 
	\end{flalign} 
	
\end{lem}

To show the above lemma, we require the below one that evaluates $f_0 (\xi_x)$ and $f_0' (\xi_x)$ as explicit rational functions of $\mathfrak{q}^x$, $\mathfrak{q}^{\xi_x}$, and $\mathfrak{q}^{\gamma_x}$.

\begin{lem} 
	
	\label{gammatderivative0} 
	
	For any $x \in [0,a]$, we have 
	\begin{flalign}
		\label{f0xi}
		f_0 ( \xi_x ) = \displaystyle\frac{\mathfrak{q}^{\xi_x - \gamma_x} - 1}{\mathfrak{q}^{x + \xi_x - \gamma_x} - 1}; \qquad f_0' (\xi_x) = \displaystyle\frac{\mathfrak{q}^{\xi_x - \gamma_x}}{(\mathfrak{q}^{x + \xi_x - \gamma_x} - 1)^2} \cdot (\mathfrak{q}^x - 1) \log \mathfrak{q}.
	\end{flalign}
	
\end{lem} 

\begin{proof} 
  The first statement of the lemma follows from \eqref{zxuf0} at $(z_x, u) = (\gamma_x, \xi_x)$. The second follows from the equality
 	\begin{flalign*}
		f_0' ( \xi_x ) & = \displaystyle\lim_{s\rightarrow x} \displaystyle\frac{f_0 ( \xi_x ) - f_0 (\xi_{s} )}{\xi_x - \xi_{s}} 
		 = (\xi_x')^{-1} \cdot \partial_x (f_0 (\xi_x)) 
		 = (\xi_x')^{-1} \cdot \partial_x \Big( \displaystyle\frac{\mathfrak{q}^{\xi_x - \gamma_x} - 1}{\mathfrak{q}^{x + \xi_x - \gamma_x} - 1} \Big) \\
		& = \displaystyle\frac{\mathfrak{q}^{x + \xi_x - \gamma_x} \log \mathfrak{q}}{(\mathfrak{q}^{x + \xi_x - \gamma_x}-1)^2 \xi_x'} ( 1 - \mathfrak{q}^{\xi_x - \gamma_x} + (1 - \mathfrak{q}^{-x}) \cdot \xi_x' + (\mathfrak{q}^{-x} - 1) \cdot \gamma_x'  ) \\
		& = \displaystyle\frac{\mathfrak{q}^{x + \xi_x - \gamma_x}}{(\mathfrak{q}^{x + \xi_x - \gamma_x}-1)^2} \cdot (1 - \mathfrak{q}^{-x}) \log \mathfrak{q},
	\end{flalign*} 
	
	\noindent where we used the first statement of \eqref{f0xi} and the fact that 
	\begin{flalign}
		\label{gammaxderivative}
		\gamma_x' = \displaystyle\frac{f_x (\gamma_x )}{f_x ( \gamma_x ) - 1} = \displaystyle\frac{\mathfrak{q}^x f_0 ( \xi_x )}{\mathfrak{q}^x f_0 ( \xi_x ) - 1} = \displaystyle\frac{\mathfrak{q}^{x + \xi_x - \gamma_x} - \mathfrak{q}^x}{1 - \mathfrak{q}^x}.
	\end{flalign}
\end{proof} 

We additionally require the below lemma evaluating the derivative of $\xi$; it can be viewed as the source of the quantity $A = A_x$ from \eqref{af0}.

\begin{lem}
	
	\label{fxiderivative} 
	
	Let $x \in [0, a]$. Recalling $A = A_x$ from \eqref{af0}, we have 
	\begin{flalign*}
		\xi_x' = -\displaystyle\frac{\mathfrak{q}^{x+2\xi_x-2\gamma_x} (\log \mathfrak{q})^2}{A (\mathfrak{q}^{x+\xi_x-\gamma_x}-1)^2}.
	\end{flalign*} 
	
\end{lem}

\begin{proof} 
	
	By \Cref{gammatderivative}, we have
	\begin{flalign*}
		\xi_x' f_0'' (\xi_x) & = \partial_x ( f_0' (\xi_x)) \\
		& = \displaystyle\frac{\mathfrak{q}^{\xi_x-\gamma_x} (\log \mathfrak{q})^2}{(\mathfrak{q}^{x+\xi_x-\gamma_x}-1)^2} +  \displaystyle\frac{\mathfrak{q}^{x+\xi_x-\gamma_x} (\mathfrak{q}^{-x}-1)( \gamma_x'-\xi_x'-1) (\log \mathfrak{q})^2}{(\mathfrak{q}^{x+\xi_x-\gamma_x}-1)^2} \\
		& \qquad - \displaystyle\frac{2\mathfrak{q}^{2x+2\xi_x-2\gamma_x} (\mathfrak{q}^{-x}-1)(\gamma_x'-\xi_x'-1) (\log \mathfrak{q})^2}{(\mathfrak{q}^{x+\xi_x-\gamma_x}-1)^3} \\
		& = \displaystyle\frac{\mathfrak{q}^{\xi_x-\gamma_x} (\log \mathfrak{q})^2}{(\mathfrak{q}^{x+\xi_x-\gamma_x}-1)^2} + \mathfrak{q}^{x+\xi_x-\gamma_x} \cdot \displaystyle\frac{\mathfrak{q}^{x+\xi_x-\gamma_x}+1}{(\mathfrak{q}^{x+\xi_x-\gamma_x}-1)^3} (\xi_x'-\gamma_x'+1) (\mathfrak{q}^{-x}-1) (\log \mathfrak{q})^2.
	\end{flalign*}
	
	\noindent Thus,
	\begin{flalign*}
		A \xi_x' & = \xi_x' \bigg( f_0'' (\xi_x) - \displaystyle\frac{\mathfrak{q}^{x+\xi_x-\gamma_x} (\mathfrak{q}^{x+\xi_x-\gamma_x}+1)}{(\mathfrak{q}^{x+\xi_x-\gamma_x}-1)^3} (\mathfrak{q}^{-x}-1) (\log \mathfrak{q})^2 \bigg) \\
		& = \displaystyle\frac{\mathfrak{q}^{\xi_x-\gamma_x} (\log \mathfrak{q})^2}{(\mathfrak{q}^{x+\xi_x-\gamma_x}-1)^2} + \displaystyle\frac{\mathfrak{q}^{x+\xi_x-\gamma_x} (\mathfrak{q}^{x+\xi_x-\gamma_x}+1) (1-\gamma_x')}{(\mathfrak{q}^{x+\xi_x-\gamma_x}-1)^3} (\mathfrak{q}^{-x}-1) (\log \mathfrak{q})^2 \\
		& = \displaystyle\frac{\mathfrak{q}^{\xi_x-\gamma_x} (\log \mathfrak{q})^2}{(\mathfrak{q}^{x+\xi_x-\gamma_x}-1)^2} - \displaystyle\frac{\mathfrak{q}^{\xi_x-\gamma_x}(\mathfrak{q}^{x+\xi_x-\gamma_x}+1)}{(\mathfrak{q}^{x+\xi_x-\gamma_x}-1)^2} (\log \mathfrak{q})^2 = -\displaystyle\frac{\mathfrak{q}^{x+2\xi_x-2\gamma_x} (\log \mathfrak{q})^2}{(\mathfrak{q}^{x+\xi_x-\gamma_x}-1)^2},
	\end{flalign*}
	
	\noindent where we used the fact from \eqref{gammaxderivative} that
	\begin{flalign*}
		1 - \gamma_x' = \displaystyle\frac{\mathfrak{q}^{-x} (1-\mathfrak{q}^{x+\xi_x-\gamma_x})}{\mathfrak{q}^{-x}-1}.
	\end{flalign*}
	 This yields the lemma.
	\end{proof} 
	
	Now we can establish \Cref{gammatderivative2}.

	\begin{proof}[Proof of \Cref{gammatderivative2}]
	
	By \Cref{gammatderivative}, \eqref{fxf0}, \Cref{fxiderivative}, the second statement of \eqref{f0xi}, we obtain
	\begin{flalign*}
		\gamma_x''  = \partial_x \Big( 1 + \displaystyle\frac{f_x (\gamma_x)}{f_x (\gamma_x) - 1} \Big) 	& = \partial_x \bigg( 1 + \displaystyle\frac{\mathfrak{q}^{-x}}{f_0 (\xi_x) - \mathfrak{q}^{-x}} \bigg) \\
		& = -\displaystyle\frac{\mathfrak{q}^{-x} f_0 (\xi_x) \log \mathfrak{q}}{ (f_0 (\xi_x) - \mathfrak{q}^{-x})^2} - \displaystyle\frac{\mathfrak{q}^{-x}  \xi_x' f_0' (\xi_x)}{( f_0 (\xi_x) - \mathfrak{q}^{-x} )^2} \\
		& = \displaystyle\frac{\mathfrak{q}^{x+3\xi_x-3\gamma_x} (1-\mathfrak{q}^{-x}) (\log \mathfrak{q})^3}{A (\mathfrak{q}^{x+\xi_x-\gamma_x} - 1)^4 ( f_0 (\xi_x) - \mathfrak{q}^{-x})^2} - \displaystyle\frac{f_x (\gamma_x) \log \mathfrak{q}}{( f_x (\gamma_x) - 1 )^2} \\
		& = \displaystyle\frac{\mathfrak{q}^{x+3\xi_x-3\gamma_x} ( f_0 (\xi_x) - \mathfrak{q}^{-x} )^2 (\log \mathfrak{q})^3}{A (1-\mathfrak{q}^{-x})^3} - \displaystyle\frac{f_x (\gamma_x) \log \mathfrak{q}}{( f_x (\gamma_x) - 1 )^2},
	\end{flalign*}

	\noindent from which \eqref{gammaderivative2} follows by the fact from \eqref{fxf0} that $f_0 (\xi_x) = \mathfrak{q}^{-x} f_x (\gamma_x)$. By \Cref{gammatderivative}, we have that $\Gamma_{x;1}= \gamma_x' (1 - \gamma_x') \log \mathfrak{q}$. Therefore, \eqref{gamma2gamma} follows from \eqref{gammaderivative2}.
\end{proof} 

\begin{lem} {If $z=(x,\gamma_x)$ with $x\in[x^{SW}_q,a]$, then
   $\xi_x\asymp \md_z$.}
 \label{xixd}
\end{lem}
\begin{proof}
   First recall the first of \eqref{f0xi}.
Since $x$ is bounded away from $0$, we have that 
\begin{eqnarray}
f_0(\xi_x)=\frac{{\mathfrak q}^{\xi_x-\gamma_x} - 1}{{\mathfrak q}^{x+\xi_x-\gamma_x} - 1} \asymp \xi_x-\gamma_x.  
\end{eqnarray}
Since   $\gamma_x \asymp \md_z^2$ {(because the arctic boundary has positive and bounded curvature)},
 it suffices to show that $f_0 (\xi_x) \asymp \md_z$. To that end, observe that 
\begin{eqnarray}
f_0 (\xi_x) = {\mathfrak q}^{-x} f_x (z_x (\xi_x)) = {\mathfrak q}^{-x} f_x (\gamma_x),  
\end{eqnarray}
where the first statement follows from \eqref{fxf0} and the second from Definition \ref{gammaxi}. Hence, it suffices to show that $f_x (\gamma_x) \asymp \md_z$. By Lemma \ref{gammatderivative} and $\gamma_x' \asymp \md_z$, it follows that $f_x (\gamma_x) \asymp \md_z$, showing the bound.
\end{proof}

\subsubsection{Arctic curve  at $q =0$}

\label{Functionsqabc1}

Here, we derive the values of the functions from the previous sections at $\mathfrak{q} = 1$. We have from \eqref{zxequation} (and the fact that $f_x (z_x)$ is constant in $x$ for $\mathfrak{q}=1$, by \eqref{fxf0}) that
\begin{flalign}
	\label{f01} 
	\gamma_x = \xi_x + x \gamma_x'.
\end{flalign} 

\noindent We also have from \Cref{gammatderivative} that 
\begin{flalign}
	\label{f021}  
  \gamma_x' = \displaystyle\frac{f_0 (\xi_x)}{f_0 (\xi_x)-1}=
\frac{(b+c-\xi_x)\xi_x}{c\xi_x+a(\xi_x-b)}.
\end{flalign}

\noindent By \eqref{f0xi}, we have
\begin{flalign*}
	f_0 (\xi_x) = \displaystyle\frac{\xi_x - \gamma_x}{x + \xi_x - \gamma_x}; \qquad f_0' (\xi_x) = \displaystyle\frac{x}{(x + \xi_x - \gamma_x)^2},
\end{flalign*}

\noindent from which it follows using \eqref{f01} and \eqref{f021} that 
\begin{multline}
	\label{xgammaxq1} 
	x = \displaystyle\frac{(f_0(\xi_x)-1)^2}{f_0' (\xi_x)};\\ \gamma_x = \xi_x + f_0 (\xi_x) (f_0 (\xi_x) - 1) f_0' (\xi_x)^{-1}=  \frac{c(a+b+c)\xi_x^2}{c \xi_x^2+a(b^2+b(c-2\xi_x)+\xi_x^2)}.
\end{multline}
{From the first line of \eqref{xgammaxq1}, together with the fact that
  \begin{eqnarray}
    \label{eq:fuq0}
    f_0(u)=\frac{u(u-b-c)}{(u+a)(u-b)}
  \end{eqnarray}
  for $\mathfrak q=1$, that follows from \eqref{f0uv}, we can express $\xi'_x$ in terms of $\xi_x$ itself, and then through \eqref{f021} we can compute the derivatives (of any order) of $\gamma_x$ with respect to $x$, as function of $\xi_x$. In particular, one finds the explicit expressions
  \begin{multline}
    \label{eq:od}
    \gamma''_x=\frac{(c\xi_x^2+a(b^2+b(c-2\xi_x)+\xi_x^2))^3}{2 abc (a+b+c)(a(b-\xi_x)-c\xi_x)^3},\\\gamma'''_x=\frac{3(c\xi_x^2+a(b^2+b(c-2\xi_x)+\xi_x^2))^4(c^2 \xi_x^2 + a^2 (b^2 + \xi_x^2 - b (c + 2 \xi_x)) - 
  a c (b^2 - 2 \xi_x^2 + b (c + 2 \xi_x)))}{4 (abc)^2 (a+b+c)^2(a(b-\xi_x)-c\xi_x)^5}.
  \end{multline}
}

\noindent The above holds for arbitrary $(a, b, c)$. Now let us set $(a, b, c) = (1, 1, 1)$. Then, we have 
\begin{flalign}
  x = \displaystyle\frac{(1 - 2 \xi_x)^2}{2 (\xi_x^2 - \xi_x+1)}; \qquad \gamma_x = \displaystyle\frac{3 \xi_x^2}{2 (\xi_x^2 - \xi_x + 1)}.
  \label{eq:xiesplicita}
\end{flalign}

\noindent This prescribes an arc of ellipse as we vary $\xi_x \in \mathbb{R} \cup \{ \infty \}$.

\subsubsection{Scaling of the complex slope near $\mathfrak A_{q;a,b,c}$}

\label{Slopede}

	The following lemma approximates the complex slope $f_x (y)$ for $(x,y)$ in the liquid region but $y$ close to $\gamma_x$. Underlying its proof is the fact that, while the function $f_x(y)$ is not differentiable in $y$ around $y=\gamma_x$ for general $x > 0$ (as is indeed implied by the below lemma), the function $f_0$ from \eqref{f0uv} is. Thus, we will use the fact from \eqref{fxf0} that $f_x (y) = \mathfrak{q}^x f_0 (u)$, for $u$ defined by setting $y = z_x (u)$, observing that $u \approx \xi_x$ as $y \approx \gamma_x$; Taylor expand $f_0$ around $\xi_x$; and use the results from \Cref{CurvatureDensity2} to evaluate the resulting coefficients.

	\begin{lem}
		
		\label{imaginaryftx0}
		{Let $C>1$ be a constant.}
		Letting $x \in [1/4, 5/4]$, the following two statements hold. 
		
		\begin{enumerate} 
			\item Recalling $\Gamma_{x;2}$ from \eqref{gammax12}, we have for $\gamma_x \le y$ that
		\begin{flalign}
			\label{fxyygammax} 
			f_x (y) = f_x (\gamma_x) + \mathrm{i} ( f_x (\gamma_x) - 1 )^2 \cdot ( 2 \Gamma_{x;2} (y - \gamma_x) )^{1/2} + (y-\gamma_x) \cdot \varsigma_1 ( x, (y-\gamma_x)^{1/2}),  
		\end{flalign}
		
		\noindent where $\varsigma_1: [1/4,5/4] \times \mathbb{R}_{\ge 0} \rightarrow \mathbb{C}$ is a function dependent on $(a,b,c,q)$ that satisfies the following. For any derivative $\partial_{\ell}$ in $(x,z,a,b,c,q)$, there exists a constant $C_{\ell} > 1$ such that $|\partial_{\ell} \varsigma_1 (x,z)| \le C_{\ell}$, for all $(x,z) \in [1/4,5/4] \times [0, C]$; $q \in [-\varepsilon, \varepsilon]$; and $a,b,c \in [1-\varepsilon, 1+ \varepsilon]$. 
	
		\item  The asymptotics can be refined close to $p^{SW}_{q;a,b,c}$. That is, assume that $\gamma_x \le y \le C \mathfrak{d}^2$ and denote
		\begin{flalign}\label{DAB}
			D' & = f_0' (0) = \displaystyle\frac{\log \mathfrak{q}}{\mathfrak{q}^x - 1}; \\ 
			A' & = f_0'' (0) + \displaystyle\frac{\mathfrak{q}^{x} + 1}{(\mathfrak{q}^{x} - 1)^2} (\log \mathfrak{q})^2; \\
			B' & = \displaystyle\frac{1}{3} \cdot \Bigg( f_0''' (0) + \displaystyle\frac{2 (\mathfrak{q}^{2x} + \mathfrak{q}^{x}+1)}{(\mathfrak{q}^{x}-1)^3} (\log \mathfrak{q})^3 +  3 f_0'' (0)   \cdot \displaystyle\frac{ (\mathfrak{q}^{x} + 1)}{\mathfrak{q}^{x} - 1} \log \mathfrak{q} \Bigg).
		\end{flalign}
		Then, we have 
		\begin{flalign}
			\begin{aligned} 
                  f_x (y) - f_x (\gamma_x)  & = \mathrm{i}  (f_x (\gamma_x) - 1)^2 \cdot ( 2 \Gamma_{x;2}  (y-\gamma_x) )^{1/2} \\ 
                 	& \qquad + \mathfrak{q}^x A'^{-1} D' (A'^{-1} B' D' - f_0'' (0)) \cdot (y - \gamma_x) + \mathfrak{d}^3 \cdot \varsigma_2 ( x, (y-\gamma_x)^{1/2}), 
                  \label{e:vm}
                 \end{aligned} 
		\end{flalign}
	
		\noindent where $\varsigma_2: \mathbb{R}_{\ge 0} \rightarrow \mathbb{C}$ is a function dependent on $(a,b,c,q)$ that satisfies the following. For any derivative $\partial_{\ell}$ in $(x,z,a,b,c,q)$, there exists a constant $C_{\ell} > 1$ such that $|\partial_{\ell} \varsigma_2 (x, z)| \le C_{\ell}$, for all $(x,z) \in [1/4,5/4] \times [0, C \mathfrak{d}]$; $q \in [-\varepsilon, \varepsilon]$; and $a,b,c \in [1-\varepsilon, 1+ \varepsilon]$. 
		
	\end{enumerate}
                
	\end{lem} 
		{In the second claim, one should have in mind the case of $\md$ small, so that $(x,\gamma_x)$ is indeed close to $x^{SW}_{q;a,b,c}$.}

	\begin{proof}		
		Throughout this proof, for any parameter $\varkappa$, we let $\mathcal{O}(\varkappa)$ denote a quantity of the form $\kappa \cdot \varsigma(x, (y-\gamma_x)^{1/2})$, where $\varsigma$ is a bounded function (that can depend on $\kappa$), that also smoothly depends on $(a,b,c,q)$, such that the following holds. For any derivative $\partial_{\gamma}$ in the variables $(x, z, a, b, c, q)$, there exists a constant $C_{\gamma} > 1$ such that $|\partial_{\gamma} \varsigma(z)| \le C_{\gamma}$.
		
		Now let $y = z_x (u)$. Since by \eqref{zxuf0} we have
		\begin{flalign*} 
			f_0 (u) = \mathfrak{q}^{-x} + \displaystyle\frac{\mathfrak{q}^{-x}-1}{\mathfrak{q}^{u+x-y} - 1},  
		\end{flalign*} 
		
		\noindent it follows that 
		\begin{flalign}
			\label{f02} 
			\begin{aligned} 
			f_0 & (\xi_x) + (u - \xi_x) \cdot f_0' (\xi_x) + \displaystyle\frac{(u - \xi_x)^2}{2} \cdot f_0'' (\xi_x) \\
			& = f_0 (u) + \mathcal{O} ( |u - \xi_x|^3 ) 
			= \mathfrak{q}^{-x} + \displaystyle\frac{\mathfrak{q}^{-x}-1}{\mathfrak{q}^{x+u-y} - 1} + \mathcal{O} ( |u - \xi_x|^3 ) \\
			& = \mathfrak{q}^{-x} + \displaystyle\frac{\mathfrak{q}^{-x}}{\mathfrak{q}^{x+\gamma_x-\xi_x} - 1} + (y - \gamma_x) \cdot \displaystyle\frac{\mathfrak{q}^{x + \xi_x - \gamma_x}}{(\mathfrak{q}^{x + \xi_x - \gamma_x} - 1)^2} \cdot (\mathfrak{q}^{-x} - 1) \log \mathfrak{q} \\
			& \qquad + (u - \xi_x) \cdot \displaystyle\frac{\mathfrak{q}^{x + \xi_x - \gamma_x}}{(\mathfrak{q}^{x + \xi_x - \gamma_x}-1)^2} (1 - \mathfrak{q}^{-x}) \log \mathfrak{q} \\
			& \qquad + (u - \xi_x)^2 \cdot \displaystyle\frac{\mathfrak{q}^{x+\xi_x - \gamma_x} (\mathfrak{q}^{x + \xi_x - \gamma_x} + 1)}{2 (\mathfrak{q}^{x + \xi_x - \gamma_x} - 1)^3} (\mathfrak{q}^{-x} - 1) (\log \mathfrak{q})^2 + \mathcal{O} ( |y-\gamma_x|^2 +  |u - \xi_x|^3),
			\end{aligned} 	
	\end{flalign} 
		
		\noindent where we used the facts that 
		\begin{flalign}
			& \partial_y (\mathfrak{q}^{x+u-y}-1)^{-1} |_{(u, y) = (\xi_x, \gamma_x)} = \mathfrak{q}^{x + \xi_x - \gamma_x} (\mathfrak{q}^{x + \xi_x - \gamma_x} - 1)^{-2} \log \mathfrak{q}; \\
			& \partial_u (\mathfrak{q}^{x+u-y}-1)^{-1} |_{(u,y) = (\xi_x, \gamma_x)} = -\mathfrak{q}^{x + \xi_x - \gamma_x} (\mathfrak{q}^{x + \xi_x - \gamma_x}-1)^{-2} \log \mathfrak{q}; \\ 
			& \partial_u^2 (\mathfrak{q}^{x + u - y}-1)^{-1} |_{(u,y) = (\xi_x, \gamma_x)} = \mathfrak{q}^{x + \xi_x - \gamma_x} (\mathfrak{q}^{x + \xi_x - \gamma_x}+1) (\mathfrak{q}^{x + \xi_x - \gamma_x}-1)^{-3} (\log \mathfrak{q})^2.
                          \label{e:vm2}
		\end{flalign}
		
		\noindent Thus, by \Cref{gammatderivative}, we have
		\begin{flalign*}
			(u - \xi_x)^2 &  \bigg( f_0'' (\xi_x) - \displaystyle\frac{\mathfrak{q}^{x + \xi_x - \gamma_x} (\mathfrak{q}^{x + \xi_x - \gamma_x} + 1)}{(\mathfrak{q}^{x + \xi_x - \gamma_x}-1)^3} (\mathfrak{q}^{-x}-1) (\log \mathfrak{q})^2 \bigg) \\
			& \qquad \qquad = \displaystyle\frac{2(\gamma_x-y) \mathfrak{q}^{x + \xi_x - \gamma_x}}{(\mathfrak{q}^{x + \xi_x - \gamma_x}-1)^2} (1 - \mathfrak{q}^{-x}) \log \mathfrak{q} + \mathcal{O} (|y - \gamma_x|^2 + |u - \xi_x|^3 ),
		\end{flalign*} 
		
		\noindent and so $|u-\xi_x| = \mathcal{O} ( |y-\gamma_x|^{1/2})$ and
		\begin{flalign*}
			u - \xi_x & = \mathrm{i} \Bigg( \displaystyle\frac{2 \mathfrak{q}^{x + \xi_x - \gamma_x} (y - \gamma_x) (1 - \mathfrak{q}^{-x})  \log \mathfrak{q}}{(\mathfrak{q}^{x + \xi_x - \gamma_x} - 1)^2} \\
			& \qquad \qquad \times \bigg( f_0'' (\xi_x) - \displaystyle\frac{\mathfrak{q}^{x+\xi_x-\gamma_x} (\mathfrak{q}^{x + \xi_x - \gamma_x}+1)}{(\mathfrak{q}^{x + \xi_x - \gamma_x}-1)^3}  (\mathfrak{q}^{-x}-1) (\log \mathfrak{q})^2 \bigg)^{-1} \Bigg)^{1/2} + \mathcal{O} (y-\gamma_x).
		\end{flalign*}
		
		\noindent Therefore, recalling the quantity $A = A_x$ from \eqref{af0}, we obtain
		\begin{flalign*}
			f_x (y) - f_x (\gamma_x) & = \mathfrak{q}^{x} \big( f_0 (u) - f_0 (\xi_x) \big) \\
			& = \mathfrak{q}^{x} f_0' (\xi_x) \cdot  (u - \xi_x) + \mathcal{O} ( |u - \xi_x|^2 ) \\
			& = \mathrm{i} \cdot \displaystyle\frac{\mathfrak{q}^{2x + \xi_x - \gamma_x} (1 - \mathfrak{q}^{-x}) \log \mathfrak{q}}{(\mathfrak{q}^{x + \xi_x - \gamma_x} - 1)^2} \Bigg( \displaystyle\frac{2\mathfrak{q}^{x + \xi_x - \gamma_x} (y-\gamma_x) (1 - \mathfrak{q}^{-x}) \log \mathfrak{q}}{A(\mathfrak{q}^{x+\xi_x - \gamma_x}-1)^2} \Bigg)^{1/2} + \mathcal{O} ( y - \gamma_x) \\
			& = \mathrm{i}  ( f_x (\gamma_x) - 1 )^3 \Bigg( \displaystyle\frac{2(y-\gamma_x) \mathfrak{q}^{3\xi_x - 3\gamma_x-x} (\log \mathfrak{q})^3}{(1-\mathfrak{q}^{-x})^3 A} \Bigg)^{1/2} + \mathcal{O} ( y - \gamma_x).
		\end{flalign*}
		
		\noindent where we used the fact that 
		\begin{flalign*}
			\displaystyle\frac{1}{(\mathfrak{q}^{x+\xi_x-\gamma_x}-1)^3} = \displaystyle\frac{(f_0 (\xi_x) - \mathfrak{q}^{-x})^3}{(\mathfrak{q}^{-x}-1)^3} = \mathfrak{q}^{-3x} \displaystyle\frac{(f_x (\gamma_x) - 1)^3}{(\mathfrak{q}^{-x}-1)^3}.
		\end{flalign*}
		
                \noindent We therefore obtain \eqref{fxyygammax} by recalling the definition \eqref{gammax12} of $\Gamma_{x;2}$.

                Next, we prove \eqref{e:vm}. Since from the argument above it follows that $|u-\xi_x|=O(\md)$, we can improve \eqref{f02} to
		\begin{flalign*}
			f_0 & (\xi_x) + (u - \xi_x) \cdot f_0' (\xi_x) + \displaystyle\frac{(u - \xi_x)^2}{2} \cdot f_0'' (\xi_x) + \displaystyle\frac{(u - \xi_x)^3}{6} \cdot f_0''' (\xi_x) \\
			& = f_0 (u) + \mathcal{O} ( \mathfrak{d}^4 ) 
			= \mathfrak{q}^{-x} + \displaystyle\frac{\mathfrak{q}^{-x}-1}{\mathfrak{q}^{x+u-y} - 1} + \mathcal{O} ( \mathfrak{d}^4 ) \\
			& = \mathfrak{q}^{-x} + \displaystyle\frac{\mathfrak{q}^{-x} - 1}{\mathfrak{q}^{x-\gamma_x+\xi_x} - 1} + (y - \gamma_x) \cdot \displaystyle\frac{\mathfrak{q}^{x + \xi_x - \gamma_x}}{(\mathfrak{q}^{x + \xi_x - \gamma_x} - 1)^2} \cdot (\mathfrak{q}^{-x} - 1) \log \mathfrak{q} \\
			& \qquad + (u - \xi_x) \cdot \displaystyle\frac{\mathfrak{q}^{x + \xi_x - \gamma_x}}{(\mathfrak{q}^{x + \xi_x - \gamma_x}-1)^2} (1 - \mathfrak{q}^{-x}) \log \mathfrak{q} \\
			& \qquad + (u - \xi_x)^2 \cdot \displaystyle\frac{\mathfrak{q}^{x+\xi_x - \gamma_x} (\mathfrak{q}^{x + \xi_x - \gamma_x} + 1)}{2 (\mathfrak{q}^{x + \xi_x - \gamma_x} - 1)^3} (\mathfrak{q}^{-x} - 1) (\log \mathfrak{q})^2 \\
			& \qquad + (u - \xi_x) (y - \gamma_x) \cdot \displaystyle\frac{\mathfrak{q}^{x + \xi_x - \gamma_x} (\mathfrak{q}^{x + \xi_x - \gamma_x} + 1)}{(\mathfrak{q}^{x + \xi_x - \gamma_x} - 1)^3}	(1 - \mathfrak{q}^{-x}) (\log \mathfrak{q})^2 \\
			& \qquad + (u - \xi_x)^3 \cdot \displaystyle\frac{\mathfrak{q}^{x + \xi_x - \gamma_x} (\mathfrak{q}^{2x + 2 \xi_x - 2\gamma_x} + 4 \mathfrak{q}^{x + \xi_x - \gamma_x}+1)}{6 (\mathfrak{q}^{x + \xi_x - \gamma_x}-1)^4} (1 - \mathfrak{q}^{-x}) (\log \mathfrak{q})^3  + \mathcal{O} ( \mathfrak{d}^4),
	\end{flalign*} 
		
		\noindent where we used \eqref{e:vm2} together with
		\begin{flalign*}
			& \partial_u \partial_y ( \mathfrak{q}^{x+u-y} - 1)^{-1} |_{(u,y) = (\xi_x, \gamma_x)} = -\mathfrak{q}^{x + \xi_x - \gamma_x} (\mathfrak{q}^{x + \xi_x - \gamma_x}+1) (\mathfrak{q}^{x + \xi_x - \gamma_x}-1)^{-3} (\log \mathfrak{q})^2; \\
			& \partial_u^3 (\mathfrak{q}^{x+u-y} - 1)^{-1} |_{(u,y) = (\xi_x, \gamma_x)} =  -\displaystyle\frac{\mathfrak{q}^{x + \xi_x - \gamma_x} (\mathfrak{q}^{2x + 2 \xi_x - 2\gamma_x} + 4 \mathfrak{q}^{x + \xi_x - \gamma_x}+1)}{(\mathfrak{q}^{x + \xi_x - \gamma_x}-1)^4} \cdot (\log \mathfrak{q})^3.
		\end{flalign*}
		\noindent Thus, by \Cref{gammatderivative0}, it follows that 
		\begin{flalign*}
			 (u - \xi_x & )^2 \cdot \Bigg( f_0'' (\xi_x) - \displaystyle\frac{\mathfrak{q}^{x+\xi_x - \gamma_x} (\mathfrak{q}^{x + \xi_x - \gamma_x} + 1)}{(\mathfrak{q}^{x + \xi_x - \gamma_x} - 1)^3} (\mathfrak{q}^{-x} - 1) (\log \mathfrak{q})^2 \Bigg)   + \displaystyle\frac{(u - \xi_x)^3}{3} \cdot f_0''' (0) \\
			& = 2 (\gamma_x - y) \cdot \displaystyle\frac{\mathfrak{q}^{x + \xi_x - \gamma_x}}{(\mathfrak{q}^{x + \xi_x - \gamma_x} - 1)^2} \cdot (1 - \mathfrak{q}^{-x}) \log \mathfrak{q} \\
			& \qquad + 2 (u - \xi_x) (y - \gamma_x) \cdot \displaystyle\frac{\mathfrak{q}^{x + \xi_x - \gamma_x} (\mathfrak{q}^{x + \xi_x - \gamma_x} + 1)}{(\mathfrak{q}^{x + \xi_x - \gamma_x} - 1)^3}	(1 - \mathfrak{q}^{-x}) (\log \mathfrak{q})^2 \\
			& \qquad + (u - \xi_x)^3 \cdot \displaystyle\frac{\mathfrak{q}^{x + \xi_x - \gamma_x} (\mathfrak{q}^{2x + 2 \xi_x - 2\gamma_x} + 4 \mathfrak{q}^{x + \xi_x - \gamma_x}+1)}{3 (\mathfrak{q}^{x + \xi_x - \gamma_x}-1)^4} (1 - \mathfrak{q}^{-x}) (\log \mathfrak{q})^3  + \mathcal{O} ( \mathfrak{d}^4) \\
			& = 2 (\gamma_x - y) \cdot \displaystyle\frac{\mathfrak{q}^{x + \xi_x - \gamma_x}}{(\mathfrak{q}^{x + \xi_x - \gamma_x} - 1)^2} \cdot (1 - \mathfrak{q}^{-x}) \log \mathfrak{q} \\
			& \qquad + 2 (u - \xi_x) (y - \gamma_x) \cdot \displaystyle\frac{ (\mathfrak{q}^{x } + 1)}{(\mathfrak{q}^{x } - 1)^2} (\log \mathfrak{q})^2 \\
			& \qquad + (u - \xi_x)^3 \cdot \displaystyle\frac{(\mathfrak{q}^{2x } + 4 \mathfrak{q}^{x}+1)}{3 (\mathfrak{q}^{x}-1)^3} (\log \mathfrak{q})^3  + \mathcal{O} ( \mathfrak{d}^4).
		\end{flalign*} 
We have replaced $\xi_x,\gamma_x$ with $0$ in the terms proportional to $(u-\xi_x)^3$ and $(u-\xi_x)^2(y-\gamma_x)$, that are $O(\md^3)$. In fact, recall from Lemma \ref{xixd} that $\xi_x\asymp \md$, and that $\gamma_x\asymp \md^2$. 		
		\noindent So, we have
		\begin{flalign*}
			\gamma_x - y = (u - \xi_x)^2 \cdot \bigg( f_0'' (0) \cdot \displaystyle\frac{\mathfrak{q}^{x} - 1}{2 \log \mathfrak{q}} + \displaystyle\frac{(\mathfrak{q}^{x} + 1)}{2 (\mathfrak{q}^{x}-1)}  \log \mathfrak{q} \bigg) + \mathcal{O} (\mathfrak{d}^3),
		\end{flalign*} 
		\noindent meaning that 
			\begin{flalign*}
			& A (u - \xi_x )^2  +B  (u - \xi_x)^3   = 2 D (\gamma_x - y)  + \mathcal{O} ( \mathfrak{d}^4).
		\end{flalign*} 
		
		\noindent where we have denoted 
		\begin{flalign*}
			A & = f_0'' (\xi_x) - \displaystyle\frac{\mathfrak{q}^{x+\xi_x - \gamma_x} (\mathfrak{q}^{x + \xi_x - \gamma_x} + 1)}{(\mathfrak{q}^{x + \xi_x - \gamma_x} - 1)^3} (\mathfrak{q}^{-x} - 1) (\log \mathfrak{q})^2; \\
			B & = \displaystyle\frac{1}{3} \cdot \Bigg( f_0''' (0) + \displaystyle\frac{2(\mathfrak{q}^{2x} +  \mathfrak{q}^{x}+1)}{(\mathfrak{q}^{x}-1)^3} (\log \mathfrak{q})^3 +  3 f_0'' (0)   \cdot \displaystyle\frac{ (\mathfrak{q}^{x
                            } + 1)}{\mathfrak{q}^{x} - 1} \log \mathfrak{q} \Bigg); \\
			D & = f_0' (\xi_x) = \displaystyle\frac{\mathfrak{q}^{x + \xi_x - \gamma_x}}{(\mathfrak{q}^{x + \xi_x - \gamma_x} - 1)^2} \cdot (1 - \mathfrak{q}^{-x}) \log \mathfrak{q} = (f_x (\gamma_x) - 1)^2 \cdot \displaystyle\frac{\mathfrak{q}^{\xi_x - \gamma_x - x}}{1 - \mathfrak{q}^{-x}} \log \mathfrak{q}.
		\end{flalign*}
		\noindent So, $|u-\xi_x| = \mathcal{O} ( |y-\gamma_x|^{1/2})$ and
		\begin{flalign*}
			u - \xi_x & = \mathrm{i} (2A^{-1} D (y - \gamma_x))^{1/2} + A^{-2} BD (y - \gamma_x) + \mathcal{O} (\mathfrak{d}^3).
		\end{flalign*}
		\noindent Therefore, recalling the quantity $A = A_x$ from \eqref{af0} and $\Gamma_{x;2}$ from \eqref{gammax12}, and using the identity ${\mathfrak q}^xf'_0(\xi_x)\sqrt{2D/A}=
                  (f_x(\gamma_x)-1)^2\sqrt{2\Gamma_{x;2}}$ that follows from
                  Lemma \ref{gammatderivative0} and \eqref{fxf0},
                   we obtain
		\begin{flalign}
			f_x (y) - f_x (\gamma_x) & = \mathfrak{q}^{x} \big( f_0 (u) - f_0 (\xi_x) \big) \\
			& = \mathfrak{q}^{x} f_0' (\xi_x) \cdot   (u - \xi_x)
                          - A^{-1} D \mathfrak{q}^x f_0'' (\xi_x) \cdot  (y - \gamma_x) +  \mathcal{O} ( \mathfrak{d}^3 ) \\
                                                 & =  \mathrm{i}  (f_x (\gamma_x) - 1)^2 \cdot ( 2 \Gamma_{x;2}  (y-\gamma_x) )^{1/2}+{\mathfrak q}^x\frac{D'}{A'}(y-\gamma_x)\Big(-f''_0(0)+\frac{D' B'}{A'}\Big)+  \mathcal{O} ( \mathfrak{d}^3 ),
                                                   \label{alittledifferent}
		\end{flalign}
		\noindent where $A',B',D'$ are as in the statement of the lemma.
	\end{proof}

        From \Cref{imaginaryftx0} we can get an approximated expression for the complex slope, in the rescaled coordinates $\hat x,\hat y$ of \eqref{e:funzrescaled}.
        In particular, \Cref{cor:cella1} is relevant when we look at points $z$ in the liquid region with  $\me_{z;q}\ll \md_z$, while \Cref{cor:cella2} is relevant when $\me_{z;q}\asymp\md_z\ll1$, which is the case very close to a tangency location.

		\begin{rem}
			\label{rem:whatissmooth}
			Throughout the remainder of this section, for any parameter $\varkappa$, we let $\mathcal{O}(\varkappa)$ denote a smooth function $\varsigma$, often in two arguments\footnote{Observe that our convention here is slightly different from that in the proof of \Cref{imaginaryftx0}, where $\varsigma$ was a function of $x$ and $(y-\gamma_x)^{1/2}$.} $(\hat{x}, \hat{y})$ (on some specified compact domain) and also dependent on $(a, b, c, q)$, such that $|\partial_{\ell} \varsigma(\hat{x},\hat{y})| \le C_{\ell} \varkappa$, where $\partial_{\ell}$ denotes a derivative (of any order) of $\varsigma$ in the variables $(\hat{x}, \hat{y}, a, b, c, q)$ (and $C_{\ell}$ is a finite constant depending only on $\ell$).
		\end{rem}

\begin{cor} 
	\label{cor:cella1}
    Let $\md,\me\lesssim1$ be positive numbers. There exists a constant
    $C>1$ such that the following holds. 
	Fix $(x', y') = (x', \gamma_{x'})$ on the arctic boundary $\mathfrak A_{q;a,b,c}$ with $x' \in [1/4, 5/4]$, and denote $f = f_{x'} (\gamma_{x'})$; $\gamma = \gamma_{x'}$; $\gamma' = \gamma_{x'}'$; and $\gamma'' = \gamma_{x'}''$. 
	 Let $z=(x,y)$ with 
	\begin{flalign}
		\label{xyxyde} 
		x = x' + (\mathfrak{d} \mathfrak{e})^{1/2} \hat{x} - \gamma' \mathfrak{d} \mathfrak{e} \hat{y}; \qquad y = \gamma_{x'} + \gamma' (\mathfrak{d} \mathfrak{e})^{1/2} \hat{x} + \mathfrak{d} \mathfrak{e} \hat{y}
	\end{flalign}
        and assume that
        \begin{equation}
          \label{eq:megaassunzione}
       \frac{ \md}{|\gamma'|},\frac{\me_{z;q}}{\me},\frac{\md_{z}}\md\in   [C^{-1},C].
        \end{equation}
        
	\noindent Denoting  
	\begin{flalign}
		\label{fxy2} 
		\hat{f}(\hat{x}, \hat{y}) = (\mathfrak{d} \mathfrak{e})^{-1/2} ( f(x,y) - f); \qquad \Gamma =  \gamma''- \gamma'(1 - \gamma') \log \mathfrak{q}, 
	\end{flalign}
	\noindent we have 
	\begin{flalign}
		\label{fapproximate} 
		\hat{f} (\hat{x}, \hat{y}) =  \displaystyle\frac{\mathrm{i} \Gamma^{1/2}}{(\gamma'-1)^2}   \big( 2 (\gamma'^2 + 1) \hat{y} - \gamma'' \hat{x}^2 \big)^{1/2}  - \displaystyle\frac{\gamma'' \hat{x}}{(\gamma'-1)^2} + \mathcal{O} ( (\mathfrak{d} \mathfrak{e})^{1/2}).
	\end{flalign}
      \end{cor}
      \begin{rem}
        \label{rem:sparrow}       Note that the assumption \eqref{eq:megaassunzione} implies in particular that the argument of the square root in \eqref{fapproximate} is positive and of order $1$.  
      \end{rem}

\begin{proof} 
  By \eqref{gamma2gamma}, we have that $\Gamma = \Gamma_{x';2}$. By \Cref{imaginaryftx0} and  \Cref{gammatderivative} it follows
    \begin{eqnarray}
    f_x(\gamma_x)-1=(\gamma' - 1)^{-1}-\frac1{(\gamma' - 1)^2}\gamma''\hat x\sqrt{\md \me}+\mathcal{O}(\md \me)\\=f-1 -\frac1{(\gamma' - 1)^2}\gamma''\hat x\sqrt{\md \me}+\mathcal{O}(\md \me) 
    \end{eqnarray}
    and $\Gamma_{x';2} = \Gamma_{x;2} + \mathcal{O} ( (\mathfrak{d} \mathfrak{e})^{1/2})$. Therefore,
\begin{flalign*}
	f(x, y) = f + \displaystyle\frac{(\mathfrak{d} \mathfrak{e})^{1/2}}{(\gamma'-1)^2} \Bigg( \mathrm{i} \cdot \bigg( 2 \Gamma \Big( (\gamma'^2 + 1) \hat{y} - \displaystyle\frac{ \gamma'' \hat{x}^2}{2} \Big) \bigg)^{1/2} - \gamma'' \hat{x} \Bigg) 
	+ \mathcal{O} (\mathfrak{d} \mathfrak{e})
\end{flalign*}
\noindent where we used the fact that 
\begin{flalign*}
	y - \gamma_x = \gamma_{x'} - \gamma_x + \gamma' (\mathfrak{d} \mathfrak{e})^{1/2} \hat{x} + \mathfrak{d} \mathfrak{e} \hat{y} & = \gamma' \big(x' - x + (\mathfrak{d} \mathfrak{e})^{1/2} \hat{x} \big) - \mathfrak{d} \mathfrak{e} \cdot  \displaystyle\frac{\gamma'' \hat{x}^2}{2} + \mathfrak{d} \mathfrak{e} \hat{y} + \mathcal{O} ( (\mathfrak{d} \mathfrak{e})^{3/2} ) \\
	& = \displaystyle\frac{ \mathfrak{d} \mathfrak{e}}{2} \cdot \big(  2\hat{y} (\gamma'^2 + 1) - \gamma'' \hat{x}^2 \big) + \mathcal{O} \big( (\mathfrak{d} \mathfrak{e})^{3/2} ). 
\end{flalign*}

\noindent Thus, the lemma follows from the definition \eqref{fxy2} of $\hat{f}$.
\end{proof}

\begin{cor} 
\label{cor:cella2}	
Let $C>1$ be a fixed constant. 
	Denote $\gamma'' = \gamma_{x^{SW}_{q;a,b,c}}'',\gamma'''=\gamma_{x^{SW}_{q;a,b,c}}'''$, and let
	\begin{flalign}
		\label{xyxydd} 
		x = x^{SW}_{q;a,b,c} + \mathfrak{d} \hat{x}; \qquad y =  \mathfrak{d}^2 \hat{y}.
	\end{flalign}	
Assume that $z=(x,y)\in \mathfrak L_{q;a,b,c}$ with {$\me_{z;q}\ge (1/C) \md, \md_z\asymp \md$.}	
	\noindent Denoting  
	\begin{flalign}
		\label{fxy2dd} 
		\hat{f}(\hat{x}, \hat{y}) = \mathfrak{d}^{-1} \cdot f(x,y); \qquad \Gamma_{\hat{x}} = \gamma''+\md\hat x(\gamma'''-\gamma''\log {\mathfrak q}+4 \gamma''^2).          
	\end{flalign}

	\noindent we have 
	\begin{flalign}
		\hat{f} (\hat{x}, \hat{y}) & =  \mathrm{i} \Gamma_{\hat{x}}^{1/2}  \Big( 2 \hat{y} - \gamma'' \hat{x}^2 - \displaystyle\frac{\gamma'''}{3} \mathfrak{d} \hat{x}^3 \Big)^{1/2}  - \gamma'' \hat{x}  - \mathfrak{d} \Big( \gamma''^2 + \displaystyle\frac{\gamma'''}{2} \Big) \hat{x}^2         + \mathfrak{d} g  \Big(\hat{y}- \gamma''\displaystyle\frac{\hat{x}^2}{2} \Big)  + \mathcal{O} ( \mathfrak{d}^2),\label{fhatcell2}
	\end{flalign}
        where $A',B',D'$ are as in \eqref{DAB} with $x=x^{SW}_{q;a,b,c}$ and
        \begin{eqnarray}
          \label{eq:theg}
    g=g(q,a,b,c):= \mathfrak{q}^{x_{q;a,b,c}^{SW}} A'^{-1} D' (D' A'^{-1} B' - f_0'' (0)).
  \end{eqnarray}
\end{cor} 
\begin{rem}
  \label{rem:nosubtract} Note that, with respect to \Cref{cor:cella1}, in the definition of $\hat f$ we do not subtract $f=f_{x'}(\gamma_{x'})$ since in the present case we have $x'=x_{q;a,b,c}^{SW}$, where $f_{x'}(\gamma_{x'})=0$. Also, the assumption $\me_{z;q}\ge (1/C) \md$ implies that imaginary part of $\hat f$ is positive and bounded away from zero {and that $\hat y$ is bounded away from zero and infinity}.
\end{rem}

\begin{proof} 
By \Cref{gammatderivative} it follows
   \begin{eqnarray}
    f_x(\gamma_x)-1=-1-\gamma''\hat x\md -(\gamma''^2+\gamma'''/2)\md^2\hat x^2+\mathcal{O}(\md^3).
   \end{eqnarray}
   We also have by \eqref{gamma2gamma}
   \begin{equation}
     \Gamma_{x;2}=\gamma''+\hat x \md( \gamma'''- \gamma'' \log {\mathfrak q}) +\mathcal{O}(\md^2)
   \end{equation}
   as well as
   \begin{eqnarray}
     y-\gamma_x=\md^2\Big(\hat y-\frac{\gamma''}2\hat x^2-\md\frac{\gamma'''}6 \hat x^3\Big)+\mathcal{O}(\md^4).
   \end{eqnarray}
The results then follows from \Cref{imaginaryftx0} {with $x=x_{q;a,b,c}^{SW}$}
  and  the definition \eqref{fxy2} of $\hat{f}$.
\end{proof}

The following two lemmas use \Cref{cor:cella1} and \Cref{cor:cella2} to approximate derivatives of the (rescaled) height function in terms of the (rescaled) complex slope $\hat{f}$ from \eqref{fxy2}.
\begin{lem} 
	\label{fh} 
	Adopt the notation and assumptions of \Cref{cor:cella1}, and denote
	\begin{flalign}\label{twocases}
          \widetilde{\mathcal{H}} (\hat{x}, \hat{y}) = \mathfrak{d}^{-1/2} \mathfrak{e}^{-3/2} \mathcal{H}_{q;a,b,c} (x, y) \quad \text{if}\quad \gamma'>0,\\
          \widetilde{\mathcal{H}} (\hat{x}, \hat{y}) = \mathfrak{d}^{-1/2} \mathfrak{e}^{-3/2}(y- \mathcal{H}_{q;a,b,c} (x, y)) \quad \text{if}\quad \gamma'<0. 
	\end{flalign}
	\noindent Then,  we have
	\begin{flalign}
		\label{derivativeshf}
		\begin{aligned} 
		& \partial_{\hat{x}} \widetilde{\mathcal{H}} (\hat{x}, \hat{y}) = \displaystyle\frac{\mathfrak{d} (1-\gamma')^3}{\pi |\gamma'|} \cdot \Real \hat{f}(\hat{x}, \hat{y}) \cdot \Imaginary \hat{f} (\hat{x}, \hat{y}) + \mathcal{O} ( (\mathfrak{d}^{-1} \mathfrak{e})^{1/2}); \\  
		& \partial_{\hat{y}} \widetilde{\mathcal{H}} (\hat{x}, \hat{y}) = \displaystyle\frac{\mathfrak{d} (1-\gamma') (\gamma'^2+1)}{\pi |\gamma'|} \cdot \Imaginary \hat{f} (\hat{x}, \hat{y}) + \mathcal{O} ( (\mathfrak{d}^{-1} \mathfrak{e})^{1/2} ).
		\end{aligned} 
	\end{flalign}
\end{lem}

\begin{proof} 
  We give the proof in the case $\gamma'>0$, the other case being analogous (one just needs to replace \eqref{eq:argv} below with $\arg v = |v|^{-1} \Imaginary v + \mathcal O ( (\Imaginary v)^3 / |v|^3)$ for $\Real v > 0$). Throughout, we abbreviate $\hat{f} = \hat{f}(\hat{x},\hat{y})$. {Recall from \Cref{cor:cella1} that $f:=f_{x'}(\gamma_{x'})$ and in particular from \Cref{rem:gamma'} that
    \begin{eqnarray}
      \label{eq:statementsf}
      f\asymp -\md<0, 1-f\asymp 1.
    \end{eqnarray}
  }
	We begin by showing the second statement in \eqref{derivativeshf}. To that end, by \eqref{xyxyde}, we have 
	\begin{flalign}
		\label{yh}
		\partial_{\hat{y}} \widetilde{\mathcal{H}} (\hat{x}, \hat{y}) =  (\mathfrak{d} \mathfrak{e}^{-1})^{1/2} \cdot (  \partial_y \mathcal{H}_{q;a,b,c} (x, y) - \gamma' \partial_x \mathcal{H}_{q;a,b,c} (x, y)). 
	\end{flalign}
	
	\noindent Using \eqref{fxygradient} and \eqref{fxy2}, so that
        \begin{eqnarray}
          \label{eq:statementsfxy}
        \Real(f(x,y))<0,\qquad\Real(f(x,y)-1)<0 ,
        \end{eqnarray}
         we obtain that 
	\begin{flalign}
		\label{xhyh1} 
		\begin{aligned} 
			\partial_x \mathcal{H}_{q;a,b,c} (x, y) & = \pi^{-1} \big( \arg (f(x,y) - 1) - \pi \big) 
			 = -\displaystyle\frac{\Imaginary f (x,y) }{\pi |f(x,y)-1|}  + {\mathcal{O} ( (\mathfrak{d} \mathfrak{e})^{3/2})} \\
			& = (\mathfrak{d} \mathfrak{e})^{1/2} \cdot \displaystyle\frac{\Imaginary \hat{f}}{\pi (f-1)^2} \cdot (f-1-(\mathfrak{d} \mathfrak{e})^{1/2} \Real \hat{f} ) + {\mathcal{O} ( (\mathfrak{d} \mathfrak{e})^{3/2})},
		\end{aligned} 
	\end{flalign} 
	
	\noindent and 
	\begin{flalign} 
		\label{xhyh2} 
		\begin{aligned} 
			\partial_y \mathcal{H}_{q;a,b,c} (x, y) & = \pi^{-1} (\pi - \arg f(x, y)) \\
			& = \pi^{-1} \sin^{-1} \bigg( (\mathfrak{d} \mathfrak{e})^{1/2} \cdot \displaystyle\frac{\Imaginary \hat{f}}{|f(x,y)|} \bigg) \\
			& =  \pi^{-1} \sin^{-1} \big( (\mathfrak{d} \mathfrak{e} f^{-2} \Real \hat{f} -  (\mathfrak{d}\mathfrak{e})^{1/2} f^{-1}) \cdot \Imaginary \hat{f} \big) + \mathcal{O} ( \mathfrak{d}^{-1/2} \mathfrak{e}^{3/2}).
		\end{aligned}
	\end{flalign}
	
	\noindent Here, in the first statement of \eqref{xhyh1}, we used \eqref{fxygradient}; in the second we used \eqref{eq:statementsf} and \eqref{eq:statementsfxy} and the fact that, for any $v \in \overline{\mathbb{H}}$ with $\Real v<0$, 
	\begin{flalign} 
		\arg v = \pi-|v|^{-1} \cdot \Imaginary v + \mathcal{O} ((\Imaginary v/|v|)^3);
		\label{eq:argv}
	\end{flalign} 
	
	\noindent and in the third we used the definition \eqref{fxy2} of $\hat{f}$, with the fact that 
	\begin{flalign*}
		|f(x,y) - 1|^{-1} & = \big| 1-f-(\mathfrak{d}\mathfrak{e})^{1/2} \Real \hat{f} - \mathrm{i} (\mathfrak{d} \mathfrak{e})^{1/2} \Imaginary \hat{f} \big|^{-1} \\ 
		& = \big( 1-f-(\mathfrak{d}\mathfrak{e})^{1/2} \Real \hat{f} \big)^{-1}	+ \mathcal{O}(\mathfrak{d}\mathfrak{e}) \\
		& = (f-1)^{-2} \cdot \big(1 - f + (\mathfrak{d} \mathfrak{e})^{1/2} \Real \hat{f} \big) + \mathcal{O} (\mathfrak{d} \mathfrak{e}),
	\end{flalign*}
	
	\noindent again by \eqref{eq:statementsf}. In the first statement of \eqref{xhyh2}, we used \eqref{fxygradient}; in the second we used \eqref{fxy2} and \eqref{eq:statementsf}, \eqref{eq:statementsfxy}; and in the third we used the facts that $-f(x,y) \asymp \mathfrak{d}$ and that 
	\begin{flalign*}
		|f(x,y)|^{-1} = \big| f + (\mathfrak{d} \mathfrak{e})^{1/2} \Real \hat{f} + \mathrm{i} (\mathfrak{d}\mathfrak{e})^{1/2} \Imaginary \hat{f} \big|^{-1} & = \big| f + (\mathfrak{d} \mathfrak{e})^{1/2} \Real \hat{f} \big|^{-1} + \mathcal{O} (\mathfrak{d}^{-1} \mathfrak{e}) \\
		& = f^{-2} \cdot ((\mathfrak{d} \mathfrak{e})^{1/2} \Real \hat{f} - f) + \mathcal{O}(\mathfrak{d}^{-1} \mathfrak{e}).
	\end{flalign*}

	\noindent Inserting \eqref{xhyh1} and \eqref{xhyh2} into \eqref{yh}, using \Cref{gammatderivative}, denoting 
	\begin{flalign}
		\label{fde} 
		\Upsilon = (\mathfrak{d} \mathfrak{e} f^{-2} \Real \hat{f} -  (\mathfrak{d}\mathfrak{e})^{1/2} f^{-1}) \cdot \Imaginary \hat{f},
	\end{flalign}
	and since by \Cref{gammatderivative} we have 
	\begin{flalign*}
		- \displaystyle\frac{\mathfrak{d}}{f} \cdot \displaystyle\frac{f^2+(f-1)^2}{\pi (f-1)^2}  \cdot  \Imaginary \hat{f} = \displaystyle\frac{\mathfrak{d} (1-\gamma') (\gamma'^2+1)}{\pi \gamma'} \cdot \Imaginary \hat{f},\quad \frac{f^3 + (f-1)^3}{f^2(f-1)^3} =(1+\gamma'^3)\frac{(1-\gamma')^2}{\gamma'^2}
	\end{flalign*}
	\noindent we obtain   
	\begin{flalign}
		\label{hy0} 
		\begin{aligned} 
                  \partial_{\hat{y}} \widetilde{\mathcal{H}} & = \displaystyle\sqrt{\frac\md\me} \frac{ \sin^{-1} (\Upsilon) - \Upsilon }\pi+   \displaystyle\frac{\mathfrak{d} (1-\gamma') (\gamma'^2+1)}{\pi \gamma'} \cdot \Imaginary \hat{f}  \\
 & \qquad + \displaystyle\frac{\mathfrak{d}^2}{\pi} \cdot \displaystyle\sqrt{\frac{\mathfrak{e}}{\mathfrak{d}}} \cdot(1+\gamma'^3)\frac{(1-\gamma')^2}{\gamma'^2}  \cdot \Real \hat{f} \cdot \Imaginary \hat{f}+ \mathcal{O} ( (\mathfrak{d} \mathfrak{e})^{1/2} ),         
		\end{aligned} 
	\end{flalign}
	
	\noindent where we have abbreviate $\partial_{\hat{y}} \widetilde{\mathcal{H}} (\hat{x},\hat{y}) = \partial_{\hat{y}} \widetilde{\mathcal{H}}$. 
	We then deduce the second bound in \eqref{derivativeshf} by combining \eqref{hy0} with 
        {$\sin^{-1}(\Upsilon)-\Upsilon=O(\Upsilon^3)$}
        (since $|\Upsilon| = \mathcal{O} ( (\mathfrak{e}/\mathfrak{d})^{1/2})$, as $|f| \asymp \mathfrak{d}$). 
	
	To address the $\hat{x}$ derivative, observe again by \eqref{xyxyde} that 
	\begin{flalign} 
		\label{xh} 
		\partial_{\hat{x}} \widetilde{\mathcal{H}} (\hat{x}, \hat{y}) =  \mathfrak{e}^{-1} \cdot (  \partial_x \mathcal{H}_q (x, y) + \gamma' \partial_y \mathcal{H}_q (x, y)).
	\end{flalign} 
	
	\noindent Thus, inserting \eqref{xhyh1} and \eqref{xhyh2} into \eqref{xh}, using \Cref{gammatderivative}, and recalling \eqref{fde}, we obtain
	\begin{flalign}
		\label{hx0}
		\begin{aligned}  
			\partial_{\hat{x}} \widetilde{\mathcal{H}} & =\displaystyle\frac{\gamma'}{\pi \mathfrak{e}} \big(  \sin^{-1} (\Upsilon) - \Upsilon \big) + \displaystyle\frac{\mathfrak{d} (1-\gamma')^3}{\pi \gamma'} \cdot \Real \hat{f} \cdot  \Imaginary \hat{f} + \mathcal{O} ( (\mathfrak{d} \mathfrak{e})^{1/2}),
		\end{aligned}
	\end{flalign}
	
	\noindent where we have abbreviated $\partial_{\hat{x}} \widetilde{\mathcal{H}} (\hat{x},\hat{y}) = \partial_{\hat{x}} \widetilde{\mathcal{H}}$. We then deduce the first bound in \eqref{derivativeshf} by combining \eqref{hx0} with  {$\sin^{-1}(\Upsilon)-\Upsilon=O(\Upsilon^3)$}
        (since $|\Upsilon| = \mathcal{O} ( (\mathfrak{e}/\mathfrak{d})^{1/2})$, as $|f| \asymp \mathfrak{d}$).
\end{proof}

\begin{lem}
  \label{lem:cella2H}
  Adopt the notation and assumptions of \Cref{cor:cella2}, and denote
  	\begin{flalign}\label{eq:HtildeH}
		\widetilde{\mathcal{H}} (\hat{x}, \hat{y}) = \mathfrak{d}^{-2} \mathcal{H}_{q;a,b,c} (x, y). 
	\end{flalign}
	\noindent Then, we have $\widetilde\cH(\hat{x},\hat{y})=\frak H_{q,\gamma'',\gamma'''}(\hat x,\hat y)+\mathcal O(\md^2)$, where
  {  \begin{multline}\label{approxG2}
     \frak H_{q,\gamma'',\gamma'''}(\hat x,\hat y)=\mathtt H_{\gamma''}(\hat x,\hat y)        + \frac{\mathfrak d}{\pi\sqrt{\gamma''}}
        \left[
            \left(\frac{\gamma''}{3}-\frac{q}{6}+\frac{2\gamma'''}{9\gamma''}\right) (2\hat{y}-\gamma''\hat{x}^2)^{3/2}
            - \frac{\gamma'''}{3\gamma''} \hat{y}\sqrt{2\hat{y}-\gamma''\hat{x}^2}
        \right],\\
\mathtt H_{\gamma''}(\hat x,\hat y)
        = \displaystyle\frac{\hat{y}}{2}
        - \frac{1}{\pi}\left[
            \frac{\hat{x}}{2}\sqrt{\gamma''}\,\sqrt{2\hat{y}-\gamma''\hat{x}^2}
            + \hat{y}\,\arctan\!\left(\frac{\sqrt{\gamma''} \hat{x}}{\sqrt{2\hat{y}-\gamma''\hat{x}^2}}\right)
        \right].
    \end{multline}
}
\end{lem}
Some useful facts about the functions $\frak H_{q,\gamma'',\gamma'''}(\hat x,\hat y)$ and $\mathtt H_{\gamma''}$ are collected in \Cref{sec:AI}.
\begin{proof}
 It follows immediately from
  \eqref{xyxydd}, \eqref{fxy2dd} and \eqref{fxygradient} that
	\begin{flalign}
		\label{derivativeshfcell2}
		\begin{aligned} 
	&	 \partial_{\hat{x}} \widetilde{\mathcal{H}} (\hat{x}, \hat{y}) =  \frac{1}{\mathfrak{d}} \left(
  - 1 + \frac{1}{\pi} \arg (\mathfrak{d} \hat{f} (\hat{x}, \hat{y}) - 1)
  \right)= - \frac{1}{\pi} \Imaginary
  (\hat{f})[ 1+\md \Real                                                                                    (\hat{f})] + {\mathcal O (\mathfrak{d}^2)}\\
                  &= - \frac{1}{\pi} \Imaginary
  (\hat{f})[ 1-\md \gamma''\hat x] + {\mathcal O (\mathfrak{d}^2)} \\  
		 &\partial_{\hat{y}} \widetilde{\mathcal{H}} (\hat{x}, \hat{y}) =	1 - \frac{1}{\pi} \arg
  (\hat{f} (\hat{x}, \hat{y}))=1-\frac1\pi\arctan\Big(\frac{\Imaginary(\hat f)}{\Real(\hat f)}\Big).	\end{aligned} 
	\end{flalign}
        From \eqref{fhatcell2} and \eqref{derivativeshfcell2} one deduces, letting $S=2\hat y-\gamma''\hat x^2$,
          \begin{flalign*}
        \partial_{\hat{x}} \widetilde\cH(\hat{x},\hat{y})
        = -\frac{1}{\pi}\sqrt{\gamma''}\,\sqrt{S}
        - \frac{\mathfrak d}{\pi}
        \left[
            \sqrt{\gamma''}\,\hat{x}\left(\gamma''-\frac q2+\frac{\gamma'''}{2\gamma''}\right)\sqrt{S}
            - \frac{\gamma'''\sqrt{\gamma''}}{6}\frac{\hat{x}^3}{\sqrt{S}}
        \right] + \mathcal{O}(\mathfrak d^2).
          \end{flalign*}
            Integrating with respect to $\hat{x}$ gives
    \begin{flalign*}
        \widetilde \cH(\hat{x},\hat{y})
        &= C(\hat{y})- \displaystyle\frac{\hat{y}}{2}
        - \frac{1}{\pi}\left[
            \frac{\hat{x}}{2}\sqrt{\gamma''}\,\sqrt{S}
            + \hat{y}\,\arctan\!\left(\frac{\sqrt{\gamma''} \hat{x}}{\sqrt{S}}\right)
          \right]
 \\
        &\qquad
        + \frac{\mathfrak d}{\pi\sqrt{\gamma''}}
        \left[
            \left(\frac{\gamma''}{3}-\frac q6+\frac{2\gamma'''}{9\gamma''}\right) S^{3/2}
            - \frac{\gamma'''}{3\gamma''} \hat{y}\sqrt{S}
        \right]
        + \mathcal{O}(\mathfrak d^2)
    \end{flalign*}
    where $C$ depends only on $\hat{y}$. Differentiating with respect to $\hat y$ and comparing with the second of \eqref{derivativeshfcell2} implies that $\partial_{\hat y}C(\hat y)=1+\mathcal O(\md^2)$ and \eqref{approxG2} follows.
\end{proof}

\subsubsection{Proof of \Cref{prop:formerassumption}, \Cref{hde0} and \Cref{prop:last}}

\label{Proof00} 

	\begin{proof}[Proof of \Cref{hde0}]
		
		The proof of the fourth statement is very similar to that of the first three and is therefore omitted. The third statement of the lemma follows from integrating the first from $z_0=(x_0, y_0)$ to $(x_1, y_0)$, where $x_1 > x_0$ is chosen so that $(x_1, y_0) \in \mathfrak{A}_{q;a,b,c}$ (which satisfies $x_1-x_0 =\me_{z_0;q}$). It thus suffices to prove the first and second statements of the lemma.

          When $\mathfrak{e}_{z_0;q} \ge \omega$ for a fixed constant $\omega>0$, these follow from \Cref{omegaf}.
          Therefore, we can assume that  $\mathfrak{e}_{z_0;q} \le \omega$ for some sufficiently small but fixed constant $\omega > 0$.
            We distinguish two situations: either $\mathfrak{e}_{z_0;q} \le (1/ {\bf K})\md_{z_0}$ with $\bfK$ a sufficiently large constant, or $\md_{z_0}\le (1/\bfK)\wedge (\bfK\mathfrak{e}_{{z_0};q})$.
            In the former case, we apply \Cref{cor:cella1} and  \eqref{fxygradient}.
            In particular, we have
            \begin{eqnarray}
              \partial_y\cH_{q;a,b,c}=1-\frac1\pi\arg(f+\sqrt{\md_{z_0}\me_{{z_0};q}}\hat f);
            \end{eqnarray}
            together with \eqref{eq:statementsf} and $|\Real(\sqrt{\md_{z_0}\me_{{z_0};q}}\hat f)|\lesssim \md_{z_0}/\sqrt{\bfK}$ that is negligible with respect to $f$ (provided that $\bfK$ is large enough), item (2) of \Cref{hde0} follows. Item (1) follows, similarly, from
            \begin{eqnarray}
               \partial_x\cH_{q;a,b,c}=-1+\frac1\pi\arg(-1+f+\sqrt{\md_{z_0}\me_{{z_0};q}}\hat f)
            \end{eqnarray}
            and the fact that $-1+f\asymp -1$ (again by \eqref{eq:statementsf}).
            
            In the case $\md_{z_0}\le (1/\bfK)\wedge (\bfK\mathfrak{e}_{{z_0};q})$, we use \Cref{cor:cella2} instead. For item (2), note that we need just to prove that  $\partial_y\cH_{q;a,b,c}$ is bounded away from zero, which follows from  $\partial_y\cH_{q;a,b,c}=1-\pi^{-1}\arg(\hat f)$ together with \Cref{rem:nosubtract}. Finally, for item (1) we need to prove that $-\partial_x \cH_{q;a,b,c}\asymp -\md_{z_0}$, which quickly follows from  $\partial_x\cH_{q;a,b,c}=-1+\frac1\pi\arg(-1+\md_{z_0}\hat f)$, \Cref{cor:cella2} {and \Cref{rem:nosubtract}}. 
	\end{proof}

	\begin{proof}[Proof of \Cref{prop:formerassumption}]
       By \Cref{omegaf}, we can restrict to the case where $\me_{z_0;q}\le \omega$ for some sufficiently small but fixed constant $\omega > 0$. 

            When $\me_{z_0;q}\le (1/\bfK)\md_{z_0}$, the statement is a consequence of \Cref{fh}, together with \Cref{rem:whatissmooth}, which guarantees that the error terms in \eqref{derivativeshf} are themselves smooth functions of the coordinates, and \Cref{rem:sparrow}, that guarantees that there is no singularity related to the square root in \eqref{fapproximate}. In fact, note that the function $\widetilde \cH$ in \Cref{fh} differs from the function $\widehat \cH$ in \Cref{prop:formerassumption} only by a space translation (coordinates are centered at $z_0$ for $\widehat \cH$ and at a point on the arctic boundary for $\widetilde \cH$) and by a bounded linear transformation of coordinates (with bounded inverse). Indeed, the latter follows from recalling the  definitions of the vectors $\bu_{z_0;q},\bw_{z_0;q}$, so that the change of coordinates in \eqref{e:funzrescaled} is
            \begin{eqnarray}
              x=x_0+\sqrt{\md\me}\sqrt{1-\gamma'^2}\hat x+\gamma'\md\me \hat y, \quad          y=y_0+\gamma'\sqrt{\md\me}\hat x-\sqrt{1-\gamma'^2}\md\me \hat y,
            \end{eqnarray}
            to be compared with \eqref{xyxyde}.
            
When $\md_{z_0}\le (1/\bfK)\wedge (\bfK\mathfrak{e}_{z_0;q})$, the statement follows a similar reasoning, using this time \Cref{lem:cella2H}, \Cref{rem:nosubtract} and \Cref{rem:whatissmooth}.    
	
	\end{proof} 
	
        \begin{proof}[Proof of \Cref{prop:last}]
      We rename $q$ in the statement as $q_0$, because we will use $q$ to indicate either the value $q_0$ or $0$.
          By \Cref{omegaf}, we can restrict to the case where $\me_{z_0;q_0}\le \omega$ for some sufficiently small but fixed constant $\omega > 0$. 
          As in the proof of \Cref{hde0} and \Cref{prop:formerassumption}, we  distinguish the case $\me_{z_0;q_0}<(1/\bfK)\md_{z_0}$ and the case $\md_{z_0}\le (1/\bfK)\wedge (\bfK\mathfrak{e}_{z_0;q})$.  $\bfK$ is large but independent of any other parameter.

Case I:  $\me_{z_0;q_0}<(1/\bfK)\md_{z_0}$. We seek to locally compare the limit shape $\cH_{q_0;1,1,1}=\cH_{q_0}$ around some point $z_0 = (x_0, y_0)$ close to its arctic boundary  (call $z'_0=(x_0', y_0')$ the point of the boundary closest to $z_0$) to the limit shape $\cH_{a,b,c}$ of the $\mathfrak{q}=1$ measure (with suitable  side-lengths $(a,b,c)$) around some point $z_1 = (x_1, y_1)$ whose closest point on the associated arctic boundary is denoted by $z'_1=(x_1', y_1')$.

Let $\md:=\md_{z_0}$ and $\me:=\me_{z_0;q_0}$.
 For each $i \in \{ 0, 1 \}$, we denote the derivatives of the arctic boundary (for $\cH_{q_0}$ if $i=0$ and $\cH_{a,b,c}$ if $i=1$) at $(x_i', y_i')$ by $(\gamma_i, \gamma_i', \gamma_i'')$. We will later impose that $\gamma':=\gamma_1' = \gamma_0'$, which fixes $(x_1', y_1')$ on the arctic boundary of the $\mathfrak{q} = 1$ limit shape.  
We write $z_0,z_1$ as
	\begin{flalign}
		\label{xyxyde0} 
          x_i = x'_i + (\mathfrak{d} \mathfrak{e})^{1/2} \hat{x_i} - \gamma' \mathfrak{d} \mathfrak{e} \hat{y_i}; \qquad y_i = \gamma_{x'_i} + \gamma' (\mathfrak{d} \mathfrak{e})^{1/2} \hat{x_i} + \mathfrak{d} \mathfrak{e} \hat{y_i}, \quad i=0,1.
	\end{flalign}	
        Note that we use the same $\md,\me$ for both, and that the rescaling of variables is the same as in \Cref{cor:cella1}. Thanks to \Cref{rem:prop}, the coordinates $\hat x_0,\hat y_0$ are of order $1$ (with $\hat y_0>0$ since $z_0$ is inside the liquid region). Our construction (see \Cref{rem:B82} below) will make sure that the same claims are true for $i=1$.

        Define, for any $q\ge 0$, $\widetilde\cH$ using the notations of Lemma \ref{fh}.
        \begin{rem}
        For $q=q_0$ and $a=b=c=1$, $\widetilde\cH$ is nothing but $\widehat \cH_{q_0;1,1,1}$ in \eqref{e:funzrescaled} (up to a shift and a bounded, linear transformation of coordinates, as remarked in the proof of \Cref{prop:formerassumption}), and we denote it $\widetilde\cH_0$.
        For $q=0$, instead, $\widetilde\cH$ is nothing but ${\widecheck\cH}_{a,b,c}$ in \eqref{eq:checkH}, again up to a shift and to the same linear change of coordinates, and we denote it $\widetilde\cH_1$. Since we are aiming at matching derivatives of the two rescaled height functions, these bounded affine transformations are inessential and it suffices to show that the statements of \Cref{prop:last} hold with the functions $\widehat{\mathcal{H}}_q, \widecheck{\mathcal{H}}_{a,b,c}$ appearing there replaced by $\widetilde{\mathcal{H}}_0$ and $\widetilde{\mathcal{H}}_{1}$ here.
  \label{rem:HinsteadofH}
        \end{rem}

Next, we verify that the map from  $(a,x)$ to $(\gamma'_x,\gamma''_x)$ is a bijection with bounded Jacobian for $a$ around $1$, when $b,c$ are fixed as $b=c=(3-a)/2$ (so that $a+b+c=3$).  It suffices to do so at $\mathfrak{q} =1$; this is permitted since is quickly seen from \Cref{equationv} that the arctic boundary  is uniformly smooth in $\mathfrak{q}.$

         From \eqref{f021} and  \eqref{eq:od} (with $b=c=(3-a)/2$)  we have $(\gamma'_x,\gamma''_x)$ as an explicit function of $(a,\xi_x)$.      The  Jacobian matrix $\mathcal J $ of the map $(a,\xi_x)\mapsto (\gamma'_x,\gamma''_x)$ and its determinant can be computed explicitly with the result that  (for $a=1$)
         \begin{eqnarray}
           \label{eq:detJ2}
           \det(\mathcal J)=\frac{2 (1 + (-1 + \xi_x) \xi_x)^2 (-2 + \xi_x (4 + \xi_x (6 + (-8 + \xi_x) \xi_x)))}{(-1 + 2 \xi_x)^5}.
         \end{eqnarray}
   The denominator vanishes for $\xi_x=1/2$,  but by the geometric characterization of Remark \ref{xiq1} we know that $\xi_x<0$ if {$x>x^{SW}>1/2$}.  Recall from \Cref{rem:rightoftangency} that we have $1/2-O(\varepsilon)\le x_0\le 1+O(\varepsilon)$. Therefore, the denominator is non-zero in the range of parameters we are interested in.          As for the numerator,
note 
          from \eqref{eq:xiesplicita}  and from $x\in[1/2,1]$ that $\xi_x$ ranges between $\frac12(1-\sqrt3)$ and $0$. On the other hand, one can quickly verify that the numerator does not vanish in this range.

         Since the map $(a,x)\mapsto (\gamma'_x,\gamma''_x)$ is a local bijection, we fix $(a,x=:x'_1)$ such that
         \begin{equation}
           \label{eq:fixed}
         \gamma'_1=\gamma'_0,\gamma_1'' = \gamma_0'' - \mathfrak{c}_0 q_0 \mathfrak{d}  
         \end{equation}
         for some small constant $\mathfrak{c}_0 > 0$ to be fixed later.
As a consequence, defining   (according to  \eqref{fxy2}) $\Gamma_0=\gamma''_0-\gamma'(1-\gamma')q_0,\Gamma_1=\gamma''_0-\mathfrak c_0 q_0\md$, we have that  $\Gamma_1 - \Gamma_0 =\gamma'(1-\gamma')q_0-\mathfrak c_0q_0\md\asymp q_0 \mathfrak{d}$,  as $\gamma' \asymp \mathfrak{d}$, provided that $\mathfrak c_0$ is small enough. Since $z_1'$ lies on the arctic boundary, fixing $x_1'$ fixes also its second coordinate.

Up to now, we have fixed $z'_1$. Now we want to fix $z_1$ (that is, $\hat x_1,\hat y_1$). To this purpose, we make the following observation.
        \begin{rem}
  \label{rem:B80} A look at the proof of  \Cref{fh} and \Cref{rem:whatissmooth} shows that, under the assumptions of \Cref{cor:cella1} we have  the following: 
if  $a,b\in\mathbb N_0$, then 
\begin{align}
\label{approxG} \partial^a_{\hat x}\partial^b_{\hat y} \widetilde{\mathcal{H}} (\hat{x},\hat{y})&=\partial^a_{\hat x}\partial^b_{\hat y} \mathfrak{G} (\hat x, \hat y) +\mathcal{O} ( \sqrt{\md\me})+\hat{\mathcal{O} }( \sqrt{\me/\md}), \\ 	\mathfrak{G} (\hat x, \hat y)& =\Xi \cdot \big( 2 (\gamma'^2 + 1) \hat y - \gamma'' \hat x^2 \big)^{3/2}, \qquad \Xi:= \displaystyle\frac{\mathfrak{d} \Gamma^{1/2}}{3 \pi |\gamma'| (1 - \gamma')}
\end{align}
where the error terms are   smooth in $q,a,b,c,\hat x,\hat y$ {(hence also in $\Gamma,\gamma',\gamma''$)} and, in addition,  $\hat{\mathcal{O}} ( \sqrt{\me/\md})$ are smooth functions of the arguments $\md,\me,\hat f,\gamma'$ only. This is the reason why we distinguish the two types of error, despite the fact that $\sqrt{\md\me}\ll\sqrt{\me/\md}$.

This should be compared with \Cref{lem:cella2H}, that holds under the assumptions of \Cref{cor:cella2} instead (an aspect to be noted is that, while level lines of $\frak G$ are simply parabolas, those of $\frak H$ are more intricate).

\end{rem}

We claim that the map from $(\hat{x}, \hat{y})$ to $\nabla \widetilde{\mathcal{H}} (\hat{x}, \hat{y})$ is a bijection with bounded Jacobian. Indeed, in view of \Cref{rem:B80}, it suffices to do this for $\mathfrak{G}$, provided that $\bfK$ is chosen large enough (recall that we are working under the assumption 
 that $\me<(1/\bfK) \md$).
            The Jacobian of the map $(\hat{x}, \hat{y})\mapsto \nabla \mathfrak G(\hat{x}, \hat{y})$ is
            \begin{equation}
              \label{eq:J}
              J(\hat x,\hat y)=\displaystyle\frac{\mathfrak{d} \Gamma^{1/2}}{ \pi |\gamma'| (\gamma'-1)}\frac1{\sqrt{ 2 (\gamma'^2 + 1) \hat y - \gamma'' \hat x^2 }}
              \begin{pmatrix}
               2\gamma''(-\gamma'' \hat x^2+(1+\gamma'^2)\hat y)&(1+\gamma'^2)\gamma''\hat x \\(1+\gamma'^2)\gamma''\hat x &-(1+\gamma'^2)^2
              \end{pmatrix}
            \end{equation}
          with
            \begin{equation}
              \label{eq:detJ}
              \det(J)=-\frac{\md^2\Gamma(1+\gamma'^2)^2\gamma''}{\pi^2\gamma'^2(1-\gamma')^2}
            \end{equation}
            Note that this expression is strictly {negative} and bounded away from zero and minus infinity since the curvature $\gamma''$ and $\Gamma$ are positive and bounded   uniformly if $q_0$ is small and $a,b,c$ are close to $1$. Recall also from Remark \ref{rem:gamma'} that $\gamma'$ is bounded away from $1$, and that $\md/\gamma'\asymp 1$.

            \begin{rem}
              \label{rem:B82}            Since the map from $(\hat{x}, \hat{y})$ to $\nabla \widetilde{\mathcal{H}} (\hat{x}, \hat{y})$ is a bijection with bounded Jacobian, we fix $(\hat x_1,\hat y_1)$ by imposing $\nabla \widetilde{\mathcal{H}}_1 (\hat{x}_1, \hat{y}_1) = \nabla \widetilde{\mathcal{H}}_0 (\hat{x}_0, \hat{y}_0)$.
Since $\Gamma_1-\Gamma_0\asymp \gamma''_0-\gamma''_1\asymp q_0\md$,  it follows from \eqref{approxG} that
$|\hat{x}_1 - \hat{x}_0| + |\hat{y}_1 - \hat{y}_0| = \mathcal{O} ( q_0 \mathfrak{d})$.
This (together with Remark \ref{rem:B80}) implies that $\partial_{\gamma} \hat{\mathcal{H}}_1 (\hat{x}_1, \hat{y}_1) = \partial_{\gamma} \hat{\mathcal{H}}_0 (\hat{x}_0, \hat{y}_0) + \mathcal{O} (q_0 \mathfrak{d})$ for any differential operator $\partial_\gamma$, which verifies the last statement in the proposition.
            \end{rem}

To complete the proof, it then suffices to show that $\nabla^2 \widetilde{\mathcal{H}}_1 (\hat{x}_1, \hat{y}_1) - \nabla^2 \widetilde{\mathcal{H}}_0 (\hat{x}_0, \hat{y}_0) \gtrsim q_0 \mathfrak{d} \cdot \Id$.  To do so, by \Cref{rem:B80}, we can replace $\widetilde \cH$ with $\mathfrak{G}$, provided that $\bfK$ is chosen large enough. Recall that going from $(\hat x_0,\hat y_0)$ to $(\hat x_1,\hat y_1)$, $\Gamma$ increases by $\asymp q_0\md$ while $\gamma''$ decreases by $\mathfrak c_0q_0\md$.
It then suffices to show that
\begin{eqnarray}
  \label{eq:instead}
  \partial_{\Gamma} \nabla^2 \mathfrak{G} (\tilde{x}(\Gamma,\gamma''), \tilde{y}(\Gamma,\gamma''))-\mathfrak c_0\partial_{\gamma''} \nabla^2 \mathfrak{G} (\tilde{x}(\Gamma,\gamma''), \tilde{y}(\Gamma,\gamma'')) \gtrsim \Id,
\end{eqnarray}
 where $(\tilde{x}(\Gamma,\gamma''), \tilde{y}(\Gamma,\gamma''))$ is chosen so that $(\gamma'',\Gamma)\mapsto \nabla \mathfrak{G} (\tilde{x}(\Gamma,\gamma''), \tilde{y}(\Gamma,\gamma'')) $ is constant. 
 Write $\tilde x=\tilde{x} (\Gamma, \gamma''), \tilde y=\tilde{y} (\Gamma, \gamma'')$ for lightness and let $\eta=2 (\gamma'^2 + 1) , \psi=\md/(3\pi\gamma'(1-\gamma'))$. Note that $\eta,\psi$ are positive and bounded away from zero and infinity, since  $\gamma'/\md\asymp 1$.
From the explicit formula for $\frak G$, it is immediate to see that 
\begin{eqnarray}
  \tilde{x} &=& \frac{\kappa_1}{\gamma''},\quad  \kappa_1 := - \frac{\eta \partial_{\hat x}\mathfrak{G}}{2
  \partial_{\hat y}\mathfrak{G}}\\
 \sqrt{\eta \tilde{y} - \gamma''  \tilde{x}^2} &=& \frac{\kappa_2}{\sqrt{\Gamma}}, \quad \kappa_2 :=
  \frac{2\partial_{\hat y}\mathfrak{G}}{3\psi \eta} > 0.
\end{eqnarray}
Since $\kappa_1,\kappa_2$ are independent of $\Gamma,\gamma''$ (because $(\gamma'',\Gamma)\mapsto \nabla\frak G(\tilde{x} (\Gamma, \gamma''), \tilde{y} (\Gamma,
\gamma''))$ is constant), it follows that
\begin{equation}
  \partial_{\Gamma} \nabla^2 \mathfrak{G}(\tilde x,\tilde y) = \psi\left(\begin{array}{cc}
    \frac{3 \kappa_1^2}{\kappa_2} & - \frac{3 \eta \kappa_1}{2 \kappa_2}\\
    - \frac{3 \eta \kappa_1}{2 \kappa_2} &  \frac{3 \eta^2}{4 \kappa_2}
  \end{array}\right), \quad   \partial_{\gamma''} \nabla^2 \mathfrak{G} (\tilde x,\tilde y)= \psi\left(\begin{array}{cc}
    - 3 \kappa_2 & 0\\
    0 & 0
  \end{array}\right).
\end{equation}
The former has eigenvalues  $0$ and $3\psi (\eta + 4 \kappa_1^2) / (4 \kappa_2) > 0$
and $\partial_{\Gamma} \nabla^2 \mathfrak{G}(\tilde x,\tilde y) - \mathfrak c_0
\partial_{\gamma''} \nabla^2 \mathfrak{G}(\tilde x,\tilde y)$ has determinant $9 \eta^2\psi^2 \mathfrak c_0/ 4>0$. Therefore, \eqref{eq:instead} follows since we have taken
$\mathfrak c_0 > 0$.

{ Case II:  $\md_z\le (1/\bfK)\wedge (\bfK\mathfrak{e}_{z;q_0})$}.
  We let as usual $\mathfrak q=e^{q}$ and $\gamma''(q,a,b,c)$ and $\gamma'''(q,a,b,c)$ denote the second and third derivative of the arctic boundary at the tangency point $p_{q;a,b,c}^{SW}$ (given by \eqref{tangpoint}).  Using \Cref{axy} to express $\gamma_x$ as a function of $(x,q,a,b,c)$, then computing its second and third derivatives at $p^{SW}_{q;a,b,c}$  and finally expanding in $q$ and in $a-1,b-1$ (with $a+b+c$ fixed to $3$), one gets 
  \begin{align}
\label{eq:theresultis}
    \gamma''(q,1+\bar a,1+\bar b,{1-(\bar a+\bar b)})=\frac43 - 2  \bar a  +O(\bar a^2+\bar b^2+q^2),\\\label{eq:theresultis2}
    \gamma'''(q,1+\bar a,1+\bar b,{1-(\bar a+\bar b)})=6\Big(-\frac49 + \frac8{9} (q + 2 \bar a + \bar b)\Big)  +O(\bar a^2+\bar b^2+q^2).
  \end{align}
  In particular,  the map $(a,b)\mapsto (\gamma'',\gamma''')$ (with $c$ fixed to be $3-a-b$) is a local bijection with uniformly bounded Jacobian. Given $q>0$, we fix $a_q,b_q,c_q=3-a_q-b_q$ so that
    \begin{eqnarray}
      \label{eq:sothat}
    \gamma''(q,1,1,1)-\gamma''(0,a_q,b_q,c_q)=0      , \quad
  \gamma'''(q,1,1,1)- \gamma'''(0,a_q,b_q,c_q)=r q+O(q^2)
    \end{eqnarray}
    where ${r\in\mathbb R }$ will be chosen later. It is quickly verified from \eqref{eq:theresultis}-\eqref{eq:theresultis2} that
    \begin{eqnarray}
     a_q-1 =      O(q^2),\quad b_q-1=q  \Big(1 - \frac{3}{16}r\Big)+O(q^2).
    \end{eqnarray}
   \begin{rem}
      \label{rem:gamma'''}
    Note that   $\gamma'''$ changes  by order $q$, but   $\gamma'''$ appears in the expressions of \Cref{cor:cella2} always multiplied by $\md$, so the effect on  $D^2\widetilde \cH$ will be of order $\md q$.
    \end{rem}
    
Recall the definition \eqref{eq:theg} of $g$. 
 We have
  \begin{eqnarray}
    \label{eq:gis}
    g(q,1+\bar a,1+\bar b,1-(\bar a+\bar b))=2+\frac13(q-7 \bar a+4 \bar b) +O(\bar a^2+\bar b^2+q^2)
  \end{eqnarray}
  To see this, just express $A',B',D'$ through the rational function $f_0$ of \eqref{f0uv2}, note that
  \begin{eqnarray}
    \label{tangpoint}
    \mathfrak{q}^{x_{q;a,b,c}^{SW}}=\frac{(\mathfrak q^a-1)(\mathfrak q^{-b}-1)}{\mathfrak q^{-b-c}-1}+1,
  \end{eqnarray}
  which follows directly from \Cref{xabc}, and expand the resulting explicit function in  $q,a-1,b-1$.

  One deduces  from \eqref{eq:gis} and \eqref{eq:theresultis}-\eqref{eq:theresultis2} that
    \begin{eqnarray}
      \label{eq:gq}
      g(q,1,1,1)=2+\frac q3+O(q^2),\quad g(0,a_q,b_q,c_q)=2+q\Big(\frac43-\frac{r}4\Big)+O(q^2).
    \end{eqnarray}

    We define $\widetilde \cH$ as in \Cref{lem:cella2H} and we call it $\widetilde \cH_0$ if $(q,a,b,c)=(q_0,1,1,1)$ and $\widetilde\cH_1$ if $(q,a,b,c)=(0,a_q,b_q,c_q)$ with $a_q,b_q,c_q$ chosen as in \eqref{eq:sothat} (with $q$ there equal to $q_0$ here). Let again $\md:=\md_{z_0}$ and
     call  $\gamma''_0=\gamma''(q,1,1,1),\gamma''_1:=\gamma''(0,a_q,b_q,c_q)$ (and similarly for $\gamma'''_0,\gamma'''_1$) the values of $\gamma''$ (resp. $\gamma'''$) appearing in $\widetilde \cH_0,\widetilde \cH_1$.
Recall that our goal is to prove \eqref{e:thereexist} for a suitable
choice of $z_1$. Reasoning similarly as in Case I, this is
equivalent to proving \eqref{e:thereexist} with
$\widehat \cH_q,\widecheck \cH_{a,b,c}$ there replaced by
$\widetilde \cH_0,\widetilde \cH_1$ here.

 \begin{rem}
   \label{rem:B82bis}  Since we want to bound the difference $D^2\widetilde \cH_0-D^2\widetilde \cH_1$ rather than the two Hessians individually, we can and will replace the expression  for  $\nabla\widetilde \cH$ in  \Cref{lem:cella2H} by the  expression $\nabla\breve\cH$ below, where we neglect both the terms of order $\md^2$ (since in the difference \eqref{e:thereexist} they would give a contribution of order $\md^2 q_0\le (1/\bfK)\md q_0$ and we are assuming that  $\bfK$ is large),  as well as the terms of order $\md$, except those that depend linearly in  $q_0$ (where the $q$-dependence can be hidden in $\gamma'',\gamma''',g$ because of \eqref{eq:sothat} and \eqref{eq:gq}):
\begin{eqnarray}
  \label{eq:breve}
  & \breve f(\hat x,\hat y)=  \mathrm{i} \sqrt{D_1}-D_2\\
  \label{eq:brevegrad}
  &\partial_{\hat x}\breve\cH(\hat x,\hat y)= - \frac{1}{\pi} \sqrt{D_1},\qquad
  \partial_{\hat y}\breve\cH(\hat x,\hat y)=\frac{1}{\pi} \arctan \left( \frac{\sqrt{D_1}}{D_2} \right)\\
  & D_1=(\gamma''+\md\hat x( \gamma'''-q\gamma'')) \Big(2 \hat{y} - \gamma'' \hat{x}^2 - \displaystyle\frac{\gamma'''}{3} \mathfrak{d} \hat{x}^3 \Big),\\ &D_2=   \gamma'' \hat{x}   
		- \mathfrak{d} \breve g\,\Big(\hat{y}- \gamma''\displaystyle\frac{\hat{x}^2}{2} \Big)+\md \hat x^2\frac{\gamma'''}{2}.
\end{eqnarray}
{We then call $\breve \cH_i,i=0,1$ the expression obtained replacing $(\gamma'',\gamma''')$ with $\gamma''_i,\gamma'''_i$, $q$ with either $q_0$ (if $i=0$) or with $0$ (if $i=1$) and 
  $\breve g$  with $\breve g_i$, where (compare with \eqref{eq:gq})}
\begin{eqnarray}
  \label{deltagi}
  \breve g_0=2+\frac q3, \quad \breve g_1=2+q\Big(\frac43-\frac{r}4\Big).
\end{eqnarray}
{Note that the derivatives in \eqref{eq:brevegrad} are not quite consistent (they do not verify $\partial_{\hat x}(\partial_{\hat y}\breve\cH)=\partial_{\hat y}(\partial_{\hat x}\breve\cH)$). This is due to neglecting the error of ${\mathcal{O}} (\mathfrak{d} \hat{f})$ in \eqref{derivativeshfcell2}. Despite that,  we will  formally replace the derivative of $\widetilde{\cH}$ by \eqref{eq:brevegrad} and its second derivatives by \eqref{eq:Hbreve11}-\eqref{eq:Hbreve22} below, since (as explained at the beginning of this remark) the error involved is negligible, for the purposes of the present proof.
}
\end{rem}

Note that $D_1$ is positive and bounded away from zero; see \Cref{rem:nosubtract}. 
The map $(\hat x,\hat y)\mapsto \nabla \breve \cH$ is a local bijection; this can be seen by computing its Jacobian matrix (that is, the Hessian of $\breve \cH$) and checking that it is not degenerate. In fact,
the matrix elements of the Hessian of $\breve \cH$ are
\begin{align}
  \label{eq:Hbreve11}
  &  \partial^2_{\hat{x}} \breve{\mathcal{H}} =-\frac{1}{2\pi\sqrt{D_1}}\partial_{\hat x} D_1  \\
  \label{eq:Hbreve12}
& \partial^2_{\hat{x}\hat y} \breve{\mathcal{H}}=-\frac{1}{\pi\sqrt{D_1}}(\gamma''+\md \hat x(\gamma'''-q\gamma''))\\
    \label{eq:Hbreve22}
&\partial^2_{\hat{y}} \breve{\mathcal{H}}=\frac1\pi\frac{D_2(\gamma''+\md \hat x(\gamma'''-q\gamma''))+D_1\md\breve g}{(D_1+D_2^2)\sqrt{D_1}} .
\end{align}

With manipulations using 
\begin{eqnarray*}
    2\hat{y} - \gamma'' {\hat{x}^2} = \frac{D_1}{\gamma''+ \mathfrak{d} \hat{x}(\gamma'''-q\gamma'')}+\frac{\gamma'''}3\md \hat x^3, \qquad \hat{x} \gamma'' = D_2 + \mathfrak{d} \breve g \frac{D_1}{2\gamma''}-\md \hat x^2\frac{\gamma'''}2+\mathcal O(\md^2),
\end{eqnarray*}
we can write
\begin{align}
\label{breveH11simp}\partial^2_{\hat{x}} \breve{\mathcal{H}} =\frac{1}{\pi\sqrt{D_1}} \left\{ \frac{\mathfrak{d}}{2\gamma''}[(\breve g +q)\gamma''-\gamma'''] D_1+D_2[\gamma''+\md\hat x(\gamma'''-q\gamma'')]  \right\}+\mathcal O(\md^2).
\end{align}
At zeroth order in $\md,q$, the Hessian of $\breve \cH$ is then (using $D_1+D_2^2=2\hat y\gamma''+O(\md)$)
\begin{eqnarray}
  \frac{\gamma''}{\pi\sqrt{D_1}}\begin{pmatrix}
                                  \gamma''\hat x & -1\\
                                  -1& \frac{\hat x}{2\hat y}
  \end{pmatrix}
\end{eqnarray}
which is non-singular (recall {from \Cref{rem:nosubtract}} that $D_1$, and therefore also $\hat y$, is bounded away from zero).

Since $(\hat x,\hat y)\mapsto \nabla \breve \cH$ is a local bijection, we fix $(\hat x_1(q_0),\hat y_1(q_0))$ in such a way that $\nabla \breve\cH_1(\hat x_1(q_0),\hat y_1(q_0))=\nabla \breve\cH_0(\hat x,\hat y)$. In particular, $D_1,D_2$ are constant with respect to $q_0$.

Note that
\begin{eqnarray}
  \label{eq:makesure}
|\hat x_1(q_0)-\hat x|+|\hat y_1(q_0)-\hat y|=O(\md q_0).
\end{eqnarray}
This follows noting that between $\breve \cH_0$ and $\breve \cH_1$, $\gamma'''$ changes by $O(q)$, and $\gamma'''$ is multiplied by $\md$ in the expression  \eqref{eq:breve} for $D_1,D_2$. It remains to check that
\begin{eqnarray}
  \label{eq:itre}
  D^2\breve\cH_1-D^2\breve \cH_0\gtrsim \md q_0 \Id.
\end{eqnarray}

From \eqref{breveH11simp} and \eqref{eq:makesure} we deduce that 
\begin{eqnarray}
  \label{eq:H11diff}
  &\partial^2_{\hat{x}} \breve{\mathcal{H}}_1-\partial^2_{\hat{x}} \breve{\mathcal{H}}_0=\frac{q_0\md}{\pi\sqrt{D_1}}\left\{D_2^2
    +r \Big[D_1(4-\gamma'')/(8\gamma'')-D_2 \hat x\Big]\right\}  +O(\md^2q_0) .
\end{eqnarray}

\noindent Similarly, from \eqref{eq:Hbreve12}, \eqref{eq:Hbreve22} we have
\begin{eqnarray}
  \label{eq:H12diff}
  &  \partial^2_{\hat{x}\hat y} \breve{\mathcal{H}}_1-\partial^2_{\hat{x}\hat y} \breve{\mathcal{H}}_0=\frac{q_0\md}{\pi\sqrt{D_1}}\left[-D_2
    +r\hat x\right] +O(\md^2q_0) \\
  &  \partial^2_{\hat y} \breve{\mathcal{H}}_1-\partial^2_{\hat y} \breve{\mathcal{H}}_0=\frac{q_0\md}{\pi\sqrt{D_1}}\left[1   -\frac{r}{2\hat y\gamma''}(\frac{D_1}4+D_2\hat x)
    \right]+O(\md^2q_0) .
\end{eqnarray}

We conclude that
\begin{eqnarray}
  \label{eq:weco}
  D^2\breve\cH_1-D^2\breve \cH_0= \frac{q_0\md}{\pi\sqrt{D_1}} M+O(\md^2q_0) ,\\
 M= \begin{pmatrix}
    D_2^2 +r \Big[{D_1}\frac{4-\gamma''}{8\gamma''}-D_2 \hat x\Big] & -D_2+r\hat x\\
    -D_2+r\hat x & 1    -\frac{r}{2\hat y\gamma''}(\frac{D_1}4+D_2\hat x)
  \end{pmatrix}.
\end{eqnarray}
Using $\gamma''=4/3+O(q_0)$,
{$D_1=(4/3)(2\hat y-(4/3)\hat x^2)+O(q_0+\md)$ and $D_2=(4/3)\hat x+O(q_0+\md)$ we get that  
\begin{eqnarray}
  \label{eq:detM}
  \det(M)=\frac r{3\hat y}\,(2\hat y-\frac43\hat x^2)(\hat y+\frac43 \hat x^2)+O(r(r+q_0+\md)).
\end{eqnarray}}
Recall from Remark \ref{rem:nosubtract} that under our assumptions, $(2\hat y-\frac43\hat x^2)$ and $\hat y$ are bounded away from zero and infinity. Then, we deduce from \eqref{eq:weco} and \eqref{eq:detM} (together with the positivity of the trace of $M$ for $r$ small) that, provided that $r$ is positive and small enough, the  inequality \eqref{eq:itre} holds.
\label{Heightf} 
        \end{proof}
	
	\subsection{Comparison between $\mathcal{H}_q$ at different $q$}
	
	\label{HqEstimate}

	In this section we show \Cref{boundarydistance} and \Cref{qlineq}.
	\begin{lem} 
		
		\label{hqomegaq}
		
		Fix real numbers $q, q', \omega, \kappa \in \mathbb{R}$ with $q' \le q$. If $\mathcal{H}_{q'} (z) \le \mathcal{H}_{q} (z) - \omega$ for all $z \in \mathfrak{L}_q$ with $\dist (z, \mathfrak{A}_{q}) = \kappa$, then $\mathcal{H}_{q'} (z) \le \mathcal{H}_{q} (z) - \omega$ for all $z \in \mathfrak{L}_q$ with $\dist (z, \mathfrak{A}_{q}) \ge \kappa$.
	\end{lem} 
	
	\begin{proof}
          The limit shape is monotone  with respect to $q$ (by \Cref{rem:stochdom}) and with respect to boundary conditions (this follows from monotonicity of the solution of the elliptic partial differential equation \eqref{eq:pdeq} with respect to the boundary data). The lemma then follows applying these monotonicities in the domain $\Lambda = \{ z \in \mathfrak{L}_q : \dist (z, \mathfrak{A}_{q}) \ge \kappa \}$.
	\end{proof}
	
	Fix $y \in [0, 49/100]$ (recall \eqref{x0y0}), and let $\mathfrak{r}_q (y):=\max\{x\in\mathbb R: (x, y) \in \mathfrak{A}_q\}$. 
	\begin{lem}
		
		\label{rqc} 
		
		There exists a constant $C > 1$ such that $C^{-1} (q-q') \le \mathfrak{r}_q (y) - \mathfrak{r}_{q'} (y) \le C (q-q')$, for any $q, q' \in [-\varepsilon, \varepsilon]$ with $q \ge q'$, if $\varepsilon$ is sufficiently small.
		
	\end{lem}

	\begin{proof} 
		We begin by setting some notation. Set $\mathfrak{q} = e^q$ and $\mathfrak{q}' = e^{q'}$; also denote $x = \mathfrak{r}_q (y)$ and $x' = \mathfrak{r}_{q'} (y)$. As in the $(a, b, c) = (1, 1, 1)$ case of \eqref{abcxy}, set $(A, B, C) = (-[-1]_{\mathfrak{q}}, [1]_{\mathfrak{q}}, [2]_{\mathfrak{q}} - [1]_{\mathfrak{q}})$ and $(A', B', C') = (-[-1]_{\mathfrak{q}'}, [1]_{\mathfrak{q}'}, [2]_{\mathfrak{q}'} - [1]_{\mathfrak{q}'})$. Also following \eqref{abcxy}, let $X = X_q (x)$, $X' = X_{q'} (x')$, $Y = Y_q (y)$, and $Y' = Y_{q'} (y)$; further denote $\Delta = q -q'$. 
		
		Now, observe that if $y = 0$ then, by \Cref{xabc} and the fact that $ |q| \le \varepsilon$, $x = 1/2 + O(\varepsilon)$ and $x' = 1/2 + O(\varepsilon)$. Therefore,  
		\begin{flalign*}
			x-x' +  \displaystyle\frac{3 \Delta}{8} + O(\Delta^2 + \varepsilon \Delta) = X - X' & = \displaystyle\frac{AB}{B+C} - \displaystyle\frac{A'B'}{B'+C'}\\
			&  = \displaystyle\frac{\mathfrak{q}-1}{1 - \mathfrak{q}^{-2}} -  \displaystyle\frac{\mathfrak{q}'-1}{1 - \mathfrak{q}'^{-2}} = \displaystyle\frac{3\Delta}{4} + O(\varepsilon \Delta).
		\end{flalign*}
		
		\noindent where the first statement follows from Taylor expansion and \eqref{abcxy} (with the fact that $x, x' = 1/2 + O(\varepsilon)$); the second follows from \Cref{xabc}; the third follows from \eqref{abcxy}; and the fourth follows from the fact that $|q|, |q'| = O(\varepsilon)$. Therefore, $x - x' = 3 \Delta/8 + O(\varepsilon \Delta)$, which for sufficiently small $\varepsilon$ yields the lemma when $y = 0$. 
		
		Next observe from the explicit form \eqref{axy2} of the arctic boundary that it is smooth in the parameter $q$. Therefore, $|x-x'| \lesssim \Delta$, so it suffices to show that $x-x' \gtrsim \Delta$. By this smoothness, we also have that the difference between the slope of the tangent line to $\mathfrak{A}_q$ at $(x, y)$ and the slope of the tangent line to $\mathfrak{A}_{q'}$ at $(x', y)$ is $O(\Delta)$. It follows since $x-  x' \gtrsim \Delta$ at $y=0$ that $x - x' \gtrsim \Delta$ for $y \le C^{-1}$, if $C > 1$ is sufficiently large. Therefore, we may assume in what follows that $y \gtrsim 1$. 
		
		Next, we have from \Cref{axy} that 
		\begin{flalign*}
			 \big( X (B+C) - Y&(A+C) - AB + (\mathfrak{q}^{-1}-1)ABY \big)^2\\& - 4ABY \big( A+C-X+(1-\mathfrak{q}^{-1})(A+C)Y \big) = 0,
		\end{flalign*} 
		and 
		\begin{flalign*} 
			\big( X' (B'+C') - Y' & (A'+C') - A'B' + (\mathfrak{q'}^{-1}-1) A'B'Y' \big)^2 \\
			& - 4A'B'Y' \big( A'+C'-X'+(1-\mathfrak{q'}^{-1})(A'+C')Y' \big) = 0. 
		\end{flalign*}
		
		\noindent Since $(a, b, c) = (1, 1, 1)$ and $q \in [-\varepsilon, \varepsilon]$, by \eqref{abcxy} we have
		\begin{flalign}
			\label{xy2}
			\begin{aligned}
				& X = X' + x - x' + \displaystyle\frac{\Delta}{2} \cdot x(x+1) + O(\Delta^2); \qquad Y = Y' + \displaystyle\frac{\Delta}{2} \cdot y(y+1)+ O (\Delta^2); \\
				& \qquad \quad  A = A' + \Delta + O(\varepsilon \Delta); \qquad B = B' = 1; \qquad C = C' - \Delta + O(\varepsilon \Delta).
			\end{aligned}
		\end{flalign}
		
		\noindent Set $\breve{X} = X' + \Delta x(x+1)/2$. It suffices to show that $-\xi (\breve{X}, Y) \gtrsim \Delta$, as $-\xi (\breve{X},Y) = \xi (X,Y) - \xi (\breve{X}, Y) \lesssim X - \breve{X} = x-x' + O(\varepsilon \Delta)$, where the first statement holds by \eqref{axy2}, the second is quickly verified using the fact that $Y \gtrsim 1$ (as this implies that $(x, y)$ is bounded away from a tangency location of $\mathfrak{A}_q$), and the third holds by \eqref{xy2} with the definition of $\breve{X}$. To see this, observe using \eqref{xy2} that
		\begin{flalign}
			\label{xi3} 
			\begin{aligned} 
				-\xi (\breve{X}, Y) & = \big( X' (B'+C') - Y' (A'+C') - A'B' + (\mathfrak{q}^{-1}-1) A'B'Y' \big)^2 \\
				& \qquad - 4A'B'Y' \big( A'+C'-X'+(1-\mathfrak{q})(A'+C')Y' \big) \\ 
				& - \big( \breve{X} (B+C) - Y(A+C) - AB + (\mathfrak{q}^{-1}-1)ABY \big)^2 \\
				& \qquad + 4ABY \big( A+C-\breve{X}+(1-\mathfrak{q})(A+C)Y \big) \\
				& = \Delta (2 X'^2 Y' + 2X' Y'^2 - 4X'^3 - 4Y'^3 + 2X'^2 + 2Y'^2 + 4X' + 4Y' - 2) + O(\varepsilon \Delta).
			\end{aligned} 
		\end{flalign}
		
		\noindent Again since $\mathfrak{A}_q$ is smooth in $q$, the distance between $(X', Y')$ to a corresponding point $(\mathfrak{r}_0 (Y'), Y') \in \mathfrak{A}$ (that is, at $q=0$) is $O(|q|) = O(\varepsilon)$. Since (by \eqref{axy2}) we have $\mathfrak{A} = \{ (x, y) \in \mathbb{R}^2 : 4x^2 - 4xy + 4y^2 - 4x - 4y + 1 = 0 \}$, we have that 
		\begin{flalign*}
			\mathfrak{r}_0 (Y') = \displaystyle\frac{Y'+1 + (6Y' - 3Y'^2)^{1/2}}{2}.
		\end{flalign*} 
		
		\noindent Inserting this, with the bound $X' = \mathfrak{r}_0 (Y') + O(\varepsilon)$, into \eqref{xi3} yields 
		\begin{flalign*}
			-\xi (\breve{X}, Y) = \displaystyle\frac{\Delta}{2} \cdot (1-2Y')(3-2Y') (6Y'-3Y'^2)^{1/2} + O(\varepsilon \Delta) \gtrsim \Delta,
		\end{flalign*} 
		
		\noindent where the last inequality follows from the fact that $1/2-Y' \gtrsim 1$ (as $1/2-y \gtrsim 1$, since $\ell_0 (x,y)$ is the south edge of $\mathfrak{X}$). As mentioned above, this yields the lemma.
              \end{proof} 

	\begin{proof}[Proof of \Cref{boundarydistance}]
		This follows from \Cref{rqc}, since $x = \mathfrak{r}_q (y)$ and $x' = \mathfrak{r}_{q'} (y)$. 
	\end{proof} 

	\begin{lem}
	
	\label{qdifferencehq} 
	
	There exists $C > 1$ such that the following holds whenever $\varepsilon$ is sufficiently small. For any $q,q' \in [-\varepsilon, \varepsilon]$, we have for any $e \in [0, \mathfrak{r}_q (y)]$ that
	\begin{flalign*}
		|\partial_x \mathcal{H}_q ( \mathfrak{r}_q (y) - e, y) - \partial_x \mathcal{H}_{q'} (\mathfrak{r}_{q'} (y) - e, y) | \le C |q-q'| y^{1/4} e^{1/2}.
	\end{flalign*}  
	
	\end{lem} 

	\begin{proof}

		Denote $x_1 = \mathfrak{r}_q (y)$ and $x = x_1 - e$, and let $z = (x, y)$; observe that $(\mathfrak{d}_z, \mathfrak{e}_{z;q}) = (y^{1/2}, e)$. Also set $X = -[-x]_q$ and $Y = -[-y]_q$. Further let $\xi_q$ and $\xi_{q'}$ denote the function defined by \eqref{xietazeta}, with the parameter equal to $q$ and $q'$, respectively (including in the definitions of $(A,B,C)$ through $(a, b, c)$). 
		
		First observe that 
		\begin{flalign}
			\label{12zeta} 
			\zeta(X,Y) \ge \displaystyle\frac{1}{2},
		\end{flalign} 
	
		\noindent since
		\begin{flalign*}
			\zeta(X,Y) = 2(A+C-X + (1 - \mathfrak{q}^{-1})(A+C)Y) \ge 2(2 - x) + \mathcal{O}(q) \ge \displaystyle\frac{1}{2},
		\end{flalign*}
	
		\noindent for $q$ sufficiently small. Further observe that 
		\begin{flalign*} 
			\Imaginary f(z) \asymp (y^{1/2} e)^{1/2}; \qquad \Real f(z) \lesssim (y^{-1/2} e)^{1/2}.
		\end{flalign*} 
	
		\noindent {Indeed, if $\me_{z;q} \le \omega$ for some sufficiently small constant $\omega > 0$, then this follows  (using the fact that $e \lesssim y^{1/2}$) from either \Cref{cor:cella1} or  \Cref{cor:cella2}, according to whether $\me_{z;q}\le (1/{\bf K}) y^{1/2}=(1/{\bf K})\md_z$ or not, with $\bf K$ a large constant}.  Otherwise, it follows from \eqref{xietazeta1}. Together with \eqref{f0uv2}, this implies that $\Imaginary v \asymp (\mathfrak{d} \mathfrak{e})^{1/2}$, which by \eqref{12zeta} and \eqref{vxietazeta} yields
		\begin{flalign}
			\label{xixyeta} 
			\xi(X,Y) \asymp -y^{1/2} e; \qquad \eta (X, Y) \lesssim (y^{-1/2} e)^{1/2}.
		\end{flalign}
	
		\noindent Further observe that 
		\begin{flalign}
			\label{qderivativexietazeta} 
			|\partial_q \xi (X,Y) | + |\partial_q \eta (X, Y)| + |\partial_q \zeta (X, Y)| \lesssim 1,
		\end{flalign} 
	
		\noindent by \eqref{abcxy} and \eqref{xietazeta}, with the fact that $\partial_q [R]  \lesssim 1$ for $R \lesssim 1$. 
		
		By \eqref{xixyeta},  \eqref{12zeta}, \eqref{qderivativexietazeta}, and the second statement in \eqref{fxygradient}, it suffices to show that $|\xi_q (X, Y)^{1/2} - \xi_{q'} (X', Y')^{1/2}| \lesssim |q-q'| y^{1/4} e^{1/2}$. Therefore, we must show that 
		\begin{flalign} 
			\label{xiq} 
			|\partial_q (\xi_q (X, Y)^{1/2})| \lesssim y^{1/4} e^{1/2}.
		\end{flalign} 
	
		\noindent  To that end, we may assume that $e = \mathfrak{e}_{z;q}$ (corresponding to $z$ is closer in the $\ell_z$ direction to the right side of $\mathfrak{A}_{q}$  than to its left side), as the proof in the alternative case is entirely analogous. Since \eqref{xixyeta} implies that $-\xi_q (X, Y) \asymp y^{1/2} e$, to verify \eqref{xiq}, we must show
		\begin{flalign}
			\label{xiq2} 
			|\partial_q \xi_q (X,Y)| \lesssim y^{1/2} e. 
		\end{flalign}
		
		Recalling that $x_1 = \mathfrak{r}_q (y)$, set $X_1 = -[-x_1]_q$. Then $\xi_q (X_1, Y) = 0$ by \eqref{axy2}, so 
		\begin{flalign*} 
			\xi_q (X, Y) & =  \xi_q ( X, Y)- \xi_q (X_1, Y)  \\
		& = (B+C)(X-X_1) \big( (B+C) (X+X_1) - 2Y(A+C) - 2AB + 2 (1 - \mathfrak{q}^{-1}) ABY \big) \\ 
		& \qquad + 4ABY (X - X_1).
		\end{flalign*} 
	
		\noindent Using the fact that for any uniformly bounded $r, r' \in \mathbb{R}$ we have $|[r] - [r']| \lesssim |r-r'|$ and $|\partial_q [r] - \partial_q [r']| \lesssim |r-r'|$, with the facts that $|x -x_1| = e$, it follows that 
		\begin{flalign}
			\label{xiq3} 
			|\partial_q \xi_q (X, Y)| \lesssim e y + e \cdot |(B+C) (X+X_1) - 2 AB|.
		\end{flalign} 
		
		\noindent Denoting {$X_{SW} = -[-x^{SW}_q]_q$}, we also have that 
		\begin{flalign*} 
			|(B+C)(X+X_1) - 2AB| \lesssim  |X+X_1 - 2X_{SW}| \lesssim  |x - x_q^{SW}| + |x_1 - x^{SW}_q| \lesssim y^{1/2},
		\end{flalign*} 
	
		\noindent where the first statement follows from \Cref{xabc}, the second from the fact that $|[r]-[r']| \lesssim |r-r'|$ for any uniformly bounded $r, r' \in \mathbb{R}$; and the third from the uniform convexity of $\mathfrak{A}_{q}$ from \Cref{convex}. Inserting this into \eqref{xiq3} yields \eqref{xiq2} and thus the lemma. 
	\end{proof}

	\begin{proof}[Proof of \Cref{qlineq}]
		 Since $\mathfrak{A}_q$ varies smoothly in $q$, we have $|\mathfrak{e}_{z;q} - \mathfrak{e}_{z;q'}| \lesssim q-q'$. If either $z \notin \mathfrak{L}_q$ or $z
            \notin \mathfrak{L}_{q'}$, we have that $\max \{ \mathfrak{e}_{z;q}, \mathfrak{e}_{z;q'} \} \le q-q'$, and so the result then follows from the third part of \Cref{hde0} (if $z = (x, y)$ with $\ell_z$ being the south edge of $\mathfrak{X}$ and $x \ge x_q^{SW}$; the other cases are entirely analogous by symmetry). So, we may assume in what follows that $z \in \mathfrak{L}_q\cap \mathfrak{L}_{q'}$. 
		
		Then, again since $|\mathfrak{e}_{z;q} - \mathfrak{e}_{z;q'}| \lesssim q-q'$, it suffices to show that 
		\begin{flalign}
			\label{hqz2} 
			\mathcal{H}_q (z) - \mathcal{H}_{q'} (z) \asymp (q-q') \mathfrak{d}_z^{1/2} \cdot \max \{ \mathfrak{e}_{z;q}^{1/2}, \mathfrak{e}_{z;q'}^{1/2} \}.
		\end{flalign} 
	
		\noindent Let us first reduce to the case when $\dist (z, \mathfrak{A}_q) \le C^{-1}$ for a sufficiently large constant $C>1$. Indeed, if \eqref{hqz2} holds under this situation, then we have $\mathcal{H}_q (z) - \mathcal{H}_{q'} (z) \gtrsim q-q'$ for $z \in \mathfrak{X}$ satisfying $\dist (z, \mathfrak{A}_q) = C^{-1}$; hence, \Cref{hqomegaq} yields $\mathcal{H}_q (z) - \mathcal{H}_{q'} (z) \gtrsim q-q'$ for all $z \in \mathfrak{X}$ satisfying $\dist (z, \mathfrak{A}_q) \ge C^{-1}$, which shows the lower bound in \eqref{hqz2}. To show the upper bound, observe from the content following \eqref{xietazeta1} that $v(z)$ is uniformly smooth in $q$ for $z \in \mathfrak{X}$ satisfying $\dist (z, \mathfrak{A}_q) \ge C^{-1}$. We also have by \eqref{xietazeta1} that $\Imaginary v(z) \gtrsim 1$ there, which implies by \eqref{fxygradient} that $\nabla \mathcal{H}_q (z)$ is uniformly smooth in $q$ for $z \in \mathfrak{X}$ satisfying $\dist (z, \mathfrak{A}_q) \ge C^{-1}$. Together with the fact that $|\mathcal{H}_q (z) - \mathcal{H}_{q'} (z)| \lesssim q-q'$ for $z \in \mathfrak{X}$ satisfying $\dist (z, \mathfrak{A}_q) \ge C^{-1}$, this implies that $|\mathcal{H}_q (z) - \mathcal{H}_{q'} (z)| \lesssim q-q'$ for all $z \in \mathfrak{X}$ satisfying {$\dist (z, \mathfrak{A}_q) \ge C^{-1}$}; this confirms the upper bound in \eqref{hqz2}.  
		
		Thus, it suffices to show \eqref{hqz2} when $\dist (z, \mathfrak{A}_q) \le C^{-1}$ for sufficiently large  $C>1$. Denote $z = (x, y)$. {Recall from \Cref{rem:rightoftangency} that we are assuming  by symmetry that $\ell_z$ is the south edge of $\mathfrak{X}$, that  $y \le 49/100$ and $x>x^{SW}_q$. Therefore,  $\mathfrak{e}_{z;q} = \mathfrak{r}_q (y) - x$; this in particular implies that $\mathfrak{e}_{z;q} \ge \mathfrak{e}_{z;q'}$.} We then have
		\begin{flalign}
			\label{h3} 
			\begin{aligned}
			\mathcal{H}_q (z) - \mathcal{H}_{q'} (z) & = \displaystyle-\int_0^{\mathfrak{e}_{z;q}} \big( \partial_x \mathcal{H}_q (\mathfrak{r}_q (y) - e, y) - \partial_x \mathcal{H}_{q'} (\mathfrak{r}_{q'} (y) - e, y) \big) de \\
			& \qquad  - \displaystyle\int_{\mathfrak{e}_{z;q'}}^{\mathfrak{e}_{z;q}} \partial_x \mathcal{H}_{q'} (\mathfrak{r}_{q'} (y) - e, y) de \\
			& = \displaystyle-\int_{\mathfrak{e}_{z;q'}}^{\mathfrak{e}_{z;q}} \partial_x \mathcal{H}_{q'} (\mathfrak{r}_{q'} (y) - e, y) dy + O( (q-q') \mathfrak{d}_z^{1/2} \mathfrak{e}_{z;q}^{3/2}),
			\end{aligned} 
		\end{flalign}
	
		\noindent where the first statement holds since $\mathcal{H}_q (\mathfrak{r}_q (y),y) = 0 = \mathcal{H}_{q'} (\mathfrak{r}_{q'} (y), y)$, and the second holds by \Cref{qdifferencehq}. Next, we have 
		\begin{flalign*}
			\displaystyle-\int_{\mathfrak{e}_{z;q'}}^{\mathfrak{e}_{z;q}} \partial_x \mathcal{H}_{q'} (\mathfrak{r}_{q'} - e, y) de \asymp \displaystyle\int_{\mathfrak{e}_{z;q'}}^{\mathfrak{e}_{z;q}} e^{1/2} \mathfrak{d}_z^{1/2} de \asymp (q-q') \mathfrak{d}_z^{1/2} \mathfrak{e}_{z;q}^{1/2},	
		\end{flalign*} 

		\noindent where in the first statement we used the first statement in \Cref{hde0} and in the second we performed the integration. Inserting this into \eqref{h3} and taking $\mathfrak{e}_{z;q}$ sufficiently small (by taking $\dist (z, \mathfrak{A}_q)$ sufficiently small) yields \eqref{hqz2} and thus the lemma.
	\end{proof}

        \subsection{Level lines close to the arctic boundary}
\label{sec:lultima}
        In this section, we prove \Cref{prop:improvedconv}.
We let
          \begin{eqnarray}
            \label{eq:xys}
z_s:=          (x_s,y_s)=(s+\mathfrak l_{q;a,b,c}(h_0),\cU^{h_0}_{q;a,b,c}(s+\mathfrak l_{q;a,b,c}(h_0)))  .
          \end{eqnarray}

        Note that the condition on
        $s_-,s_+,s_0$ in the statement of the proposition ensures that
        if $s\in[s_-,s_+]$, then $ z_s$
          (with $a=b=c=1$)
          is inside the liquid region  $\mathfrak L_{q}$ and if $s=s_0$, then $z_{s_0}$ is to  the right of $p^{SW}_q$ and bounded away from $p^{SE}_q$.
          Note also that, since $\cH_{q;a,b,c}(z_s)=h_0$ (for any $a,b,c$),
          \begin{eqnarray}
            \label{dsV}
            \partial_s \mathcal U^{h_0}_{q;a,b,c}(s+\mathfrak l_{q;a,b,c}(h_0))=-\frac{\partial_x\cH_{q;a,b,c}(z_s)}{\partial_y\cH_{q;a,b,c}(z_s)}.
          \end{eqnarray}
          Since $h_0$ is smaller than $\varepsilon$, if $\varepsilon$
          is small $z_s$ is close to the arctic boundary for every
          $s\in[s_-,s_+]$ and we can assume that either
          $\me_{z_s;q}\le (1/\bfK)\md_{z_s}$ with $\bfK$ a large constant,
          or $\md_{z_s}\le (1/\bfK)\wedge(\bfK\me_{z_s;q}) $. We will
          actually distinguish three situations:
\begin{enumerate}
\item [(I)]\label{casoI}   $\me_{z_s;q}\le (1/\bfK)\md_{z_s}$ and {$z_s$ is to the right of  $p_q^{SW}$}.

\item [(II)]  \label{casoII} $\md_{z_s}\le (1/\bfK)\wedge(\bfK\me_{z_s;q}) $.

\item [(III)] \label{casoIII} $\me_{z_s;q}\le (1/\bfK)\md_{z_s}$ and $z_s$ is {to the left of $p_q^{SW}$}.  Note that as $s\to0$, $\me_{z_s;q}\sim s$  while $\md_{z_s}$ tends to $\sqrt{h_0}$. 
As in \Cref{fh}, write
$\cH_{q}(x_s,y_s)=y_s-\md_s^{1/2}\me_s^{3/2}\widetilde\cH(\hat x_s,\hat
y_s)$, with $(\me_s,\md_s)=(\me_{z_s;q},\md_{z_s})$. We can
approximate the derivatives of 
$\widetilde\cH(\hat x_s,\hat y_s)$ as in \Cref{rem:B80}.

  \end{enumerate}
Note that in cases (I) and (II), we have $\md_{z_s}\asymp s$.
Since $s_0$ satisfies \eqref{eq:condizs}, then $z_{s_0}$ falls into either  case (I) or (II) (with $a=b=c=1$).

{We will first prove (in Sections \ref{sec:item1} to \ref{sec:item5}) the statements of \Cref{prop:improvedconv} under the assumption that $s_-\ge {\bf K}h_0^{1/2}$, with ${\bf K}$ a large constant. This means that $\me_{z_s;q}\le (1/{\bf K})\md_{z_s}$ for all $s\in[s_-,s_+]$, possibly with a different (but still large) constant ${\bf K}$. Note that in this case, \Cref{314item1bis} of \Cref{prop:improvedconv} is void, since the interval $[s_-,C^{-1}h_0^{1/2}]$ is empty for $C>1$. Afterwards, in \Cref{sec:theothercase} we briefly outline the changes needed in the case where $s_-\le {\bf K}h_0^{1/2}$.}

For the proof of \Cref{prop:improvedconv}, we need to approximate precisely the level lines $\mathcal U^{h_0}$ of the height function.
The following statement will be of use in the case where $z_s$ falls into case (I) above:
  \begin{lem}
    \label{lem:sofferto}
If $\me_{z_s;q}\le (1/\bfK)\md_{z_s}$ and $z_s$ is  to the right of $p_q^{SW}$, we have
      \begin{eqnarray}
            \label{eq:dsV1}
            \partial_s \mathcal U^{h_0}_{q;a,b,c}(s+\mathfrak l_{q;a,b,c}(h_0))=\gamma'_{s}(1+\hat{\mathcal O}((\me_{z_s;q}/\md_{z_s})^{1/2})+\mathcal O((\me_{z_s;q}\md_{z_s})^{1/2}))
      \end{eqnarray}
     and
             \begin{eqnarray}
            \label{eq:dsV2}
            \partial^2_s \mathcal U^{h_0}_{q;a,b,c}(s+\mathfrak l_{q;a,b,c}(h_0))=\gamma''_s+\hat{\mathcal O}((\me_{z_s;q}/\md_{z_s})^{1/2})+\mathcal O((\me_{z_s;q}\md_{z_s})^{1/2}),
             \end{eqnarray}
              where $\gamma'_s,\gamma''_s$ are  the first and second derivative of the arctic curve computed at the point  $\bar z_s\in\mathfrak A_{q;a,b,c}$ closest to $z_s$.
              The  error terms are as in \Cref{rem:B80}.
  \end{lem}

  \begin{proof}[Proof of \Cref{lem:sofferto}]

Given $h > 0$ let $\widetilde{\mathcal{U}} = \widetilde{\mathcal U}_{q;a,b,c}^h$ denote the $\mathfrak{d}^{-1/2} \mathfrak{e}^{-3/2} h$-th level line of $\widetilde{\mathcal{H}}_{q;a,b,c}$ adopting the notations of
\Cref{fh} where
\begin{eqnarray}
  \label{eq:adopt}
\gamma':=\gamma'_s,\gamma'':=\gamma''_s,
\md:=\md_{z_s},\me:=\me_{z_s;q},   
\end{eqnarray}
namely, 
$\widetilde{\mathcal{H}}_{q;a,b,c} (\hat{s}, \widetilde{\mathcal{U}}_{q;a,b,c}^h (\hat{s})) = \mathfrak{d}^{-1/2} \mathfrak{e}^{-3/2} h$. Recall that we abbreviate $(q;a,b,c) = q$ if $(a,b,c) = (1,1,1)$ and $(q;a,b,c) = (a,b,c)$ if $q=0$, and that $\mathcal{U} = \mathcal{U}_{q;a,b,c}^h$ similarly denotes the $h$-th level line of $\mathcal{H}_{q;a,b,c}$. Let us first evaluate $\partial_s \mathcal{U},\partial_s^2 \mathcal{U}$ in terms of $\partial_{\hat{s}}\widetilde{\mathcal{U}}, \partial_{\hat{s}}^2 \widetilde{\mathcal{U}}$.

	\noindent Recall the change of coordinates in \eqref{xyxyde}, where we will eventually take $(x,y)$ to be $z_s$, so that $(x',y')=\bar z_s$. Since  the level lines of $\mathcal{H}$ and $\widetilde{\mathcal{H}}$ coincide, we find that 
	\begin{flalign}
		\label{uw1} 
	\mathcal{U} (x' + (\mathfrak{d} \mathfrak{e})^{1/2} \hat{s} - \gamma' \mathfrak{d} \mathfrak{e} \widetilde{\mathcal{U}} (\hat{s})) = y' + \hat{s} \gamma' (\mathfrak{d} \mathfrak{e})^{1/2} + \mathfrak{d} \mathfrak{e} \widetilde{\mathcal{U}}( \hat{s}).
	\end{flalign}
         Differentiating w.r.t. $\hat s$ we get
         \begin{eqnarray}
           \label{eq:aerop1}
          \cU'(x)=\frac{\gamma'+\sqrt{\md\me}\widetilde\cU'(\hat s)}{1-\gamma'\sqrt{\md\me}\widetilde \cU'(\hat s)}  , \quad \cU''(x)=\widetilde\cU''(\hat s)\frac{1+\gamma'^2}{(1-\gamma'\sqrt{\md\me}\widetilde \cU'(\hat s))^3}         
         \end{eqnarray}

	\noindent where we have denoted 
	\begin{flalign} 
          x = x' + (\mathfrak{d} \mathfrak{e})^{1/2} \hat{s} - \gamma' \mathfrak{d} \mathfrak{e} \widetilde{\mathcal{U}} (\hat{s}).
          \label{holdunchanged}
	\end{flalign}
Next, note that
\begin{equation}
  \label{eq:fDT}
          \widetilde \cU'(\hat s)=-\left.\frac{\partial_{\hat x}\widetilde \cH(\hat x,\hat y)}{\partial_{\hat y}\widetilde \cH(\hat x,\hat y)}\right|_{\hat x=\hat s,\hat y=\widetilde {\mathcal U}(\hat s)}, \quad      \widetilde \cU''(\hat s)=-\left.\frac{(A,D^2\widetilde\cH(\hat x,\hat y)A)}{\partial_{\hat y}\widetilde\cH(\hat x,\hat y)}\right|_{\hat x=\hat s,\hat y=\widetilde {\mathcal U}(\hat s)},A=(1, \widetilde \cU'(\hat s)).
        \end{equation}
        By the first of \eqref{eq:fDT}, together with \eqref{approxG} and the bounds $\gamma'\asymp \md, \me\lesssim\md$, we set  $x$ equal to $s+\mathfrak l_{q;a,b,c}(h_0)=x_s$ (recall \eqref{eq:xys}) and we deduce \eqref{eq:dsV1}. To prove \eqref{eq:dsV2}, note first of all that
        $\cU''(x)=\widetilde\cU''(\hat s)(1+\gamma'^2)(1+\mathcal O(\md^{3/2}\me^{1/2}))$.
Then,  \eqref{eq:dsV2} follows from \eqref{eq:fDT} and \eqref{approxG}.
\end{proof}

\subsubsection{Proof of 
  \Cref{314item1} of \Cref{prop:improvedconv}}
\label{sec:item1}
{Recall that we are assuming that  $\me_{z_s;q}\le (1/{\bf K})\md_{z_s}$ for all $s\in[s_-,s_+]$, in particular for $s=s_0$.}

We wish to fix $a,b,c=3-a-b$ to meet the requirements of \Cref{prop:improvedconv}. We first note that the map $(a,b,x)\mapsto (\gamma'_x,\gamma''_x,\gamma'''_x)$ is a local bijection with bounded Jacobian. Since the arctic boundary is smooth in $\frak q,a,b,c$, it is enough to prove this for $a=b=\mathfrak q=1$. Recall from \Cref{Functionsqabc1} that the map $x\mapsto \xi_x$ is a local bijection. Next, we claim that the map $(a,b,\xi_x)\mapsto (\gamma'_x,\gamma''_x,\gamma'''_x)$ is a local bijection with bounded Jacobian, at $a=b=\mathfrak q=1$. For this, one uses the explicit formulas \eqref{f021}, \eqref{eq:od} to find that the determinant of the Jacobian matrix is
\begin{eqnarray}
  \label{eq:jacobiandet}
  \frac{((-2 + 2 \xi_x - 2 \xi_x^2)^6 (-8 + 16 \xi_x - 24 \xi_x^2 + 16 \xi_x^3 - 8 \xi_x^4)}{24(1-2\xi_x)^9}.
\end{eqnarray}
As discussed just after \eqref{eq:detJ2}, in the range of parameters we are interested in, $\xi_x$ ranges from $\frac12(1-\sqrt3)+O(\varepsilon)$ and $O(\varepsilon)$ (recall that $q,|a-1|,|b-1|=O(\varepsilon)$), and it is easily checked that in this range the ratio \eqref{eq:jacobiandet} is bounded away from zero and infinity.

Consequently, we can fix the parameters $x,a,b$ (recall that $c$ is fixed to $c=3-a-b$) in order to tune the first three derivatives of the arctic boundary $\mA_{a,b,c}$. In particular (using also \Cref{lem:sofferto} to relate the derivative of the level lines to those of the arctic boundary) we can fix 
{$a-1,b-1=O(q) $} and $\frak r$ (the latter is equivalent to fixing the coordinate $x$) in such a way that \eqref{eq:meccia} holds, and at the same time (letting $z_{s_0}:=(s_0+\mathfrak l_{q}(h_0),\cU^{h_0}_{q}(s_0+\mathfrak l_{q}(h_0)))$, $z'_{s_0}:=(s_0+\mathfrak l_{q}(h_0)+\frak r,\cU^{h_0}_{a,b,c}(s_0+\mathfrak l_{q}(h_0)+\frak r)) $)
\begin{eqnarray}\label{eq:g'''''}
  \gamma''(q,1,1,1)-\gamma''(0,a,b,c)=3 C^{-1}q {s_0} ,  \qquad \gamma'''(q,1,1,1)=\gamma'''(a,b,c)
\end{eqnarray}
where $\gamma''(q,1,1,1),\gamma'''(q,1,1,1)$ denote the derivatives of the arctic boundary $\mA_q$ at the point closest to $z_{s_0}$, and $\gamma''(0,a,b,c),\gamma'''(0,a,b,c)$  denote the derivatives of the arctic boundary $\mA_{a,b,c}$ at the point closest to $z'_{s_0}$. Here, $C$ is the  positive constant (independent of
 $\bfK$) that appears in the statement of \Cref{prop:improvedconv}.

  For later convenience, we make the following observation.
  \begin{rem}\label{remarkone}
    The values $a,b,c$ satisfy $(a-1),(b-1),(c-1)=O(q)$. Moreover, as in \Cref{rem:B82}, the rescaled coordinates $(\hat x_0,\hat y_0),(\hat x_1,\hat y_1)$ of the points $z_{s_0},z'_{s_0}$ satisfy $|\hat x_0-\hat x_1|,|\hat y_0-\hat y_1|=O(\md q)=O(s_0 q)$. Finally,  since $z_{s_0},z'_{s_0}$ are on two level lines of the limit shapes, corresponding to  different parameters $(q,1,1,1)$ and $(0,a,b,c)$, but to the same height value $h_0$, from \eqref{approxG} we have that
    \begin{eqnarray}
      \frac{2(\gamma'^2+1)\hat y_0-\gamma''(q,1,1,1) \hat x^2_0}{2(\gamma'^2+1)\hat y_1-\gamma''(0,a,b,c) \hat x^2_1}=\frac{\Gamma_1^{1/3}}{\Gamma_0^{1/3}}+\mathcal O({\bf K}^{-1/2}s_0),
    \end{eqnarray}
    where (cf. \eqref{fxy2}) $\Gamma_0=\gamma''(q,1,1,1)-\gamma'(1-\gamma')q$, $\Gamma_1=\gamma''(0,a,b,c)$.
  \end{rem}

 Next, we claim that with this choice of $a,b,c$,  
\begin{eqnarray}
   \label{eq:d2Vvar}
   \partial^2_s\mathcal U^{h_0}_q(s_0+\mathfrak l_q(h_0))-\partial^2_s\mathcal U^{h_0}_{a,b,c}(s_0+\mathfrak r+\mathfrak l_{q}(h_0))\ge 2C^{-1} q s_0.
 \end{eqnarray}
This follows from \eqref{eq:g'''''} and 
 \eqref{eq:dsV2}, provided that $\bfK$ is large
 enough. Indeed, on one
 hand, the error term $\mathcal O(\sqrt{\md\me})$ in  \eqref{eq:dsV2} is smooth with
 respect to $a,b,c,q$, so that it contributes order
 $\sqrt{\md\me}O(\max(q,|a-1|,|b-1|,|c-1|)=O(\sqrt{\md\me} q)=O(s_0
 q/\bfK^{1/2})\le (10C)^{-1} q {s_0}$ to the l.h.s. of \eqref{eq:d2Vvar}, provided that $\bf K$ is large enough.   On the other hand, also the error term
 $\hat{\mathcal O}(\sqrt{\me/\md})$ contributes order
 $O(s_0 q/\bfK^{1/2})$. In fact,  by \Cref{rem:B80} it is smooth in $a,b,c$ (which change by
 order $\md q_0\asymp q s_0$) and it depends on $q$ only through $\Gamma$ which (as one sees from its definition
  \eqref{fxy2}) changes by order $\gamma'q\asymp q s_0$ when   $\log \mathfrak q$ changes from $0$ to $q$.

  We have proven \eqref{eq:d2Vvar} at $s=s_0$ and now we extend the proof (with $2C^{-1}$ replaced by $C^{-1}$)  to all $s\in[s_-,s_+]$. Thus, \Cref{314item1} then follows by integration with respect to $s$. Recall that we are working under the assumption $s_-\ge {\bf K}h_0^{1/2}$, which implies that $z_s$ falls into case (I)  described  at the beginning of \Cref{sec:lultima}.
    Since the arctic curve is smooth (in $q,a,b,c$), it follows from \Cref{lem:sofferto}, \eqref{eq:g'''''} and \eqref{eq:d2Vvar} that
    \begin{eqnarray}
   \partial^2_s\mathcal U^{h_0}_q(s+\mathfrak l_q(h_0))-\partial^2_s\mathcal U^{h_0}_{a,b,c}(s+\mathfrak r+\mathfrak l_{q}(h_0))\ge 2C^{-1} q s_0-O(q (s_+-s_-)^2)\ge C^{-1} q s_0
    \end{eqnarray}
   for all $s\in[s_-,s_+]$, where the last bound comes from the assumption that
    $s_+ - s_- \le  C^{-1}s_0^{1/2}$, taking $C$  large enough.
    \Cref{314item1} of \Cref{prop:improvedconv} is proven, under the assumption that $s_-\ge {\bf K}h_0^{1/2}$.

\subsubsection{Proof of 
  \Cref{314item2} of \Cref{prop:improvedconv}}
\label{sec:item3case1}
{Recall that we are assuming that  $\me_{z_s;q}\le (1/{\bf K})\md_{z_s}$ for all $s\in[s_-,s_+]$.}
Note first of all that
\begin{eqnarray}
  \label{eq:dhV}
  \partial_h \mathcal U^h_{q;a,b,c}(s+\mathfrak l_{q;a,b,c}(h))=\Big(\partial_y\cH_{q;a,b,c}(\mathfrak l_{q;a,b,c}(h)+s,\mathcal U^h_{q;a,b,c}(\mathfrak l_{q;a,b,c}(h)+s))\Big)^{-1},
\end{eqnarray}
since $\cH_{q;a,b,c}(\mathfrak l_{q;a,b,c}(h)+s,\mathcal U^h_{q;a,b,c}(\mathfrak l_{q;a,b,c}(h)+s))=h$. Therefore, \eqref{e:primamail} follows if we  prove
\begin{multline}
  \partial_y\cH_{q}(\mathfrak l_{q}(h_0)+s_0,\mathcal U^h_{q}(\mathfrak l_{q}(h_0) +s_0))\\\le
  \partial_y\cH_{a,b,c}(\mathfrak l_{q}(h_0)+s_0+\mathfrak r,\mathcal U^h_{a,b,c}(\mathfrak l_{q}(h_0) +s_0+\mathfrak r)).
\end{multline}

 Observe from \eqref{xyxyde} that 
\begin{flalign*} 
	\widetilde{\cH}_{q;a,b,c} \bigg( \displaystyle\frac{\overline{x} + \gamma' \overline{y}^h}{(\mathfrak{d} \mathfrak{e})^{1/2} (\gamma'^2+1)}, \displaystyle\frac{\overline{y}^h - \gamma' \overline{x}}{\mathfrak{d} \mathfrak{e} (\gamma'^2+1)} \bigg) = \displaystyle\frac{h}{\mathfrak{d}^{1/2} \mathfrak{e}^{3/2}},
\end{flalign*} 
\noindent where 
\begin{flalign*}
	\overline{x} = s - x'; \qquad \overline{y}^h = \mathcal{U}_{q;a,b,c}^h (s) - y'.
\end{flalign*}

\noindent Differentiating with respect to $h$ yields 
\begin{flalign*} 
	\displaystyle\frac{\partial_h \mathcal{U}_{q;a,b,c}^h}{\gamma'^2+1} \cdot \Bigg( & (\mathfrak{d} \mathfrak{e})^{-1/2} \gamma' \partial_{\hat{x}} \widetilde{\mathcal{H}}_{q;a,b,c} \bigg( \displaystyle\frac{\overline{x} + \gamma' \overline{y}^h}{(\mathfrak{d} \mathfrak{e})^{1/2} (\gamma'^2+1)}, \displaystyle\frac{\overline{y}^h - \gamma' \overline{x}}{\mathfrak{d} \mathfrak{e} (\gamma'^2+1)} \bigg) \\
	& \qquad + (\mathfrak{d} \mathfrak{e})^{-1} \partial_{\hat{y}} \widetilde{\mathcal{H}}_{q;a,b,c} \bigg( \displaystyle\frac{\overline{x} + \gamma' \overline{y}^h}{(\mathfrak{d} \mathfrak{e})^{1/2} (\gamma'^2+1)}, \displaystyle\frac{\overline{y}^h - \gamma' \overline{x}}{\mathfrak{d} \mathfrak{e} (\gamma'^2+1)} \bigg) \Bigg)  = \mathfrak{d}^{-1/2} \mathfrak{e}^{-3/2}.
\end{flalign*} 

\noindent Thus, 
\begin{flalign}
	\label{vhh} 
	\partial_h \mathcal{U}_{q;a,b,c}^h  = \mathfrak{d}^{1/2} \mathfrak{e}^{-1/2} (\gamma'^2+1) \cdot \big( (\mathfrak{d} \mathfrak{e})^{1/2} \gamma' \partial_{\hat{x}} \widetilde{\mathcal{H}}_{q;a,b,c} (\hat{x}, \hat{y}) + \partial_{\hat{y}} \widetilde{\mathcal{H}}_{q;a,b,c} (\hat{x}, \hat{y}) \big)^{-1}.
\end{flalign}

  Now we specialize to $s=s_0$. To prove \Cref{314item2} of \Cref{prop:improvedconv}, it is enough to show that, when changing parameters from $(0;a,b,c)$ to $(q,1,1,1)$ (with $a,b,c$ the values chosen in the proof of \Cref{314item1}), the denominator of \eqref{vhh} decreases. The derivative $\partial_{\hat y}\widetilde\cH$ is related to $\hat f$ via \eqref{derivativeshf}, where the ratio $\md/|\gamma'|\asymp 1$,  and \eqref{fapproximate} gives an approximation for $\hat f$. Recalling \Cref{remarkone} and in particular the fact that the parameter
  $\Gamma$ in \eqref{fxy2} increases by $\Gamma_1-\Gamma_0\asymp\gamma' q\asymp \md q\asymp q s_0$, we see that  $\partial_{\hat y}\widetilde\cH$ decreases by order $q s_0$. The derivative $\partial_{\hat x}\widetilde\cH$ also changes by order $q s_0$ but, since it is multiplied by $(\md\me)^{1/2}\gamma'\lesssim {\bf K}^{-1/2}$ in \eqref{vhh}, we do not need to determine the sign of this change. Altogether the denominator of \eqref{vhh} decreases, and \eqref{e:primamail} follows,  under the assumption that $s_-\ge {\bf K}h_0^{1/2}$.

\subsubsection{Proof of 
  \Cref{314item3} of \Cref{prop:improvedconv}}
 {Recall that we are assuming that  $\me_{z_s;q}\le (1/{\bf K})\md_{z_s}$ for all $s\in[s_-,s_+]$.}
        	The denominator on the right side of \eqref{vhh} is bounded below and changes by $O(\md_{z_s} q)$ upon changing $(q;a,b,c)$ from $(q;1,1,1)$ to $(0;a,b,c)$. Therefore, it follows that 
		\begin{flalign} 		\label{eq:nogdt}
			|\partial_h \mathcal{U}_q^h(\mathfrak l_q(h_0)+s) - \partial_h \mathcal{U}_{a,b,c}^h(\mathfrak l_q(h_0)+s+\frak r)| \lesssim  \md_{z_s}^{3/2} \me_{z_s;q}^{-1/2} q \asymp  q\md_{z_s}^{5/3} h^{-1/3}\lesssim q s_+^{5/3}h^{-1/3},
		\end{flalign} 
		where in the second step we use the third part of  \Cref{hde0}.	
	The bound \eqref{eq:primamail'} then follows by integrating with respect to $h$. 
                \subsubsection{Proof of \Cref{heighthabcq2} of \Cref{prop:improvedconv}}
                \label{sec:item5}
By symmetry, we may show \eqref{heighthabcq} when $x=(x,y)$ is such that $x \ge x_q^{SW}$, so that $s \asymp \mathfrak{d}_z$. {Let $z_0=(x_0,y_0)$ with $x_0 = \mathfrak{l}_q (h_0) + s$ and $y_0 = \mathcal{U}_q^{h_0} (x_0)$.}
By \eqref{eq:meccia}, \eqref{eq:meccia2}, and \eqref{eq:weaker}, we have that $\mathcal{U}_q^{h_0} (x_0) \ge \mathcal{U}_{a,b,c}^{h_0} (x_0 + \mathfrak{r}) + \mathfrak{r}'$. {In particular, the left-hand side of \eqref{heighthabcq} is non-negative at $z_0$. Hence, it is enough to prove the claim at $z=(x_0,y)$ with $y<y_0$. Let $x_0\ge x^{SW}_q$ and $y,y'\le y_0$ be such that  $\cH_q(x_0,y')=\cH_{a,b,c}(x_0+\frak r,y+\frak r')=h \in [0, h_0]$. } For any $h \in [0, h_0]$, we have 
\begin{flalign*}
y'-y=	\mathcal{U}_q^h (x_0) - (\mathcal{U}_{a,b,c}^h (x_0 + \mathfrak{r}) + \mathfrak{r}') \gtrsim -qs^{5/3} \displaystyle\int_0^{h_0} v^{-1/3} dv \asymp -q \mathfrak{d}_{z_0}^{5/3} h_0^{2/3},
\end{flalign*}
\noindent where the first inequality follows {from integrating \eqref{eq:nogdt} (with $\md_{z_s}$ not yet replaced by $s_+$)} (and using the bound $\mathcal{U}_q^{h_0} (x_0) \ge \mathcal{U}_{a,b,c}^{h_0} (x_0 + \mathfrak{r}) + \mathfrak{r}'$), and the second follows from the fact that $s \asymp \mathfrak{d}_{z_0}$.
{Write
\begin{eqnarray}
 \cH_{a,b,c}(x_0+\frak r,y+\frak r')= \cH_q(x_0,y')=\cH_q(x_0,y)+\int_y^{y'}\partial_u \cH_q(x_0,u)d u
\end{eqnarray}
and recall \eqref{eq:lpa2}. We have
$\mathcal{H}_q (w)^{1/3} \mathfrak{d}_{w}^{-2/3} \asymp (\mathfrak{d}_{w}^{-1} \mathfrak{e}_{w;q})^{1/2}$ (by the third item of \Cref{hde0}) and we prove below that
\begin{eqnarray}
  \label{slightlyoff}
(\md_z^{-1} \me_{z;q})^{1/2} \lesssim (\md_{z_0}^{-1} \me_{z_0;q})^{1/2}.  
\end{eqnarray}
for $z=(x_0,y),y\le y_0$. Then, we deduce
\begin{eqnarray}
 \cH_{a,b,c}(x_0+\frak r,y+\frak r')- \cH_q(x_0,y)\gtrsim-q\md_{z_0}^{5/3}h_0^{2/3}(\md_{z_0}^{-1} \me_{z_0;q})^{1/2}\asymp -q \md_{z_0} h_0
\end{eqnarray}
which confirms \eqref{heighthabcq}.  It remains to prove
\eqref{slightlyoff}. To that purpose, first observe that since $\me_{z;q}\lesssim \md_{z}$ for all $z \in \mathfrak{L}_q$, we can assume that the r.h.s. of \eqref{slightlyoff} is smaller than some small constant $\epsilon$. Next, note that the slope $\rho_y$ of the arctic boundary with vertical coordinate $y$ satisfies $\rho_y\le C y^{1/2}$ for some constant $C$ (as  the arctic boundary has uniformly positive and bounded curvature). Increasing {$y=\md_w^2$} by $\delta$
increases $\me_{z;q}$ by $\delta \rho_y^{-1} (1+o(\delta))$. Hence,
$\partial_{y} (\me_{w;q}/\md_w) =  (\frac1{\rho_y\md_{w;q}} - \frac{\me_{w;q}}{2\md_w^3})\ge \frac1{\md_w^2}(\frac1{C}-\frac{\me_{w;q}}{2\md_w})$. For $w=z_0$, the quantity in parenthesis is larger than $1/C-\epsilon^2/2>0$, where positivity holds choosing $\epsilon$ small enough. Therefore, the ratio $\me_{w;q}/\md_w$ decreases when moving $w$ from $z_0$ vertically downward. As a consequence, the quantity in parenthesis is positive for all $w$ between $z$ and $z_0$, and \eqref{slightlyoff} follows.

\subsubsection{Proof of \Cref{prop:improvedconv} in the case $s_-\lesssim h_0^{1/2}$}
\label{sec:theothercase}
  The proof of \Cref{prop:improvedconv} requires some additional arguments in the case $s_-\lesssim h_0^{1/2}$. The reason is that, in that case, the interval $[s_-,s_+]$ contains points for which  $\md_{z_s}$ is of the same order as $\me_{z_s;q}$. In this regime, the limit shape and its level lines cannot be approximated via \Cref{rem:B80} and we will rely on \Cref{lem:cella2H} instead.
In the present section, we will  show how the proof of \Cref{314item1} is modified with respect to the case  $s_-\ge {\bf K}h_0^{1/2}$ (similar considerations hold also  for \Cref{314item3}--\Cref{heighthabcq2}, but we will skip details there) and we will prove \Cref{314item1bis} (recall that in the case $s_-\ge {\bf K}h_0^{1/2}$ this statement is void) as well as \Cref{314item2}.
Let us additionally assume that $\md_{z_ {s_0}}\le (1/\bfK)\wedge(\bfK\me_{z_{s_0};q}) $, since in the case $\me_{z_{s_0};q}\le (1/{\bf K})\md_{z_ {s_0}}$ the proof is closer to that in \Cref{sec:item1}.

We need to prove that we can find parameters $a,b,c\in[1-C\varepsilon,1+C\varepsilon]$, as well as $\frak r,\frak r'$, such  that \eqref{eq:meccia2}, \eqref{eq:weaker} and \eqref{e:primamail} hold, as well as \eqref{eq:meccia} and \eqref{eq:r'}. The choice of $\frak r$ will correspond to the choice of $\hat x_1$ below, while $\frak r'$ is defined by \eqref{eq:r'} itself.
We start with \eqref{eq:meccia2} and we restrict our attention to the case where $z_s$ is either to the left of $p^{SW}_q$, or it is to the right and $\md_{z_ {s}}\le (1/\bfK)\wedge(\bfK\me_{z_{s};q}) $, since in the opposite case the proof is similar to that in \Cref{sec:item1}.  We use the notation of \Cref{lem:cella2H}, and write
\begin{eqnarray}
  \label{eq:cnfs}
x_s=x_q^{SW}+\md \hat x_s, y_s=\md^2 \hat y_s,\qquad \widetilde\cH(\hat x_s,\hat y_s)=\md^{-2}\cH_{q;a,b,c}(z_s),  \qquad \md=\md_{z_{s_0}}\asymp s_0\asymp h_0^{1/2}.
\end{eqnarray}
As shown in \Cref{Proof00}, the map $(a,b)\mapsto (\gamma''(q,a,b,c),\gamma'''(q,a,b,c))$ (with $c$ fixed to $3-a-b$ and $\gamma''(q,a,b,c),\gamma'''(q,a,b,c)$ denoting the second and third derivatives of the arctic boundary at $p^{SW}_{q;a,b,c}$) is a local bijection with uniformly bounded Jacobian. 
We then choose $a,b,c$ such that $a-1,b-1,c-1=O(q)$ and
\begin{eqnarray}
  \label{eq:TheChoice}
  \gamma''(q,1,1,1)-\gamma''(0,a,b,c)=-q \md (2\hat x_{s_0}+\eta)\\
  \gamma'''(q,1,1,1)-\gamma'''(0,a,b,c)=3q\gamma''(q,1,1,1)
\end{eqnarray}
for some small $\eta>0$, and we wish to show that with this choice \eqref{eq:meccia2} holds.
Thanks to \Cref{lem:cella2H}, we can replace the level line $\mathcal U^{h_0}_{q;a,b,c}$ (or equivalently, the corresponding level line of the rescaled height function $\widetilde \cH$) with the $h:=h_0/\md^2\asymp 1$ level line of $\frak H_{\gamma''(q,a,b,c),\gamma'''(q,a,b,c)}$. In fact, when comparing the second derivatives of the level lines with parameters $(q,1,1,1)$ and $(0,a,b,c)$, the error term $\mathcal O(\md^2)$ in \Cref{lem:cella2H} translates into an error term of order $q\md^2$, where the factor $q$ comes from $a-1,b-1,c-1=O(q)$. Note that $q\md^2\asymp q s_0\md$ is by an order $\md\le (1/\bf K)$ smaller than the gain in curvature of order $qs_0$ that we are after in \eqref{eq:meccia2}, and therefore we can safely neglect it.

To proceed, we need to state a few facts about the explicit functions $\frak H_{q,\gamma'',\gamma'''}$ and $\mathtt H_{\gamma''}$ defined in  \Cref{lem:cella2H}. The proofs are based on elementary  calculus, and can be found in \Cref{sec:AI}. We begin with the following:
\begin{lem}\label{prop:AI} Fix $\gamma''>0,\gamma'''\in\mathbb R$. Let $\mu\in\mathbb R$, $\mathfrak H_{q,\gamma'',\gamma'''}$ be as in \eqref{approxG2} and let
  \begin{eqnarray}\label{eq:C138}
 		X=\hat x+\frac{\mathfrak d q}{2\gamma''} (\mu \hat x+\gamma'' \hat x^2);
		\qquad
		Y={\hat y}.
  \end{eqnarray}
  Then,
  \begin{eqnarray}
    \label{eq:AI1}
    \frak H_{q,\gamma''+\mu q \md,\gamma'''+3\gamma'' q}(\hat x,{\hat y})=    \frak H_{0,\gamma'',\gamma'''}(X,Y)+\mathcal O(q \md^2+q^2\md).
  \end{eqnarray}
\end{lem}

          \begin{definition}\label{def:PhiMQ}
            With the notations  of \Cref{prop:AI}, set  $Q = \mathfrak{d} q / 2$ and $M = \mu/\gamma''$. Define 
                \begin{equation}
                \Phi_{M,Q}({\hat x}):={\hat x}+Q(M{\hat x}+ {\hat x}^2),
              \end{equation}
              so that the change of variables \eqref{eq:C138} from $(X,Y)$ to $({\hat x},{\hat y})$ is $(X,Y)=(\Phi_{M,Q}({\hat x}),{\hat y})$. Define $G_{M,Q}({\hat x},{\hat y}):=\mathtt H_{\gamma''}\bigl(\Phi_{M,Q}({\hat x}),{\hat y}\bigr).$ Given $h>0$ and ${\hat x}$ such that ${\hat x},\Phi_{M,Q}({\hat x})>-\sqrt{2h/\gamma''}$, let $U_h({\hat x})$ and $U_{h;M,Q}({\hat x})$ denote the unique values such that
              \begin{equation}
                \mathtt H_{\gamma''}({\hat x},U_h({\hat x}))=
                G_{M,Q}({\hat x},U_{h;M,Q}({\hat x}))\Bigl(=\mathtt H_{\gamma''}(\Phi_{M,Q}({\hat x}),U_{h;M,Q}({\hat x}))\Bigr)=h.
                \label{gaqu}
              \end{equation}

          \end{definition}
             \begin{rem}\label{rem:HhH} From  \Cref{prop:AI} it follows that if we wish to estimate (neglecting error terms $\mathcal O(q^2 \md+q \md^2)$) the difference between the derivatives of the level lines at level $h$  of $ \frak H_{q,\gamma''_0,\gamma'''+3\gamma''_0 q}$ and of $ \frak H_{0,\gamma'',\gamma'''}$ with $\gamma''_0-\gamma''=\mu q \md$, it is enough to compare the corresponding level lines of  $\mathtt H_{\gamma''}(\Phi_{M,Q}({\hat x}),{\hat y})$ and of $\mathtt H_{\gamma''}({\hat x},{\hat y})$, that is, $U_{h;M,Q}$ and $U_{h;0,0}$.
              \end{rem}

	The next lemma compares the derivatives with respect to $h$ of $U_{h;M,Q}$ and $U_{h;0,0}$ at points where the slopes coincide.
	\begin{lem}
		\label{lem:first-order-height-comparison}Fix $h>0$ and a base point ${\hat x}_0\ge 0$. Let ${\hat x}_1={\hat x}_1^{[M]}$ be chosen so that
	\begin{align*}
		U_{h;M,Q}'({\hat x}_0)=U_{h;0,0}'({\hat x}_1)
	\end{align*}
	(suppressing from the notation the dependence of $\hat{x}_1$ on $\hat{x}_0$ and $Q$.)
		The choice $M=-2{\hat x}_0-\eta$, for sufficiently small $\eta>0$, has the following property: for each fixed $h'\in(0,h]$, there is a constant
		$c=c(h',h,{\hat x}_0)>0$, such that
		\begin{align}
			\partial_hU_{h';M,Q}({\hat x}_0)-\partial_hU_{h';0,0}({\hat x}_1^{[-2 {\hat x}_0-\eta]})
			\ge cQ\eta.
			\label{eq:common-A-h-derivative-quantitative-claim}
		\end{align}
                For $0<\varepsilon<h$ one has $\inf\{ c(h',h,{\hat x}_0):h'\in [\varepsilon,h]\}>0$.
	\end{lem}

Finally, the next lemma allows to compare the second derivatives with respect to ${\hat x}$ of $U_{h;M,Q}$ and $U_{h;0,0}$ in the neighborhood of points where the slopes coincide.

	\begin{lem}
          \label{lem:3.1}
		Let $h>0,\delta>0,\eta>0,{\hat x}_0\ge 0,w\in\mathbb R$, write ${\hat x}'={\hat x}_0+w$, and choose 
		$M=-2{\hat x}_0-\eta$. 
                There exists a constant $c>0$ such that
		 for fixed $h$, ${\hat x}_0$, and $w$ satisfying $w+{\hat x}_0>\delta-\sqrt{2h/\gamma''}$, it holds for $\eta$ small enough
		\begin{align}
			&\bigl(U_{h;M,Q}({\hat x}')-U_{h;M,Q}({\hat x}_0)\bigr)
			-\bigl(U_{h;0,0}({\hat x}_1+w)-U_{h;0,0}({\hat x}_1)\bigr)
			\ge c \delta^{1/2} Q w^2
                          +\mathcal{O}(Q^2),
			\label{eq:shifted-height-final-statement}
		\end{align}
                as $Q\to0$,
		where ${\hat x}_1={\hat x}_1^{[-2{\hat x}_0-\eta]}$. 
	\end{lem}
        
        We go back to the proof of \eqref{eq:meccia2}.
        Thanks to \Cref{rem:HhH}, within the level of approximation we seek in \Cref{prop:improvedconv}, we can replace the $h$-level lines of $\frak H_{\gamma''(q,1,1,1),\gamma'''(q,1,1,1)}$ and $\frak H_{\gamma''(0,a,b,c),\gamma'''(0,a,b,c)}$
with $U_{h;M,Q}$ and $U_{h;0;0}$ respectively, defined in \Cref{def:PhiMQ}. Then, recalling the rescaling of variables \eqref{eq:cnfs}, \eqref{eq:meccia2} follows directly from \Cref{lem:3.1}, where the assumption $w+\hat x_0>\delta-\sqrt{2h/\gamma''}$ there (with $\delta\asymp C^{-1}$) is guaranteed by the assumption $s\ge C^{-1}h_0^{1/2}$ in \Cref{314item1} (indeed, recall from \eqref{eq:cnfs} that the variable $\hat x$ is rescaled by a factor $\md\asymp h_0^{1/2}$  with respect to $s$).

Next, we prove
    \Cref{314item2}. Here we only discuss the case when
    $h \ge \mathbf{K}^{-1} h_0$, as the case $h < \mathbf{K}^{-1} h_0$
    can be addressed very similarly to as in \Cref{sec:item3case1} (as
    then we are in the regime where
    $\mathfrak{e} \le \mathbf{K}^{-1} \mathfrak{d}$, where the
    approximations used there apply). In the case $h \ge \mathbf{K}^{-1} h_0$, \Cref{314item2}
    follows directly from \Cref{lem:first-order-height-comparison}
    (keeping \Cref{rem:HhH} in mind). Note that, with respect to
    \eqref{e:primamail},
    \eqref{eq:common-A-h-derivative-quantitative-claim} gives a
    quantitative lower bound of order $Q\asymp q \md$ on the
    difference of $h$-derivatives of the $(q,1,1,1)$ and $(0,a,b,c)$
    level lines. This positive margin allows to neglect the errors
    (that are of higher order in $q$ and/or $\md$, thanks to
    \Cref{lem:cella2H} and \Cref{rem:HhH}) involved in replacing this
    difference of level lines of $\widetilde\cH$ with the difference
    of the lines $U_{h;M,Q},U_{h;0;0}$ defined in
    \Cref{def:PhiMQ}. 

Finally, we prove \Cref{314item1bis}. 
Let    $s\in[s_-,C^{-1} h_0^{1/2}]$, that is, 
 $z_s$ satisfies condition (III) described  at the beginning of \Cref{sec:lultima}: $\me_{z_s;q} \lesssim \md_{z_s}/C$ for some large $C>1$  and $z_s$ is to the left of $p^{SW}_q$ (for the rest of this proof, we will write $\me,\md$ instead of $\me_{z_s;q},\md_{z_s}$). In this case,  the level lines of $\mathcal{H}$ and $\widetilde{\mathcal{H}}$ no longer coincide, since we are in the case $\gamma'<0$ of \eqref{twocases}. In fact, {the $h$-level line $\mathcal U=\mathcal U^h_{q;a,b,c}$ of $\cH_{q;a,b,c}$ coincides with the $\mathfrak{d}^{-1/2} \mathfrak{e}^{-3/2} h$ level line   $\widetilde{\mathcal{U}} =  \widetilde{\mathcal{U}}_{q;a,b,c}^h$  of }
	\begin{flalign} \label{nolonger}
		\mathfrak{d}^{-1/2} \mathfrak{e}^{-3/2} y - \widetilde{\mathcal{H}}_{q;a,b,c} (\hat{x}, \hat{y}) =  \mathfrak{d}^{-1/2} \mathfrak{e}^{-3/2} (y' + (\mathfrak{d} \mathfrak{e})^{1/2} \gamma' \hat{x} + \mathfrak{d} \mathfrak{e} \hat{y}) - \widetilde{\mathcal{H}}_{q;a,b,c} (\hat{x}, \hat{y}).
	\end{flalign}

	\noindent Formulas \eqref{uw1} to \eqref{holdunchanged} hold unchanged.

        To determine $\widetilde{\mathcal U}$, we set $\hat y=\widetilde{\mathcal U}(\hat x)$ in \eqref{nolonger}, differentiate with respect to $\hat x$ and impose the derivative to equal zero. Using
 \eqref{approxG} to approximate the derivatives of $\widetilde\cH$ with those of $\mathfrak G$, we get immediately
 (letting $\Delta:=2(\gamma'^2+1)\widetilde {\mathcal U}(\hat s)-\gamma''\hat s$)
        \begin{eqnarray}
          \label{eq:Uy'}
          \widetilde {\mathcal U}'(\hat s)=-\frac{\frac{\gamma'}{\sqrt{\md \me}}+3 \Xi\Delta^{1/2}\gamma''\hat s\sqrt{\frac\me\md}+\mathcal O(\me)+\hat{\mathcal O}(\me/\md)}{1-3 \Xi\Delta^{1/2}(1+\gamma'^2)\sqrt{\frac\me\md}}.
        \end{eqnarray}
        Plugging this into the first of \eqref{eq:aerop1}, we deduce
        \begin{eqnarray}
          \label{plu}
          \mathcal U'(x)=3\Xi \Delta^{1/2}
\frac{  |\gamma'|\sqrt{\frac\me\md}-\frac{\gamma''}{(1+\gamma'^2)}\hat s\me }{1-3\Xi \Delta^{1/2}\sqrt{\frac\me\md} }\Big(1+\mathcal O(\sqrt{\me\md})+\hat{\mathcal O}((\me\md)^{1/2})\Big).
        \end{eqnarray}
This holds for general $(q,a,b,c)$. Now we compare the cases $(q,a,b,c)=(q,1,1,1)$ and $(q,a,b,c)=(0,a,b,c)$, 
where the values $a,b,c$ are those fixed so that \eqref{eq:TheChoice} holds.

          Taking $x=s+\mathfrak l_q(h_0)$, we claim that
  \begin{eqnarray}
              \label{argue}
\left|            \partial_s\mathcal U_q^{h_0}(s+\mathfrak l_q(h_0))-\partial_s\mathcal U_{a,b,c}^{h_0}(s+\mathfrak l_q(h_0)+\mathfrak r)\right| \lesssim \sqrt{\me \md} q s_0
           \lesssim C^{-1/2} q s_0^2 \end{eqnarray}
            for all { $s\le s_1:= C^{-1}h_0^{1/2}$}, where the second bound follows from \eqref{eq:cnfs} and $\me\lesssim \md/C$. The first estimate in \eqref{argue} follows quickly from \eqref{plu} together with the facts that $|\gamma'| \asymp \mathfrak{d}$, and that $\gamma'$ and $\gamma''$ change by $\mathcal{O}(q\mathfrak{d}^2)$ and $\mathcal{O}(q\mathfrak{d})$, respectively, from passing from $\mathfrak{A}_q$ to $\mathfrak{A}_{a,b,c}$ (where that $\gamma'$ changes by $\mathcal{O}(q\mathfrak{d}^2)$ follows from the facts that $\gamma' = 0$ in both $\mathfrak{A}_q$ and $\mathfrak{A}_{a,b,c}$ at their SW tangency locations, and that $\gamma''$ changes by $\mathcal{O}(q\mathfrak{d})$ when passing from $\mathfrak{A}_q$ to $\mathfrak{A}_{a,b,c}$).

            We have then \eqref{argue} for $s\le s_1$ and, on the other hand, for $s= s_1$ we deduce from \eqref{eq:shifted-height-final-statement}
	 (where in our case $\delta\asymp C^{-1}$, as indicated above)  that
         \begin{eqnarray}
              \label{450}
              \mathcal U^{h_0}_q(\mathfrak l_q(h_0)+s_1)- (\mathfrak r'+ \mathcal U^{h_0}_{a,b,c}(\mathfrak l_q(h_0)+s_1+\frak r))\gtrsim  C^{-1/2}q s_0^3,
            \end{eqnarray}
            where we used the fact that $s_0-s_1\gtrsim s_0$, that follows from the assumptions of \Cref{prop:improvedconv}, as well as $s\le C^{-1}h_0^{1/2}$.
            Integrating \eqref{argue} from $s$ to  $s_1$, noting that $s_1-s\lesssim  C^{-1} s_0$ and using  \eqref{450} we deduce
             \begin{flalign*}
              \mathcal U^{h_0}_q(\mathfrak l_q(h_0)+s)- (\mathfrak r'+ \mathcal U^{h_0}_{a,b,c}(\mathfrak l_q(h_0)+s+\frak r))\gtrsim C^{-1/2}q s_0^3-q s^3_0C^{-3/2} \gtrsim C^{-1/2}q s_0^3>0,
             \end{flalign*}
             where in the last bound we  assumed that $C$ is large enough. \Cref{314item1bis} is proven.

\subsubsection{Properties of $\frak H_{q,\gamma'',\gamma'''}$, $\mathtt H_{\gamma''}$ and their level lines}
\label{sec:AI}
In this section, we prove the statements of \Cref{prop:AI}, \Cref{lem:first-order-height-comparison} and \Cref{lem:3.1} about the function $\mathfrak H$ and its level lines.

	\begin{proof}[Proof of \Cref{prop:AI}]
Write for lightness of notation $\psi:=\gamma''+\mu q\mathfrak d$ and $
		S(\hat x,{\hat y}):=\sqrt{2{\hat y}-\gamma'' \hat x^2}$.
	For $\mathtt H_{\gamma''}$ one computes 
	\begin{align}\label{derivativesttH}
		\hat{x} \partial_{x} \mathtt H_{\gamma''}({\hat x},{\hat y})
		&=-\frac{\hat{x} \sqrt{\gamma''}}{\pi}S({\hat x},{\hat y})=2\gamma''	\partial_{\gamma''}\mathtt H_{\gamma''}({\hat x},{\hat y});    \quad \partial_y\mathtt H_{\gamma''}({\hat x},{\hat y})= \frac12-\frac1\pi\arctan(\frac{\sqrt{\gamma''}{\hat x}}{S({\hat x},{\hat y})}).
	\end{align}
A Taylor expansion  gives
	\begin{align*}
		\mathtt H_{\psi}({\hat x},{\hat y})
		&=\mathtt H_{\gamma''}(X,Y)-\frac{\mathfrak d q}{2\gamma''}
		\Bigl(\mu X+  \gamma'' X^2   \Bigr)
                  \,\partial_x\mathtt H_{\gamma''}(X,Y)
		+\mu q\mathfrak d\,\partial_{\gamma''}\mathtt H_{\gamma''}(X,Y)
		+\mathcal{O}(\mathfrak d^2q^2)
		\\
		&=\mathtt H_{\gamma''}(X,Y)
		+\frac{\mathfrak d q}{\pi\sqrt\gamma''}\frac{\gamma'' X^2}{2}   S(X,Y)
		+\mathcal{O}(\mathfrak d^2q^2).
	\end{align*}
	
	\noindent Also, in the term in  $\frak H_{q,\psi,\gamma'''+3\gamma'' q}$ proportional to  $\mathfrak d$, we may replace
	${\hat x},{\hat y},\gamma''+\mu q \md$ there by $X,Y,\gamma''$ up to an error $\mathcal{O}(\mathfrak d^2q)$, obtaining
	\begin{multline*}
		\frac{\mathfrak d}{\pi\sqrt{\gamma''}}
		\left[
		\left(\frac{\gamma''}3+\frac{2\gamma'''}{9\gamma''}\right)S(X,Y)^3
		-\frac{\gamma'''}{3\gamma''}Y S(X,Y)
		\right] \\
		 +\frac{\mathfrak d q}{\pi\sqrt{\gamma''}}
		\left[  \frac12 S(X,Y)^3	-    Y S(X,Y)
		\right]
		+\mathcal{O}(\mathfrak d^2 q).
	\end{multline*}
        
	Using $S(X,Y)^3=(2Y-\gamma''X^2)S(X,Y)$, the part proportional to  $\mathfrak d q$ becomes
	\begin{align*}
	-	\frac{\mathfrak d q}{\pi\sqrt{\gamma''}}\frac{\gamma'' X^2}2     S(X,Y),
	\end{align*}
	which cancels exactly with the $\mathfrak d q$ contribution from $\mathtt H_{\psi}({\hat x},{\hat y})$. Therefore, the claim of the lemma follows.	
	\end{proof}

\begin{proof}[Proof of \Cref{lem:first-order-height-comparison}]
		For notational convenience, write $U_{h'}:=U_{h';0,0}$ and
		$\widetilde U_{h'}:=U_{h';M,Q}$. By \eqref{gaqu},
		\begin{align}\label{eq:UU}
			\widetilde U_{h'}(t)=U_{h'}\bigl(\Phi_{M,Q}(t)\bigr).
		\end{align} Take first $\eta=0$. 
		Evaluating at the fixed base point ${\hat x}_0$ and using $M=-2{\hat x}_0$ gives for every $h'\in(0,h]$
		\begin{align}
			\widetilde U_{h'}({\hat x}_0)=U_{h'}\bigl({\hat x}_0+Q(M {\hat x}_0+{\hat x}_0^2)\bigr)
                  =U_{h'}({\hat x}_0-Q{\hat x}_0^2).
                  			\label{eq:h-derivative-perturbed-exact}
		\end{align}
		Differentiating \eqref{eq:UU} with respect to $t$ and then
		evaluating at $t={\hat x}_0,h'=h$ gives
		\begin{align}
			\widetilde U_h'({\hat x}_0)
			&=U_h'({\hat x}_0-Q{\hat x}_0^2)\bigl(1+Q(M+2{\hat x}_0)\bigr)
			=U_h'({\hat x}_0-Q{\hat x}_0^2),
			\label{eq:xprime-exact-slope-match}
		\end{align}
		because $M+2{\hat x}_0=0$. By the definition of ${\hat x}_1$, this shows that
		\begin{align}
			{\hat x}_1^{[-2{\hat x}_0]}={\hat x}_0-Q{\hat x}_0^2.
			\label{eq:xprime-common-choice}
		\end{align}
	
		\noindent Substituting \eqref{eq:xprime-common-choice} into
		\eqref{eq:h-derivative-perturbed-exact} yields $\widetilde U_{h'}({\hat x}_0)=U_{h'}({\hat x}_1^{[-2 {\hat x}_0]})$.
		Differentiating this identity with respect to the level parameter $h'$ gives
                	\begin{align}
			\partial_hU_{h';M,Q}({\hat x}_0)-\partial_hU_{h';0,0}({\hat x}_1^{[-2 {\hat x}_0]})
			=0.
			\label{eq:common-A-h-derivative-claim}
		\end{align}

		For the strict comparison, assume $Q>0$ and choose $M=-2{\hat x}_0-\eta$. Put
		\begin{align*}
			\Phi:=\Phi_{M,Q}({\hat x}_0)={\hat x}_0+Q(M{\hat x}_0+{\hat x}_0^2)={\hat x}_0-Q{\hat x}_0^2-Q\eta {\hat x}_0.
		\end{align*}
		Then for every $h'\in(0,h]$,
		\begin{align*}
			\widetilde U_{h'}({\hat x}_0)=U_{h'}(\Phi),
		\end{align*}
		so it remains to compare $\partial_hU_{h'}(\Phi)$ with $\partial_hU_{h'}({\hat x}_1)$, with ${\hat x}_1={\hat x}_1^{[-2{\hat x}_0-\eta]}$.
		
		We first record two monotonicity facts.
		Let
                \begin{equation}
                  \label{eq:thetayh}
                \theta=\theta({\hat x},h'):=\pi\left.\partial_y\mathtt H_{\gamma''}({\hat x},{\hat y})\right|_{{\hat y}=U_{h'}({\hat x})}=\frac{\pi}{\partial_h U_{h'}({\hat x})}\in(0,\pi)  .
                \end{equation}
                 From \eqref{derivativesttH} we have $S({\hat x},U_{h'}({\hat x}))=\sqrt{\gamma''}{\hat x} \tan(\theta)$ and therefore
		\begin{align}
			U_{h'}'({\hat x})=-\frac{\partial_x \mathtt H_{\gamma''}({\hat x},U_{h'}({\hat x}))}{\partial_y \mathtt H_{\gamma''}({\hat x},U_{h'}({\hat x}))}=\sqrt{\gamma''}\frac{S({\hat x},U_{h'}({\hat x}))}{\theta}=\gamma'' {\hat x}\,\frac{\tan\theta}{\theta}.\label{uprimo}
		\end{align}
		Substituting $\sqrt{2U_{h'}({\hat x})-\gamma'' {\hat x}^2}=\sqrt{\gamma''}{\hat x}\tan\theta$ into the identity
		$\mathtt H_{\gamma''}({\hat x},U_{h'}({\hat x}))=h'$ yields
		\begin{align*}
			h'=\frac{\gamma'' {\hat x}^2}{2\pi}\,f(\theta),\qquad f(\theta)=\frac{\theta-\sin\theta\cos\theta}{\cos^2\theta}.
		\end{align*}
Since $f'(\theta)=\frac{2\theta\sin\theta}{\cos^3\theta}>0$, it follows that $\theta({\hat x},h')$ is strictly decreasing as ${\hat x}$ increases. Therefore
		$\partial_hU_{h'}({\hat x})=\pi/\theta({\hat x},h')$ is strictly increasing in ${\hat x}$.
		Let us also show that $U_{h'}'({\hat x})$ is strictly increasing in ${\hat x}$. Using
		${\hat x}=\sqrt{2\pi h'/[f(\theta)\gamma'']}$, we may rewrite
		\begin{align}\label{Uhprimo}
			U_{h'}'({\hat x})=\sqrt{2\pi h'/\gamma''}\,\psi(\theta),
			\qquad
			\psi(\theta):=\frac{\tan\theta}{\theta\sqrt{f(\theta)}}.
		\end{align}
		A direct computation shows that $\psi'(\theta)<0$, so $U_{h'}'({\hat x})$ is strictly decreasing as a function of
		$\theta$, hence strictly increasing as a function of ${\hat x}$.
		Differentiating $\widetilde U_h(t)=U_h(\Phi_{M,Q}(t))$ with respect to $t$ and
		evaluating at $t={\hat x}_0$ gives
		\begin{align*}
			\widetilde U_h'({\hat x}_0)
			&=U_h'(\Phi)\bigl(1+Q(M+2{\hat x}_0)\bigr) =(1-Q\eta)U_h'(\Phi).
		\end{align*}
		Hence, by the defining property of ${\hat x}_1={\hat x}_1^{[-2{\hat x}_0-\eta]}$,
		\begin{align}
			U_h'({\hat x}_1)=(1-Q\eta)U_h'(\Phi)<U_h'(\Phi).
			\label{eq:xprime-tilted-slope-match}
		\end{align}
              Note that $U_h' (\Phi) \asymp 1$, thanks to \eqref{Uhprimo} and the fact that $\theta=\pi\partial_y\mathtt H_{\gamma''}$ is bounded away from zero, which follows from \eqref{derivativesttH} together with the fact that the $S({\hat x},U_{h}({\hat x}))$ is bounded away from zero if $h$ is.  
		Since $U_h'$ is strictly increasing on $(0,\infty)$ and $U_h' (\Phi) \asymp 1$,  \eqref{eq:xprime-tilted-slope-match}
		implies $\Phi - {\hat x}_1^{[-2{\hat x}_0-\eta]} \asymp \eta Q$. Hence, since $\partial_x\partial_hU_h \gtrsim 1$ as shown above,
		\begin{align*}
			\partial_hU_{h'}(\Phi)-\partial_hU_{h'}({\hat x}_1^{[-2{\hat x}_0-\eta]} ) \gtrsim (\Phi - {\hat x}_1^{[-2{\hat x}_0-\eta ]}) \gtrsim \eta Q, 
		\end{align*}
	
		\noindent which shows \eqref{eq:common-A-h-derivative-quantitative-claim}. The claim $\inf\{ c(h',h,{\hat x}_0):h'\in [\varepsilon,h]\}>0$ follows from the fact that for $h'$ bounded away from zero, the argument of the square root in the definition of $\mathtt H_{\gamma''}$ is also bounded away from zero.
              \end{proof}

	\begin{proof}[Proof of \Cref{lem:3.1}]
          We will prove the statement with $\eta=0$. The claim then follows by continuity (indeed, the presence of $\eta$ cannot produce a linear term in $w$ in the r.h.s. of \eqref{eq:shifted-height-final-statement}, since the derivatives of $U_{h;M,Q}$ and $U_{h;0;0}$ at ${\hat x}_1$ and ${\hat x}_0$ coincide). The assumption ${\hat x}_0+w> \delta-\sqrt{2h/\gamma''}$ (and therefore also ${\hat x}_1+w> \delta-\sqrt{2h/\gamma''}+O(Q)$)  implies that the function $U_h$ in the proof is computed at values of ${\hat x}$ where it is smooth.
		Write
		\begin{align*}
			U(t):=U_{h;0,0}(t),
			\qquad
			\widetilde U(t):=U_{h;M,Q}(t),
			\qquad
			\Phi(t):=\Phi_{M,Q}(t).
		\end{align*}
		By \eqref{eq:xprime-common-choice}, we have ${\hat x}_1={\hat x}_0-Q{\hat x}_0^2=\Phi({\hat x}_0)$. Since $\widetilde U(t)=U(\Phi(t))$, we have $\widetilde U({\hat x}_0)=U({\hat x}_1)$. At the shifted point ${\hat x}'={\hat x}_0+w$,
		\begin{align*}
			\Phi({\hat x}')
			&={\hat x}'+Q\bigl(M{\hat x}'+({\hat x}')^2\bigr)
			={\hat x}_0+w-Q{\hat x}_0^2+Qw^2
			={\hat x}_1+w+Qw^2.
		\end{align*}
		Therefore
		\begin{align*}
			\widetilde U({\hat x}')=U({\hat x}_1+w+Qw^2).
		\end{align*}
		Subtracting the unperturbed shifted difference gives
		\begin{align*}
			&\bigl(\widetilde U({\hat x}')-\widetilde U({\hat x}_0)\bigr)-\bigl(U({\hat x}_1+w)-U({\hat x}_1)\bigr)=U({\hat x}_1+w+Qw^2)-U({\hat x}_1+w).
		\end{align*}
		Because ${\hat x}'> -\sqrt{2h/\gamma''}$ and ${\hat x}_1+w={\hat x}'+\mathcal{O}(Q)$, the points ${\hat x}_1+w$ and
		${\hat x}_1+w+Qw^2$ remain in a compact subinterval of $(-\sqrt{2h/\gamma''},\infty)$ for all
		sufficiently small $Q$.
		Taylor expansion at ${\hat x}_1+w$ yields
		\begin{align*}
			U({\hat x}_1+w+Qw^2)-U({\hat x}_1+w)
			&=Qw^2U'({\hat x}_1+w)+\mathcal{O}(Q^2).
		\end{align*} As proven above, $U'$ is positive and strictly increasing in ${\hat x}$. Moreover, $U'(-\sqrt{2h/\gamma''}+\delta)\gtrsim \delta^{1/2}$ holds by \eqref{uprimo}, since $\theta \gtrsim 1$ (as for $\delta \ll 1$ we have $|\theta - \pi| \ll 1$) and $S(\hat x,U_h(\hat x))\asymp \delta^{1/2}$, as is easily seen from the definition $S(\hat x,\hat y)=\sqrt{2\hat y-\gamma''\hat x^2}$. This implies \eqref{eq:shifted-height-final-statement}.
	\end{proof}

\bibliographystyle{alpha}
\bibliography{References}

\end{document}